\documentclass[11pt,a4,onesided]{amsart}
 \usepackage{epsfig}
\usepackage{amssymb}
\usepackage{color}
\usepackage{mathrsfs}
\usepackage{hyperref}
\newcommand\N{{\mathbb N}}
\newcommand\R{{\mathbb R}}
\newcommand\T{{\mathbb T}}
\newcommand\C{{\mathbb C}}

\def\AA{{\mathcal A}}
\def\BB{{\mathcal B}}
\def\CC{{\mathcal C}}
\def\DD{{\mathcal D}}
\def\EE{{\mathcal E}}

\def\GG{{\mathcal G}}
\def\HH{{\mathcal H}}

\def\JJ{{\mathcal J}}
\def\KK{{\mathcal K}}
\def\LL{{\mathcal L}}

\def\NN{{\mathcal N}}
\def\OO{{\mathcal O}}

\def\RR{{\mathcal R}}
\def\SS{{\mathcal S}}
\def\TT{{\mathcal T}}
\def\UU{{\mathcal U}}

\def\BBB{{\mathscr B}}
\def\CCC{{\mathscr C}}

\def\QQQ{{\mathscr Q}}

\newcommand{\init}[1]{#1_{\mbox{{\scriptsize in}}}}
\newcommand{\initem}[1]{#1_{\mbox{{\scriptsize {\em in}}}}}

\newcommand\dt{{\frac{\mathrm d}{\mathrm dt}}}

\newcommand{\dd}{{\, \mathrm d}}

\newcommand{\td}[2]{\frac{\mathrm d #1}{\mathrm d #2}}
\newcommand{\n}[1]{{\left\| #1 \right\|}}

\newcommand{\vs}{\vspace{0.3cm}}
\newcommand{\fa}{\forall \,}
\newcommand{\sk}{\smallskip}
\newcommand{\mk}{\medskip}

\newcommand{\var}{\varepsilon}
\newcommand{\Id}{\hbox{Id}}

\let\oldmarginpar\marginpar
\renewcommand\marginpar[1]{\-\oldmarginpar[\raggedleft\footnotesize #1]%
{\raggedright\footnotesize #1}}
\newcommand{\Black}{\color{black}}

\def\eps{{\varepsilon}}

\newtheorem{theo}{Theorem}
\newtheorem{prop}[theo]{Proposition}
\newtheorem{lem}[theo]{Lemma}
\newtheorem{cor}[theo]{Corollary}

\theoremstyle{definition}
\newtheorem{defin}[theo]{Definition}

\theoremstyle{remark}
\newtheorem{rem}[theo]{Remark}
\newtheorem{rems}[theo]{Remarks}

\newcommand{\beqn}{\begin{equation}}
\newcommand{\eeqn}{\end{equation}}
\newcommand{\bear}{\begin{eqnarray}}
\newcommand{\eear}{\end{eqnarray}}
\newcommand{\bean}{\begin{eqnarray*}}
\newcommand{\eean}{\end{eqnarray*}}

\newcommand{\ds}{\displaystyle}

\def\Nt{|\hskip-0.04cm|\hskip-0.04cm|}
\def\signsm{\bigskip \begin{center} {\sc St\'ephane
      Mischler\par\vspace{1mm} Universit\'e Paris-Dauphine \& IUF\par
      CEREMADE, UMR CNRS 7534\par
      Place du Mar\'echal de Lattre de Tassigny \par
      75775 Paris Cedex 16 FRANCE\par\vspace{1mm} e-mail:}
    \tt{mischler@ceremade.dauphine.fr} \end{center}}

\def\signcm{\bigskip \begin{center} {\sc Cl\'ement
      Mouhot\par\vspace{1mm}
      University of Cambridge\par
      DPMMS, Centre for Mathematical Sciences\par
      Wilberforce road, Cambridge CB3 0WA, UK
      \par\vspace{1mm} e-mail:}
    \tt{C.Mouhot@dpmms.cam.ac.uk} \end{center}}
    
\def\signmpg{\bigskip \begin{center} {\sc Maria Pia
      Gualdani\par\vspace{1mm} George Washington University \par
      Mathematics Department \par
      2115 G Street NW \par
      Washington DC, 20052, USA\par\vspace{1mm} e-mail:}
    \tt{gualdani@gwu.edu} \end{center}}

\begin{document}

\title[Factorization of Non-Symmetric Operators\dots]{Factorization of
  Non-Symmetric Operators and Exponential $H$-Theorem}

\author{M.P. Gualdani, S. Mischler, C. Mouhot}

% \footnotetext[1]{Department of Mathematics, The University of Texas
%   at Austin, 1 University Station C1200, Texas 78712, USA.  E-mail:
%   \texttt{gualdani@math.utexas.edu}}

% \footnotetext[2]{CEREMADE, Universit\'e Paris IX-Dauphine, Place du
%   Mar\'echal de Lattre de Tassigny, 75775 Paris, France.  E-mail:
%   \texttt{mischler@ceremade.dauphine.fr}}

% \footnotetext[3]{CEREMADE, Universit\'e Paris IX-Dauphine, Place du
%   Mar\'echal de Lattre de Tassigny, 75775 Paris, France.  E-mail:
%   \texttt{cmouhot@ceremade.dauphine.fr}}

\maketitle

\begin{center} {\bf Preliminary version of \today}
\end{center}

\begin{abstract} 
  We present an abstract method for deriving decay estimates on the
  resolvents and semigroups of non-symmetric operators in Banach
  spaces in terms of estimates in another smaller reference Banach
  space. This applies to a class of operators writing as a
  regularizing part, plus a dissipative part. The core of the method
  is a high-order quantitative factorization argument on the
  resolvents and semigroups.  We then apply this approach to the
  Fokker-Planck equation, to the kinetic Fokker-Planck equation in the
  torus, and to the linearized Boltzmann equation in the
  torus. % It requires each time significant
  % efforts and new estimates in order to prove the assumptions of the
  % abstract method. {\color{red}{MG: i'd not put this last sentence. it
  %     gives a negative tone to the general method}}

  We finally use this information on the linearized Boltzmann
  semigroup to study perturbative solutions for the nonlinear
  Boltzmann equation. We introduce a non-symmetric energy method to
  prove nonlinear stability in this context in $L^1_v L^\infty_x
  (1+|v|^k)$, $k >2$, with sharp rate of decay in time.

  As a consequence of these results we obtain the first constructive
  proof of exponential decay, with sharp rate, towards global
  equilibrium for the full nonlinear Boltzmann equation for hard
  spheres, conditionally to some smoothness and (polynomial) moment
  estimates. This improves the result in
  \cite{Desvillettes-Villani-EB} where polynomial rates at any order
  were obtained, and solves the conjecture raised in
  \cite{MR554086,Ce82,Reza-Villani} about the optimal decay rate of
  the relative entropy in the $H$-theorem.
 %  We apply these functional analysis results to several PDEs: the (kinetic)
%   Fokker-Planck equations and the linear(ized) Boltzmann equation in the
%   space homogeneous setting or in the confined case. As an application, we
%   deduce from our result the first result of exponential decay towards
%   global equilibrium for the full Boltzmann equation conditionnally to some
%   smoothness and moments estimates (improving on the result
%   \cite{Desvillettes-Villani-EB} where the rate was ``almost exponential'',
%   that is polynomial with exponent as high as wanted). % A quantitative
% %   close-to-homogeneous theory for the Boltzmann equation can also be
% %   steadily obtained from our results, improving on the non-constructive
% %   results \cite{Arkeryd-Esposito-Pulvirenti}.
\end{abstract}

\vspace{1.2cm}

One-parameter semigroups and linear evolution equations [See also
34G10, 34K30], 35P15 Estimation of eigenvalues, upper and lower
bounds, 47H20 Semigroups of nonlinear operators [See also 37L05,
47J35, 54H15, 58D07], 35Q84 Fokker-Planck equations, 76P05 Rarefied
gas flows, Boltzmann equation [See also 82B40, 82C40, 82D05].

\vspace{1cm}


\newpage

\tableofcontents

\newpage

%%%%%%%%%%%%%%%%%%%% Introduction %%%%%%%%%%%%%%%%%%%%%%%%%%%%%%%%%%%

\bigskip
\section{Introduction}
\label{sec:intro}
\setcounter{equation}{0}
\setcounter{theo}{0}

\bigskip
\Black

\subsection{The problem at hand}
\label{sec:problem-at-hand}

This paper deals with (1) the study of resolvent estimates and decay
properties for a class of linear operators and semigroups, and (2) the
study of relaxation to equilibrium for some kinetic evolution
equations, which makes use of the previous abstract tools.

Let us give a brief sketch of the first problem. Consider two Banach
spaces $E \subset \EE$, and two unbounded closed linear operators $L$
and $\LL$ respectively on $E$ and $\EE$ with spectrum $\Sigma(L),
\Sigma(\LL) \subsetneq \C$. They generate two $C_0$-semigroups $S(t)$
and $\SS(t)$ respectively in $E$ and $\EE$. Further assume that
$\LL_{|E} = L$, and $E$ is dense in $\EE$.  The theoretical question
we address in this work is the following:
\begin{center}
  \emph{Can one deduce quantitative informations on $\Sigma(\LL)$
    and $\SS(t)$ in terms of informations on $\Sigma(L)$ and $S(t)$?}
\end{center}

% Assuming that
% some informations on the localization of the spectrum $\Sigma(L)$ and
% on the asymptotic behavior of the semigroup $S(t)$ are known, can we
% obtain similar informations on the localization of the spectrum
% $\Sigma(\LL)$ and on the asymptotic behavior of the semigroup
% $\SS(t)$? In other words, can we carry over the larger space $\HH$
% some information we have on $\LL$ in the smaller space $H$?
We provide here an answer for a class of operators $\LL$ which split
as $\LL = \AA + \BB$, where the spectrum of $\BB$ is well localized
and the iterated convolution $(\AA \SS_\BB )^{*n}$ maps $\EE$ to $E$
with proper time-decay control for some $n \in \N^*$. We then prove that
(1) $\LL$ inherits most of the spectral gap properties of $L$, (2)
explicit estimates on the rate of decay of the semigroup $\SS(t)$ can
be computed from the ones on $S(t)$. The core of the proposed method
is a quantitative and robust factorization argument on the resolvents
and semigroups, reminiscent of the Dyson series.

In a second part of this paper, we then show that the kinetic
Fokker-Planck operator and the linearized Boltzmann operator for hard
sphere interactions satisfy the above abstract assumptions, and we
thus extend their spectral-gap properties from the linearization space
(a $L^2$ space with Gaussian weight prescribed by the equilibrium) to
larger Banach spaces (for example $L^p$ with polynomial decay). It is
worth mentioning
%that we provide results in the non spatially homogeneous
%situation, and
that the proposed method provides optimal rate of decay and there is
no loss of accuracy in the extension process from $E$ to $\EE$ (as
would be the case in, say, interpolation approaches).

Proving the abstract assumption requires significant technical efforts
for the Boltzmann equation and leads to the introduction of new tools:
some specific estimates on the collision operator, some iterated averaging
lemma and a nonlinear non-symmetric energy method. As a conclusion we
obtain a set of new stability results for the Boltzmann equation for
hard spheres interactions in the torus as discussed in the next
section.

\subsection{Motivation}
\label{sec:math-phys-motiv}

The motivation for the abstract part of this paper, i.e. enlarging the
functional space where spectral properties are known to hold for a
linear operator, comes from nonlinear PDE analysis. 

The first motivation is when the linearized stability theory of a
nonlinear PDE is \emph{not} compatible with the nonlinear theory. More
precisely the natural function space where the linearized equation is
well-posed and stable, with nice symmetric or skew-symmetric
properties for instance, is ``too small'' for the nonlinear PDE in the
sense that no well-posedness theorem is known even locally in time (or
even conjectured to be false) in such a small space. This is
the case for the classical Boltzmann equation and therefore it is a
key obstacle in obtaining perturbative result in natural physical
spaces and connecting the nonlinear results to the perturbative
theory.

This is related to the famous $H$-theorem of Boltzmann. The natural
question of understanding mathematically the $H$-theorem was
emphasized by Truesdell and Muncaster \cite[pp 560-561]{MR554086}
thirty years ago: ``{\it Much effort has been spent toward proof that
  place-dependent solutions exist for all time. [\dots] The main
  problem is really to discover and specify the circumstances that
  give rise to solutions which persist forever.  Only after having
  done that can we expect to construct proofs that such solutions
  exist, are unique, and are regular.}''  

The precise issue of the rate of convergence in the $H$-theorem was
then put forward by Cercignani \cite{Ce82} (see also \cite{Ce88}) when
he conjectured a linear relationship between the entropy production
functional and the relative entropy functional, in the spatially
homogeneous case.  % This ``functional inequality'' version of the
% $H$-theorem
% of it was conjectured by Cercignani~\cite{Ce82} in the spatially
% homogeneous case, in the form of a linear control of the relative
% entropy in terms of the entropy production functional. This conjecture
%has motivated important works from the early 1990's on. 
While this conjecture has been shown to be false in general
\cite{BC99}, it gave a formidable impulse to the works on the
Boltzmann equation in the last two decades
\cite{CC94,CC92,TV99,BC99,Vi03a}. It has been shown to be almost true
in \cite{Vi03a}, in the sense that polynomial inequalities relating
the relative entropy and the entropy production hold for powers close
to $1$, and it was an important inspiration for the work
\cite{Desvillettes-Villani-EB} in the spatially inhomogeneous case.

However, due to the fact that Cercignani's conjecture is false for
physical models \cite{BC99}, these important progresses in the far
from equilibrium regime were unable to answer the natural conjecture
about the correct timescale in the $H$-theorem, and to
prove the {\em exponential decay in time of the relative
  entropy}. Proving this exponential rate of relaxation was thus
pointed out as a key open problem in the lecture notes
\cite[Subsection~1.8, page 62]{Reza-Villani}. % about his breakthrough
% result together with Desvillettes already mentioned of convergence to
% equilibrium with polynomial rates.
This has motivated the work \cite{Mcmp} which answers this question,
but only in the spatially homogeneous case. % The core idea was to
% connect the nonlinear estimates with the linearized theory (providing
% exponential decay estimates) through the enlargement of sectorial
% decay estimates on the linearized semigroup in order to overcome the
% functional spaces incompatibility above-mentioned.  

In the present paper we answer this question for the full Boltzmann
equation for hard spheres in the torus.  We work in the same setting
as in~\cite{Desvillettes-Villani-EB}, that is under some {\it a
  priori} regularity assumptions (Sobolev norms and polynomial moments
bounds). We are able to connect the nonlinear theory in
\cite{Desvillettes-Villani-EB} with the perturbative stability theory
first discovered in \cite{MR0363332} and then revisited with
quantitative energy estimates in several works including
\cite{MR2013332} and \cite{MNeu}. This connexion relies on the
development of a perturbative stability theory in natural physical
spaces thanks to the abstract extension method. Let us mention here
the important papers \cite{MR946973,MR900501,MR1215007,MR1338453}
which proved for instance nonlinear stability in spaces of the form
$L^1_v W^{s,p}_x(1+|v|^k)$ with $s>3/p$ and $k>0$ large enough, by
non-constructive methods.

We emphasize the dramatic gap between the spatially homogeneous
situation and the spatially inhomogeneous one. In the first case the
linearized equation is coercive and the linearized semigroup is
self-adjoint or sectorial, whereas in the second case the equation is
{\em hypocoercive} and the linearized semigroup is neither sectorial,
nor even hypoelliptic.

The second main motivation for the abstract method developed here is
considered in other papers \cite{MMcmp,AGGMMS}. It concerns the
existence, uniqueness and stability of stationary solutions for
degenerate perturbations of a known reference equation, when the
perturbation makes the steady solutions leave the natural
linearization space of the reference equation. Further works
concerning spatially inhomogeneous granular gases are in
progress. % (thanks to
% perturbative or implicit function argument, see for instance) as well
% as the proof of long time convergence to the equilibrium (see for
% instance \cite{Mcmp,MMcmp}).  In this paper we construct a general
% abstract approach including the examples of these previous papers as
% particular cases.  Within this general framework we give new results
% for the Fokker-Planck equation,
% %the (hypoelliptic) kinetic Fokker-Planck equation, 
% and most
% importantly the full (inhomogeneous) Boltzmann equation.
%{\Green\tiny We also refer to the note \cite{Mouhot-ISAAC} where some results of this paper were announced.}

%\smallskip
%{\Red Fausse piste, vraie reponse ici et depuis [Mcmp2006,MMopus]}
% The main outcome of this paper is thus {\bf a proof of the exponential
%   convergence to equilibrium for the full Boltzmann equation in the
%   torus, under {\it a priori} smoothness assumptions}. 

\subsection{Main results}
\label{sec:main-results}

We can summarize the main results established in this paper as
follows: 
\sk 

\noindent {\em Section~2.} We prove an abstract theory for enlarging
(Theorem~\ref{theo:factor:enlarg}) the space where the spectral gap
and the discrete part of the spectrum is known for a certain class of
unbounded closed operators. We then prove a corresponding abstract
% (although more ``concrete'' when applications to PDE are considered) 
theory for enlarging (Theorem~\ref{theo:EnlargingSGdecay}) the space
where explicit decay semigroup estimates are known, for this class of
operators. This can also be seen as a theory for obtaining
quantitative spectral mapping theorems in this setting, and it works
in Banach spaces.  \sk

\noindent {\em Section~3.} We prove a set of results concerning
Fokker-Planck equations.  The main outcome is the proof of an explicit
spectral gap estimate on the semigroup in $L^1_{x,v}(1+|v|^k)$, $k>0$ as
small as wanted, for the kinetic Fokker-Planck equation in the torus
with super-harmonic potential (see Theorems~\ref{theo:FP1}
and~\ref{theo:KFP2}).  \sk

\noindent {\em Section~4.} We prove a set of results concerning the
linearized Boltzmann equation. The main outcome is the proof of
explicit spectral gap estimates on the linearized semigroup in $L^1$
and $L^\infty$ with polynomial moments (see
Theorem~\ref{theo:LIBE2}). More generally we prove explicit spectral
gap estimates in any space of the form $W^{\sigma,q}_vW^{s,p}_x(m)$,
$\sigma \le s$, with polynomial or stretched exponential weight $m$,
including the borderline cases $L^\infty_{x,v}(1+|v|^{5+0})$ and
$L^1_vL^\infty_x(1+|v|^{2+0})$. We also make use of the factorization
method in order to study the structure of singularities of the
linearized flow (see Subsection~\ref{sec:study-sing-struct}).

\sk 

\noindent {\em Section~5.} We finally prove a set of results concerning the
nonlinear Boltzmann equation in perturbative setting. The main
outcomes of this section are: (1) The construction of perturbative
solutions close to the equilibrium or close to the spatially
homogeneous case in $W^{\sigma,q}_vW^{s,p}_x(m)$, $s>6/p$ with
polynomial or stretched exponential weight $m$, including the borderline
cases $L^\infty_{x,v}(1+|v|^{5+0})$ and $L^1_vL^\infty_x(1+|v|^{2+0})$ without
assumption on the derivatives: see Theorem~\ref{theo:NLBE} in a
close-to-equilibrium setting, and Theorem~\ref{theo:NLBE+} in a
close-to-spatially-homogeneous setting. (2) We give a proof of the
exponential $H$-theorem: we show exponential decay in time of the
relative entropy of solutions to the fully nonlinear Boltzmann
equation, conditionnally to some regularity and moment bounds. Such
rate is proven to be sharp.  This answers the conjecture in
\cite{Desvillettes-Villani-EB,Reza-Villani} (see
Theorem~\ref{theo:expEB}). We also finally apply the factorization
method and the Duhamel principle to study the structure
of singularities of the nonlinear flow in perturbative regime (see
Subsection~\ref{sec:struct-sing-nonl}).  \bigskip

Below we give a precise statement of what seems to us the main result
established in this paper.
\begin{theo}\label{theo:BEglobal}
  The Boltzmann equation 
  \begin{equation*}
\left\{ 
\begin{array}{l} \ds
  \partial_t f + v \cdot \nabla_x f = Q(f,f), \quad t \ge 0, \ x \in
  \T^3, \ v \in \R^3, \vs \\ \ds
  Q(f,f)  :=  
  \int_{\R^3} \int_{\mathbb{S}^{2}}
  \Big[f(x,v') \, f(x,v'_*) - f(x,v) \, f(x,v_*) \Big] \, |v-v_*| \dd v_* \dd \sigma \vs \\ \ds
  v' = \frac{v+v_*}2 + \sigma \, \frac{|v-v_*|}2, \qquad v'_* =
  \frac{v+v_*}2 - \sigma \, \frac{|v-v_*|}2
\end{array}
\right.
\end{equation*}
with hard spheres collision kernel and periodic boundary conditions is
globally well-posed for non-negative initial data close enough to the
Maxwellian equilibrium $\mu$ or to a spatially homogeneous profile in
$L^1_vL^\infty_x(1+|v|^k)$, $k >2$. 

The corresponding solutions decay exponentially fast in time with
constructive estimates and with the same rate as the linearized flow
in the space $L^1_vL^\infty_x(1+|v|^k)$. For $k$ large enough (with
explicit threshold) this rate is the sharp rate $\lambda>0$ given by
the spectral gap of the linearized flow in $L^2(\mu^{-1/2})$.

Moreover any solution that is a priori bounded uniformly in time in
$H^s_{x,v}(1+|v|^k)$ with some large $s,k$ satisfies the exponential
decay in time with sharp rate $O(e^{-\lambda \, t})$ in $L^1$ norm, as
well as in relative entropy. 
\end{theo}

\subsection{Acknowledgments} 
We thank Claude Bardos, Jos\'e Alfr\'edo Ca\~nizo, Miguel Escobedo,
Bertrand Lods, Mustapha Mokhtar-Kharroubi, Robert Strain for fruitful
comments and discussions.  The third author also wishes to thank
Thierry Gallay for numerous stimulating discussions about the spectral
theory of non-self-adjoint operators, and also for pointing out the
recent preprint \cite{Helffer-Sjoestrand}. The authors wish to thank
the funding of the ANR project MADCOF for the visit of MPG in
Universit\'e Paris-Dauphine in spring 2009 where this work was
started. The third author's work is supported by the ERC starting
grant MATKIT. The first author is supported by NSF-DMS
1109682. Support from IPAM (University of California Los Angeles) and
ICES (The University of Texas at Austin) is also gratefully
acknowledged.

%%%%%%%%%%%%%%%%%%%%%%%%%%%%%%%%%%%%%%%%%%%%
%\newpage

%\bigskip
\section{Factorization and quantitative spectral mapping theorems}
\label{sec:factorization}
\setcounter{equation}{0}
\setcounter{theo}{0}
%\bigskip\bigskip

% \footnotesize
% \Green{\bf 
% A FAIRE:

% - Si on veut pouvoir \^etre citer il nous faut un th\'eor\`eme (ou
% corollaire) qui r\'esume la situation usuelle dans laquelle on peut
% appliquer notre m\'ethode.

% - section 5. Dans les deux cas "triviaux" (Boltzmann lineaire homogene et
% Fokker-Planck homogene) pr\'esenter les d\'etails de la preuve du bon
% comportement du semi-groupe $e^{t L}$. Cela ne rajoutera pas plus de deux
% pages et rendra l'article vraiment abordable ``pour les nuls".

% - Question 1. Introduire tot (section 3) la notion de projecteur au sens du
% Kato (bien citer). Quand a-t-on $(I-\Pi) E = $ orthogonal de $E_1$ au sens
% du produit scalaire de $E$. Je suppose que cela est ok lorsque $E$ est un
% Hilbert.  Dans ce cas les deux notions de projecteurs coincident. NON il
% faut que $L$ soit auto-adjoint. On utilise ici autre chose $\LL^* \phi =
% 0$?

% - Question 3. Revoir Boltzmann cinetique linearise qui me pose peut-\^etre
% un probleme. Celui-ci d'abord o\`u on doit passer de $L^1_{x,v}$ \`a
% $L^2_{x,v}$. De plus, dans ce cas sait-on d\'emontrer $N(\LL) = N(L)$? Ce
% n'est pas clair. Mais ce n'est pas clair que l'on en ait besoin

% }\Black\normalsize

% \bigskip

\subsection{Notation and definitions}
\label{sec:notation-definitions}

For a given real number $a \in \R$, we define the half complex plane
$$
\Delta_a := \left\{ z \in \C, \, \Re e \, z > a \right\}.
$$
% and more generally, for $\theta \in [\pi/2,\pi)$, we also consider the (larger) sectorial set 
% \[
% \Delta_{a,\theta} := \left\{ z \in \C, \ | \arg(z-b)| \le \theta \mbox{ for some } b >a \right\}.
% \]
% Note that it includes the previous case of a half-plane when
% $\theta=\pi/2$: $\Delta_a = \Delta_{a,\pi/2}$. 

\smallskip
For some given Banach spaces $(E,\|\cdot \|_E)$ and $(\EE,\| \cdot
\|_\EE)$ we denote by $\mathscr{B}(E, \EE)$ the space of bounded linear
operators from $E$ to $\EE$ and we denote by $\| \cdot
\|_{\mathscr{B}(E,\EE)}$ or $\| \cdot \|_{E \to \EE}$ the associated norm
operator. We write $\mathscr{B}(E) = \mathscr{B}(E,E)$ when $E=\EE$.
We denote by $\mathscr{C}(E,\EE)$ the space of closed unbounded linear
operators from $E$ to $\EE$ with dense domain, and $\mathscr{C}(E)=
\mathscr{C}(E,E)$ in the case $E=\EE$. 

\smallskip For a Banach space $X$ and $\Lambda \in \mathscr{C}(X)$ we
denote by $S_\Lambda(t)$, $t \ge 0$, its semigroup, by
$\textrm{Dom}(\Lambda)$ its domain, by $N(\Lambda)$ its null space and
by $\mbox{R}(\Lambda)$ its range. We also denote by $\Sigma(\Lambda)$
its spectrum, so that for any $z$ belonging to the resolvent set
$\rho(\Lambda) := \C \backslash \Sigma(\Lambda)$ the operator $\Lambda
- z$ is invertible and the resolvent operator
$$
\RR_\Lambda(z) := (\Lambda -z)^{-1}
$$
is well-defined, belongs to $\mathscr{B}(X)$ and has range equal to
$D(\Lambda)$.  We recall that $\xi \in \Sigma(\Lambda)$ is said to be
an eigenvalue if $N(\Lambda - \xi) \neq \{ 0 \}$. Moreover an
eigenvalue $\xi \in \Sigma(\Lambda)$ is said to be \emph{isolated} if
\[
\Sigma(\Lambda) \cap \left\{ z \in \C, \,\, |z - \xi| \le r \right\} =
\{ \xi \} \ \mbox{ for some } r >0.
\]
In the case when $\xi$ is an isolated eigenvalue we may define
$\Pi_{\Lambda,\xi} \in \mathscr{B}(X)$ the associated spectral projector by 
\begin{equation}
\label{def:SpectralProjection} 
\Pi_{\Lambda,\xi} := - {1 \over 2i\pi} \int_{ |z - \xi| = r' } \RR_\Lambda(z) \dd z
\end{equation}
with $0<r'<r$. Note that this definition is independent of the value
of $r'$ as the application $ \C \setminus \Sigma(\Lambda) \to
\mathscr{B}(X)$, $z \to \RR_{\Lambda}(z)$ is holomorphic.  For any
$\xi \in \Sigma(\Lambda)$ isolated, it is well-known (see
\cite[III-(6.19)]{Kato}) that $\Pi_{\Lambda,\xi}^2=\Pi_{\Lambda,\xi}$,
so that $\Pi_{\Lambda,\xi}$ is indeed a projector, and that the
associated projected semigroup
\[
S_{\Lambda,\xi}(t) := -\frac{1}{2i \pi} \int_{|z-\xi|=r'} e^{zt}
\RR_\Lambda (z) \dd z , \quad t \ge 0,
\]
satisfies 
\begin{equation}
\label{def:SpectralPSg} 
S_{\Lambda,\xi}(t)  
= \Pi_{\Lambda,\xi}  S_{\Lambda}(t) =
S_\Lambda(t)  \Pi_{\Lambda,\xi}, \quad t \ge 0.
\end{equation}

\smallskip When moreover the \emph{algebraic eigenspace}
$\mbox{R}(\Pi_{\Lambda,\xi})$ is finite dimensional we say that $\xi$
is a \emph{discrete eigenvalue}, written as $\xi \in
\Sigma_d(\Lambda)$. In that case, $\RR_\Lambda$ is a meromorphic
function on a neighborhood of $\xi$, with non-removable finite-order
pole $\xi$, and there exists $\alpha_0 \in \N^*$ such that
\[
\mbox{R}(\Pi_{\Lambda,\xi}) = N(\Lambda -\xi)^{\alpha_0} = N(\Lambda
-\xi)^\alpha \ \mbox{ for any } \ \alpha \ge \alpha_0.
\]
On the other hand, for any $\xi \in \C$ we may also define the
``classical algebraic eigenspace"
$$
M(\Lambda-\xi) := \lim_{\alpha \to \infty} N(\Lambda-\xi)^\alpha.
$$
We have then $M(\Lambda-\xi) \not= \{0\}$ if $\xi \in \Sigma(\Lambda)$
is an eigenvalue and $M(\Lambda-\xi) = \mbox{R}(\Pi_{\Lambda,\xi})$ if
$\xi$ is an isolated eigenvalue. %$\xi \in \Sigma_d(\Lambda)$.
  
\smallskip
Finally for any $a \in \R$ %, \theta \in [\pi/2,\pi)$ 
such that
\[
\Sigma(\Lambda) \cap \Delta_{a} = \left\{ \xi_1, \dots, \xi_k\right\}
\]
where $\xi_1, \dots, \xi_k$ are distinct discrete eigenvalues, we define
without any risk of ambiguity 
\[
\Pi_{\Lambda,a} := \Pi_{\Lambda,\xi_1} + \dots + \Pi_{\Lambda,\xi_k}.
\]

\subsection{Factorization and spectral analysis}
\label{sec:spectr-analys-fact}

The main abstract factorization and enlargement result is:
% In both case, note that it is required that the
% \emph{smaller} functional space is a Hilbert space.
%% CM: Hilbert pas necessaire a ce niveau, mais seulement dans la
%% suite avec les semigroupes. 
\begin{theo}[Enlargement of the functional space]\label{theo:factor:enlarg}
  Consider two Banach spaces $E$ and $\EE$ such that $E \subset \EE$
  with continuous embedding and $E$ is dense in $\EE$.  Consider an
  operator $\LL \in \mathscr{C}(\EE)$ such that $L := (\LL) _{|E} \in
  \mathscr{C}(E)$.  Finally consider a set $\Delta_{a}$ as
  defined above.

We assume:
  \begin{itemize}
  \item[{\bf (H1)}] {\bf Localization of the spectrum in $E$.}  There
    are some distinct complex numbers $\xi_1, \dots, \xi_k \in
    \Delta_a$, $k \in \N$ (with the convention $\{\xi_1, \dots, \xi_k
    \}= \emptyset$ if $k=0$) such that
\[
\Sigma(L) \cap \Delta_a = \left\{ \xi_1, \dots,
    \xi_k \right\} \subset \Sigma_d(L) \quad \mbox{(distinct discrete eigenvalues).}
\]

\item[{\bf (H2)}] {\bf Decomposition.} There exist  $\AA,\BB $ some operators  defined on $\EE$ such that  
$\LL = \AA + \BB$ and 
    \begin{itemize}
    \item[(i)] $\BB \in \CCC(\EE) $ is such that  $\RR_\BB(z)$ is  bounded in $\BBB(\EE)$ uniformly on $z \in \Delta_a$
    and $\| \RR_\BB(z) \|_{\BBB(\EE)} \to 0$ as $\Re e \, z \to \infty$, 
    %a closed unbounded operator on $\EE$   with domain containing $\mbox{\textrm{Dom}}(\LL)$  and
      in particular 
      \[
      \Sigma(\BB) \cap \Delta_a = \emptyset; 
      \]
     \item[(ii)] $\AA \in \mathscr{B}(\EE)$ is a bounded operator on $\EE$; 
     \item[(iii)]    There is $n \ge 1$ such that the operator $(\AA \RR_{\BB}(z))^n$
    is  bounded in $\BBB(\EE,E)$ uniformly on $z \in \Delta_a$.
      
%    There is $n \ge 1$ such that the operator $(\AA
%       \RR_{\BB}(z))^n$ is bounded from $\EE$ to $E$ locally uniformly in $z \in \Delta_a \backslash \{ \xi_1, ... , \xi_k \} $. 
%       
%     { \Green or (but that is the same!) 
%     Moreover the application
%       \[
%       \Delta_a  \backslash \{ \xi_1, ... , \xi_k \} \to \mathscr{B}(\EE,E), \quad z \mapsto (\AA R_{\BB}(z))^n 
%       \]
%       is continuous. }
%       % on $\Delta_a$ with values in $\mathscr{B}(\EE,E)$.
% Moreover the application
       % \[
       % \Delta_a \to \mathscr{B}(\EE,E), \quad z \mapsto (\AA R_{\BB}(z))^n 
       % \]
       % is holomorphic on $\Delta_a$ with values in $\mathscr{B}(\EE,E)$.
     \end{itemize}
\smallskip
%
%   \item[{\bf (H3)}] 
%   %{\bf \Green\footnotesize Technical assumption for avoiding pathological situations.} 
%       There is some $z \in \Delta_a$ such
%     that $(\LL - z)$ is invertible in $\EE$. 
%     %     {\Green\footnotesize and the application $z \mapsto (\AA
%     %     R_{\BB}(z))^n$ is holomorphic on $\Delta_a$ with
     %values in $\mathscr{B}(\EE,E)$.}
  \end{itemize}
\smallskip

  Then we have in $\EE$: 
  \begin{itemize}
  \item[(i)] The spectrum satisfies: 
    $\Sigma(\LL) \cap \Delta_a = \{ \xi_1, \dots, \xi_k \}$.
\smallskip

\item[(ii)] For any $z \in \Delta_a \setminus \{ \xi_1, \dots, \xi_k
  \}$ the resolvent satisfies: 
\begin{equation}
    \label{eq:factor-enlarg}
    \RR_{\LL}(z) = \sum_{\ell=0} ^{n-1} (-1)^\ell \RR_{\BB}(z) \left( \AA
      \RR_{\BB}(z) \right)^\ell + (-1)^n \RR_{L}(z) \left( \AA
      \RR_{\BB}(z) \right)^n.
  \end{equation}

\item[(iii)] For any $\xi_i \in \Sigma(L) \cap \Delta_a = \Sigma(\LL)
  \cap \Delta_a$, $i=1, \dots,
  k$, we have 
\[
\forall \, m \ge 1, \quad N(L -\xi_i)^m = N(\LL
-\xi_i)^m 
\ \mbox{ and } \ 
M(L -\xi_i) = M(\LL -\xi_i)
\]
and at the level of the spectral projectors
\[
\left\{ 
\begin{array}{l}
(\Pi_{\LL,\xi_i})_{|E} = \Pi_{L,\xi_i}
\vspace{0.2cm} \\
S_{\LL,\xi_i}(t) = S_{L,\xi_i}(t) \Pi_{\LL,\xi_i} = S_{L}(t)
\Pi_{\LL,\xi_i}.
\end{array}
\right.
\]
 \end{itemize}
\end{theo}

\begin{rems}\label{rem:factorization} 
  \begin{enumerate}
   % \item Let us comment on the global meaning of this theorem: as
   %  $E \subset \EE$, and the assumptions on the spectral gap are within
   %  the space $E$, this theorem is a recipee for \emph{enlarging} the
   %  functional space of a spectral gap.
    
  \item In words, assumption {\bf (H1)} is a weak formulation of a
    spectral gap in the initial functional space $E$. The assumption
    {\bf (H2)} is better understood in the simplest case $n=1$, where
    it means that one may decompose $\LL$ into a regularizing part
    $\AA$ (in the generalized sense of the ``change of space'' $\AA
    \in \mathscr{B}(\EE,E)$) and another part $\BB$ whose spectrum is
    ``well localized'' in $\EE$: 
  % \item In the simplest case where $\AA \in \mathscr{B}(\EE,E)$ is
  %   ``regularizing'', the assumption {\bf (H2)-(iii)} is trivially
  %   satisfied with $n=1$. 
    for instance when $\BB-a'$ is dissipative with $a'<a$ then the
    assumption {\bf (H2)-(i)} is satisfied.
  \item There are many variants of sets of hypothesis for the
    decomposition assumption.  In particular, assumptions
    \textbf{(H2)-(i)} and \textbf{(H2)-(iii)} could be
    weakened. However, (1) these assumptions are always fulfilled by
    the operators we have in mind, (2) when we weaken
    \textbf{(H2)-(i)} and/or \textbf{(H2)-(iii)} we have to compensate
    them by making other structure assumptions.  We present below
    after the proof a possible variant of
    Theorem~\ref{theo:factor:enlarg}.
 \item One may relax {\bf (H2)-(i)} into $\Sigma(\BB) \cap \Delta_a
    \subset \{ \xi_1, \dots, \xi_k \}$ and the bound in
    \textbf{(H2)-(iii)} could be asked merely locally uniformly in $z
    \in \Delta_a \backslash \{ \xi_1, \dots , \xi_k \} $. 
  \item One may replace $\Delta_a \setminus \{ \xi_1, \dots, \xi_k \}$
    by any nonempty open connected set $\Omega \subset \C$. 
  \item This theorem and the next ones in this section can also be
    extended to the case where $E$ is not necessarily included in
    $\EE$. This will be studied and applied to some PDE problems in
    future works.
  \end{enumerate}
\end{rems}

\begin{proof}[Proof of Theorem~\ref{theo:factor:enlarg}]
  Let us denote $\Omega := \Delta_a \setminus \{ \xi_1, \dots, \xi_k \}$
  and let us define for $z \in \Omega$
\[
\UU(z) := \sum_{\ell=0} ^{n-1} (-1)^\ell \RR_\BB(z) \left( \AA \RR_{\BB}(z) 
    \right)^\ell  + (-1)^n \RR_L(z) \left( \AA \RR_{\BB}(z)
    \right)^n.
\]
Observe that thanks to the assumption {\bf (H2)}, the operator $\UU(z)$
is well-defined and bounded on $\EE$.

\smallskip
\noindent {\it Step 1. $\UU(z)$ is a right-inverse of $(\LL-z)$ on
  $\Omega$.}  
For any $z \in \Omega$, we compute
 \begin{eqnarray*} 
   (\LL-z) \UU(z)
   &=& \sum_{\ell = 0}^{n-1} (-1)^\ell \, (\AA + (\BB-z)) \, \RR_\BB(z)
   \, \left(\AA \RR_\BB(z) \right)^\ell \\
   &&
   \quad + (-1)^n \, (\LL-z) \, \RR_\LL(z) \,  \left(\AA \RR_\BB(z)\right)^n  \\
   &=&  \sum_{\ell = 0}^{n-1} (-1)^\ell \, \left( \AA \RR_\BB(z) \right)^{\ell+1} 
   \quad +  \sum_{\ell = 0}^{n-1} (-1)^\ell \, \left( \AA \RR_\BB(z)\right)^\ell \\
   &&
   \quad + (-1)^n \,  \left(\AA \RR_\BB(z) \right)^n  =  \mbox{Id}_{\EE}. 
\end{eqnarray*} 

\smallskip
\noindent {\it Step 2. $(\LL-z)$ is invertible on $\Omega$.} First we
observe that there exists $z_0\in \Omega$ such that $(\LL-z_0)$ is
invertible in $\EE$. Indeed, we write
$$
\LL - z_0 = (\AA \, \RR_\BB(z_0) +  \mbox{Id}_{\EE} ) \, (\BB-z_0)
$$
with $\| \AA \, \RR_\BB(z_0) \| < 1$ for $z_0 \in \Omega$, $\Re e z_0$
large enough, thanks to assumption {\bf(H2)-(i)}. As a consequence
$(\AA \, \RR_\BB(z_0) +  \mbox{Id}_{\EE})$ is invertible and so is $\LL - z_0 $ as
the product of two invertible operators.

\Black

Since we assume that $(\LL-z_0)$ is invertible in $\EE$ for some $z_0
\in \Omega$, we have $\RR_\LL(z_0) = \UU(z_0)$. And if
\[
\| \RR_\LL(z_0) \|_{\mathscr{B}(E)} = \| \UU(z_0) \|_{\mathscr{B}(E)} \le C
\]
for some $C \in (0,\infty)$, then $(\LL-z)$ is invertible on the disc
$B(z_0,1/C)$ with
\begin{equation}
\label{eq:devResolvante} 
\forall \, z \in B(z_0,1/C), \quad 
\RR_\LL(z) = \RR_\LL(z_0)  \, \sum_{n=0}^\infty (z_0 - z)^n \, {%\Blue
 \RR_\LL( z_0)^{n}} ,
\end{equation} 
and then again, arguing as before, $\RR_\LL(z) = \UU(z)$ on $B(z_0,1/C)$
since $\UU(z)$ is a left-inverse of $(\LL-z)$ for any $z \in \Omega$.
Then in order to prove that $(\LL - z)$ is invertible for any $z \in
\Omega$, we argue as follows. For a given $z_1 \in \Omega$ we consider
a continuous path $\Gamma$ from $z_0$ to $z_1$ included in $\Omega$,
i.e. a continuous function $\Gamma : [0,1] \to \Omega$ such that
$\Gamma(0) = z_0$, $\Gamma(1) = z_1$.
Because of assumption {\bf
  (H2)} we know that $(\AA \RR_\BB(z))^\ell$, $1 \le \ell \le n-1$, and
${%\Blue 
\RR_L (z) } (\AA \RR_\BB(z))^n$ are locally uniformly bounded in $\mathscr{B}(\EE)$ on
$\Omega$, which implies
$$
\sup_{z \in \Gamma([0,1])} \| \UU(z) \|_{ \mathscr{B}(\EE)} := C_0 <
\infty.
$$
Since $(\LL - z_0)$ is invertible we deduce that $(\LL - z)$ is
invertible with $\RR_\LL(z)$ locally bounded around $z_0$ with a bound
$C_0$ which is uniform along $\Gamma$ (and a similar series expansion
as in \eqref{eq:devResolvante}). By a continuation argument we hence
obtain that $(\LL-z)$ is invertible in $\EE$ all along the path
$\Gamma$ with
\[
\RR_\LL(z) = \UU(z) \ \mbox{ and } \ \| \RR_\LL(z) \|_{\mathscr{B}(\EE)} =
\| \UU(z) \|_{\mathscr{B}(\EE)} \le C_0.
\]
Hence we conclude that $(\LL-z_1)$ is invertible with $\RR_\LL(z_1) =
\UU(z_1)$. 

This completes the proof of this step and proves $\Sigma(\LL) \cap
\Delta_a \subset \{ \xi_1, \dots, \xi_k \}$ together with the point (ii)
of the conclusion.

\smallskip
\noindent
{\it Step 3. Spectrum, eigenspaces and spectral projectors.} 
On the one hand, we have 
\[
N(L-\xi_j)^\alpha \subset
N(\LL-\xi_j)^\alpha, \quad j = 1, \dots, k, \ \alpha \in \N,
\]
so that $ \{ \xi_1, \dots, \xi_k \} \subset \Sigma(\LL) \cap \Delta_a$.
The other inclusion was proved in the previous step, so that these two
sets are equals.  We have proved
\[
\Sigma(\LL) \cap \Delta_a = \Sigma(L) \cap \Delta_a.
\]

Now, we consider a given eigenvalue $\xi_j$ of $L$ in $E$.  We know
(see \cite[paragraph I.3]{Kato}) that in $E$ the following Laurent
series holds
$$
\RR_L(z) = %z^{-1} \, R_{-1} +
\sum_{\ell=-\ell_0}^{+\infty} (z-\xi_j)^{\ell} \, \CC_\ell, \quad
\CC_{\ell} = (\LL-\xi_j)^{|\ell|-1} \Pi_{\LL,\xi_j}, \ \ell_0 \le \ell
\le -1, 
%L^{\ell-1}_{-2}+\sum_{\ell=0}^{\infty} z^{\ell} \, L^\ell_0
$$
for $z$ close to $\xi_j$ and for some bounded operators $\CC_\ell \in
\mathscr{B}(E)$, $\ell \ge 0$. The operators $\CC_{-1}, \dots,
\CC_{-\ell_0}$ satisfy the range inclusions
\[
\mbox{R}(\CC_{-2}), \dots, \mbox{R}(\CC_{-\ell_0})
\subset \mbox{R}(\CC_{-1}).
\]
This Laurent series is convergent on $B(\xi_j,r) \backslash \{\xi_j\}
\subset \Delta_a$. The Cauchy formula for meromorphic functions applied
to the circle $\{z, \,\, |z-\xi_j|=r \}$ with $r$ small enough thus
implies that
\[
\Pi_{L,\xi_j} = \CC_{-1} \ \mbox{ so that } \ \CC_{-1} \not = 0
\]
since $\xi_j$ is a discrete eigenvalue. 

Using the
definition of the spectral projection operator
\eqref{def:SpectralProjection}, the above expansions and the Cauchy
theorem we get for any small $r > 0$
\begin{eqnarray*}
 \Pi_{\LL,\xi_j} &:=& {(-1)^{n+1} \over 2i\pi}
\int_{|z-\xi_j| = r} \RR_L(z) \, (\AA \, \RR_\BB(z))^n \dd z
\\
&=& \int_{|z-\xi_j| = r} \sum_{\ell=\ell_0}^{-1} {\CC_\ell \, (z-\xi_j)^\ell} \, (\AA \, \RR_\BB(z))^n \dd z
\\
&&+ \int_{|z-\xi_j| = r} \sum_{\ell=0}^{\infty} {\CC_\ell \, (z-\xi_j)^\ell} \, (\AA \, \RR_\BB(z))^n \dd z,
\end{eqnarray*}
where the first integral has range included in $\mbox{R}(\CC_{-1})$ and the second
integral vanishes in the limit $r \to 0$. 
%
%Assumption {\bf (H3)} then implies that
%$$
%(\AA \, \RR_\BB(z))^n = \sum_{\beta=0}^{\infty} z^\beta \, \CC_\beta,
%\quad \CC_\beta \in \BBB(\EE,E),
%$$
%which is convergent (in norm) on the disc $B(\xi_j,r)$.  
We deduce  that 
\[
M(\LL-\xi_j) = \mbox{R}(\Pi_{\LL,\xi_j}) \subset \mbox{R} (\CC_{-1}) =
\mbox{R}(\Pi_{L,\xi_j}) = M(L-\xi_j).
\]
Together with 
\[
M(L-\xi_j) = N (L-\xi_j)^{\alpha_0} \subset N(\LL-\xi_j)^{\alpha_0}
\subset M(\LL-\xi_j) \ \mbox{ for some } \alpha_0 \ge 1
\]
we conclude that $M(\LL - \xi_j) = M(L - \xi_j)$ and
$N((\LL-\xi_j)^\alpha) = N((L-\xi_j)^\alpha)$ for any $j= 1, \dots, k$
and $\alpha \ge 1$.

Finally, the proof of $\Pi_{\LL,\xi_j| E}= \Pi_{L,\xi_j}$ is
straightfoward from the equality 
\[
\RR_\LL(z) f = \RR_L(z) f \ \mbox{ when } \ f \in E
\]
and the integral formula \eqref{def:SpectralProjection} defining the
projection operator.
\end{proof}

Let us shortly present a variant of the latter result where the
assumption {\bf (H2)} is replaced by a more algebraic one. The proof
is then purely based on the factorization method and somehow
simpler. The drawback is that it requires some additional assumption
on $B$ at the level of the \emph{small} space (which however is not so
restrictive for a PDE's application perspective but can be painful to
check).
%{\Red Non, absolument pas, (H3') n'est pas une hypoth\`ese ``alg\'ebrique",
%ni, et encore moins!, l'hypoth\`ese de la remarque puisqu'on y fait une 
%hypoth\`ese de surjectivit\'e qui est un r\'esultat ``d'existence" par excellence,
%donc c'est justement une hypoth\`ese d'analyse. En revanche, la preuve 
%du th\'eor\`eme est purement alg\'ebrique. } 

% {\color{red}{both theorems have the label " Enlargement of the functional
%   space". maybe a different label for this one?}}
\begin{theo}[Enlargement of the functional space, purely algebraic
  version]\label{theo:factor:enlarg:algeb}
  
  Consider the same setting as in Theorem~\ref{theo:factor:enlarg},
  assumption {\bf (H1)}, and where assumption {\bf (H2)} is replaced
  by
    
       \begin{itemize}
       \item[{\bf (H2')}] {\bf Decomposition.} There exist operators
         $A, B$ on $E$ such that $L = A + B$ (with corresponding
         extensions $\AA,\BB$ on $\EE$) and
    \begin{itemize}
    \item[(i$'$)] $B$ and $\BB$ are closed unbounded operators on $E$
      and $\EE$ (with domain containing $\mbox{\textrm{Dom}}(L)$ and
      $\mbox{\textrm{Dom}}(\LL)$) and
      \[
      \Sigma(B) \cap \Delta_a  =  \Sigma(\BB) \cap \Delta_a = \emptyset.
      \]
     \item[(ii)] $\AA \in \mathscr{B}(\EE)$ is a bounded operator on
       $\EE$. 
     \item[(iii)] There is $n \ge 1$ such that the operator $(\AA
       \RR_{\BB}(z))^n$ is bounded from $\EE$ to $E$ for any $z \in
       \Delta_a$.
     \end{itemize} 
   \end{itemize} 
   
   Then the same conclusions as in Theorem~\ref{theo:factor:enlarg}
   hold.
\end{theo}

\begin{rem}
  %\begin{enumerate}
  Actually there is no need in the proof that $(B-z)^{-1}$ for
    $z \in \Delta_a$ is a bounded operator, and therefore assumption
    {\bf (H2')} could be further relaxed to assuming only
    $(B-z)^{-1}(E) \subset \mbox{Dom}(L) \subset E$ (bijectivity is
    already known in $\EE$ from the invertibility of
    $(\BB-z)$). However these subtleties are not used at the level of
    the applications we have in mind.
\end{rem}

\begin{proof}[Proof of Theorem~\ref{theo:factor:enlarg:algeb}]
The Step 1 is unchanged, only the proofs of Steps 2 and 3 are
modified: 

\smallskip
\noindent
{\it Step 2. $(\LL-z)$ is invertible on $\Omega$.}  
Consider $z_0 \in \Omega$. First observe that if the
operator $(\LL - z_0)$ is bijective, then composing to the left the
equation 
\[
(\LL-z_0) \UU(z_0) = \mbox{Id}_{\EE}
\]
by $(\LL-z_0)^{-1}=\RR_\LL(z_0)$ yields $\RR_\LL(z_0) = \UU(z_0)$ and we
deduce that the inverse map is bounded (i.e. $(\LL-z_0)$ is an
invertible operator in $\EE$) together with the desired formula for the
resolvent. Since $(\LL-z_0)$ has a right-inverse it is surjective.

Let us prove that it is injective. Consider $f \in N(\LL-z_0) \subset
\EE$: 
\[
(\LL-z_0)f = 0 \ \mbox{ and thus } (\mbox{Id} + \GG(z_0)) (\BB-z_0) f
=0 \quad \mbox{with } \GG(z_0) := \AA \RR_\BB(z_0).
\]
We denote $\bar f := (\BB -z_0) f \in \EE$ and obtain 
\[
\bar f = - \GG(z_0) \bar f \ \Rightarrow \ \bar f = (-1)^n \,
\GG(z_0)^n \bar f
\]
and therefore, from assumption {\bf (H2')}, we deduce that $\bar f \in
E$. Finally $f = \RR_\BB(z_0) \bar f = \RR_B(z_0) \bar f \in \mbox{Dom}(L)
\subset E$. Since $(L-z_0)$ is injective we conclude that $f=0$. 

This completes the proof of this step and proves $\Sigma(\LL) \cap
\Delta_a \subset \{ \xi_1, \dots, \xi_k \}$ together with the point (ii)
of the conclusion.

\smallskip
\noindent
{\it Step 3. Spectrum, eigenspaces and spectral projectors.} 
On the one hand, 
\[
N (L-\xi_j)^\alpha \subset N(\LL-\xi_j)^\alpha, \quad j = 1, \dots, k,
\quad \alpha \in \N, 
\]
so that $\Sigma(\LL) \cap \Delta_a \supset \{ \xi_1, \dots, \xi_k \}$.
Since the other inclusion was proved in the previous step, we conclude that 
\[
\Sigma(L) \cap \Delta_a = \Sigma(\LL) \cap \Delta_a.
\]

On the other hand, let us consider an eigenvalue $\xi_j$, $j=1, \dots,
k$ for $\LL$, some integer $\alpha \ge 1$ and some $f \in
N(\LL-\xi_j)^\alpha$:
\[
\left( \LL - \xi_j \right)^\alpha (f) = 0.
\]
Using the decomposition of {\bf (H2)} and denoting 
$\bar f=(\BB-\xi_j)^\alpha f$ we deduce 
\[
\left( \mbox{Id} + \GG(\xi_j)\right)^\alpha \bar f  = 0 \ \mbox{ with } \
\GG(\xi_j) := \AA \RR_\BB(\xi_j). 
\]
By expanding this identity we obtain 
\[
\bar f = \GG(\xi_j) \OO_\alpha(\xi_j) (\bar f) 
\]
where $\OO_\alpha(\xi_j)$ is a finite sum of powers of $\GG(\xi_j)$
(with $\alpha$ terms and exponents between $0$ and $\alpha-1$). By
iterating this equality ($\GG(\xi_j)$ and $\OO_\alpha(\xi_j)$
commute), we get
\[
\bar f = \GG(\xi_j)^n \OO_\alpha(\xi_j)^n \bar f.
\]
This implies, arguing as in the previous step, that $\bar f \in
E$ and finally $f \in \mbox{Dom}(L) \subset E$. This proves that 
\[
N(\LL-\xi_j)^\alpha = N(L-\xi_j)^\alpha
\]
and since the eigenvalues are discrete, it straightforwardly completes
the proof of the conclusions (i) and (ii).
% Hence we conclude that $\Sigma(L) \cap \Delta_a = \{ \xi_1, \dots,
% \xi_k \}$ with 
% \[
% N (L-\xi_j)^\alpha = N(\LL-\xi_j)^\alpha, 
% \ \alpha \in \N \ \mbox{ and thus } \ M( \LL - \xi_j) =
% M(L-\xi_j), \quad i = 1, \dots, k.
% \]
Finally, the fact that $\Pi_{\LL,\xi_j| E} = \Pi_{L,\xi_j}$ is a
straightforward consequence of $\RR_\LL(z) (f) = \RR_L(z)(f)$ when $f
\in E$ and of the formula \eqref{def:SpectralProjection} for the
projector operator.
\end{proof}

\subsection{Hypodissipativity}
\label{sec:hypo-diss-coerc}

Let us first introduce the notion of  {\em hypodissipative}
operators and discuss its relation with the classical notions of {\em dissipative}
operators and  {\em coercive} operators as well as its relation with the recent terminology
of {\em hypocoercive} operators % used in the so-called ``hypocoercivity''
% theory 
(see mainly~\cite{MR2562709} and then
\cite{MNeu,MR2215889,DMS} for related references).
% We shall briefly
% recall statement of Hille-Yosida together with the statement of a
% variant of the
% Hille-Yosida theorem (e.g. \cite[Chap I, Theorem 3.1]{Pazy}) and the
% {\em Lumer-Phillips
%   theorem}. 
% (e.g. \cite[Chap I, Theorem 4.4]{Pazy}).

\begin{defin}[Hypodissipativity]\label{def:hypodissipative} 
  Consider a Banach space $(X,\| \cdot \|_X)$ and some operator
  $\Lambda \in \mathscr{C}(X)$. 
  We say that $(\Lambda-a)$ is \emph{hypodissipative} on $X$
   if there exists some norm $\Nt  \cdot \Nt_X$ on $X$   equivalent to the initial norm $\| \cdot \|_X$ such that
\begin{equation}\label{eq:def-dissipative} 
  \forall \, f \in  D(\Lambda) , \,\, \exists \, \varphi \in F(f) \ \mbox{
    s.t. } \  \Re e
  \, \langle \varphi, (\Lambda-a) \, f \rangle \le 0, 
\end{equation} 
where $\langle \cdot , \cdot \rangle$ is the duality bracket for the duality in $X$ and $X^*$
and $F(f)  \subset X^*$ is the dual set of $f$
defined by
$$
F(f) = F_{\Nt \cdot \Nt} (f) := \left\{\varphi \in X^*; \,\, \langle \varphi , f \rangle = \Nt
  f \Nt_X^2 = \Nt \varphi \Nt^2_{X^*} \right\}.
$$
\end{defin} 

\begin{rems}
\begin{enumerate} 
\item An hypodissipative operator $\Lambda$ such that $\Nt \cdot \Nt_X
  = \|\cdot \|_X$ in the above definition is nothing but a
  \emph{dissipative} operator, or in other words, $-\Lambda$ is an
  accretive operator.
\item When $\Nt \cdot \Nt_X$ is an Hilbert norm on $X$, we have $F(f)
  = \{ f \}$ and \eqref{eq:def-dissipative} writes
\begin{equation}\label{eq:def-dissipativeHilbert} 
\forall \, f \in D(\Lambda), \quad \Re e \, (\!( \Lambda f,f )\!)_X \le a \ \Nt f \Nt_X^2,
\end{equation}
where $(\!( \cdot , \cdot )\!)_X$ is the scalar product associated to
$\Nt \cdot \Nt_X$.  In this Hilbert setting such a hypodissipative
operator shall be called equivalently \emph{hypocoercive}.
\item When $\Nt \cdot \Nt_X = \| \cdot \|_X$ is an Hilbert norm on
  $X$, the above definition corresponds to the classical definition of
  a \emph{coercive} operator.
\item In other words, in a Banach space (resp. an Hilbert space) $X$,
  an operator $\Lambda \in \CCC(X)$ is hypodissipative
  (resp. hypocoercive) on $X$ if $\Lambda$ is dissipative
  (resp. coercive) on $X$ endowed with a norm (resp. an Hilbert norm)
  equivalent to the initial one. Therefore the notions of
  hypodissipativity and hypocoercivity are invariant under change of
  equivalent norm.
\end{enumerate}
\end{rems}

\smallskip The concept of \emph{hypodissipativity} seems to us
interesting since it clarifies the terminology and draws a bridge
between works in the PDE community, in the semigroup community and in
the spectral analysis community. For convenience such links are
summarized in the theorem below. This theorem is a non standard
formulation of the classical Hille-Yosida theorem on $m$-dissipative
operators and semigroups, and therefore we omit the proof.

% {\color{red}{ , as one can figure out thanks to the
% following result.}}

% \vspace{0.5cm}
% {\color{red}{MG: From the sentence above in red it seems that the following theorem is our result. I would rather say the following and then delete (1) in the Remark below: 

% " For convenience such links are summarized in the theorem below. The theorem is a non standard formulation of the classical Hille-Yosida theorem on $m$-dissipative operators and semigroups, and therefore we omit the proof." }}

\begin{theo}\label{theo:4equivDissipatif} Consider  $X$ a  Banach space 
  and $\Lambda$ the generator of a $C_0$-semigroup $S_\Lambda$. We
  denote by $\RR_\Lambda$ its resolvent.  For given constants $a \in
  \R$, $M > 0$ the following assertions are equivalent:
\begin{itemize}
\item[(i)] $\Lambda-a$ is hypodissipative; 

\item[(ii)] the semigroup satisfies the growth estimate
  \[
  \forall \, t \ge 0, \quad \|S_\Lambda (t) \|_{\BBB(X)} \le M \, e^{a \, t};
  \]

\item[(iii)] $\Sigma(\Lambda) \cap \Delta_a = \emptyset$ and
  \[
  \forall \, z \in \Delta_a, \quad \|\RR_\Lambda (z)^n \| \le
  \frac{M}{(\Re e \, z - a)^n};
\]

\item[(iv)] $\Sigma(\Lambda) \cap (a,\infty)= \emptyset$ and there
  exists some norm $\Nt \cdot \Nt$ on $X$ equivalent to the norm $\|
  \cdot \|$:
$$
\forall \, f \in X \qquad \| f \| \le \Nt   f \Nt \le M \, \|f \|,
$$
such that 
$$ 
\forall \, \lambda > a , \,\, \forall \, f \in D(\Lambda) ,\quad \Nt
(\Lambda - \lambda) \, f \Nt \ge (\lambda - a) \, \Nt f \Nt.
$$
\end{itemize}
\end{theo}

\begin{rems}\label{rem:dissipativeResolvante} 
\begin{enumerate}
% \item We omit the proof of Theorem~\ref{theo:4equivDissipatif} since
%   it is a non standard formulation of the classical Hille-Yosida
%   theorem on $m$-dissipative operators and semigroups.
\item We recall that $\Lambda-a$ is maximal if
$$
\mbox{R}(\Lambda - a) = X.
$$
This further condition leads to the notion of $m$-hypodissipative,
$m$-dissipative, $m$-hypocoercive, $m$-coercive operators.

\item The Hille-Yosida theorem is classically presented as the
  necessary and sufficient conditions for an operator to be the
  generator of a semigroup. Then one assumes, additionally to the
  above conditions, that $\Lambda - b$ is maximal for some given $b
  \in \R$. Here in our statement, the existence of the semigroup being
  assumed, the maximality condition is automatic, and
  Theorem~\ref{theo:4equivDissipatif} details how the operator's,
  resolvent's and the associated semigroup's estimates are linked.
\item In other words, the notion of hypodissipativity is just another
  formulation of the minimal assumption for estimating the growth of a
  semigroup. Its advantage is that it is arguably more natural from a
  PDE viewpoint.
\item The equivalence (i) $\Leftrightarrow$ (iv) is for instance a
  consequence of \cite[Chap 1, Theorem 4.2]{Pazy} and \cite[Chap 1,
  Theorem 5.3]{Pazy}. All the other implications are also proved in
  \cite[Chap 1]{Pazy}.
\end{enumerate}
\end{rems}

% \begin{defin}\label{def:hypocoercif} Consider  $X$ a Banach space, $\Lambda \in \mathscr{C}(X)$, $a \in \R$. 
%   We say that $\Lambda$ is $a$-hypocoercive if $\Lambda \in \GGG(X)$,
%   for some closed subspace $X_0\subset X$ of finite codimension we
%   have that $\Lambda_{|X_0} \in \mathscr{C}(X_0)$, $\Lambda_{|X_0}-a$
%   is hypodissipative on $X_0$ and there is $X'_0$ some (finite
%   dimensional) supplementary vector space $X'_0 \oplus X_0 = X$ which
%   is invariant under the action of $\Lambda$: $\Lambda(X'_0) \subset
%   X'_0$.
% \end{defin}

Let us now give a synthetic statement adapted to our purpose. We omit
the proof which is a straightforward consequence of the Lumer-Philipps
or Hille-Yosida theorems together with basic matrix linear algebra on
the finite-dimensional eigenspaces. The classical reference for this
topic is \cite{Kato}.

\begin{theo}\label{theo:sect3thHC} 
  Consider a Banach space $X$, a generator $\Lambda \in \CCC(X)$ of a
  $C_0$-semigroup $S_\Lambda$, $a \in \R$ and distinct $\xi_1, \dots,
  \xi_k \in \Delta_a$, $k \ge 1$.  The following assertions are
  equivalent:
\begin{itemize}
\item[(i)] There exist $g_1, \dots,g_m$ linearly independent
  vectors so that the subspace $\hbox{\rm Span}\{g_1, \dots, g_m\}$ is
  invariant under the action of $\Lambda$, and
$$
\forall \, i \in\{1 , \dots, m \}, \,\,\, 
\exists \, j \in \{1 , \dots, k \}, \quad g_i \in M(\Lambda-\xi_j).
$$
Moreover there exist $\varphi_1, \dots,\varphi_m$ linearly independent
vectors so that the subspace $\hbox{\rm Span}\{\varphi_1, \dots,
\varphi_m\}$ is invariant under the action of $\Lambda^*$. These two
families satisfy the orthogonality conditions $\langle \varphi_i,g_j
\rangle = \delta_{ij}$ and the operator $\Lambda - a$ is
hypodissipative on $\hbox{\rm Span}\{\varphi_1, \dots,
\varphi_m\}^\perp$:
$$
\forall \, f \in \bigcap_{n=1}^m \hbox{\rm Ker} (\varphi_i) \cap
D(\Lambda), \,\, \exists \, f^* \in F_{\Nt\cdot \Nt}(f) , \quad
\Re e \, \langle f^*, (\Lambda-a) f \rangle \le 0.
$$
\item[(ii)] There exists a decomposition $X = X_0 \oplus \dots \oplus
  X_k$ where (1) $X_0$ and $(X_1 + \dots + X_k)$ are invariant by the
  action of $\Lambda$, (2) for any $j=1, \dots, k$ $X_j$ is a
  finite-dimensional space included in $M(\Lambda - \xi_j)$, and (3)
  $\Lambda - a$ is hypodissipative on $X_0$:
$$
\forall \, f \in D(\Lambda) \cap X_0, \,\, \exists \, f^* \in
F_{\Nt\cdot \Nt}(f) , \quad \Re e \, \langle f^*, (\Lambda-a) f \rangle
\le 0.
$$
\item[(iii)] There exist some finite-dimensional projection operators
  $\Pi_1, \dots, \Pi_k$ which commute with $\Lambda$ and such that
  $\Pi_i \, \Pi_j = 0$ if $i \not = j$, and some operators $T_j =
  \xi_j \, \mbox{{\rm Id}}_{Y_j} + N_j$ with $Y_j := \mbox{{\rm
      R}}(\Pi_j )$, $N_j \in \mathscr{B}(Y_j)$ nilpotent, so that the
  following estimate holds
   \begin{equation}\label{bddSlambda} 
     \forall \, t \ge 0, \quad 
     \Bigl\| S_\Lambda(t) - \sum_{j=1}^k e^{t \, T_j}\, \Pi_j  \Bigr\|_{\mathscr{B}(X)} 
     \le C_a \, e^{ a \, t} , 
   \end{equation} 
   for some constant $C_a \ge 1$.   
 \item[(iv)] The spectrum of $\Lambda$ satisfies
\[
\Sigma(\Lambda) \cap \Delta_a = \{\xi_1, \dots , \xi_k  \} \in
\Sigma_d(\Lambda) \quad \mbox{(distinct discrete eigenvalues)}
\]
and $\Lambda-a$ is hypodissipative on $\mbox{{\rm
    R}}(I-\Pi_{\Lambda,a})$.
\end{itemize}

\smallskip Moreover, if one (and then all) of these assertions is
true, we have 
\[
\left\{
\begin{array}{l} 
X_0 = \mbox{{\rm R}}(I-\Pi_{\Lambda,a}), \vs \\ 
X_j = Y_j = M(\Lambda - \xi_j), \vs \\ 
\Pi_{\Lambda,\xi_j} = \Pi_j, \vs \\
T_j = \Lambda \Pi_{\Lambda,\xi_j}.
\end{array}
\right.
\]

As a consequence, we may write
$$
\RR_\Lambda(z) = \RR_0(z) + \RR_1(z),
$$
where $\RR_0$ is holomorphic and bounded on $\Delta_{a'}$ for any $a'
> a$ and
$$
\RR_1(z) = \sum_{j=1}^k \left( {\Pi_j \over z-\xi_j} +
\sum_{n=2}^{\beta_j} {N_j^n \over (z-\xi_j)^n} \, \Pi_j \right).
 $$
\end{theo}

\begin{rem} When $X$ is a Hilbert space and $\Lambda $ is a
  self-adjoint operator, the assumption (i) is satisfied with $k=1$,
  $\xi_1 = 0$, as soon as there exist $g_1, \dots, g_k \in X$
  normalized such that $g_i \perp g_j$ if $i \not= j$, $\Lambda g_i =
  0$ for all $i=1,\dots,k$, and
$$
\forall \, f \in X_0 := \hbox{\rm Span} \{g_1, \dots, g_k\}^\perp,
\quad \langle \Lambda f , f \rangle \le a \, \langle f,f \rangle. %, \quad a < 0.
$$
\end{rem}

\subsection{Factorization and quantitative spectral mapping theorems}
\label{sec:fact-quant-part}

The goal of this subsection is to establish quantitative decay
estimates on the semigroup in the larger space $\EE$. Let us recall
the key notions of \emph{spectral bound} of an operator $\LL$ on
$\EE$: 
\begin{equation*}
  s(\LL) := \sup\{ \Re e \, \xi \ : \ \xi \in \Sigma(\LL) \}
\end{equation*}
and of \emph{growth bound} of its associated semigroup
\begin{equation*}
  w(\LL) := \inf_{t >0} \frac1t \left\| S_\LL(t) \right\| = \lim_{t \to +\infty}
  \frac1t \left\| S_\LL(t) \right\|.
\end{equation*}
It is always true that $s(\LL) \le w(\LL)$ but we are interested in
proving the equality with quantitative estimates, in the larger space
$\EE$. Proving such a result is a particular case of a \emph{spectral
  mapping theorem}.

Let us first observe that in view of our previous factorization result
the natural control obtained straightforwardly on the resolvent in the
larger functional space $\EE$ is a \emph{uniform control on vertical
  lines}. It is a classical fact that this kind of control is not
sufficient in general for inverting the Laplace transform and
recovering spectral gap estimates on a semigroup from it.  

Indeed for semigroups in Banach spaces the equality between the
spectral bound and the growth bound is false in general when assuming
solely that the resolvent is uniformly bounded in any $\Delta_a$ with
$a>s(\LL)$ (with bound depending on $a$). A classical counterexample
\cite[Chap.~5, 1.26]{MR1721989} is the derivation operator $\LL f =
f'$ on the Banach space $C_0(\R_+) \cap L^1(\R_+, e^s \dd s)$ of
continuous functions that vanish at infinity and are integrable for
$e^s \dd s$ endowed with the norm
\begin{equation*}
  \| f \|= \sup_{s \ge 0} |f(s)| + \int_0 ^{+\infty} |f(s)| e^s \dd
  s. 
\end{equation*}
Another simple counterexample can be found in \cite{MR1300417}:
consider $1 \le p < q < \infty$ and the $C_0$-semigroup on
$L^p(1,\infty) \cap L^q (1,\infty)$ defined by 
\[
(T(t)f)(s) = e^{t/q} f(se^t), \quad t > 0, \ s >1. 
\]

However for semigroups in Hilbert spaces, the
Gerhart-Herbst-Pr\"uss-Greiner theorem
\cite{MR0461206,MR715559,MR743749,MR839450} (see also
\cite{MR1721989}) asserts that the expected semigroup decay $w(\LL) =
s(\LL)$ is in fact true, under this sole pointwise control on the
resolvent.  While the constants seem to be non-constructive in the
first versions of this theorem, Engel and Nagel gave a comprehensive
and elementary proof with constructive constant in \cite[Theorem 1.10;
chapter V]{MR1721989}. Let us also mention on the same subject
subsequent works like Yao \cite{Yao1995} and Blake \cite{MR1823064},
and more recently \cite{Helffer-Sjoestrand}.

The main idea in the proof of \cite[Theorem 1.10, chapter
V]{MR1721989}, which is also used in \cite{Helffer-Sjoestrand}, is to
use a Plancherel identity on the resolvent in Hilbert spaces in order
to obtain explicit rates of decay on the semigroup in terms of bounds
on the resolvent. We will present in a remark how this interesting
argument can be used in our case, but instead our proof will use a
more robust argument valid in \emph{Banach} spaces, which is made
possible by the additional factorization structure we have. The key
idea is to translate the factorization structure at the level of the
semigroups.

% Now we present our main \emph{enlargement}  of the functional space of
% a quantitative spectral mapping theorem (in the sense of a semigroup
% decay estimate as above). 
% The focus is on the combination of the above mentioned method
% with our factorization theorem in order to modify the space for the
% decay estimate on the semigroup.

\smallskip
 
We shall need the following definition on the convolution of semigroup
(corresponding to \emph{composition} at the level of the resolvent
operators).

\begin{defin}[Convolution of semigroups] 
  Consider some Banach spaces $X_1$, $X_2$, $X_3$. For two
  one-parameter families of operators
  \[
  S_1 \in L^1(\R_+; \BBB(X_1,X_2)) \ \mbox{ and } \ S_2 \in L^1(\R_+;
  \BBB(X_2,X_3)),
  \]
  we define the convolution $S_2 \ast S_1 \in L^1(\R_+; \BBB(X_1,X_3))$ by
  $$
  \forall \, t \ge 0, \quad (S_2 * S_1)(t) := \int_0^t S_2(s) \, S_1 (t-s) \dd s. 
  $$
  When $S_1=S_2$ and $X_1=X_2=X_3$, we define recursively $S^{(*0)} =
  \Id$ and $S^{(*\ell)} = S * S^{(*(\ell-1))} $ for any $\ell \ge 1$.
\end{defin}

\begin{rems}
  \begin{enumerate}
  \item Note that this product law is in general not commutative.
%  \item If $S_i \in C^{\alpha_i}$ with $\alpha_i \in \N$ we deduce
%    $S_2 * S_1 \in C^{\alpha_1 + \alpha_2 + 1}$. \noteC{Is it this
%      remark true? I think there are convolutions of two continuous
%      functions which are nowhere differentiable (see
%      e.g. Bogdanowicz, Proc. AMS, 1965)?}
  \item A simple calculation shows that that if $S_i$ satisfies
    \[
    \fa t \ge 0, \quad \|S_i(t) \|_{\BBB(X_i,X_{i+1})} \le C_i \, t^{\alpha_i} \, e^{a_i
      \, t}
    \]
    for some $a_i \in \R$, $\alpha_i \in \N$, $C_i \in (0,\infty)$,
    then 
    $$
    \forall \, t \ge 0, \quad \|S_1 * S_2 (t) \|_{\BBB(X_1,X_2)} \le
    C_1 \, C_2 \, {\alpha_1! \, \alpha_2! \over (\alpha_1 +
      \alpha_2)!} \, t^{\alpha_1+\alpha_2+1} \, e^{\max(a_1,a_2) \,
      t}.
$$ 
  \end{enumerate}
\end{rems}

\begin{theo}[Enlargement of the functional space of the semigroup
  decay]\label{theo:EnlargingSGdecay}
  Let $E, \EE$ be two Banach spaces with $E \subset \EE$ dense with
  continuous embedding, and consider $L\in~\CCC(E)$, $\LL \in
  \mathscr{C}(\EE)$ with $\LL_{|E} = L$ and $a \in \R$.

We assume the following:
\begin{itemize}
\item[{\bf (A1)}] %the small space $E$ is a Hilbert space,
  $L$ generates a semigroup $e^{tL}$ on $E$,  $L-a$ is
  hypodissipative on $\mbox{{\rm R}}(\mbox{{\rm Id}}-\Pi_{L,a})$ and
  \[
  \Sigma(L) \cap \Delta_{a} := \{\xi_1, \dots, \xi_k \} \subset
  \Sigma_d(L) \quad \mbox{(distinct discrete eigenvalues)}
  \]
  (with $\{\xi_1, \dots, \xi_k\} = \emptyset$ if $k=0$).
  \smallskip

\item[{\bf (A2)}] There exist $\AA, \BB \in \mathscr{C}(\EE)$
  such that $\LL = \AA + \BB$ (with corresponding restrictions $A,B$ on
  $E$), some $n \ge 1$ and some constant $C_a > 0$ so
  that% , $\alpha \in [0,1)$ {\color{green}{(What is
      % $\alpha$?)}}, such that
  \begin{itemize}
  \item[{\bf (i)}]  $(\BB-a)$ is hypodissipative on $\EE$;
  \item[{\bf (ii)}]  $\AA \in \BBB( \EE)$ and $A \in \BBB(E)$;
  \item[{\bf (iii)}] $T_n := (\AA \, S_\BB)^{(*n)}$ satisfies $\|T_n
    (t) \|_{\BBB(\EE,E)} \le C_a \, e^{a \, t} $.
   \end{itemize}
 \end{itemize}

 Then $\LL$ is hypodissipative in $\EE$ with 
 \begin{equation}\label{eq:estimSGfort} 
   \forall \, t \ge 0, \quad \left\|  S_\LL(t) - \sum_{j=1}^k S_L(t) \, \Pi_{\LL,\xi_j} \right\|_{\BBB(\EE)}
   \le C' _{a} \, t^{n} \, e^{a \, t} % \quad \ \forall
   % \, a' > a
 \end{equation}
 for some explicit constant $C' _{a}>0$ depending on the constants in
 the assumptions. Moreover we have the following \emph{factorization
   formula} at the level of semigroups on $\EE$:
 \begin{multline}\label{eq:factor-semigroup}
   S_\LL(t) = \sum_{j=1}^k S_L(t) \, \Pi_{\LL,\xi_j} + \sum_{\ell=0}
   ^{n-1} (-1)^\ell \, (\mbox{{\em Id}} - \Pi_{\LL,a}) \, S_\BB \ast \left( \AA S_\BB \right)^{\ast \ell}(t) \\
   + (-1)^n \left[ (\mbox{{\em Id}} - \Pi_{L,a}) S_L\right] \ast
   \left( \AA S_\BB \right)^{\ast n}(t).
 \end{multline}
\end{theo}

\begin{rems}
\begin{enumerate}
\item It is part of the result that $\BB$ generates a~semigroup on
  $\EE$ so that {\bf (A2)-(iii)} makes sense. Except for the
  assumption that $L$ generates a semigroup, all the other assumptions
  are pure functional, either on the discrete eigenvalues of $L$ or on
  $L$, $\BB$, $\AA$, $A$ and $T_n$, and do not require maximality
  conditions.

% \item Observe that we require in this theorem that $E$, i.e. the small
%   target space, be a Hilbert space, which was not required in the
%   previous enlargement theorem dealing with spectral analysis.

\item Assumption {\bf (A1)} could be alternatively formulated
  by mean of any of the equivalent assertions listed in
  Theorem~\ref{theo:sect3thHC}.
\end{enumerate}
\end{rems}

\begin{proof}[Proof of Theorem~\ref{theo:EnlargingSGdecay}.] 
  We split the proof into four steps. 
  
\smallskip
\noindent {\sl Step 1. } First remark that since $B=L-A$, $A \in
\BBB(E)$, and $L$ is $m$-hypodissipative then $B$ is
$m$-hypodissipative and generates a strongly continuous semigroup
$S_B$ on $E$.

Because of the hypodissipativity of $\BB$, we can extend this
semigroup from $E$ to $\EE$ and we obtain that $\BB$ generates a
semigroup $S_\BB$ on $\EE$. To see this, we may argue as follows.  We
denote by $\Nt \cdot \Nt_E$ a norm equivalent to $\| \cdot \|_E$ so
that $B-b$ is dissipative in $(E,\Nt \cdot \Nt_E)$ and $\Nt \cdot
\Nt_\EE$ a norm equivalent to $\| \cdot \|_\EE$ so that $\BB-b$ is
dissipative in $(\EE,\Nt \cdot \Nt_\EE)$, for some $b \in \R$ large
enough.  We introduce the new norm 
\[
\Nt f \Nt_\epsilon := \Nt f \Nt_\EE + \epsilon \, \Nt f \Nt_E \quad
\mbox{ on } \quad E
\]
so that $\Nt \cdot \Nt_\epsilon$ is equivalent to $\Nt \cdot \Nt_E$
for any $\epsilon >0$.  Since $B-b$ is $m$-dissipative in $(E,\Nt
\cdot \Nt_\epsilon)$, the Lumer-Phillips theorem shows that the
operator $B-b$ generates a semigroups of contractions on $(E,\Nt \cdot
\Nt_\epsilon)$, and in particular %for any $f \in E$, $t \ge 0$,
$$
\fa f \in E, \ \fa t \ge 0, \quad \Nt S_{(B-b)}(t) f \Nt_\EE +
\epsilon \, \Nt S_{(B-b)}(t) f \Nt_E \le \Nt f \Nt_\EE + \epsilon \,
\Nt f \Nt_E.
$$
Letting $\epsilon$ going to zero, we obtain
$$
\forall \, f \in E, \ \forall \, t \ge 0, \quad \Nt S_{B}(t) f
\Nt_\EE \le e^{t \, b} \, \Nt f \Nt_\EE.
$$
Because of the continuous and dense embedding $E \subset \EE$ we
deduce that we may extend $S_B(t)$ from $E$ to $\EE$ as a family of
operators $S(t)$ which satisfies the same estimate. We easily conclude
that $S(t)$ is a semigroup with generator $\BB$, or in
other words, $\BB$ generates a semigroup $S_\BB = S$ on $\EE$.

Finally, since $\LL = \AA + \BB$ and $\AA \in \BBB(\EE)$, we deduce
that $\LL$ generates a semigroup.

\smallskip
\noindent
{\sl Step  2. } We have from {\bf (A2)-(i)} that \begin{equation}\label{bdd:T0}
\forall \, t \ge 0, \quad \|S_\BB(t)\|_{\EE \to \EE} \le C \, e^{at}
\end{equation} and we easily deduce (by iteration) that $T_\ell := (\AA \,
S_\BB)^{(*\ell)}$, $\ell \ge 1$, satisfies 
\begin{equation}\label{bdd:T1} 
  \forall \, t \ge 0, \,\, \forall \, \ell \ge 1, \quad \| T_\ell (t)
  \|_{\BBB(\EE)} \le C_\ell \, t^{\ell -1}\, e^{at} 
\end{equation}
for some constants $C_\ell>0$, $\ell \ge 1$. 

Let us define
$$
\UU_\ell := (-1)^\ell \, \left(\Id_\EE - \Pi_{\LL,a}\right) \, S_{\BB} \ast (\AA \,
S_\BB )^{(\ast \ell)}, \quad 0 \le \ell \le n-1.
$$
From \eqref{bdd:T0} and \eqref{bdd:T1} and the boundedness of
$\Pi_{\LL,a}$ we get
\begin{equation}\label{bdd:T2} 
  \forall \, t \ge 0,   \quad \| \UU_\ell(t) \|_{\BBB(\EE)} \le C_\ell \, t^{\ell}\, e^{at}, \quad
  \,\, 0 \le \ell \le n-1.
\end{equation}
By applying standard results on Laplace transform, we have for any $f
\in \EE$
$$
\forall \, z \in \Delta_a, \quad \int_0^{+\infty} e^{zt}\, \UU_\ell(t)
f \dd t = (-1)^{\ell} \, \left(\Id_\EE-\Pi_{\LL,a}\right) \, \RR_\BB(z) \, (\AA \,
\RR_\BB(z))^\ell f.
$$
Then the inverse Laplace theorem implies, for $\ell=0, \dots, n-1$,
that
\begin{multline}\label{eq:RiemannUell} 
  \forall \, a' > a, \quad \UU_\ell (t) f \\ = {(-1)^\ell \over 2i\pi}
  \, \left(\Id_\EE - \Pi_{\LL,a}\right) \, \int_{a'-i\infty}^{a'+i\infty} e^{zt} \,
  \RR_\BB(z) \, (\AA \, \RR_\BB(z))^\ell \, f \dd z
  \\
  := \lim_{M \to \infty}{(-1)^n \over 2i\pi} \, \left(\Id_\EE - \Pi_{\LL,a}\right) \,
 \int_{a'-iM}^{a'+iM}
  e^{zt} \, \RR_\BB(z) \, (\AA \, \RR_\BB(z))^\ell \, f \dd z,
\end{multline}
where the integral along the complex line $\{a' + i y, \,\, y \in \R
\}$ may not be absolutely convergent, but is defined as the above
limit.

Let us now consider the case $\ell =n$ and define 
\begin{multline*}
  \UU_n(t) = (-1)^n \, \left(\Id_\EE - \Pi_{\LL,a}\right) \, \left[ S_L \ast (\AA \,
  S_\BB )^{(\ast n)} \right]\\ = (-1)^n \left[\left(\Id_E -
    \Pi_{L,a}\right) S_L\right] 
\ast (\AA \, S_\BB )^{(\ast n)}. 
\end{multline*}
Observe that this one-parameter family of operators is well-defined
and bounded on $\EE$ since $(\AA \, S_\BB )^{(\ast n)}$ is bounded
from $\EE$ to $E$ by the assumption {\bf (A2)-(iii)}. Moreover for $f
\in \EE$ we have
\begin{equation*}
  \left\| (\AA \, S_\BB )^{(\ast n)} (t) f\right\|_E \le C_a \,
  t^{n-1} \, e^{at} \, \| f \|_\EE
\end{equation*}
and since from {\bf (A1)}
\begin{equation*}
  \left\| \left[\left(\Id_E - \Pi_{L,a}\right) S_L\right](t)  g \right\| \le C'_a
  \, e^{at} \, \| g \|_E
\end{equation*}
for $g \in E$ we deduce by convolution that 
\begin{equation}\label{estim:Un}
  \left\|\UU^n (t) f\right\|_E \le C'' _a \, t^{n} \, e^{at} \, \| f \|_\EE
\end{equation}
(for some constants $C_a, C'_a, C''_a >0$). Observe finally that 
\begin{equation*}
\forall \, z \in \Delta_a, \quad  %\mbox{ with } \Re e \, z < \min\{\Re e \,
%\xi_1, \dots, \Re e \, \xi_k \}, \\ 
\int_0^{+\infty} e^{zt}\, (\Id_E - \Pi_{L,a}) S_L(t)
\dd t = (\Id_E - \Pi_{L,a}) \RR_L(z)
\end{equation*}
by classical results of spectral decomposition. 

Therefore the inverse Laplace theorem implies that for
    any $a'>a$ close enough to $a$ (so that $a' < \min\{\Re e \,
    \xi_1, \dots, \Re e \, \xi_k \}$) it holds
\begin{equation*}
 \UU_n(t) f := \lim_{M \to \infty} \UU_{n,M}(t) f
\end{equation*}
with 
\begin{equation*}
  \UU_{n,M}(t) f := {(-1)^n \over
    2i\pi} \, \left(\Id_E - \Pi_{L,a}\right) \, \int_{a'-iM}^{a'+iM} e^{zt} \, \RR_L(z) \, (\AA \,
  \RR_\BB(z))^n f \dd z. 
\end{equation*}
  
 \smallskip
\noindent
{\sl Step 3. } Let us prove that the following representation formula
holds
\begin{equation}\label{eq:representation} 
  \forall \, f \in \EE, \,\,
  \forall \, t \ge 0, \quad S_\LL(t) f = \sum_{j=1}^k S_{\LL,\xi_j}
  (t) \, f + \sum_{\ell = 0}^n \UU_\ell (t) \, f,
\end{equation}
where $S_{\LL,\xi_j} (t) = S_\LL(t)  \Pi_{\LL,\xi_j}$ and
$\Pi_{\LL,\xi_j}$ is the spectral projection as defined in
\eqref{def:SpectralProjection}. 
% We compute 
% \begin{equation*}
%  \UU_n(t) :=
% \lim_{M \to \infty} \UU_{n,M}(t)
% \end{equation*}
% with 
% \begin{equation*}
% \UU_{n,M}(t) := {(-1)^n \over
%   2i\pi}\int_{a-iM}^{a+iM} e^{zt} \, \RR_L(z) \, (\AA \, \RR_\BB(z))^n
% \dd z. 
% \end{equation*}

%We argue as in the step 2. 
Consider $f \in D(\LL)$ and define $f_t = S_\LL(t) f$. From {\bf (A2)}
there exists $b \in \R$ and $C_b \in (0,\infty)$ so that
\begin{equation}\label{eq:RiemannL1} 
t \mapsto f_t \in C^1(\R_+;\EE) \qquad\hbox{and}\qquad 
\| f_t \|_\EE \le C_b \, e^{b \, t}\, \| f \|_\EE
\end{equation}
and therefore the inverse Laplace theorem implies
for $b' > b$ 
\begin{equation}\label{eq:RiemannL2} 
  \forall \, z \in \Delta_{b'}, \quad
  r(z) := \int_0^{+\infty} f_t \, e^{-z \, t} \dd t = - \RR_\LL (z) \, f
\end{equation} 
is well-defined as an element of $\EE$, and
\begin{equation}\label{eq:RiemannL3} 
\forall \, t \ge 0, \quad f_t
  = {1 \over 2 i \pi} \int_{b'-i\infty}^{b'+i\infty} e^{z t} \, r(z)
  \dd z := \lim_{M\to\infty} {1 \over 2 i \pi} \int_{b'-iM}^{b'+iM}
  e^{z t} \, r(z) \dd z . 
\end{equation} 
% It is worth emphasizing that in \eqref{eq:RiemannL3} the right hand
% side integral may be only \Red ``a semi-convergent" integral (it is
% not absolutely convergent) \Black Moreover, inequality (b) in
% \eqref{eq:RiemannL1} classically implies (see e.g. \cite[section
% 1.7]{Pazy}) that for $z \in \Delta_{b'}$ \bear\nonumber r(z) &=&
% \int_0^\infty e^{t \, \LL} \, e^{-z \, t}\dd t f_0
% \\ \label{eq:RiemannL4} &=& (z-\LL)^{-1} \, f_0 = - R_\LL (z) \, f_0.
% \eear 

Combining the definition of $f_t$ together with \eqref{eq:RiemannL3}
and \eqref{eq:RiemannL2} we get
\begin{equation}\label{eq:RiemannL5} 
  - S_\LL(t) f = \lim_{M \to \infty}
  I_{b',M} 
\end{equation} 
where
$$
\forall \, c \in \R \setminus \Re e (\Sigma(\LL)), \quad I_{c,M} := {1 \over 2 i
  \pi} \int_{c-iM}^{c+iM} e^{z t} \, \RR_\LL (z) \, f \dd z .
$$

Now from {\bf (A2)-(iii)} we have that $(\AA \,
    \RR_\BB(z))^n$ defined as
$$
(\AA \, \RR_\BB(z))^n = \int_0^\infty e^{z \, t}\, T_n (t) \dd t 
$$
is holomorphic on $\Delta_a$ with values in $\BBB(\EE,E)$. Hence the
assumptions {\bf (H1)-(H2)} of Theorem~\ref{theo:factor:enlarg} are
satisfied. We deduce that
\[
\Sigma(\LL) \cap \Delta_a = \Sigma(L) \cap \Delta_a
\]
with the same eigenspaces for the discrete eigenvalues $\xi_1, \dots,
\xi_k$. 

Moreover, thanks to {\bf (A1)} and {\bf (A2)-(i)} we have
$$
\forall \, a' > a, \,\, \forall \, \eps > 0, \quad 
\left\{ 
\begin{array}{l} \ds
  \sup_{z \in K_{a',\eps}} \|\RR_L(z) \|_{\BBB(E)} \le C_{a',\eps}, \vs
  \\ \ds 
  \sup_{z \in \Delta_{a'}} \|\RR_\BB(z) \|_{\BBB(\EE)} \le C_{a'},
\end{array}
\right.
$$
with 
\[
K_{a',\eps} := \Delta_{a'} \setminus  \big( B(\xi_1,\eps) \cup \, \dots \,
\cup B(\xi_k,\eps)\big).
\] 

As a consequence of the factorization formula \eqref{eq:factor-enlarg}
we get
$$
\forall \, a' > a, \,\, \forall \, \eps > 0, \quad \sup_{z \in
  K_{a',\eps}} \left\|\RR_\LL(z) \right\|_{\BBB(\EE)} \le C_{a',\eps}.
$$
Thanks to the identity 
\[
\forall \, z \notin \Sigma(\LL), \quad \RR_\LL(z) = z^{-1} \, [- \Id + \RR_\LL(z) \,
\LL] 
\]
and the above bound, we have (remember that $f \in D(\LL)$)
\begin{equation}\label{estim:Mtoinfty}
  \sup_{z; \, |\Im m \, z|\ge M, \, \Re e \, z \ge a'} \| \RR_\LL(z) \,
  f \|_{\BBB(\EE)} \mathop{\longrightarrow}_{M \to \infty} 0.
\end{equation}

We then choose $a' > a$ close enough to $a$ and $\eps > 0$ small
enough so that 
\[
B(\xi_1,\eps) \cup \, \dots \, \cup B(\xi_k,\eps) \subset \Delta_{a'}.
\]
Since $\RR_\LL$ is a meromorphic function on $\Delta_a$ with poles
$\xi_1$, \dots, $\xi_k$, we compute by the Cauchy theorem on path
integral
\begin{equation}
I_{b',M} = I_{a',M} + \sum_{j=1}^k S_{\LL,\xi_j} \, f+ \eps_1(M) 
\end{equation}
with
$$
\eps_1(M) = \left[ {1 \over 2 i \pi} \int_{a}^{b'} \, e^{(x+iy) \, t}
  \, \RR_\LL(x + i y) \, f \dd x \right]_{y=-M}^{y =M} \longrightarrow 0
$$
as $M \to 0$ thanks to \eqref{estim:Mtoinfty}.

On the other hand, because of Theorem~\ref{theo:factor:enlarg}, we may
decompose
\begin{multline}\label{eq:RiemannL7} 
  I_{a',M} = {1 \over 2 i \pi}
  \int_{a'-iM}^{a'+iM} e^{z t} \, \sum_{\ell=0}^{n-1} (-1)^\ell \RR_\BB(z)
  \, (\AA \, \RR_\BB (z))^\ell \, f \dd z \\ + { (-1)^n
    \over 2 i \pi} \int_{a'-iM}^{a'+iM} e^{z t} \, \RR_L(z) \, (\AA \, \RR_\BB (z))^n\, f \dd z.
\end{multline}
Note that the limit in (\ref{eq:RiemannL7}) as $M$ goes to infinity is
well defined. Hence \eqref{eq:RiemannL5} and \eqref{eq:RiemannL7}
yields to
\begin{multline*}
  S_\LL(t) f = \sum_{j=1}^k S_{\LL,\xi_j}
  (t) \, f + {1 \over 2 i \pi}
  \int_{a'-i\infty}^{a'+i\infty} e^{z t} \, \sum_{\ell=0}^{n-1} (-1)^\ell \RR_\BB(z)
  \, (\AA \, \RR_\BB (z))^\ell \, f \dd z \\ + { (-1)^n
    \over 2 i \pi} \int_{a'-i\infty}^{a'+i\infty} e^{z t} \, \RR_L(z) \, (\AA \,
  \RR_\BB (z))^n\, f \dd z. 
\end{multline*}
But since 
\begin{equation*}
  \sum_{j=1}^k S_{\LL,\xi_j}(t) =  \Pi_{\LL,a}  S_{\LL} (t)
\end{equation*}
we deduce that the sum of the two last terms in the previous equation
belongs to $R(\Id_\EE-\Pi_{\LL,a})$, and finally we have
\begin{multline*}
  S_\LL(t) f = \sum_{j=1}^k S_{\LL,\xi_j} (t) \, f \\ + {1 \over 2 i
    \pi} \int_{a'-i\infty}^{a'+i\infty} e^{z t} \, \sum_{\ell=0}^{n-1}
  (-1)^\ell \, \left(\Id_\EE - \Pi_{\LL,a}\right) \, \RR_\BB(z) \, (\AA \, \RR_\BB
  (z))^\ell \, f \dd z \\ + { (-1)^n \over 2 i \pi}
  \int_{a'-i\infty}^{a'+i\infty} e^{z t} \, \left(\Id_E - \Pi_{L,a}\right) \,
  \RR_L(z) \, (\AA \, \RR_\BB (z))^n\, f \dd z.
\end{multline*}

As a consequence, we deduce,  thanks to the step 2, that 
\begin{equation*}
  \forall \, f \in D(\LL), \,\,
  \forall \, t \ge 0, \quad S_\LL(t) f = \sum_{j=1}^k S_{\LL,\xi_j}
  (t) \, f + \sum_{\ell = 0}^n \UU_\ell (t) \, f.
\end{equation*}
% \begin{equation}\label{eq:decompIaM}
%   \lim_{M \to\infty} I_{a,M} =  \sum_{\ell=0}^{n-1} \UU_\ell(t) f+ \UU_n(t) f.
% \end{equation}
% and
% \eqref{eq:decompIaM} and
Then using the density of $D(\LL) \subset \EE$ we obtain the
representation formula \eqref{eq:representation}. We have thus
established \eqref{eq:factor-semigroup}.

%Alternatively, one could also
%consider $f \in D(\LL) \cap D(\LL^2)$ in a smaller subspace where the
%integral is absolutely convergent and then argue by density. 

  \smallskip
  \noindent {\sl Step 4. Conclusion.} We finally obtain the time decay
  \eqref{eq:estimSGfort} by plugging the decay estimates
  \eqref{bdd:T2} and \eqref{estim:Un} into the representation formula
  \eqref{eq:representation}.
\end{proof}

\begin{rem}
Let us explain how, in the case where $E$ is a Hilbert space, the
decay estimate on $\UU_n(t)$ can be obtained by reasoning
\emph{purely} at the level of resolvents, thanks to the Plancherel
theorem. Let us emphasize that the following argument does not require
a Hilbert space structure on the large space $\EE$. 

Consider $f \in \mbox{Dom}(\LL) \subset \EE$ and $\phi \in
\mbox{Dom}(L^*) \subset E^* = E$ ($E^*$ denotes the dual space of $E$
and $L^*$ the adjoint operator of $L$). Let us estimate
\begin{eqnarray}\nonumber 
  \left \langle \phi, \UU_n(t)  f \right \rangle
  &:=& - \frac{1}{2 i \pi} \lim_{M\to\infty} \int_{a'-iM}^{a'+iM} {e^{z \,
      t} \over t} \, \left\langle \RR_{L^*}(z) \phi , ( \AA \, \RR_\BB(z))^n  f \right \rangle \dd z.
\end{eqnarray} 
Applying the Cauchy-Schwarz inequality, we get 
\begin{eqnarray*} 
&&\left|   \left \langle \phi, \UU_n(t)  f \right \rangle \right| \le \frac{e^{a't}}{2\pi} \,
  \int_{-\infty} ^\infty \| \RR_{L^*} (a' + iy) \phi \|_{E} \, \|  (\AA \, \RR_\BB(a' + iy))^n f \|_E \dd y \\
   &&\le  \frac{e^{a't}}{2\pi} \,
  \left( \int_\R  \| \RR_{L^*} (a'+iy) \phi \|^2 _E \dd y \right)^{1/2}  
   \left( \int_\R \|  (\AA \, \RR_\BB(a' + iy))^n f  \|^2 _E \dd s \right)^{1/2}.
\end{eqnarray*}

For the first term, we make use of (1) the identity 
$$
\RR_{L^*} (a'+iy) = \left(\mbox{Id}_{E^*} + (a'-b) \, \RR_{L^*}(a'+i
  y)\right) \, 
\RR_{L^*}(b+ i y ),
$$
and (2) the fact that
\[
\| \RR_{L^*}(a'+i y) \|_{\BBB(E)} = \| \RR_{L}(a'+i y) \|_{\BBB(E)} 
\]
is uniformly bounded for $y \in \R$, and (3) the Plancherel theorem in
the Hilbert space $E$ and (4) the fact that 
\[
\| S_{L^*}(t) \|_{\BBB(E^*)} = \| S_{L}(t) \|_{\BBB(E)} \le C_b \, e^{ b
  \, t} \ \mbox{ for } b > \max (\xi_j),
\]
and we get
\begin{eqnarray*}
  \int_{\R} \| \RR_{L^*} (a'+iy) \phi \|_{E}^2 \dd s &\le &
 C_1  \int_{\R} \| \RR_{L^*} (b+iy) \phi \|_{E}^2 \dd s \\
   &\le & 2 \pi \, C_1 \, \int_0 ^{+\infty} 
  \big\| e^{- b \, t} \, e^{t \, L^*} \phi \big\|_{E}^2 \dd t \\
  &\le&2 \pi \, C_1  \left( \int_0 ^{+\infty} \big\| e^{-b \, t} \, e^{t \, L^*}
    \big\|_{\BBB(E)}^2 \dd t \right) \, \| \phi \|_{E}^2 \\
  &\le& C_2 \, \| \phi \|^2 _{E}. 
\end{eqnarray*} 

As for the second term, we identify the Laplace transform of
$(\AA \, \RR_\BB(z))^nf$ on $\Delta_a$ as $T_n(t)f$ arguing as before,
and % efine
% \[
% \varphi : \R_+ \to E, \quad \varphi(t) := T_n(t) \, f
% \]
% so that
% \begin{equation}\label{bdd:varphi}
% \|\varphi(t) \|_E \le C \, e^{at} \, \|f \|_\EE. 
% \end{equation}
% We have $\varphi \in C^1$ since $f \in \mbox{Dom}(\LL)$ and its
% Laplace transform function $r$ satisfies
% $$
% \left\{ 
% \begin{array}{l} \ds
% \forall \, z \in \Delta_a, \quad r(z) = \int_0^\infty e^{- z \, t} \, \varphi(t) \dd t, 
% \vs \\ \ds
% \forall \, t \ge 0, \quad \varphi (t) = {1 \over 2i\pi}
% \int_{a-i\infty}^{a+i\infty} e^{zt} \, r(z) \dd z
% \end{array}
% \right.
% $$
% and a straightforward computation yields
% $$
% r(z) =  (\AA \, \RR_\BB(z))^n f.
% $$
% The
the Plancherel theorem in $E$ gives 
\begin{eqnarray*}
  \int_{\R} \|  (\AA \, \RR_\BB)^n(a'+iy)  \, f  \|_E^2 \dd y 
  &=&
  \int_{\R}  \|   r(a'+iy)  \|_E^2 \dd y
  \\
  &=&
  2\pi \int_0^\infty  \| \varphi (t) \, e^{-a't}   \|_E^2 \dd t
  \\
  &\le& 2\pi  \, \| f \|_\EE^2 \, \int_0 ^{+\infty}  C_a^2\, e^{2 (a-a') \, t}  \dd t
  \\
  &\le&C_3 \, \|f \|_\EE^2.
\end{eqnarray*} 

Putting together these three estimates, we obtain
\[
\qquad
 \forall \, \phi \in D(L^*), \,\, \forall \, f \in D(\LL), \quad
\left|  \left \langle \phi, \UU_n(t)  \, f \right \rangle \right|
\le e^{at} \, (C_1 \, C_2)^{1/2} \,  \|f \|_\EE \, \|\phi \|_{E^*},
\]
so that, from the fact that $D(\LL)$ is dense in $\EE$,
\begin{equation}
 \,\, \forall \, f_0 \in \EE, \quad 
\| \UU_n(t)  \,  f \|_{\EE} \le C_{E\subset\EE} \,  \| \UU_n(t)  \, f  \|_{E} \le 
C_4 \,  e^{a' \, t}  \,  \|f \|_\EE
\end{equation}
where $C_{E \subset \EE}$ is the bound of the continuous embedding
from $E$ to $\EE$. 
\end{rem}

\begin{rem}
  There is another way to interpret the factorization formula at the
  level of semigroups. Consider the evolution equation $\partial_t f =
  \LL f$ and introduce the spliting 
\[
f = \sum_{i=1} ^k S_{\LL,\xi_i} \init{f} + f^1 + \dots + f^{n+2},
\]
 with
\begin{equation*}
\left\{ 
\begin{array}{lcl}\ds
\partial_t f^1 &=& \BB f^1  , \quad \init{f^1} = \left( \Id -\Pi_{L,a}\right) \init{f}, 
\vs \\ \ds 
\partial_t f^\ell &=&  \BB f^\ell  + \AA f^{\ell-1}, \quad \init{f^\ell} = 0, \quad 2 \le \ell \le n, 
\vs \\ \ds 
\partial_t f^{n+1} &=&  \LL f^{n+1}  + (\Id-\Pi_{\LL,a}) \AA f^{n}, \quad \init{f^{n+1}} = 0,
\vs \\ \ds
\partial_t f^{n+2} &=& \LL f^{n+2} + \Pi_{\LL,a} \AA f^{n} , \quad \init{f^{n+2}} =0.
\end{array}
\right.
\end{equation*}
This system of equations on $(f^\ell)_{1 \le \ell \le n+2}$ is
compatible with the equation satisfied by $f$, and it is possible to
estimate the decay in time inductively for $f^\ell$ (for the last
equation one uses $f^{n+2} = \Pi_{\LL,a} f^{n+2} = -
\Pi_{\LL,a}(f^1+\dots+f^{n+1})$ and the decay of the previous terms).

We made the choice to present the factorization theory from the
viewpoint of the resolvents as it reveals the algebraic structure in a
much clearer way, and also is more convenient for obtaining properties
of the spectrum and precise controls on the resolvent in the large
space.
\end{rem}

Let us finally give a lemma which provides a practical criterion for
proving assumptions {\bf (A2)-(iii)} in the enlargement 
theorem~\ref{theo:EnlargingSGdecay}:

\begin{lem}\label{lem:Tn} 
  % We consider the same setting as in
  % Theorem~\ref{theo:EnlargingSGdecay} and we assume {\bf (A1)} and
  % {\bf (A2)-(i)-(ii)}.
  
    Let $E, \EE$ be two Banach spaces with $E \subset \EE$ dense with
  continuous embedding, and consider $L\in~\CCC(E)$, $\LL \in
  \mathscr{C}(\EE)$ with $\LL_{|E} = L$ and $a \in \R$.

We assume:

\begin{itemize}
\item[{\bf (A3)}] there exist some ``intermediate spaces" (not
  necessarily ordered)
\[
E = \EE_J, \  \EE_{J-1},  \ \dots, \ \EE_2, \  \EE_1 = \EE, \quad J \ge 2,
\]
such that, still denoting  $\BB := \BB_{|\EE_j}$, $\AA := \AA_{|\EE_j}$, 

\begin{itemize}
\item[{\bf (i)}] $(\BB-a)$ is hypodissipative and $\AA$ is bounded on
  $\EE_j$ for $1 \le j \le J$.  
%\smallskip

 % \item[{\bf (ii)}] $\AA \in \BBB(\EE_j)$ for $1 \le j \le J-1$.
% \smallskip

 \item[{\bf (ii)}] There are some constants $\ell_0 \in \N^*$, $C
   \ge 1$, $K \in \R$, $\alpha \in [0,1)$ such that
   $$
   \forall \, t \ge 0, \quad 
   \|T_{\ell_0}(t) \|_{\BBB(\EE_{j},\EE_{j+1})} \le C \, {e^{Kt}\over t^\alpha},
   $$ 
   for $1 \le j \le J-1$, with the notation $T_\ell := (\AA \,
   S_\BB)^{(*\ell)}$.
  \end{itemize}  \end{itemize}

Then for any $a' > a$, there exist some explicit constants $n \in \N$,
$C_{a'} \ge 1$ such that
$$
\forall \, t \ge 0, \quad 
\|T_n(t) \|_{\BBB(\EE,E)} \le C_{a'}  \, e^{ a' \, t}.
$$
\end{lem}

\begin{proof}[Proof of Lemma~\ref{lem:Tn}.] 
On the one hand, {\bf (i)-(ii)} imply for $1 \le j \le J-1$ that 
\begin{equation}\label{lem:Tn-1}
\|T_1(t) \|_{\BBB(\EE_j)} \le C_a \, e^{at}
\end{equation}
and next
\begin{equation}\label{lem:Tn-2}
\|T_\ell \|_{\BBB(\EE_j)} \le C_a \, t^{\ell} \, e^{at}.
\end{equation}

On the other hand, for $n = p \, \ell_0$, $p \in \N^*$, we write 
\begin{eqnarray*}
  T_n(t) &=& \underbrace{ (T_{\ell_0} * \dots * T_{\ell_0} )(t) }_{\mbox{p
      times}}
  \\
  &=& \int_0^t \dd t_{p-1} \, \int_0^{t_{p-1}} \dd t_{p-2} \, \dots \,
  \int_0^{t_2} \dd t_1 \, T_{\ell_0} (\delta_p) \, \dots \, T_{\ell_0} (\delta_1)
\end{eqnarray*} 
with 
\[
\delta_1 = t_1, \quad \delta_2 = t_2 - t_1, \ \dots, \  
\delta_{p-1} = t_{p-1} - t_{p-2} \mbox{ and } \ \delta_p = t-t_{p-1}.
\]

For $p > J$, there exist at least $J-1$ increments $\delta_{r_1}$,
\dots, $\delta_{r_{J-1}}$ such that $\delta_{r_j} \le t/(p-J)$ for any
$1 \le j \le J-1$, otherwise there exist $\delta_{q_1}$, \dots,
$\delta_{q_{p-J}}$ such that $\delta_{q_j} > t/(p-J)$, and
$$
t = \delta_1 + \dots + \delta_p \ge \delta_{q_1} + \dots +
\delta_{q_{p-J}} > (p-J) \, { t \over p-J} =t
$$
which is absurd. 

Now, using {\bf (A3)-(ii)} in order to estimate $\|T_\ell (\delta_{r_j})
\|_{\BBB(\EE_j,\EE_{j+1})}$ and \eqref{lem:Tn-2} in order to bound the
other terms $\| T_\ell (\delta_r) \|_{\BBB(\EE_r)}$ in the appropriate
space, we have with $\QQQ := \{r_1, \dots, r_{J-1} \}$,
\begin{multline*}
\|T_n(t)  \|_{\BBB(\EE,E)} \\ \le    
   \int_0^t \!\! dt_{p-1} \int_0^{t_{p-1}} \!\!\!\! 
   \dd t_{p-2} \, \dots \int_0^{t_2} \!\! \dd t_1 \, \prod_{r \notin \QQQ  } C_a \, \delta_r^\ell \, e^{a \, \delta_r} 
    \, \prod_{q \in \QQQ} C \, { e^{ K \, \delta_q} \over \delta_q^\alpha}
    \\
  \le (C_a \, t)^{p-J} \, C^J \,   e^{a \, t}  \,  e^{ K \, {J \,  t \over p-J}}
   \int_0^t \!\! \dd t_{p-1} \int_0^{t_{p-1}} \!\!\!\! \dd t_{p-2} \,
   \dots \int_0^{t_2} \!\! \dd t_1  
    \, \prod_{j= 1}^J {1  \over \delta_{r_j}^\alpha} 
     \\
 \le  C' \, e^{ (a +  {K \, J \over p-J} ) \, t}\, t^{2p-1-J-J\alpha} 
   \int_0^1 \!\! \dd u_{p-1} \int_0^{u_{p-1}} \!\!\!\! \dd u_{p-2} \,
   \dots \int_0^{u_2} \!\! \dd u_1  
    \, \prod_{j= 1}^{p-1} {1  \over (u_{j+1} - u_j)^\alpha} ,
\end{multline*}
with the convention $u_p = 1$. Since the last integral is finite for any $p \in \N$, we easily conclude by just
taking $p$ (and then $n$) large enough so that $a+KJ/(p-J) < a'$. 
\end{proof}

%%%%%%%%%%%%%%%%%%%   SECTION 3 :   Fokker-Planck   %%%%%%%%%%%%%%%%%%%%%%%%%%%%%%%%

\bigskip
\section{The Fokker-Planck equation}
\label{sec:FP}
\setcounter{equation}{0}
\setcounter{theo}{0}

%%%%%%%%%%%%%%%%%%%%%%%%%%%%%%%%%%%%%%%%%%%%%%%%%%%%%%%%%%%%%%%%%%%%%%%%%%%%%

Consider the Fokker-Planck equation
\begin{align}\label{eq:FP}
  \partial_t f = L f := \nabla_v \cdot \left( \nabla_v f + F \, f\right),
  \quad f_0(\cdot) = \init{f}(\cdot),
\end{align}
on the density $f = f_t(v)$, $t \ge 0$, $v \in \R^d$ and where the
(exterior) force field $F = F(v) \in \R^d$ takes the form
\begin{equation}
  \label{eq:F}
  F = \nabla_v \phi + U,
\end{equation}
with confinement potential $\phi : \R^d \to \R$ of class $C^2$ and non gradient force field perturbation $U : \R^d \to \R^d$ of class $C^1$ so that 
\beqn\label{eq:FPstructure}
  \forall \, v \in \R^d, \quad \nabla_v \cdot (  U(v) \, e^{-\phi(v)}  )  = 0. 
\end{equation}
It is then clear that a stationary solution is 
$$
\mu (v) :=  e^{-\phi(v)}. 
$$

\smallskip In order for $\mu$ to be the global equilibrium we make
the following additional classical assumptions on the $\phi$ and $U$:

\begin{itemize}
\item[{\bf (FP1)}] The Borel measure associated to the function $\mu$
  and denoted in the same way, $\mu ({\rm d} v) := e^{-\phi(v)} \dd v$, is a
  probability measure and the function $\phi$ is $C^2$ and satisfies
  one of the two following large velocity asymptotic conditions
  \begin{equation}\label{eq:Cond1Pointcare} 
    \liminf_{|v|\to \infty} \, \left(
  {v \over |v|} \cdot \nabla_v \phi(v) \right) > 0 
\end{equation}
or 
\begin{equation}\label{eq:Cond2Pointcare}
  \exists \, \nu \in (0,1) \quad 
  \liminf_{|v|\to \infty} \left(  \nu \, |\nabla_v \phi|^2 - \Delta_v \phi \right)  > 0
\end{equation}
while the force field $U$ satisfies the growth condition
\begin{equation*}\label{FP-H1-2}
  \forall \, v \in \R^d,  
  \quad |U(v)| \le C \, \left(1 + \left|\nabla_v \phi(v)\right|\right).
\end{equation*}
\end{itemize}

It is crucial to observe that {\bf (FP1)} implies  that the measure $\mu$
satisfies the Poincar\'e inequality
\begin{align}\label{eq:FPPoincare}
  \int_{\R^d} \left| \nabla_v \left( \frac{f}{\mu} \right) \right|^2 \, \mu(dv)
  \ge 2 \, \lambda_P \, \int_{\R^d} f^2 \, \mu^{-1}(dv) \quad \mbox{ for } \quad
  \int_{\R^d} f \dd v =0,
\end{align}
for some constant $\lambda_P >0$.  We refer to the recent paper
\cite{MR2386063}  for an introduction to this important subject 
as well as to the references therein for further developments.  Actually
the above hypothesis {\bf (FP1)} could be replaced by assuming
directly that \eqref{eq:FPPoincare} holds. However, the conditions
\eqref{eq:Cond1Pointcare} and \eqref{eq:Cond2Pointcare} are more
concrete and yield criterion that can be checked for a given
potential. 

\smallskip The fundamental example of a suitable confinement potential
$\phi \in C^2(\R^d)$ which satisfies our assumptions is when
\begin{equation}\label{eq:FPexple} 
  \phi(v) \approx  \alpha \, |v|^\gamma \,\,\,\,\hbox{and}\,\,\,\, \nabla\phi(v) \approx
  \alpha \, \gamma \, v \, |v|^{\gamma-2} \quad\hbox{as}\quad |v| \to +\infty
\end{equation}
for some constants $ \alpha > 0$ and $\gamma \ge 1$. For instance, the
harmonic potential $\phi (v) = |v|^2/2 -(d/2) \, \ln (2\pi)$
corresponds to the normalised Maxwellian equilibrium $\mu(v) =
(2\pi)^{-d/2} \, \exp(-|v|^2/2)$.

\subsection{The Fokker-Planck equation: model and results}
For some given Borel weight function $m=m(v) > 0$ on $\R^d$, let us define
$L^p(m)$, $1 \le p \le 2$, as the Lebesgue space associated to
the norm
$$
 \| f \|_{L^p(m)} := \|f \, m \|_{L^p} = \left( \int_{\R^d} f^p(v)
   \, m(v)^p \dd v \right)^{1/p}.
$$
For any given positive weight, we define the \emph{defect weight
  function} 
\begin{equation}\label{eq:defpsimp} 
  \psi_{m,p}  := (p-1) \, {  |\nabla m|^2 \over m^2} + { \Delta m \over m} + \left(1 - {1 \over
    p} \right) \mbox{div} \, F - F \cdot {\nabla m \over m}.
\end{equation}
Observe that $\psi_{\mu^{-1/2},2}=0$: $\psi_{m,p}$ quantifies some
error to this reference case. 

%\[
%\psi_{m,p} := {1 \over m^p} \, \Bigl\{m^{p-2}\, |\nabla m|^2 + (1 -
% {1 \over p}) (\hbox{div} F) \, m^p - {1 \over p} F \cdot \nabla m^p
% \Bigr\}.
%\]
Let us enounce two more assumptions: 
\begin{itemize}
\item[{\bf (FP2)}] The weight $m$ satisfies $L^2(\mu^{-1/2}) \subset
  L^p(m)$ (recall $p \in [1,2]$) and the condition
\[
\limsup_{|v| \to \infty} \, \psi_{m,p} = a_{m,p}< 0.
\] 
 
\item[{\bf (FP3)}] There exists a positive Borel weight $m_0$ such
  that $L^2(\mu^{-1/2}) \subset L^q(m_0)$ for any $q \in [1,2]$ and
  there exists $b \in\R$ so that
\[
\left\{ 
\begin{array}{l} \ds
\sup_{q \in [1,2], \ v \in \R^d}
  \psi_{m_0,q} \le b, \vs \\ \ds 
  \sup_{x \in \R^d} \left( \frac{\Delta
    m_0}{m_0} - \frac{|\nabla m_0|^2}{m_0^2} \right)
\le b. 
\end{array}
\right.
\] 
\end{itemize}
The typical weights $m$ satisfying these assumptions are $m(v) \approx
e^{\kappa \, \phi}$ with $\kappa \in [0,1/2]$, $m(v) = e^{\kappa \,
  |v|^\beta}$ with $\beta \in [0,1]$ and $\kappa>0$ appropriately
chosen, or $m(v) \approx \langle v \rangle^k$, at large velocities.
\medskip

Here is our main  result on the Fokker-Planck equation. 

\begin{theo}\label{theo:FP1} 
  Assume that $F$ satisfies {\bf (FP1)} and consider a $C^2$ weight
  function $m > 0$ and an exponent $p \in [1,2]$ so that {\bf (FP2)}
  holds if $p=2$ and {\bf (FP2)}-{\bf (FP3)} holds if $p \in [1,2)$.
  
  Then for any initial datum $\initem{f} \in L^p(m)$, the associated
  solution $f_t$ to~\eqref{eq:FP} satisfies the following decay
  estimate
  \begin{equation}\label{eq:FPconclusionThFP} 
    \forall \, t \ge 0, \quad 
    \left\| f_t - \mu \, \langle \initem{f} \rangle \right\|_{L^p(m)} \le C \,
  e^{-\lambda_{m,p} \, t} \, \left\| \initem{f} - \mu \, \langle \initem{f} \rangle
  \right\|_{L^p(m)},
\end{equation}
with $\lambda_{m,p} := \lambda_P$ if $\lambda_P < |a_{m,p}|$, and
$\lambda_{m,p} < |a_{m,p}|$  as close as wanted to $|a_{m,p}|$ else,
and where we use the notation
\[
\langle \initem{f} \rangle := \int_{\R^d} \initem{f} \dd v.
\]
\end{theo}

\begin{rems}
  \begin{enumerate} 
  \item Note that this statement implies in particular that the
    spectrum of $\LL$ in $L^p(m)$ satisfies for $a$ as above:
\[
\Sigma(\LL) \subset \{ z\in \C \; | \; \Re e (z) \le a \}\cup \{ 0\} ,
\]
%  $L^2(\mu^{-1/2})$ that
%\[
%\Sigma(\LL) \subset \{ z\in \C \; | \; \Re e (z) \le
%  -\lambda_P \}\cup \{ 0\} 
%\]
%where $\lambda_P$ is the Poincar\'e constant in (\ref{eq:FPPoincare}), 
and that
the null space of $\LL$ is exactly $\R \mu$. 
\item When $m=\tilde m (\phi)$ and $\hbox{div} \, U = U \cdot \nabla
  \phi = 0$, an alternative choice for the \emph{defect weight
    function} associated to the weight $m$ and $p \in [1,2]$ could be
  $\psi_{m,p} =: \psi^1_{m,p} + \psi^2_{m,p}$ with
\begin{eqnarray*}
\psi^1_{m,p} &=& \frac{1}{p \,m^2 \mu^p} \, \nabla_v \cdot \left[
{\mu^p \, m^{2p-2}} \, \nabla_v \left( \frac{1}{m^{2p-4}} \right) \right] 
  \\
\psi^2_{m,p} &=&
\frac{(p-1)}{p} \,  m^{2p-2} \, \nabla_v \cdot \left[  \frac{\mu}{m^{2p-4}}\,
  \nabla_v \cdot \left( \frac{1}{m^2 \, \mu} \right) \right].
\end{eqnarray*}
Notice that again $\psi_{\mu^{-1/2},2} = 0$. The first part $\psi^1_{m,p}$
is related to the change in the Lebesgue exponent from $2$ to $p$, and
the second part $\psi^2_{m,p}$ is related to the change of weight from
$\mu^{-1/2}$ to $m$.
\item Concerning the weight function $m$, other technical assumptions
  could have been chosen for the function $m(v)$, however the
  formulation {\bf (FP2)-(FP3)} seems to us the most natural one since
  it is based on the comparison of the Fokker-Planck operators for two
  different force field. In the case $U=0$, $p=2$ and $m = e^{\phi/2}$
  the condition {\bf (FP2)} is nothing but the classical condition
  \eqref{eq:Cond2Pointcare} with $\nu = 1/2$.  In any case, the core
  idea in the decomposition is that a coercive $\BB$ in $\EE$ is
  obtained by a negative local perturbation of the whole operator.
\item By mollification the $C^2$ smoothness assumption of $m$ could be
  relaxed: if $m(v)$ is not smooth but $\tilde m(v)$ is smooth,
  satisfies {\bf (FP2)-(FP3)} and is such that $c_1\, m(v) < \tilde
  m(v) \le c_2 \, m(v)$, then it holds
    \begin{multline*}
      \| f_t -\mu \|_{\EE} \le C \,  
      \| f_t -\mu \|_{L^p(\tilde m)} \\ \le  C' \;
      e^{-\tilde\lambda \, t} \, \left\| \init{f}-\mu
      \right\|_{L^p(\tilde m)} \; \le\; C'' \, e^{-\tilde\lambda
        \, t} \, \left\| \init{f} - \mu \right\|_{\EE}.
    \end{multline*}
\item It is easy to extend the well-posedness of the Fokker-Planck
  equation to measure solutions, and using the case $p=1$ in the
  previous theorem (under appropriate assumptions on the weight) we
  deduce the following decay estimate 
\begin{equation*}
  \forall \, t \ge 0, \quad  \left\| f_t -  \mu \, \langle \init{f} \rangle  \right\|_{M^1(m^{-1})} \le 
  C \, e^{ - \lambda_{m,1} \, t} \, \left\| \init{f} -  \mu \langle \init{f} \rangle \right\|_{M^1(m^{-1})}
\end{equation*}
where $M^1(m^{-1})$ denotes the weighted space of measures with finite
mass. 
\end{enumerate}
\end{rems}

For concrete applications, for $\phi$ satisfying the power-law
asymptotic condition \eqref{eq:FPexple}, we have the following decay
rates depending on the weight $m$ and the exponent $\gamma$ in
\eqref{eq:FPexple}:
\begin{prop}\label{prop:FP1}
Assume that $\phi$ satisfies \eqref{eq:FPexple} with exponent $\gamma
\ge 1$, then: 
\begin{itemize}
\item[\textbf{(W1)}] {\bf Exponential energy weight.} For all $\gamma
  \ge 1$, the weight $m = e^{\kappa \, \phi}$ is allowed, where
  $\kappa$ satisfies $\kappa \in (0,1/2]$ when $p=2$ and $\kappa \in
  (0,1/2)$ when $p \in [1,2)$.

  Moreover, in these spaces the estimate we obtain on the exponential
  decay rate is the optimal Poincar\'e constant
\[
\lambda_{m,p}:= \lambda_P \ \mbox{ when }  \ \gamma > 1
\]
%or $p=2$ 
 while in the critical case $\gamma=1$ it is given by $\lambda_{m,p} =
 \lambda_P$ when $\lambda_P < \kappa \, (1-p\kappa)$, and by any
 $0 \le \lambda_{m,p} < \kappa \, (1-p\kappa)$ else 
% \[
% \lambda_{m,p} := \min\left\{\lambda_P ; \kappa \, (1-p\kappa)+0\right\}
% \]
% for the critical case $\gamma = 1$
 (which degenerates to zero as $\kappa \to 0$). The constant in front
 of the exponentially decaying term in \eqref{eq:FPconclusionThFP}
 blows-up as $\lambda_{m,p} \to \kappa \, (1-p\kappa)$ in the last
 case. %, $p \in [1,1/2)$.

% (since $a_{m,p} = -\infty$). 
\smallskip

\item[\textbf{(W2)}] {\bf Stretched exponential weight.} For all
  $\gamma > 1$, the weight $m = e^{\kappa \, |v|^\beta}$ is allowed
  for any $\kappa > 0$, $p \in [1,2]$ and $ 2-\gamma \le \beta <
  \gamma$.

Moreover, in these spaces the estimate we obtain on the exponential
  decay rate is the optimal Poincar\'e constant 
\[
\lambda_{m,p}:= \lambda_P \ \mbox{ when } \ \gamma + \beta > 2,
\]
while in the critical case $\beta = 2 - \gamma$ it is given by
$\lambda_{m,p} = \lambda_P$ if $\lambda_P < \kappa \beta \gamma$, and
by any $0 \le \lambda_{m,p} < \kappa \beta \gamma$ else % . 
% \[
% \lambda_{m,p} := \max\{ \lambda_P; \kappa \beta \gamma+0\} 
% \]
% for the critical case $\beta = 2 - \gamma$ (with $\gamma \in (1,2)$),
(which degenerates to zero as $\kappa$ goes to zero).  The constant in
front of the exponentially decaying term in
\eqref{eq:FPconclusionThFP} blows-up as $\lambda_{m,p} \to \kappa
\beta \gamma$ in the last case.
% The spectral
%gap is again the (optimal) Poincar\'e constant $\lambda_P$ when $\gamma + s > 2$ while 
%  (since $\lambda_{p,m} = +\infty$), however when $s\beta
%+2 =2$, we prove a rate of decay given by the minimum
%of the Poincar\'e constant and $\lambda_{p,m} = C_\theta s^2 \beta/p$.
 \smallskip

\item[\textbf{(W3)}] {\bf Polynomial weight.}  For all $\gamma \ge 2$,
  the weight $m = \langle v \rangle^k$ is allowed for the Lebesgue
  exponent $p \in [1,2]$ under the condition
\begin{equation}\label{eq:FPcondpolyk}
(\gamma-2+d)\left(1-\frac1p\right) < k.
\end{equation}

Moreover, in these spaces the estimate we obtain on the exponential
  decay rate is the optimal Poincar\'e constant 
\[
\lambda_{m,p}:= \lambda_P \ \mbox{ when } \ \gamma > 2,
\]
while in the critical case $\gamma=2$ it is given by $\lambda_{m,p} =
\lambda_P$ if $\lambda_P < 2k - 2d (1-1/p)$, and by any $0 \le
\lambda_{m,p} < \lambda_P$ else (which degenerates to zero as $\kappa$
goes to zero).  The constant in front of the exponentially decaying
term in \eqref{eq:FPconclusionThFP} blows-up as $\lambda_{m,p} \to 2k
- 2d (1-1/p)$ in the last case.
% \[
% \lambda_{m,p} := \min\left\{ \lambda_P; 2k - 2 d \, \left(1
%     -\frac1p\right) +0 \right\}
% \]
% for the critical case $\gamma = 2$, which degenerates to zero when
% $p=1$ and $k$ goes to $0$. 
%of the Poincar\'e constant and $\lambda_{p,m} = 2 (2k-d(p-1))$.
\end{itemize}
\end{prop}

\begin{rems}
\begin{enumerate}
% \item When $a_{p,m} < - \lambda_P$, one can improve the statement with
%   the sharper rate of convergence $a = -\lambda_P$, by considering a
%   more detailed expansion of the semigroup for the first non-zero
%   eigenvalue.
  %
\item Observe how the polynomial weights are sensitive to the Lebesgue
  exponent $p$ in the condition \eqref{eq:FPcondpolyk}. We believe the
  restriction on the polynomial weight (depending on $p$, $\gamma$ and
  $d$) to be optimal. Accordingly we expect that in the case
  $\gamma=2$ the optimal value of the spectral gap is given by 
  \[
  \lambda_{m,2} := \max\left\{\lambda_P;2k - 2 d \, \left(1
      -\frac1p\right) \right\}.
  \]
  This is still an open question that needs to be proven, or disproven. 
   However we can give a partial positive answer: for potentials $\phi$
  satisfying \eqref{eq:FPexple} with $\gamma=2$, and polynomial
  weights $m = \langle v \rangle^k$, then the constant $\lambda_{m,2}
  = 2k-d$, $k >d/2$, coincides with the value of the spectral gap
  explicitly computed by Gallay and Wayne in \cite[Appendix
  A]{MR1912106}.
\item Observe furthermore that in the case of a polynomial weight we
  require the confinement potential to be quadratic or
  over-quadratic. This is reminiscent of the logarithmic Sobolev
  inequality, however this is strictly weaker than asking the
  confinement potential to satisfy the logarithmic Sobolev
  inequality. It is an open question to know whether a spectral gap
  still exists when the potential is subquadratic ($\gamma \in [1,2)$)
  and the weight is polynomial.
%  Again we believe this condition to be optimal and it is
%  another interesting open question to prove or disprove the existence
%  of a spectral in Lebesgue spaces with polynomial weights when $s \in
%  [1,2)$.
%
\item When $\gamma \ge 2$, $p=1$ and the weight is polynomial any $k
  >0$ is allowed, which means that it almost includes $L^1$ without
  weight. We expect that in the limit case $L^1$ there is no spectral
  gap and the continuous spectrum touches zero in the complex plane.
% \item Let us emphasize the three thresholds obtained for $\gamma =1$,
%   $\gamma \in (1,2)$ and $\gamma = 2$.
  %
 
%
%\item Here discuss Gallay-Wayne maybe in more details, and the result of
%Villani (hypocoercivity memoir) based on entropic hypocoercivity. 
%Maybe gather here also all the open questions raised in the remarks. 
 %
 \item 
  %Let us make a last important remark. 
   Another strategy for proving the decay of the semigroup could have
   been the use of interpolation between the exponential relaxation in
   $E$ together and a uniform bound in $L^1$ (provided by mass
   conservation and preservation of non-negativity). However, first,
   it would not recover optimal rates of decay, and second, most
   importantly, it would not apply to semigroups which do not preserve
   non-negativity (and consequently do not preserve the $L^1$ norm),
   such as those obtained by linearization of a bilinear operator that
   we consider see later in this paper.
 \end{enumerate}
\end{rems}

We give a simple application of our main result, related to the remark
(2) above. 

\begin{cor}\label{corFP} Assume that $\phi$ satisfies
  \eqref{eq:FPexple} with exponent $\gamma \in [1,2)$. Then for any $k
  > 0$, there exists $C = C(k,\gamma,d) \in (0,\infty)$ such that for
  any initial datum $\initem{f}\in L^1(\langle v \rangle^k)$, the
  solution to the initial value problem (\ref{eq:FP})-(\ref{eq:F})
  satisfies the decay estimate
\begin{equation}\label{eq:FPdecaypoly}
  \forall \, t \ge 0, \quad 
  \left\| f_t -  \mu \, \langle \initem{f} \rangle  \right\|_{L^1} \le 
 C \,  t^{- \frac{k}{2-\gamma}}\,  
 \left\| \initem{f} -  \mu \, \langle \initem{f} \rangle \right\|_{L^1(\langle v \rangle^k)}.
\end{equation}
\end{cor} 

\begin{rem} 
  A similar result has been proved in \cite[Theorem 3]{MR1751701}
  under the additional and fundamental assumptions that $\init{f}$ is
  non negative and has finite energy and entropy. Moreover the decay
  rate obtained in \cite{MR1751701} was only of order
  $t^{-(k-2)/(2(2-\gamma))}$ and remains valid for $\gamma \in (0,1)$.
\end{rem}

%\subsection{On the proof of Theorem~\ref{theo:FP1}. }
\subsection{Proof of the main results}
The proof of Theorem~\ref{theo:FP1} is based on the combination of the
spectral gap in the space $L^2(\mu^{-1/2})$ given by Poincar\'e's
inequality together with the extension to functional spaces of the form
$L^p(m)$, by applying Theorem~\ref{theo:EnlargingSGdecay}.

\smallskip Before going into the proof of Theorem~\ref{theo:FP1}, let
us remark that most of the interesting external forces and weights do
satisfy our assumptions, as detailed below.

\begin{lem}\label{lem:admissible} When $\phi$ satisfies
  \eqref{eq:FPexple} and $U \equiv 0$, conditions {\bf (FP1)}-{\bf
    (FP2)}-{\bf (FP3)} are met under conditions {\bf (W1)}, {\bf (W2)}
  and {\bf (W3)} in the statement of Proposition~\ref{prop:FP1}.
\end{lem}

\begin{proof}[Proof of Lemma~\ref{lem:admissible}.] For the sake of
  simplicity we assume $\phi(v) = |v|^\gamma$, $\gamma >
  0$, for $|v|$ large enough, and we show that the large velocity
  behavior properties in {\bf (FP1)}-{\bf (FP2)}-{\bf (FP3)} hold
  under the suitable conditions.  The proof in the general case
  \eqref{eq:FPexple} is exactly similar.

 \smallskip
 First we compute for large velocities
 $$
 \nabla \phi = \gamma \, v \, |v|^{\gamma-2}, \quad
 \hbox{div} \, F = \Delta \phi = \gamma \, ( d +\gamma-2) \,
 |v|^{\gamma-2},
$$
and we observe that both conditions \eqref{eq:Cond1Pointcare} and
\eqref{eq:Cond2Pointcare} (for any $\nu \in (0,1)$) are satisfied when
$\gamma \ge 1$, so that condition {\bf (FP1)} holds.

\smallskip\noindent {\sl Step 1. Exponential weight.} We consider $m
:= \exp ( \kappa \, |v|^\beta)$, $\kappa, \, \beta > 0$, and we
compute for large velocities
\begin{eqnarray*}
  &&
  \nabla m = \kappa \, \beta \,  v \, |v|^{\beta-2} \, m,
  \quad
  \Delta m = \kappa \, \beta \, (\beta-1)   |v|^{\beta-2} \, m +  \kappa^2 \, \beta^2 \,   |v|^{2\beta-2} \, m.
\end{eqnarray*}
We observe that in that case
\begin{eqnarray*}
  \psi_{m,p}  % &=& (p-1) { |\nabla m|^2 \over m^2} +\frac{\Delta m}{m}
  % + \left( 1- \frac1p \right) \Delta \phi-
  % \nabla \phi  \cdot {\nabla m \over m} \\
  &\approx&  (p-1) { |\nabla m|^2 \over m^2} +\frac{\Delta m}{m} -  \nabla \phi  \cdot {\nabla m \over m}
  \\
  &\approx&  (p-1) \, \kappa^2 \, \beta^2 \, |v|^{2\beta-2} + \kappa^2 \, \beta^2 \, |v|^{2\beta-2}
  -  \kappa \, \beta \, \gamma \, |v|^{\beta+\gamma-2} \\ 
  &\approx& p \, \kappa^2 \, \beta^2 \, |v|^{2\beta-2} -  \kappa \, \beta \, \gamma \, |v|^{\beta+\gamma-2}
\end{eqnarray*}
since the third term is always smaller that the fourth term when
$\beta>0$ and using the asymptotic estimates.  The condition $2-\gamma
\le \beta$ comes from (and is equivalent to) the fact that the last
term does not vanish in the large velocity asymptotic and the
condition $\beta \le \gamma$ comes from (and is equivalent to) the
fact that the last term is not negligible with respect to the first
term in the large velocity asymptotic.

 \smallskip
 When $\beta=\gamma$, we find 
 $$
 \psi_{m,p}  \, \approx \, \kappa \gamma^2 \, (p \kappa -1) \, |v|^{2\gamma-2},
 $$ 
 from which we get the condition $ p\kappa < 1$, and we conclude to
 $a_{m,p} = -\infty$ when $\gamma > 1$ while $a_{m,p} = \kappa \, (
 p\kappa - 1)$ when $\gamma=1$. However in order to have
 $L^2(\mu^{-1/2}) \subset L^p(m)$, we find the additional condition
 $\kappa \in (0,1/2)$.

 \smallskip When $\beta < \gamma$, we find
 $$
 \psi_{m,p} \, \approx \, - \kappa \, \beta \, \gamma \,
 |v|^{\beta+\gamma-2},
 $$ 
 so that $a_{m,p} = - \infty$ when $\beta> \gamma - 2$ and $a_{m,p} =
 - (\kappa \beta \gamma)$ when $\beta = \gamma - 2$.
 
 Finally, condition {\bf (FP3)} is always satisfied for $\gamma \ge 1$
 with $m_0 := e^{\kappa \, \phi}$, $\kappa \in (0,1/2)$.
 
 \smallskip\noindent {\sl Step 2. Polynomial weight.} We consider $m
 := \langle v \rangle^k$, $k > 0$, and we compute for large velocities
\begin{equation*}
\left\{ 
  \begin{array}{ll} \ds
    \nabla m = k \,  v \, \langle v \rangle^{k-2}, \quad 
    & \Delta m \approx k \, (d + k-2) \, \langle v \rangle^{k-2}, \vs \\ \ds 
    \nabla \phi = \gamma \,  v \, \langle v \rangle^{\gamma-2}, \quad 
    & \Delta \phi \approx \gamma \, (d+\gamma-2) \, \langle v \rangle^{\gamma-2}. 
\end{array}
\right.
\end{equation*}

It holds
\begin{eqnarray*}
  \psi_{m,p}  % &=& (p-1) { |\nabla m|^2 \over m^2} +\frac{\Delta m}{m}
  % + \left( 1- \frac1p \right) \Delta \phi-
  % \nabla \phi  \cdot {\nabla m \over m} \\
  &\approx&  \left(1 - {1 \over p}\right) \, \Delta \phi  
  -  \nabla \phi  \cdot {\nabla m \over m}
  \\
  &\approx&  \left(1 - {1 \over p}\right) \,  
  \gamma \, ( d +\gamma-2) \, \langle v\rangle^{\gamma-2}  -  \gamma
  \, k \, \langle v\rangle^{\gamma - 2},
\end{eqnarray*}
since the first and second terms are negligible as soon as $\gamma
>-1$. We assume $\gamma \ge 2$ so that the limit is non-zero. We
easily deduce the condition \eqref{eq:FPcondpolyk} and $a_{m,p} =
-\infty$ when furthermore $\gamma > 2$ while $a_{m,p}:= 2 d \, (1
-1/p) - 2k$ when $\gamma =2$.
 \end{proof}

\begin{lem}\label{lem:FP-DissipBB}
  Under the assumptions {\bf (FP1)-(FP2)}, there exists $M,R$ such
  that
  $$
  \BB := \LL - \AA, \quad \AA f := M \, \chi_R \, f
  $$
  satisfies the dissipativity estimate 
\begin{equation}\label{eq:FP-DissipBB}
  \forall \, t \ge 0, \quad  \| S_\BB(t) f \|_{L^p(m)} \le e^{-
    \lambda_{m,p} \, t} \,  \|  f \|_{L^p(m)}. 
\end{equation}
\end{lem}

\begin{proof}[Proof of Lemma~\ref{lem:FP-DissipBB}] 
  We calculate
\begin{multline*}
\int_{\R^d} (\LL f) \, |f|^{p-2} f \, m^p  \dd v \\
= \int_{\R^d} (\Delta f) \,  |f|^{p-2} f \, m^p \dd v + 
\int_{\R^d} \hbox{div} \, (F \, f) \, |f|^{p-2} f \, m^p \dd v =: T_1 + T_2 . 
 \end{multline*}
%For the first term, we introduce the function  $h := m^\theta \, f$, $\theta := p / (2p-2)$, the real number $\Theta := p-\theta (p-1) = p/2$, and
%observing that $\theta \, (p-1) = \Theta$, we compute 
%\begin{eqnarray*}
%T_1
%&=& - \int  \nabla ( h^{p-1} \, m^\Theta) \cdot \nabla (h \, m^{-\theta}) 
%\\
%&=& - \int  \{ \nabla h^{p-1}   \cdot \nabla h \, m^{\Theta+\theta}  + h^p \, \nabla m^\Theta \cdot  \nabla  m^{-\theta} \} 
%\\
%&=& - (p-1) \int  | \nabla h|^2 \, h^{p-2}  \, m^{\Theta+\theta} + {p^2 \over 4(p-1)} \int  f^p \, m^p \, {|\nabla m|^2 \over m^2}  
%\end{eqnarray*}
%thanks to one integration by parts. 
 For the first term $T_1$, we compute
\begin{eqnarray*}
T_1
&=& - \int_{\R^d} \nabla \left(|f|^{p-2} f \, m^p\right) \cdot \nabla  f \dd v
\\
&=& - \int_{\R^d}  \left[ \nabla \left( |f|^{p-2} f\right) \cdot \nabla f \, m^p 
 + p \, |f|^{p-2} f \, m^{p-1}\, \nabla f \cdot  \nabla   m \right] \dd v
\\
&=& - (p-1) \int_{\R^d} |\nabla f|^2 \, f^{p-2} \, m^p \dd v   
+ \int_{\R^d} |f|^p \, \hbox{div} \, \left( m^{p-1} \, \nabla m \right)
\dd v
\end{eqnarray*}
thanks to two integrations by parts. 
For the second term, we write 
\begin{eqnarray*}
T_2 &=& \int_{\R^d} (\hbox{div} \, F) \, |f|^p \, m^p \dd v + \int_{\R^d} (F
\cdot \nabla f) \, |f|^{p-2} f \, m^p \dd v
\\
&=&  \int_{\R^d} (\hbox{div} \, F) \, |f|^p \, m^p \dd v - {1 \over p}
\int_{\R^d} |f|^p \, \hbox{div} \, (F \, m^p) \dd v
\end{eqnarray*}
by integration by parts again. All together, we obtain the following
identity and estimate
\begin{eqnarray*}
  \int_{\R^d} (\LL f) \, |f|^{p-2} f \, m^p \dd v 
  &=& (1-p)  \int_{\R^d} | \nabla f|^2 \, f^{p-2}  \,
  m^p  \dd v 
  +  \int_{\R^d} |f|^p  \, m^p \, \psi_{m,p} \dd v 
  \\
  &\le&  \int_{\R^d} |f|^p  \, m^p \, \psi_{m,p} \dd v.
\end{eqnarray*}
From {\bf (FP2)}, for any $a > a_{m,p}$, we may find $M$ and $R$ large
enough so that
 $$
 \forall \, v \in \R^d, \quad \psi_{m,p} - M \, \chi_R \le a. 
 $$
 As a consequence, we deduce
 $$
 \int_{\R^d} (\BB f) \, |f|^{p-2} f \, m^p \dd v \le a \, \int_{\R^d}
 |f|^p \, m^p \dd v,
$$
from which \eqref{eq:FP-DissipBB} immediately follows. 
\end{proof}

We now shall prove a lemma about the regularization properties of the
Fokker-Planck equation. It is related to the notion of
\emph{ultracontractivity} and is well-known; we include a sketch of
its proof for clarity and in order to make the constants explicit.

\begin{lem}\label{lem:FP-RegBB}
  Under the assumptions {\bf (FP3)}, there are $b,C>0$ such that for
  any $p,q$ with $1 \le p \le q \le 2$, we have
\begin{equation}\label{lem:FP-RegBB}
\forall \, t \ge 0, \quad \|S_\BB(t) f \|_{L^q(m_0)} \le 
C \, e^{2bt} \, t^{-{d\over2} \left({1 \over p} - {1 \over q}\right)} \,  \|  f \|_{L^p(m_0)}.
\end{equation}

As a consequence, under the assumptions {\bf (FP2)}-{\bf (FP3)}, there
are $b,C >0$ such that for any $p,q$ with $1 \le p \le q \le 2$, we
have
\begin{equation}\label{lem:FP-RegAB}
  \forall \, t \ge 0, \quad \| T_\ell(t) f \|_{L^q(m)} \le 
  C \, e^{2bt} \, t^{-{d\over2} \left({1 \over p} - {1 \over q}\right)} \,  \|  f \|_{L^p(m)}
\end{equation}
with $\ell = 1$ when $L^p(m) \subset L^p(m_0)$ and with $\ell = 2$ in
the general case.
\end{lem}

\begin{proof}[Proof of Lemma~\ref{lem:FP-RegAB}] 
  From condition {\bf (FP3)} on $\psi_{m_0,p}$, by arguing as in the
  proof of Lemma~\ref{lem:FP-DissipBB} we obtain for any $p \in [1,2]$
  \begin{equation}\label{lem:FP-RegBB1} 
    \forall \, t \ge 0,
    \quad \|S_\BB(t) f \|_{L^p(m_0)} \le C_{pp} \, \| f
    \|_{L^p(m_0)}, \quad C_{pp} := e^{bt}.
\end{equation} 

\smallskip In order to establish the gain of integrability estimate we
have to use the non positive term involving the gradient in a sharper
way, i.e. not merely the fact that it is non-positive.  It is enough
to do that in the simplest case when $p=2$. Let us consider a solution
$f_t$ to the equation
$$
\partial_t f_t = \BB \, f_t, \qquad f_0 \in L^2(m_0).
$$
From the computation made in the proof of Lemma~\ref{lem:FP-DissipBB},
we have
\begin{multline*}
{{\rm d} \over {\rm d}t} \int_{\R^d} {f^2_t \over 2}  \, m_0^2 \dd v  
= - \int_{\R^d} |\nabla f_t|^2 \, m_0^2 \dd v
+ \int_{\R^d} f^2_t \, \left\{\psi_{m_0,2} - M \, \chi_R \right\} \,
m_0^2 \dd v \\
= - \int_{\R^d} |\nabla ( f_t m_0) |^2 \dd v
+ \int_{\R^d} f^2_t \, \left\{\psi_{m_0,2} - \frac{|\nabla m_0|^2}{m_0^2} +
  \frac{\Delta m_0}{m_0} - M \, \chi_R \right\} \,
m_0^2 \dd v \\ 
\le - \int_{\R^d} |\nabla ( f_t m_0) |^2 \dd v
+ 2 b \, \int_{\R^d} f^2_t \, m_0^2 \dd v.
\end{multline*}
Using Nash's inequality (\cite[Chapter~8]{MR1817225})
\[
\left( \int_{\R^d} g^2 \dd v \right) \le K_d \, \left( \int_{\R^d} |
  \nabla_v g |^2 \dd v \right)^{\frac{d}{d+2}} \, \left( \int_{\R^d} |g|
  \dd v \right)^{\frac{4}{d+2}}
\]
(for some constant $K_d>0$ depending on the dimension) applied
to $g = f_t \, m_0$, we get 
\begin{multline*}
{{\rm d} \over {\rm d}t} \int_{\R^d} {f^2_t \over 2}  \, m_0^2 \dd v  \\
\le - K_d ^{-1} \, \left( \int_{\R^d} |f_t|\, m_0 \dd v \right)^{-\frac{4}d}
\, \left(  \int_{\R^d} | f_t|^2 \, m_0^2 \dd v \right)^{\frac{d+2}{d}}
+ 2b \, \int f^2_t   \, m_0^2 \dd v. 
\end{multline*}

We then introduce the notation 
\[
X(t) := \|f_t \|_{L^2(m_0)}^2, \quad Y(t) := \| f_t \|_{L^1(m_0)}.
\]
Since $Y_t \le C \, Y_0$ for $t \in [0,1]$ by the previous step, we
end up with the differential inequality
\begin{equation}\label{eq:Nash0}
  \fa 0 \le t \le 1, \quad X'(t) \le 
  - 2 \, K Y_0^{-4/d} \, X(t)^{1+{2\over d}} + 2b \, X(t), 
\end{equation}
with $K \in (0,\infty)$.  On the one hand, if 
\[
X_0 > \left( \frac{2b}{K}\right)^{d/2} \, Y_0^2
\]
we define 
\[
\tau := \sup \left\{ t \in [0,1] \ \mbox{ s.t. } \  \forall \, s \in
[0,t], \  X(s) \ge \left( \frac{2b}{K}\right)^{d/2} \,\,
\right\}  \in (0,1],
\]
and the previous differential inequality implies
$$
\forall \, t \in (0,\tau), \quad X'(t) \le - K Y_0^{-4/d} \,
X(t)^{1+{2\over d}},
$$
which in turns implies 
\begin{equation}\label{eq:Nash2}
\forall \, t \in (0,\tau), \quad X(t) \le \left( {2 \, K \,
    Y_0^{-4/d} \, t \over d} \right)^{-d/2}.
\end{equation}

On the other hand, when $\tau < 1$ (so that $X(\tau) = (2b/K)^{d/2} \,
Y_0^2$), which includes the case $\tau=0$ and $X_0 \le (2b/K)^{d/2} \,
Y_0^2$, we simply drop the negative part in the right hand side of
\eqref{eq:Nash0} and get
\begin{equation}\label{eq:Nash4}
\forall \, t \in (\tau,1], \quad X(t) \le e^{(t-\tau) \, 2b}  \, \left( \frac{2b}{K}\right)^{d/2} \, Y_0^2.
\end{equation}

Gathering \eqref{eq:Nash2} and \eqref{eq:Nash4}, we obtain
\begin{equation}\label{eq:Nash1}
\fa t \in [0,1], \quad X(t)^{1/2} \le C \, t^{-d/4} \, e^{2b t} \, Y_0.
\end{equation}
Putting together \eqref{eq:Nash1} and the estimate
\eqref{lem:FP-RegBB1} with $p=2$ for the later times $t \ge 1$ we
conclude that
% On the one hand, if $X_0$ is such that $X_0 \le C \, Y_0^2$ then the
% first step implies that On the other hand, if $X_0 \in \CC := \{ f;
% \,\, \| f \|_{L^2(m_0)}^2 > (b/K)^{1/2} \, \| f \|_{L^1(m_0)}
% \} $ then the above differential inequality implies
%$$
%{d \over dt } X(t) \le -  K Y_0^{-4/d} \, X(t)^{1+{2\over d}} \quad\hbox{as long as}\quad X_t \in \CC,\,\, 0 \le t \le 1, 
%$$
%which in turns implies 
%\begin{equation}\label{eq:Nash2}
%X(t) \le \Bigl( {2 \over d} \, K \, Y_0^{-4/d} \, t \Bigr)^{-d/2}.
%\end{equation}
%Gathering \eqref{eq:Nash1} and \eqref{eq:Nash2}, we have proved 
$$
\forall \, t \ge 0, \quad \| S_\BB (t) f \|_{L^2(m_0)} \le C_{12} \,
\| f \|_{L^1(m_0)},  \quad C_{12}:= C \, e^{2bt} \, t^{-d/4}.
$$
Using twice the Riesz-Thorin interpolation theorem on the operator
$S_\BB(t)$ which acts in the spaces $L^1 \to L^1$, $L^2 \to L^2$ and
$L^1 \to L^2$, we obtain
 $$
 \| S_\BB (t) f \|_{L^q(m_0)} \le C_{qp} \, \|  f \|_{L^p(m_0)},
 \quad C_{qp} := C_{22}^{2 - 2/p} \, C_{11}^{2/q-1} \,
 C_{12}^{2/p-2/q},
 $$
 for any $1 \le p \le q \le 2$, which concludes the proof of
 \eqref{lem:FP-RegBB}.
  \end{proof}

\begin{proof}[Proof of Theorem~\ref{theo:FP1}.]
Let us proceed step by step. 
\smallskip
 
\noindent {\bf Step 1. The $L^2$ case for energy weight.} Let us first
review the spectral gap properties of the Fokker-Planck equation in
the space $L^2(\mu^{-1/2})$. On the one hand, performing one
integration by parts, we have
\begin{eqnarray*}
  \int_{\R^d} \hbox{div} \, (U f) \, \mu^{-1} \, f  \dd v
  &=&\int_{\R^d} \hbox{div} \, (U \mu) \,  (\mu^{-1} \, f )^2 \dd v \\
  && \qquad \qquad + {1 \over 2} \int_{\R^d} U \, \mu \cdot \nabla (\mu^{-1} f )^2 \dd v
  \\
  &=&{1 \over 2} \int_{\R^d} \hbox{div} \, (U \mu) \,  (\mu^{-1} \, f
  )^2 \dd v = 0.
\end{eqnarray*}
It is then immediate to check thanks to the Poincar\'e inequality
\eqref{eq:FPPoincare} that 
\begin{align*}
  2 \, \Re e \, \langle Lf, f \rangle_{L^2(\mu^{-1/2})}  := &\int_{\R^d} L \bar f \, f
  \, \mu^{-1}(\dd v) + \int_{\R^d} L f \, \bar f \, \mu^{-1}(\dd v) 
  \\
 = &- 2 \, \int_{\R^d} \left|\nabla_v \left(\frac{f}{\mu} \right) \right|^2
 \, \mu(\dd v)
 \\
 \le &- 2 \, \lambda_P \int_{\R^d} f^2 \, \mu^{-1} \,\dd v
\end{align*}
when $\langle f \rangle =0$. For any $\init{f} \in L^2(\mu^{-1/2})$
such that $\langle \init{f} \rangle = 0$ and then $\langle f_t \rangle
=0$ for any $t \ge 0$, we deduce that the solution $f_t$ to the
Fokker-Planck equation satisfies
$$
{{\rm d} \over {\rm d} t} \|f_t \|_{L^2(\mu^{-1/2})} \le - \lambda_P \,
\| f_t \|_{L^2(\mu^{-1/2})}
$$
from which we obtain estimate \eqref{eq:FPconclusionThFP} in the case
of the small space $E := L^2(\mu^{-1/2})$.

\smallskip\noindent {\bf Step 2. The $L^2$ case with general weight.}
Let us write $\EE = L^2(m)$ with $m$ satisfying {\bf (FP2)} and $E =
L^2(\mu^{-1/2})$, and denote by $\LL$ and $L$ the Fokker-Planck when
considered respectively in $\EE$ and $E$. 

We split the operator as $\LL = \AA + \BB$ with 
$$
\AA f := M \, \chi_R f \quad\hbox{and}\quad \BB f := \hbox{div} \,
(\nabla f + F \, f) - M \chi_R f.
$$
We then have $\AA \in \BBB(\EE,E)$ and, thanks to
Lemma~\ref{lem:FP-DissipBB}, we know that $\BB-a$ is dissipative for
any fixed $a > a_{m,2}$.  We can therefore apply
Theorem~\ref{theo:EnlargingSGdecay} with $n=1$ which yields the
conclusion.  \smallskip

\noindent {\bf Step 3. The $L^p$ case, $p \in [1,2]$.} With the same
splitting we have $\AA \in \BBB(\EE)$ as well as $T_2 (t)$ satisfies
condition {\bf (iii)} in Lemma~\ref{lem:Tn} thanks to
lemma~\ref{lem:FP-RegAB}. We can conclude by applying
Theorem~\ref{theo:EnlargingSGdecay} with $n=2$.
\end{proof}

\begin{proof}[Proof of Corollary~\ref{corFP}.]
  We proceed along the line of the proof of
  \cite[Theorem~3]{MR1751701}. Without loss of generality, we may
  assume that $\langle \init{f} \rangle = 0$.  For any $R > 0$, we
  split the initial datum as
$$
\left\{
\begin{array}{l} \ds
\init{f} := \init{f}^1 + \init{f}^2 \vs \\ \ds
\init{f}^1 := \init{f} \, {\bf 1}_{|v| \le R} - \left \langle
  \init{f} \, {\bf 1}_{|v| \le R} \right\rangle \vs \\ \ds
\init{f}^2 := \init{f} \, {\bf 1}_{|v| \ge R} - \left \langle
  \init{f} \, {\bf 1}_{|v| \ge R} \right\rangle
\end{array}
\right.
$$
and we denote by $f_t^1$ and $f_t^2$ the two solutions of the
Fokker-Planck equation respectively associated with the initial data
$\init{f}^1$ and $\init{f}^2$.  Since $\init{f}^1 \in
L^1(e^{|v|^{2-\gamma}})$ and satisfies $\langle \init{f}^1 \rangle
=0$, we may apply Theorem~\ref{theo:FP1} and we get
$$
\left\| f_t ^1 \right\|_{L^1(e^{|v|^{2-\gamma}})} \le C \, e^{-
  \lambda \, t} \, \| \init{f}^1 \|_{L^1(e^{|v|^{2-\gamma}})} \le C \,
e^{- \lambda \, t} \, {e^{R^{2-\gamma}} \over R^k} \, \left\|
  \init{f}^1 \right\|_{L^1(\langle v \rangle^k)}.
$$
On the other hand, the mass conservation for the Fokker-Planck
equation implies
$$
\| f^2_t \|_{ L^1} \le \| \init{f}^2 \|_{ L^1 } \le {1 \over R^k} \,
\| \init{f}^2 \|_{L^1(\langle v \rangle^k)}.
$$
We conclude by gathering the two estimates and choosing $R$ such that
$R^{2-\gamma} = \lambda t$.
\end{proof}

\subsection{The kinetic  Fokker-Planck equation in a periodic box}

Consider the equation
\begin{align}\label{eq:KFP}
  \partial_t f = L f := \nabla_v \cdot \left( \nabla_v f + \phi \,
    f\right) - v \cdot \nabla_x f, 
  \quad f_0(\cdot) = \init{f}(\cdot),
\end{align}
for $f = f_t(x,v)$, $t \ge 0$, $x \in \T^d$ the flat $d$-dimensional
torus, $v \in \R^d$, and for some velocity potential $\phi=\phi(v)$. 

\begin{itemize}
\item[{\bf (KFP1)}] The function $\phi$ is $C^2$ and such that $\mu
  (\dd v) := e^{-\phi(v)} \dd v$ is a probability measure and 
\[
W_\phi(v) := \frac{\Delta_v \phi}{2} - \frac{|\nabla_v \phi|^2}{4}
\xrightarrow[|v|\to+\infty]{} -\infty.
\] 
Moreover we assume that 
\begin{equation}\label{eq:reg-kfp}
\frac{\left| \nabla^s _v \phi (v)\right|}{\left| \nabla_v \phi (v)\right|}
\xrightarrow[|v|\to+\infty]{} 0 \quad \mbox { for } \ s=2,\dots, 4.
\end{equation}
This assumption is needed when deriving hypoelliptic regularization
estimates which involves taking velocity derivatives of the equation. 
\smallskip
\end{itemize}

Observe that the condition {\bf (KFP1)} is satisfied for any 
\[
\phi(v) = C_\phi \, (1+|v|^2)^{\gamma/2},  \quad \gamma>1
\]
(but does not cover the borderline case $\phi \sim |v|$ for the
Poincar\'e inequality). And as before it implies that the probability
measure $\mu$ satisfies the Poincar\'e inequality
\eqref{eq:FPPoincare} in the velocity space for some constant
$\lambda_P >0$.  It also implies the stronger inequality
  \begin{equation}\label{eq:Poincare-fort}
    \int_{\R^d} \left| \nabla_v \left( \frac{f}{\mu} \right) \right|^2
    \, \mu(\dd v)
    \ge 2 \, \bar \lambda_P \, \int \left( f - \int_{\R^d} f(v_*) \dd v_*
    \right) ^2 \, \left(1+ \left|\nabla_v \phi\right|^2\right)
    \, \mu^{-1}(\dd v) 
\end{equation}
for some constant $\bar \lambda_P >0$ (see \cite{MRS10} for a
quantitative proof).

\smallskip

For simplicity we normalize without loss of generality the volume of
the space torus to one.

% We assume on the confining potential $\Phi=\Phi(x)$:
% \begin{itemize}
% \item[{\bf (KFP1)}] The function $\Phi \in C^\infty$ is such that $\pi (dx)
%   := e^{-\Phi(x)} \, dx$ is a probability measure and $\nabla^q \Phi
%   \in L^\infty$ for all $q \ge 2$, and
% \[
% W_\Phi(x) := \frac{\Delta \Phi}{2} - \frac{|\nabla \Phi|^2}{4}
% \sim_{|x|\infty +\infty} - \frac{|\nabla \Phi|^2}{4}
% \xrightarrow[|x|\to+\infty]{} -\infty.
% \] 
% \end{itemize}
% \smallskip 

% Hence it satisfies the Poincar\'e inequality on $\R^d$:
% \begin{align}\label{4}
%   \int_{\R^d} \left| \nabla_x \left( \frac{g}{\pi} \right) \right|^2 \, \pi(dx)
%   \ge 2 \, \lambda_P \, \int g^2 \, \pi^{-1}(dx) \quad \mbox{ for } \quad
%   \int g \, dx =0,
% \end{align}
% for some constant $\lambda_P >0$. 

% Let us denote by $\gamma(x,v)$ the probability measure on $\R^d \times \R^d$ defined as
% %$$
% %\mu(v) := {1 \over (2\pi)^{d/2}} \, e^{-|v|^2/2},
% %$$
% %and 
% $$
% \gamma(x,v) := {1 \over (2\pi)^{d/2}} \, e^{-|v|^2/2 -
%   \Phi(x)} = \pi \, \mu 
% $$
% where 
% \[
% \mu(v) := {1 \over (2\pi)^{d/2}} \, e^{-|v|^2/2}.
% \]

Let us denote the probability measure $\mu(x,v)=
e^{-\phi(v)}$. Let us consider the functional space
$$
L^2(\mu^{-1/2}) := \left\{ f\in L^2(\T^d \times \R^d) \ \Big| \ \int_{\T^d
    \times \R^d} f^2 \,
  \mu^{-1} \dd x \dd v < +\infty \right\},
$$
equipped with its norm 
\[
\| f \|_{L^2(\mu^{-1/2})} := \left( \int_{\T^d \times \R^d} f^2 \, 
  \mu^{-1} \dd x \dd v \right)^{1/2}.
\]
It is immediate to check that $L(\mu)=0$ and 
\begin{align*}
  \Re e \, \langle L f, f \rangle_{L^2(\mu^{-1/2})} := &\int_{\T^d
    \times \R^d} L \bar f \, f \, \mu^{-1} \dd x \dd v +
  \int_{\T^d \times \R^d} L f \, \bar f
  \, \mu^{-1} \dd x \dd v \\
  = & \ \Re e \, \langle Lf, f \rangle_{L^2(\mu^{-1/2})} = -
  \int_{\T^d \times \R^d} \left|\nabla_v \left(\frac{f}{\mu} \right)
  \right|^2 \, \mu \dd x \dd v \le 0.
\end{align*}

We also similarly define the weighted Sobolev spaces 
$$
H^s(\mu^{-1/2}) := \left\{ f \in H^s_{\mbox{{\scriptsize loc}}} \left(\T^d
  \times \R^d\right) \ \Big| \ \forall \, |j| \le s, \ \int_{\T^d \times
    \R^d} (\partial^j f)^2 \, \mu^{-1} \dd x \dd v < +\infty \right\},
$$
for $s \in \N$ and $j \in \N^d$ multi-index (with
$|j|=j_1+\dots+j_d$), equipped with its norm
\[
\| f \|_{H^s(\mu^{-1/2})} := \left( \sum_{|j| \le s} \int_{\T^d
    \times \R^d}
  (\partial^j f)^2 \mu^{-1} \dd x \dd v \right)^{1/2}.
\]

Let us first prove an hypocoercivity result on the kinetic
Fokker-Planck equation in the torus. The proof is a variation of the
method developed in the recent works \cite{DMScras,DMS}, partly
inspired from the paper \cite{MR2215889}. In \cite{DMS} the kinetic
equation is studied in the whole space with confining potential. This
result is also related to the works \cite{MR2034753} and
\cite{MR2562709} on the kinetic Fokker-Planck equation in the whole
space with a confining potential.

\begin{theo}\label{theo:KFP1} 
  Assume that $\phi$ satisfies {\bf (FP1)-(FP2)}. Then for any initial
  datum $\initem{f} \in L^2(\mu^{-1/2})$, the solution to the initial
  value problem (\ref{eq:KFP}) satisfies
\begin{equation*}%\label{spectralGFP1}
  \forall \, t \ge 0, \quad 
  \left\| f_t -  \mu \, \langle \langle \initem{f} \rangle \rangle  \right\|_{L^2(\mu^{-1/2})} \le 
  C \, e^{ -\lambda_{KFP} \, t} \, 
  \left\| \initem{f} -  \mu \, \langle \langle \initem{f} \rangle \rangle \right\|_{L^2(\mu^{-1/2})},
\end{equation*}
for some constructive constant $C>0$ and ``hypocoercivity'' constant
$\lambda_{KFP} >0$ depending on $\phi$, with the notation
\[
\langle \langle \initem{f} \rangle \rangle := \int_{\T^d \times \R^d}  \initem{f} \dd
x \dd v.
\]
Moreover the proof below provides a quantitative estimate from below
on the optimal decay $\lambda_{KFP}$. 
\end{theo}

\begin{rems}
  \begin{enumerate}
\item More generally for $s \in \N^*$, if $\phi$ is $C^{q+2}$ and
satisfies {\bf (FP1)-(FP2)}, then for any initial datum $\init{f}\in
H^s(\mu^{-1/2})$, the solution to the initial value problem
(\ref{eq:KFP}) satisfies
\begin{equation*}%\label{spectralGFP1}
  \forall \, t \ge 0, \quad 
  \left\| f_t -  \mu \, \langle \langle \init{f} \rangle \rangle  \right\|_{H^s(\mu^{-1/2})} \le 
  C \, e^{ -\lambda_{KFP} \, t} \, \left\| \init{f} -  \mu \, \langle \langle
    \init{f} \rangle \rangle \right\|_{H^s(\mu^{-1/2})}.
\end{equation*}
\item Note that this statement implies in particular in
  $L^2(\mu^{-1/2})$ (and in fact also in $H^s(\mu^{-1/2})$) that
\[
\Sigma\left(\LL\right) \subset \{ z\in \C \; | \; \Re e (z) \le
  -\lambda_{KFP} \}\cup \{ 0\} 
\]
and that the null space of $\LL$ is exactly $\R \mu$.
  
\item Observe that, on the contrary to the previous spatially
  homogeneous case, the optimal rate of decay $\lambda_{KFP}$ is in
  general different from the Poincar\'e constant of $\Phi$. It depends
  for instance on the size of the spatial domain.
  \end{enumerate}
\end{rems}

\begin{proof}[Proof of Theorem~\ref{theo:KFP1}]
  Without loss of generality we assume that $\langle \langle
  \init{f}\rangle \rangle =0$.  Let us denote by
  \[ {\mathcal T} := v \cdot \nabla_x, \quad \bar {\mathcal L} :=
  \nabla_v \cdot \left( \nabla_v + \phi \right)
\]
and let us introduce the projection operator 
\[
\Pi f := \left( \int_{\R^d} f \dd v \right) \, \mu
\]
and the auxiliary operator 
\[
\mathcal U := \left( \Id + (\mathcal T\Pi)^* (\mathcal T\Pi) \right)^{-1}
(\mathcal T \Pi)^*.
\]

Then one can check by elementary computations that 
\[
\Pi \mathcal T \Pi = 0 \quad \mbox{ and } \quad \mathcal U = \Pi
\mathcal U
\]
and 
\begin{multline*}
  \frac{\dd}{\dd t} \left( \frac12 \, \| f\|_{L^2(\mu^{-1/2})} ^2
    + \varepsilon \, \langle \mathcal U f, f \rangle_{L^2(\mu^{-1/2})} \right) \\
  = \left\langle \bar {\mathcal L} f, f \right\rangle_{L^2(\mu^{-1/2})} +
  \varepsilon \, \langle \mathcal U \mathcal T \Pi f, f
  \rangle_{L^2(\mu^{-1/2})} +
  \varepsilon \, \langle \mathcal U \mathcal T (\Id-\Pi) f, f \rangle_{L^2(\mu^{-1/2})} \\
  - \varepsilon \, \langle \mathcal T \mathcal U f,f
  \rangle_{L^2(\mu^{-1/2})} + \varepsilon \, \langle \mathcal U \bar
  {\mathcal L} f, f
  \rangle_{L^2(\mu^{-1/2})}
\end{multline*}
(observe that $\langle \mathcal U f, \bar {\mathcal L} f
\rangle_{L^2(\mu^{-1/2})}=0$ since $\mathcal U = \Pi \mathcal U$).

By explicit computation one can show that $\mathcal U$, $\mathcal T
\mathcal U$, $\mathcal U \mathcal T$ and $\mathcal U \bar {\mathcal
  L}$ are bounded, by using that the operators
\[
\nabla_x \, \left(1- \alpha \, \Delta_x \right)^{-1} \quad \mbox{ and
} \quad \left(1- \alpha \, \Delta_x \right)^{-1} \, \nabla_x 
\quad \mbox{ with } \quad \alpha = \int_{\R^d} |v|^2 \, \mu \dd v
\]
are bounded in $L^2_x$. This implies
\begin{multline*}
  \varepsilon \, \langle \mathcal U \mathcal T (1-\Pi) f, f
  \rangle_{L^2(\mu^{-1/2})} - \varepsilon \, \langle \mathcal T \mathcal
  U f,f \rangle_{L^2(\mu^{-1/2})} - \varepsilon \, \langle \mathcal U
  \bar {\mathcal L}  f, f \rangle_{L^2(\mu^{-1/2})} \\
  \le \lambda_P \, \| (1-\Pi)f \|^2 _{L^2(\mu^{-1/2})} + C \,
  \varepsilon^2 \, \| \Pi f \|_{L^2(\mu^{-1/2})} ^2
\end{multline*}
for some constant $C>0$.  

Finally one uses the Poincar\'e inequality on the velocity variable 
\[
- \left\langle \bar {\mathcal L} f, f \right\rangle_{L^2(\mu^{-1/2})} \le - 2 \,
\lambda_P \, \| (1-\Pi) f \|_{L^2(\mu^{-1/2})} ^2
\]
and the formula 
\[
\mathcal U \mathcal T \Pi f = \left( \left(1- \alpha \, \Delta_x
  \right)^{-1} \circ (\alpha \, \Delta_x) \rho \right) \, \mu \quad
\mbox{ where } \quad \rho = \int_{\R^d} f \dd v
\]
which implies that 
\[
\left\langle \mathcal U \mathcal T \Pi f, f \right\rangle_{L^2(\mu^{-1/2})} \le -
2 \, \lambda' \, \| \Pi f \|_{L^2(\mu^{-1/2})}
\]
(we have used here $\langle \langle f_t \rangle \rangle =0$ for all
times $t \ge 0$) with
\[
\lambda' = \frac{\alpha \, \lambda' _P}{1+ \alpha \, \lambda'_P}
\]
where $\lambda' _P >0$ is the Poincar\'e constant for the
Poincar\'e-Wirtinger inequality on the torus.  This concludes the
proof of hypocoercivity by choosing some $\varepsilon$ chosen small
enough.
\end{proof}

Let us now consider some given Borel weight function $m = m(v)>0$ on
$\R^d$
% \begin{align}\label{new_weight_II}
% m(v):= e^{-\theta(\phi(v))}, %\quad e(x,v) := |v|^2/2 + \Phi(x),
% \end{align}
% with $\theta \in C^\infty$
and the associated Banach space
$L^p(m)$, $p\in [1,2]$, equipped with the norm
$$
\| f \|_{L^p(m)} = \left( \int_{\T^d \times \R^d} |f|^p \, m^p
  \dd x \dd v \right)^{1/p}.
$$ 
% We shall distinguish again between the \emph{energy weight case}
% $\theta(z) = z$, the \emph{stretched exponential weight case}
% $\theta(z) = \kappa \, z^\beta$, $\kappa >0$, $\beta \in (0,1)$,
% and the \emph{polynomial weight case} $\theta(z) = k \, \log z$.

\smallskip

We consider again the defect weight function $\psi_{p,m}$
(see~\eqref{eq:defpsimp}) and we shall assume again {\bf
  (FP2)-(FP3)}. % on which we shall assume again (observe
% that it is a function of $v$ only):
% \begin{itemize}
% \item[{\bf (KFP2)}] The function $\psi_{p,m}(v)$ defined as satisfies
% \begin{align*}
%   \psi_{p,m}(v) \xrightarrow[|v| \to \infty]{} a_{p,m} 
%   \in [-\infty,0).
% \end{align*}
% \end{itemize}
Pairs of potential-weight functions $(\phi,m)$ satisfying these
assumptions are detailed in Proposition~\ref{prop:FP1}.  \bigskip

The main result of this section is the following theorem:

\begin{theo}\label{theo:KFP2} 
    Assume that $m$, $p \in [1,2]$, $F\in C^{2}$
  satisfy {\bf (KFP1)-(FP2)-(FP3)}. Then for any initial
  data $\initem{f} \in L^p(m)$ the corresponding solution to
  \eqref{eq:KFP} satisfies
\begin{equation*}%\label{spectralGFP1}
  \forall \, t \ge 0, \quad  \left\| f_t -  \mu \, \langle \langle
    \initem{f} \rangle \rangle  \right\|_{L^p(m^{-1})} \le 
  C \; e^{ -\lambda_{m,p} \, t} \, \left\| \initem{f} -  \mu \langle
    \langle 
    \initem{f} \rangle \rangle \right\|_{L^p(m^{-1})},
\end{equation*}
with $\lambda_{m,p} := \lambda_{KFP}$ if $\lambda_{KFP} < |a_{m,p}|$, or
$\lambda_{m,p} < |a_{m,p}|$ is as close as wanted to $|a_{m,p}|$
else. 

From Proposition \ref{prop:FP1} we deduce the same estimates on the
rates $\lambda_{m,p}$ depending on the choices of $\phi$ and $m$ as in
the spatially homogeneous case, but where $\lambda_P$ is replaced by
$\lambda_{KFP}$.
% \begin{itemize} 
% \item[(i)] {\bf Energy weight.} In the case $m= \mu$
%   and $\phi=(1+|v|^2)^{s/2}$ with $s>1$, any $p \in (1,2]$ is allowed
%   (in any dimension) and the value of the spectral gap is the
%   (optimal) hypocoercivity constant $\lambda_{KFP}$ (since $a_{p,m} =
%   -\infty$).
% \item[(ii)] {\bf Stretched exponential weight.} In the case $m =
%   e^{-C_\theta \phi^\beta}$ with $C_\theta >0$, $\beta \in(0,1)$ and
%   $\phi=(1+|v|^2)^{s/2}$ with
% \[
% s\beta + s \ge 2,
% \]
% any $p\in [1,2]$ is admissible (in any dimension). The spectral gap is
% again the (optimal) hypocoercivity constant $\lambda_{KFP}$ when $s
% \beta + s >2$ (since $a_{p,m} = -\infty$), however when $s\beta +2
% =2$, our estimate on the rate of decay 
% \[
% \lambda_{p,m} = \min \left\{ \lambda_{KFP};\, C_\theta s^2 \beta/p
% \right\}. 
% \]

% \item[(iii)] {\bf Polynomial weight.} In the last case
%   $m = \phi^{-k}$ and $\phi=(1+|v|^2)^{s/2}$ with $s\ge2$, the choice
%   of $p$, $k$, $s$ and the dimension $d$ are connected by the relation
% \[
% (s-2+d)(p-1) < sk. %\ \mbox{ (see Figure~\ref{fig1}).}
% \]
% Our estimate on the spectral gap is the optimal rate $\lambda_{KFP}$
% when $s>2$, and when $s=2$, it is given by 
% \[
% \lambda_{m,p} := \min \left\{ \lambda_{KFP}; \, 2 (2k-d(p-1))
% \right\}.
% \]
% \end{itemize}
\end{theo}

\begin{rem}
Note that this statement implies in particular in
  $L^p(m)$ that
\[
\Sigma\left(\LL\right) \subset \{ z\in \C \; | \; \Re e (z) \le
  - \lambda_{m,p} \}\cup \{ 0\} 
\]
and the null space of $\LL$ is exactly $\R \mu$. All the other remarks
after Theorem~\ref{theo:FP1} and Proposition~\ref{prop:FP1} extend as
well (in particular the remark on measure solutions). However the open
questions raised in these remarks are probably harder in this
spatially inhomogeneous setting.
\end{rem}

Before going into the proof of Theorem~\ref{theo:KFP2}, let us again
prove a lemma about the regularization properties of the kinetic
Fokker-Planck equation at hand. This result is related to the notion
of \emph{hypoellipticity}, it is folklore but hard to find so we
include a sketch of proof (following closely the methods and
discussions in \cite[Section~A.21]{MR2562709}) for clarity and in
order to make explicit the estimate.

\begin{lem}\label{lem:reg-kfp}
  Under the assumptions {\bf (KFP1)-(FP2)} the semigroup of the
  equation~\eqref{eq:KFP} is well-defined in the space
  $L^1(\mu^{-1/2})$ and satisfies
\[
% \left\{ 
% \begin{array}{l}\displaystyle
\| S_{\LL}(t) f \|_{L^2(\mu^{-1/2})} \le \frac{C_L}{t^{\zeta}} \, \|
f\|_{L^1(\mu^{-1/2})} %\vspace{0.3cm} % \\ \displaystyle
% \| S_{\LL}(t) f \|_{L^\infty(\mu^{-1/2})} \le \frac{C_L}{t^{\zeta'}} \, \|
% f\|_{L^2(\mu^{-1})}
% \end{array}
% \right.
\]
for some constant $\zeta >0$. 
\end{lem}

\begin{proof}[Proof of Lemma~\ref{lem:reg-kfp}] 
The estimate 
\[
\frac{\dd}{\dd t} \int_{\T^d \times \R^d} f \, \mu^{-1/2} \dd x \dd v
= \int_{\T^d \times \R^d} f \, W_\phi \, \mu^{-1/2} \dd x \dd v \le C
\, \int_{\T^d \times \R^d} f \, \mu^{-1/2} \dd x \dd v
\]
easily ensures that the semigroup is well-defined in
$L^1(\mu^{-1/2})$. 

We rewrite the equation on $h = f/\sqrt \mu \in L^2$ (the unweighted
Lebesgue space) and we consider the functional
\[
\HH(t) := \| h \|^2_{L^2} + a^2 \, \left\| \nabla_x h
\right\|^2_{L^2} + 2b \, \left\langle \nabla_x ( D_x^{1/3} h),
\nabla_v ( D_x ^{1/3} h) \right\rangle_{L^2} + c^2 \, \left\|
  \nabla^3 _v h \right\|^2_{L^2}
\]
for some constants $a,b,c \in \R$, where $D_x :=
(1-\Delta_x)^{1/2}$. Since
\begin{multline*}
  \left\langle \nabla_x ( D_x^{1/3} h), \nabla_v ( D_x ^{1/3} h)
  \right\rangle_{L^2} = \left\langle \nabla_x h, \nabla_v (
    D_x ^{2/3} h) \right\rangle_{L^2} \\
  \le \frac\alpha2 \, \left\| \nabla_x h \right\|^2 _{L^2} +
  \frac2\alpha \, \left\| \nabla_v ( D_x
    ^{2/3} h) \right\|^2 _{L^2} \\
\le  \alpha \, \left\| \nabla_x h \right\|^2 _{L^2} + \frac{2}{\alpha} \, \left\|
  \nabla^3 _v h \right\|^2_{L^2} 
\end{multline*}
for any $\alpha >0$, it is clear that $\HH$ is equivalent to 
\[
\| h \|^2_{L^2} + \left\| \nabla_x h \right\|^2_{L^2} + \left\|
  \nabla^3 _v h \right\|^2_{L^2}
\]
as soon as $c << ab$. 

Then computations lead to
\[
\frac{\dd}{\dd t} \HH(t) \le - K \, \left( \| h \|^2_{H^{1/3}} + \left\|
    \nabla_x h \right\|^2_{H^{1/3}} + \left\| \nabla^3 _v h
  \right\|^2_{H^{1/3}} \right)
\]
for some constant $K>0$ by using the Poincar\'e inequality
\eqref{eq:Poincare-fort} in the velocity variable, the regularity
assumption \eqref{eq:reg-kfp} in {\bf (KFP1)} and the mixed-term
estimate
\[
\frac{\dd}{\dd t} \left\langle \nabla_x ( D_x^{1/3} h),
\nabla_v ( D_x ^{1/3} h) \right\rangle_{L^2} = - \left\| \nabla_x
D_x^{1/3} h \right\|^2 _{L^2} + \ \mbox{error terms\dots}.
\]

Then by interpolation with the $L^1$ norm of $h$ we deduce that 
\[
\frac{\dd}{\dd t} \HH(t) \le - K \, \frac{\HH(t)^{1+\beta}}{\| h \|_{L^1}
  ^{2\beta}}
\]
which concludes the proof of the first inequality. % The second one is
% obtained by duality as in the homogeneous case (only the transport
% operator is changed by a minus sign). 
\end{proof}

\begin{proof}[Proof of Theorem~\ref{theo:KFP2}.]
The proof is similar to that of Theorem \ref{theo:FP1}.

We first write
$$
\LL = \bar \LL + \TT, \quad \bar \LL := \nabla_v
\cdot \left(\nabla_v f + \phi \, f \right), \quad \TT := -v
\cdot \nabla_x f
$$
in $L^p(m)$ and the corresponding splitting $\LL = \bar \LL + \TT$ in
$L^2(\mu^{-1/2})$. The operator $\bar \LL$ is symmetric in $L^2(\mu^{-1/2})$
since
$$
\left\langle \bar \LL f, g \right\rangle_{L^2(\mu^{-1/2})} = - \int_{\T^d
  \times \R^d}
\nabla_v\left( \frac{f}{\mu} \right) \cdot \nabla_v \left(
  \frac{g}{\mu} \right) \, \mu \dd x \dd v.
$$
The operator $\TT$ is skew-symmetric both in $L^2(\mu^{-1/2})$
and $L^p(m)$.

Then we define the decomposition $\LL= \AA + \BB$ with
$$
\AA f := \chi_R \, M \, f \ \mbox{ and } \ \BB f:= \LL f - \chi_R\, M
\, f
$$
and $\chi_R = \chi_R(v)$ is the characteristics function of $v \in
B(0,R)$. The rest of the proof is strictly similar to that of
Theorem~\ref{theo:FP1}. 
\end{proof}

% \subsection{Comparison with existing results}
% \label{sec:comp-with-exist}

% Here is an overview of the rate of decay that we obtain
% for the admissible weights of Proposition \ref{prop:admissible}
% relative to a potential $\phi = (1+|v|^2)^{s/2}$:
% \bigskip

\subsection{Summary of the results}
\label{sec:summary-results}

Let us conclude this section with a summary of the results we have
established, both for the Fokker-Planck equation \eqref{eq:FP} or the
kinetic Fokker-Planck equation \eqref{eq:KFP} in the torus with
velocity potential $\phi(v) \approx \langle v \rangle^\gamma$ at
infinity. The constant $\lambda_*>0$ denotes either $\lambda_P$ for
the Fokker-Planck equation, or $\lambda_{KFP}$ for the kinetic
Fokker-Planck equation in the torus.  \bigskip

%\begin{center}
\hspace{-2.2cm}
    \begin{tabular}{ | c ||  c | c | c |}%p{6cm} |}
   %  \begin{tabular}{|c|c|c|}
%\hline
%\multicolumn{3}{|c|}{Team sheet} \\
      \hline  
      {\bf{ Weight }} & admissible $p$  & admissible $\gamma$ & spectral gap $\lambda$ \\ \hline\hline
       %$m = \mu$ & $1<p\le 2$ &  $s\ge 1$ &  $\lambda_P$ (optimal) \\
       %\hline
      $m = e^{\phi/2}$ & $p=2$ &  $\gamma \ge 1$ &
      $\lambda_*$ (optimal) \\  \hline
      $m = e^{\kappa \, \phi}$, $\kappa \in (0,1/2)$ & $1 \le p <2$ &
      $\gamma \ge 1$ &
      $\min\left\{ \lambda_*; \, \kappa (1-p\kappa) +0\right\}$ \\  \hline
      $m = e^{\kappa \, |v|^\beta}$, $\kappa,\beta>0$
      & $1 \le p \le 2$ &  $2-\gamma < \beta <\gamma$ &
      $\lambda_*$ (optimal) \\  \hline
      $m = e^{\kappa \, |v|^\beta}$, $\kappa,\beta>0$
      & $1 \le p \le 2$ &  $\beta+\gamma=2$ &
      $\min\left\{ \lambda_*; \, \kappa \beta \gamma +0 \right\}$ \\  \hline
       % $m = \phi^{-k}$ & $p=2$ &  $s=2$, $2k>d$ & $\min\{ \lambda_P,
       % (2k-d)/2\}$ (optimal) \\ \hline
       $m = \langle v \rangle^{k}$, $k>d(1-\frac1p)$ &  $1\le p \le2$ &  $\gamma>2$&
       $\lambda_*$ (optimal) \\ \hline
       $m = \langle v \rangle^{k}$, $k>(\gamma-2+d)(1-\frac1p)$  &  $1
       \le p \le 2$ & $\gamma=2$ 
       &  $\min\left\{ \lambda_* ; \, 2k - 2d
         (1-\frac1p)+0\right\}$ \\ \hline
      \end{tabular}
%\end{center}
\bigskip

The optimality of the estimates in the 2d, 4th and 6th line
is open. 
% Here discuss Gallay-Wayne maybe in more details, and the result of
% Villani (hypocoercivity memoir) based on entropic hypocoercivity. 

% Maybe gather here also all the open questions raised in the remarks. 

%%%%%%%%%%%%%%%%%%%% Boltzmann %%%%%%%%%%%%%%%%%%%%%%%%%%%%%%
\section{The linearized Boltzmann equation}
\label{sec:boltzmannLin}
\setcounter{equation}{0}
\setcounter{theo}{0}

%\Red ATTENTION AUX NOTATIONS: $\lambda \leftrightarrow a$ pour \^etre en coherence avec la section 2? \Black

%\Blue 
%$$
%- \lambda \ \leftrightarrow \ a ? 
%$$
%\Black

Consider the Boltzmann equation for hard spheres in the torus in
dimension $d=3$, which writes
\begin{equation}
\label{eq:NLBE}
\partial_t f = Q(f,f) - v\cdot \nabla_x f,
\end{equation}
for $f=f_t(x,v) \ge 0$, $x \in \T^3$ ($3$-dimensional flat torus), $v
\in \R^3$, and where the collision operator $Q$ is defined as
\begin{equation}
\label{eq:coll-op}
  Q(f,g)  :=  
  \int_{\R^3} \int_{\mathbb{S}^{2}}
  \Big[f(v') \, g(v'_*) - f(v) \, g(v_*) \Big] % \, b(\cos \theta) 
\, |v-v_*| \dd v_* \dd \sigma. 
\end{equation}
In \eqref{eq:coll-op} and below, we use the notations  
\begin{equation}
\label{eq:defcoll}
v' = \frac{v+v_*}2 + \sigma \, \frac{|v-v_*|}2, \qquad v'_* =
\frac{v+v_*}2 - \sigma \, \frac{|v-v_*|}2,
\end{equation}
with $\cos \theta = \sigma \cdot (v-v_*)/|v-v_*|$. % and $\gamma \in
% (0,1]$,
%and where $b$ is function on $[-1,1]$ such that 
% \begin{equation}\label{eq:hypb}
% \forall \, z \in [-1,1], \quad 0< b_0 \le b(z) \le b_1 <\infty \quad
% \mbox{ and } \quad \int_{\mathbb{S}^{d-1}} b(\cos \theta) \dd \sigma =1,
% \end{equation}
% which includes the case of  the {\em hard spheres} model  in dimension $d=3$ where
% $\gamma=1$ and  $b$ is constant.
%\Blue Because we consider the Boltzmann equation for hard spheres, we may assume without loss
%of generality that $b(\cos\theta) = 1/(4\pi)$ so that $\| b \|_{L^1} = 1$. \Black
We assume without loss of generality that the torus has volume
one. Then global equilibria are absolute Maxwell functions which
depend neither on time nor on position (see \cite[Chap. II,
sect. 7]{Ce88} for instance). By normalization of the mass, momentum
and energy, we consider the following equilibrium
\begin{equation}\label{def:secBL:mu}
\mu(v) := {1 \over (2\pi)^{3/2}} \, e^{-|v|^2/2}. 
\end{equation}

Consider the linearization $f = \mu + h$, then at first order the
linearized equation around the equilibrium is
\begin{equation}\label{eq:LIBE}
\partial_t h =\LL h : = \bar \LL h  - v\cdot \nabla_x h,
\end{equation}
for $h=h(t,x,v)=h_t(x,v)$, $x \in \T^3$, $v \in \R^3$ and
\begin{equation*}
  \bar \LL h  := \int_{\R^3} \int_{\mathbb{S}^{2}} \Big[\mu(v'_*) \, h(v')
  + \mu(v') \, h(v'_*) - \mu(v_*) \, h(v) - \mu(v) \, h(v_*) \Big]
  % \, b(\cos \theta) 
\, |v-v_*| \dd v_* \dd \sigma.
\end{equation*}

Following standard notations, we introduce the \emph{collision frequency} 
\[
\nu(v) := 4 \pi \, \int_{\R^3} \mu(v_*) \, |v-v_*| \dd v_* = 4 \pi \, \left( \mu \ast |
\cdot | \right) (v)
\]
which satisfies for some constants $\nu_0, \nu_1 > 0$ %(beacuse $0<\gamma\le 1$) 
\[
\forall \, v \in \R^3, \quad 0 < \nu_0 \le \nu_0 \, (1+|v|) \le
\nu(v) \le \nu_1 \, (1+|v|),
\]

\begin{rem}
The collision frequency satisfies in fact the explicit bounds
\[
\forall \, v \in \R^3, \quad 4\pi \max\left\{ |v|,
  \sqrt{2/(e\pi)}\right\} \le \nu(v) \le 4\pi (|v| + 2)
\]
that we shall use in the sequel. Indeed, on the one hand, the lower
bound follows from the Jensen inequality
$$ 
\nu(v) \ge 4\pi \left| \int_{\R^3} (v-v_*) \, \mu(v_*) \, \dd v_* \right| = 4\pi \, |v| 
$$
and 
\begin{eqnarray*}
(4\pi)^{-1} \nu (v) 
&\ge& 
 \int_{|v_*-v| \ge 1}  \mu(v_*) \, \dd v_*  
\\
&\ge& 
 \int_{|v_*| \ge 1}  \mu(v_*) \, \dd v_*  
\\
&\ge& \sqrt{\frac{2}{\pi}} \int_1^\infty e^{-r_*^2/2} \, r^2_* \, \dd r_*
\\
&\ge& \sqrt{\frac{2}{\pi}} \int_1^\infty e^{-r_*^2/2} \, r_* \, \dd
r_* = \sqrt{\frac{2}{e \pi}}.
% \\
% &\ge& {1 \over (\pi)^{1/2}} \int_1^\infty e^{-u} \, \dd u = {1 \over e \, (\pi)^{1/2}}. 
\end{eqnarray*}
One the other hand, we have 
\bean
(4\pi)^{-1} \nu (v) 
&\le& 
 \int_{\R^3}  |v| \mu(v_*) \, \dd v_*  + \int_{\R^3}  
\left( \frac12 + \frac{|v_*|^2}2 \right) \, \mu(v_*) \, \dd v_*  
\\
 &=& |v| + 2. 
\eean
\end{rem}

\normalsize

\subsection{Review of the decay results on the semigroup}
\label{sec:review-decay-results}

Let us briefly review the existing results concerning the decay
estimates on the semigroup of $\LL$ for hard spheres in the torus.  

In the spatially homogeneous case, the study of the linearized
collision operator $\bar \LL$ goes back to Hilbert
\cite{MR1511713,MR0056184} who computed the collisional invariant, the
linearized operator and its kernel in the hard spheres case, and
showed the boundedness and complete continuity of the non-local part
of $\bar \LL$. This operator is self-adjoint non-positive and
generates a strongly continuous contraction semigroup in the space
$L^2_v (\mu^{-1/2})$. Carleman~\cite{Carleman} then proved the existence of
a spectral gap for $\bar \LL$ by using Weyl's theorem and the
compactness of the non-local part of $\bar \LL$ proved by
Hilbert. Grad~\cite{Grad1,Grad2} then extended these results to the
so-called ``hard potentials with cutoff''.  All these results are
based on non-constructive arguments. The first constructive estimates
in the hard spheres case were obtained only recently in
\cite{Baranger-Mouhot} (see also \cite{MR2254617} for more general
interactions). Note that these spectral gap estimates can easily be
extended to the spaces $H^s _v(\mu^{-1/2})$, $s \in \N^*$, by
reasoning as in the proof of Lemma~\ref{lem:decomp-new} below when we
introduce derivatives.

Let us also mention the works
\cite{WCUh:LBE:70,MR0398379,Boby:maxw:88} for the different setting
of \emph{Maxwell molecules} where the eigenbasis and eigenvalues are
explicitely computed by Fourier transform methods. Although these
techniques do not apply here, the explicit formula computed are an
important source of inspiration for dealing with more general physical
models.
% as a consequence of
% \cite[Proposition 5.6]{MR2264623} and the approach introduced in
% \cite{Mcmp}. % (see also \cite{Grad1,Grad2}).

The complete linearized operator $\LL$ is the sum of the self-adjoint
non-positive operator $\bar \LL$ and the skew-symmetric transport
operator $- v \cdot \nabla_x$. It was first established in
\cite[Theorem 1.1]{MR0363332} that it has a spectral gap in the
Hilbert space $L^2_v H^s_x(\mu^{-1/2})$, $s \in \N$, by
non-constructive arguments. Then using an argument initially due to
Grad \cite{MR0184507} for constructing local-in-time solutions Ukai
\cite{MR0363332}, showed that the spectral property also holds in
$L^\infty_v H^s_x((1+|v|)^k \mu^{-1/2})$, $k >3/2$. In \cite[Theorems
1.1 \& 3.1]{MNeu}, quantitative spectral gap estimates are established
in $H^s_{v,x}(\mu^{-1/2})$, $s \in \N^*$, following partly ideas from
\cite{MR1946444,MR2013332,MR2095473,MR2562709}.

% It is likely that the theory developed in \cite{DMS}
% provides a constructive proof of \cite[Theorem 1.1]{MR0363332} in the
% space $L^2_{x,v}(\mu^{-1/2})$, handling the case $s = 0$ which is not
% included in \cite{MNeu}, and it is therefore an alternative approach
% to the one developed here.
  
For the spatially homogeneous case, in \cite{MR946973} the decay
estimate of $e^{t \bar \LL}$ was extended to $L^1$ with polynomial
weight by an intricate non-constructive approach: the decay bound on
the resolvent is deduced from the spectrum localization with no
constructive estimate, and then the decay of the semigroup is obtained
by some decomposition of the solution. This argument was then extended
to $L^p$ spaces in \cite{MR1215007,MR1338453}. In \cite{Mcmp}, this
decay estimate was extended to the space $L^1(m)$ for a stretched
exponential weight $m$, by constructive means, with optimal rate. Let
us also mention that in \cite{MR900501} some non-constructive decay
estimates were obtained in a Sobolev space in position combined with a
polynomially weighted $L^\infty$ space in velocity (integrating first
in $x$ and then taking the supremum in $v$, which is reminiscent of
the norms we shall use in the sequel). We also refer to the book
\cite{MR1612403} by M.  Mokhtar-Kharroubi and the more recent paper
\cite{MR2216092} for an overview of the spectral analysis and the
semigroup growth estimate available for the linear Boltzmann equation
as it appears in neutron transport.

\subsection{The main hypodissipativity results}
For some given Borel weight function $m >0$ on $\R^3$, let us define
$L^q _v L^p_x(m)$, $1 \le p,q \le \infty$, as the Lebesgue space
associated to the norm
\[
\| h \|_{L^q _v L^p_x(m)} =:  \| \| h (\cdot,v) \|_{L^p_x} \, m(v)  \|_{L^q_v}. 
\]
%with the following generalizations when $p=+\infty$ and $q \in
%[1,+\infty)$: 
%\[
%\| h \|_{L^q _v L^\infty_x(m)} =:  \left( \int_{\R^d} \left( \sup_{x
%      \in \T^d} |h(x,v)| \right)^{q} \, m(v) \dd v  \right)^{1/q} , 
%\]
%when $p \in [1,+\infty)$ and $q = +\infty$: 
%\[
%\| h \|_{L^\infty _v L^p_x(m)} =:  \sup_{v \in \R^d} \left( \left(
%    \int_{\T^d} |h(x,v)|^p \dd x \right)^{1/p} \, m(v) \right) , 
%\]
% and when $p = q = +\infty$
%\[
%\| h \|_{L^\infty _v L^\infty_x(m)} =: \sup_{v \in \R^d} \left( \left( \sup_{x
%      \in \T^d} |h(x,v)| \right) \, m(v) \right) .
%\]
We also consider the higher-order Sobolev subspaces $W^{\sigma,q}_v
W^{s,p}_x(m)$ for $\sigma, s \in \N$ defined by the norm
\begin{equation}\label{eq:def-norms}
\| h \|_{W^{\sigma,q}_v W^{s,p}_x(m)} := \sum_{i,j \in \N^d, \ |i|\le
  \sigma, |j| \le s, \ |i|+|j|\le \max\{\sigma;s\}} \left\| \left\| \partial^i_v
\partial^j _x h (\cdot,v) \right\|_{L^p_x}\right\|_{L^q_v(m)}.
\end{equation}
This definition reduces to the usual weighted Sobolev space $W^{s,p}
_{x,v}(m)$ when $q=p$ and $\sigma=s$, and we recall the shorthand
notation $H^s_\cdot :=W^{s,2}_\cdot$.

\medskip We present now our set of hypodissipativity results for the
semigroup associated to the linearized Boltzmann
equation~\eqref{eq:LIBE}.

\begin{theo}\label{theo:LIBE2} 
  Consider the space $\EE = W^{\sigma,q}_v W^{s,p}_x (m)$ with
  $s,\sigma \in \N$, $\sigma \le s$, and with one of the following
  choices of weight and Lebesgue exponents:
\begin{itemize}
\item[\textbf{(W1)}] $m=\mu^{-1/2}$, $q=p=2$; 
\item[\textbf{(W2)}] $m = e^{\kappa \, |v|^\beta}$,  $\kappa > 0$, $\beta \in (0,2)$
  and $p,q \in [1,+\infty]$; 
\item[\textbf{(W3)}] $m = \langle v \rangle^k$, $k > k_q^* $ and
  $p,q \in [1,+\infty]$, where
  \[
  k_q ^* := \frac{3 + \sqrt{49 - 48/q}}{2}. 
  \]
\end{itemize}

Then there are constructive constants $C\ge 1$, $\lambda>0$, such that
the operator $\LL$ defined in \eqref{eq:LIBE} satisfes in $\EE$:
\[
\left\{ 
\begin{array}{l}
\Sigma(\LL) \subset \left\{ z\in \C \; | \; \Re e\, (z) \le -\lambda \}\cup \{ 0\right\}
\vspace{0.2cm} \\
N(\LL) = \mbox{{\em Span}}\left\{ \mu, \, v_1 \, \mu,\, \dots, \, v_d \, \mu,\,
|v|^2 \, \mu \right\},
\end{array}
\right.
\]
and is the generator of a strongly continuous semigroup 
\[
h_t := S_\LL(t) \initem{h} %= e^{t \, \LL} \, \initem{h} 
\ \mbox{ in } \ \EE,
\]
solution to the initial value problem \eqref{eq:LIBE}, which
satisfies:
\begin{equation*}
  \forall \, t \ge 0, \quad 
  \left\| h_t -  \Pi \initem{h} \right\|_{\EE} \le 
  C \, e^{ -\lambda \, t} \, \left\| \initem{h} - \Pi \initem{h} \right\|_{\EE},
\end{equation*}
where $\Pi \initem{h}$ stands for the projection onto $N(\LL)$ defined by
\eqref{def:SpectralProjection}, or more explicitly by 
\begin{equation}
\label{eq:project-explicit}
\left\{ 
\begin{array}{l} \ds
\Pi g := \sum_{i=0}^{4} \left( \int_{\T^3 \times \R^3} g \, \varphi_i
  \dd x \dd v \right) \, \varphi_i \, \mu, \vs \\ \ds
\varphi_0 = 1, \,\, \varphi_i = v_i, \, 1 \le i \le 3, \,\,  
\varphi_{4} = \frac{\left(|v|^2- 3\right)}{18}.
\end{array}
\right.
\end{equation}

Moreover $\lambda$ can be taken equal to the spectral gap of $\LL$ in
$H^s(\mu^{-1/2})$ (with $s \in \N$ as large as wanted) in the cases
{\bf (W1)}-{\bf (W2)}. This is still true in the case {\bf (W3)} when
$k$ is big enough (with constructive threshold).
\end{theo}

\begin{rems}\label{rem:theo:LIBE2}
 \begin{enumerate}   
 \item An important aspect of this decay result is that the rate
   $\lambda$ is equal to the spectral gap in the smaller
   space $H^s(\mu^{-1/2})$. This is an {\em optimal timescale}.  For
   weights of the form {\bf (W3)} such optimality requires $k$ large
   enough.
 
%the obtaining of 
   %{\em optimal timescales}, i.e. the rate $\lambda$ is as close as
   %wanted to the spectral gap in the smaller space $H^s(\mu^{-1/2})$
   %(in the case {\bf (W3)} this requires that $k$ be large enough). 

   % We note
   % that it is likely that the ``$+0$'' loss on the rate could be
   % removed by refining the expansion of the semigroup to the second
   % smallest non-zero eigenvalue.
   % \Blue Traiter cela ??? \Black
  % In the case $L^1(1+|v|^k)$, $k >2$ on the opposite, the
  % loss on the rate as compared to Theorem~\ref{theo:LIBE1} can be made
  % as small as wanted only when $k$ is large enough (see the $L^1$
  % estimate in Lemma~\ref{lem:decomp-new}).
%\item \textcolor{red}{Possible generalization to cutoff hard
%    potentials to be added as a remark (only problem is for $L^\infty$
%  estimates, but it is nicer to give a unified statement).}
  
 \item Another important point of Theorem \ref{theo:LIBE2} is the
   spectral analysis of the linearized Boltzmann equation in Lebesgue
   spaces associated to a {\em polynomial} weight function. Apart from
   the non-constructive works \cite{MR946973,MR900501}, all the
   previous works were considering spaces with Gaussian decay in
   velocity dictated by the equilibrium $\mu$, or more recently
   stretched exponential weights in \cite{Mcmp,MMcmp,MR2476685}. We
   also refer to \cite{MR2832638,MR2821681} where polynomial weights
   are considered for a fragmentation equation.

%\item \textcolor{red}{Last paragraph added. Open question to be put
%  somewhere on the optimal spectral gap for lower polynomial weights.}

 \item Observe that we could replace $k_q^*$ by the slightly better
   exponent $k^{**}_q \le k_q ^*$ defined as the solution to the
   equation $\phi_q(k^{**}_q) = 1$ with
$$
  \phi_q(k) := \Bigl( {4 \over k+2} \Bigr)^{1/q}  \Bigl( {4 \over k-1}  \Bigr)^{1-1/q}.
$$
This last condition comes from a careful application of the
Riesz-Thorin interpolation inequality, as will be seen in the
proof. 

\item Observe that the thresholds $k^*_q, k^{**} _q$ (related to the
  decomposition of the operator) are $k_1^*= k_1^{**} =2$ in the case
  $q=1$ and $k_\infty^* = k_\infty ^{**} = 5$ in the case
  $q=+\infty$. It is remarkable that on both cases these numbers
  correspond to the threshold for the energy to be finite. For $q \in
  (1,+\infty)$ the asymptotic velocity decay suggested byx the
  finiteness of the energy is $k_q ^{***} = 5 - 3/q$ and our thesholds
  exponents $k_q^* \ge k_q ^{**} > k_q ^{***}= 5 - 3/q$ are close to
  it. There is a further loss $1-1/q$ on the threshold for the
  spectral gap (due to the fact that the reminder estimates in the
  decomposition are applied with the negative weight $\nu^{-1/q'}$,
  see later in the proofs), which leads to the conditions $k>2$ when
  $q=1$ and $k>6$ when $q=+\infty$. The optimality of these conditions
  is an open question suggested by our study.

\item As for the Fokker-Planck equation in the previous section, we
  observe a threshold condition on the polynomial degree to recover
  the optimal spectral: the weaker the growth of the weight function
  is, the more the semigroup ``ignores'' some discrete eigenvalues in
  the sense of having time decay worse than these eigenvalues, with
  eventually a time decay worse than the spectral gap and degenerating
  to zero. This suggests a ``tide'' phenomenon for the continuous
  spectrum, i.e. that depends on this weight and moves
  towards zero as the weight is weakened and approach to the critical
  ``energy space'' $L^1 _2$ in velocity. Let us also mention that
  interestingly such a phenomenon has also been observed by Bobylev in
  \cite{MR0398379} for the linearized spatially homogeneous Boltzmann
  equation associated to Maxwell molecules. In this case an explicit
  calculation (by mean of Fourier transform analysis) can be
  performed.

\item We note that even if our main goal here is to relax the tail
  decay condition on the solution, our general method is also useful
  for relaxing the \emph{regularity condition} on the solution. As a
  side result, it hence provides an alternative strategy to
  \cite{MR2215889,DMS} in order to study the linearized semigroup
  without regularity assumptions in various hypocoercive contexts. We
  refer to \cite{KCW1,KCW2} where some aspects of these works are
  revisited in this spirit, with in particular a crucial use of our
  \emph{iterated averaging lemma} (see below). In this paper we will
  give some applications of this regularity side of our method in
  order to understand the structure of propagation of the
  singularities for the Boltzmann equation.
\end{enumerate}
\end{rems}

% \subsection{Preliminary tools}
% \label{sec:preliminary-tools}

%\cite[Propositions~2.1-2.3-2.4-2.5-2.6]{Mcmp} 

\subsection{Strategy of the proof}
\label{sec:intr-trun}

\subsubsection{Methodology} The strategy is inspired from the
methodological approach in \cite[Theorem 4.2]{Mcmp}; it crucially uses
the abstract enlargement Theorem~\ref{theo:EnlargingSGdecay}. The
starting point is the quantitative hypocoercivity theorem in a small
Hilbert space setting from \cite{MNeu}, and we use % , see
% Remark~\ref{rem:theo:LIBE2} and on a convenient
a decomposition of $\LL$ found in~\cite{Mcmp}. We fundamentally
extend \cite[Theorem 4.2]{Mcmp} in several aspects: (1) we include
spatial dependency in the torus, (2) we enlarge to $L^1$ spaces with
{\em polynomial weights}, and (3) we enlarge to $L^\infty$ spaces with
polynomial or exponential weights. Extensions (2) and (3) result from
new estimates on the remaining operator $\BB^2_\delta$ in
$L^p_v(m)$, see Lemma~\ref{lem:decomp-new} below, while extension (1)
also takes advantage of the new abstract extension
Theorem~\ref{theo:EnlargingSGdecay} and a new result of smoothness for
iterated velocity averages for solutions to kinetic equations, see
Lemma~\ref{lem:regularization-lbe}.  % As a consequence, we dramatically
% improve the coercivity and hypocoercivity theory for the linearized
% Boltzmann equation obtained in
% \cite{Carleman,Grad1,Grad2,Baranger-Mouhot,Mcmp,MR0363332,MNeu}.

\subsubsection{Steps of the proof}
Consider a decomposition of the operator 
\[
\LL = \AA + \BB \quad \mbox{where} \quad \AA=\AA_\delta \quad
\mbox{and} \quad \BB = \BB^1 + \BB^2_\delta
\]
are suitable operators which are defined through an appropriate
mollification-truncation process, described later
on. %(cf. \cite{MR2081030,Mcmp}).
As a first step we estimate the remainder term $ \BB^2_\delta$ and
show that it is small in various norms. The estimate in $L^1(\langle v
\rangle^k)$, $k > 2$, is obtained by carefully exploiting a refined
version of the Povzner inequality. The estimate in $L^\infty(\langle v
\rangle^k)$ is obtained by using a representation of the gain term for
radially symmetric functions inspired from the physics literature,
which has been used for the Boltzmann equation for Bosons gas in
\cite{Semikov1995,Semikov1997,MR1958975,MR2342205}. As a second and
easier step we deduce that $\AA$ has smoothing effect in the
$v$-variable and that $\BB-a$ is dissipative with $a < 0$.  % As
% remarked in \cite{MR2832638,MR2821681} for the fragmentation equation
% we do not need here that the {\it ``remainder term"} $ \BB^2_\delta$
% tends to $0$ as $\delta$ tends to $0$, but only that $\BB-a$ is
% dissipative. While the first assertion seems to be false for a
% polynomial weight, the second can be proved. {\Red Ces deux assertions
%   ne sont effectivement pas equivalentes (cf. transport operator) mais
%   il me semble qu'on montre plus, i.e. qu'on montre bien en fait la
%   petitesse de $\BB^2$} 
In a third step, we prove some new regularity estimates on iterated
velocity averages of a solution to a kinetic transport equation and we
deduce some regularity estimates in both position and velocity
variables on the iterated time-convolutions of $\AA_\delta \,
S_{\BB_\delta}(t)$. The new feature of these regularity estimates is
that they hold for solutions merely $L^1$, whereas classical averaging
lemmas~\cite{MR923047} are well-known to degenerate in $L^1$. Finally,
the known spectral analysis of the linearized Boltzmann equation in
$H^1_{x,v}(\mu^{-1/2})$ proved in \cite{MR0363332,MNeu}, the space
extension theory developed in section~\ref{sec:factorization} and all
the preceding steps yield the full proof of
Theorem~\ref{theo:LIBE2}.

\subsubsection{The decomposition of the linearized operator} 
Let us first recall the usual decomposition 
\[
Q(g,f) = Q^+(g,f) - Q^-(g,f)
\]
of the bilinear collision operator with
\begin{equation}
  \label{eq:decQ}
  \left\{ 
    \begin{array}{l} \displaystyle
    Q^+(g,f) := 
    \int_{\R^3} \int_{\mathbb{S}^{2}} f(v') \, g(v'_*) % \, b(\cos \theta) 
    \, |v-v_*| \dd v_* \dd \sigma \vspace{0.3cm} \\ \displaystyle
    Q^-(g,f) := 
    \int_{\R^3} \int_{\mathbb{S}^{2}} f(v) \, g(v_*) % \, b(\cos \theta) 
    \, |v-v_*| \dd v_* \dd \sigma.
  \end{array}
\right.
\end{equation}

We introduce the decomposition of the linearized operator used in this
section. For any $\delta \in (0,1)$, we consider $\Theta_\delta =
\Theta_\delta(v,v_*,\sigma) \in C^\infty$ bounded by one on the set
\[
\left\{ |v| \le \delta^{-1} \ \mbox{ and } \ 2 \delta \le |v-v_*| \le
  \delta^{-1} \ \mbox{ and } \  | \cos \theta |
  \le 1 -2 \delta \right\}
\]
and whose support is included in 
\[
\left\{ |v| \le 2 \, \delta^{-1} \ \mbox{ and } \  \delta \le |v-v_*| \le
  2 \delta^{-1} \ \mbox{ and } \   | \cos \theta |
  \le 1 -\delta \right\}.
\]
We define the splitting
\[
\bar\LL h = \bar\AA_\delta h + \bar \BB_\delta h 
\]
with 
\begin{multline*}
 \bar\AA_\delta h(v)  % = \LL^{+,*}_{S,\delta}\,  h (v)
 := \int_{\R^d} \int_{\mathbb{S}^{d-1}} \Theta_\delta \, \Big[\mu(v'_*) \, h(v')
  + \mu(v') \, h(v'_*)  \\ - h(v_*) \, \mu(v) \Big]
  %\, b(\cos \theta) 
\, |v-v_*| \dd v_* \dd \sigma . 
% - \left(
%     \int_{\R^d} h(v_*) \, |v-v_*|^\gamma \dd v_* \right) \, \mu(v)
\end{multline*}

%NOTATIONS: 
%$$
%Q(\mu,f) + Q(f,\mu) = \LL^{+,*}_{S,\delta} (f) +  \LL^{+,*}_{rem} f - \nu f 
%\quad
%et
%\quad
%\tilde \AA^c_\delta \to \tilde  \LL^{+,*}_{rem}
%$$

Thanks to the truncation, we can use the so-called Carleman
representation (see \cite[Chapter~1, Section~4.4]{Villani-handbook})
and write the truncated operator $\bar \AA_\delta$ as an integral
operator
\begin{equation}\label{eq:def-barAAd}
 \bar\AA_\delta h(v)=  \int_{\R^d} k_\delta(v,v_*) \, h(v_*) \dd v_*
\end{equation}
for some smooth kernel $k_\delta \in C_c^\infty(\R^d \times \R^d)$.

Defining the corresponding remainder operator
\begin{multline}\label{eq:def-remain}
  \bar \BB^2_\delta h(v) % =  \LL^{+,*}_{R,\delta} h(v)
  := \int_{\R^d} \int_{\mathbb{S}^{d-1}} \left(
    1- \Theta_\delta\right) \, \Big[\mu(v'_*) \, h(v') + \mu(v') \,
  h(v'_*)
  \\
  - h(v_*) \, \mu(v) \Big] %\, b(\cos \theta) 
\, |v-v_*| \dd v_* \dd \sigma, 
\end{multline}
we have therefore the representation 
$
\bar \BB_\delta =  - \nu +  \bar\BB^2 _\delta.
%= \bar\BB^2 _\delta - \left(\mu \ast |
%  \cdot|^\gamma\right).
$
 We can then write a decomposition for the complete linearized operator $\LL$
\[
\LL = \AA_\delta + \BB_\delta 
\]
with 
\begin{equation*}
\begin{cases}
 \AA_\delta = \bar \AA_\delta % = \LL^{+,*}_{S,\delta}
\vspace{0.3cm} \\
 \BB_\delta = \BB^1  + \BB^2_\delta, \quad \BB^1  = -\nu - v \cdot \nabla_x, 
 \quad \BB^2_\delta = \bar \BB^2 _\delta. %\LL^{+,*}_{R,\delta}.
\end{cases}
\end{equation*}
% so that we expect (and shall prove) that $\AA_\delta$ has good velocity smoothing property, 
% $\BB^1 $ is explicit, dissipative and the associated semigroup enjoys  smoothing property for velocity averaging, 
% $\BB^1 $ is small enough.

We also define the nonnegative operator $ \tilde \BB_\delta ^2$
by
\begin{multline}\label{eq:def-remain-positive}
  \tilde \BB_\delta ^2 h(v) := \int_{\R^d} \int_{\mathbb{S}^{d-1}}
  \left( 1- \Theta_\delta\right) \, \Big[\mu(v'_*) \, h(v') + \mu(v')
  \, h(v'_*)
  \\
  + h(v_*) \, \mu(v) \Big] %\, b(\cos \theta)
  \, |v-v_*| \dd v_* \dd \sigma.
\end{multline}
It is obvious that $| (\BB^2 _\delta h) (v)| \le (\tilde \BB^2 _\delta
|h|) (v)$, and therefore any control in weighted Lebesgue space on
$\tilde \BB^2 _\delta$ implies a similar control on $\BB^2 _\delta$.

\subsection{Integral estimates with polynomial weight on the remainder}
\label{sec:integr-estim-coll}

Let us first prove some smallness estimates on the remainder term
$\BB^2 _\delta$ in the norm $L^1(\nu \langle v \rangle^k) \to
L^1(\langle v \rangle^k)$, as $\delta$ goes to zero. Since the
position $x$ is just a parameter for the operator
$\BB^2 _\delta$, we restrict the analysis to the velocity variable
only without loss of generality. This estimate improves on the estimate
\cite[Proposition~2.1]{Mcmp} since it handles \emph{polynomial
  weights} instead of stretched exponential weights. This dramatically
enlarges the functional space in which we can control the semigroup,
and it is also more natural from the perspective of the Cauchy problem
for the fully nonlinear equation. The cornerstone of
the proof is a careful use of a Povzner inequality \emph{with sharp
  constants}. % {\Red resultat legerement apparente dans
  % Arkeryd 1988 mais obtenu non par Povzner, mais par une decomposition
  % dyadique incomprehensible, pas la peine de l'evoquer ici
  % probablement \dots Dans Arkeryd-Esposito-Pulvirenti ils travaillent
  % dans $L^\infty$ en vitesse et prouvent une estimation proche de
  % celle que l'on fait par la suite mais moins optimale.}
%\Green (see \cite[Chapter~2,
%Section~2.2]{Villani-handbook} and the references therein). \Black

\begin{lem}\label{lem:Aintegral}
  For any $k > 2$ and $\delta \in (0,1)$, the remainder collision
  operator $\BB^2 _\delta$ defined in \eqref{eq:def-remain} % (as well
  % in fact as $\tilde \BB^2_\delta$ defined in
  % \eqref{eq:def-remain-positive})
  satisfies
\begin{equation}\label{eq:Aintegral}
\forall \, h \in L^1(\langle v \rangle^{k+1}), \quad \left\|
  \BB^2 _\delta h \right\|_{L^1(\langle v \rangle^k)} \le
\left(   { 4 \over k+2} +  \eps_k  (\delta)\right) \, \left\| h \right\|_{L^1(\nu \,
  \langle v \rangle^k)},
\end{equation}
where $\eps_k(\delta)$ is a constructive constant depending on $k$ and
approaching zero as $\delta$ goes to zero.
\end{lem}

Before going into the proof of Lemma~\eqref{eq:Aintegral} we shall
review a classical tool in the Boltzmann theory, i.e. a sharp version of
the \emph{Povzner (angular averaging) lemma}. The key estimate we
use was implicit in \cite{MR1461113}, \cite{MR1478067} or \cite[Lemma
2.2]{MR1697562}, and was made explicit with sharp constants in
\cite{MR2096050}, from which we adapt the following statement.

\begin{lem}[Sharp Povzner Lemma]\label{lem:Povznerbis}
  For any $k>2$, we have
\begin{multline*}
  \forall \, v, v_* \in \R^3, \quad \int_{\mathbb S^{2}} \Big[
  |v'_*|^k + |v'|^k - |v_*|^k - |v|^k \Big]  \,
  \dd\sigma \\ \le C_k \, \left( |v|^{k-1} \, |v_*| + |v| \, |v_*|^{k-1}
  \right) - (4\pi -\gamma_k) \,  |v|^k ,  
\end{multline*}
where $\gamma_k := 16 \pi / (k+2)$, so that in particular $\gamma_k
\to 0$ as $k \to \infty$, and $C_k > 0$ is a constant depending on
$k$.
  
\end{lem}

\begin{proof}[Proof of Lemma~\ref{lem:Povznerbis}] 
  We know from \cite[Corollary~3 and the remark that follows
  it]{MR2096050} that for any $k > 2$, it holds
  \begin{equation}
    \label{eq:povzner}
    \int_{{\mathbb S}^{2}} \left( |v'|^{k} + |v_*'|^{k}
    \right) \dd  \sigma \leq \gamma_k   
    \left(|v|^2 + |v_*|^2\right)^{k/2},  
  \end{equation}
from which we deduce that 
  \begin{multline*}
    \int_{\mathbb S^{2}} \Big[ |v'_*|^k + |v'|^k - |v_*|^k - |v|^k
    \Big] \dd \sigma \\
    \le \gamma_k \, \Big[ \left( |v_*|^2 + |v|^2 \right)^{k/2} 
    - |v_*|^k - |v|^k \Big]  - (4\pi -\gamma_k) \left( |v|^k + |v_*|^k  \right).
  \end{multline*}
We conclude the proof by using the elementary inequality 
$$
(y+z)^{k/2} - y^{k/2} - z^{k/2} \le 2^{k/2} \, ( y^{k/2-1/2} \, z^{1/2} + y^{1/2} \, z^{k/2-1/2}), 
$$
for any $y,z \ge 0$, in order to bound the first term. 
\end{proof}

Let us now go back to the proof of Lemma~\ref{lem:Aintegral}. 

\begin{proof}[Proof of Lemma~\ref{lem:Aintegral}]
  Since $\langle v \rangle^k \le (1+|v|^k) \le 2^{k/2} \, \langle v
  \rangle^k$, it is enough to prove the result with the weight $m:= 1
  + |v|^k$. We compute 
% (we do not recall below that all
%estimates are performed in the velocity space only)
\begin{equation*}
  \left\| \BB^2 _\delta h \right\|_{L^1(m)}  
  \le \int_{\R^3 \times \R^3 \times \mathbb S^{2}} \left( 1-
    \Theta_\delta\right) \, \Big[\mu'_* \, |h'| + \mu' \, |h'_*| + \mu
  \, |h_*|  \Big]  
  \, |v-v_*| \, m \dd v \dd v_* \dd \sigma.
\end{equation*}
We first crudely bound from above the truncation function as follows
\begin{multline*} 
  \left\| \BB^2 _\delta h \right\|_{L^1(m)} 
  \le 
   \int_{\{ |\cos \theta| \in [1-\delta,1]\}}   \mu_* \,
   |h|\Big[m'+ m'_*+ m_* \Big] 
  \, |v-v_*|  \dd v \dd v_* \dd \sigma
  \\
  + \int_{\{ |v-v_*| \le \delta\}}  \mu_* \, |h|\Big[
  m'+ m'_*+  m_* \Big] 
  \, |v-v_*|  \dd v \dd v_* \dd \sigma 
  \\
  + \int_{\{|v| \ge \delta^{-1} \mbox{ {\scriptsize or} }
    |v-v_*| \ge \delta^{-1} \}} \Big[\mu'_* \, |h'| + \mu' \,
  |h'_*| + \mu \, |h_*| \Big] 
    \, |v-v_*| \, m \dd v \dd v_* \dd \sigma,
\end{multline*}
where the change of variable $(v',v'_*,\sigma) \to (v,v_*,\sigma)$ has been used in the
two first integral terms, so that 
\begin{multline} \label{eq:AB1} 
  \left\| \BB^2 _\delta h
  \right\|_{L^1(m)} \\ \le 2^{k/2} \, \left( \int_{\{ |\cos \theta| \in [1-\delta,1]\}}
    d\sigma + \delta \right) \, \int_{\R^3 \times \R^3} \mu_*
  \, \langle v_* \rangle^{k+1} \, |h| \, \langle v \rangle^{k+1} \dd v
  \dd v_*
  \\
+ \int_{\R^3 \times \R^3 \times
    \mathbb S^2} \chi_{\delta^{-1}} \, \Big[\mu'_* \, |h'| + \mu' \,
  |h'_*| + \mu \, |h_*| \Big] \, |v-v_*| \,  m \dd v \dd v_* \dd \sigma
 \end{multline}
where $\chi_{\delta^{-1}}(v,v_*)$ is the characteristic function of
the set 
\[
\left\{\sqrt{|v|^2+|v_*|^2} \ge \delta^{-1} \mbox{ or }
|v-v_*| \ge \delta^{-1} \right\}.
\]

The first term in the right hand side of \eqref{eq:AB1} is easily controlled as
$\OO(\delta) \| h \|_{L^1(\nu m)}$.

In order to deal with the second term we write
 \begin{multline}  \label{eq:AB2}
  \int_{\R^3 \times \R^3 \times \mathbb S^2} \chi_{\delta^{-1}}
  \, \Big[\mu'_* \, |h'| + \mu' \, |h'_*| + \mu \, |h_*| \Big] \,
|v-v_*| \,  m \dd v \dd v_* \dd \sigma
\\
=
\int_{\R^3 \times \R^3 \times \mathbb S^2} \chi_{\delta^{-1}}
  \, \Big[\mu'_* \, |h'| + \mu' \, |h'_*| - \mu_*
  \, |h| - \mu \, |h_*| \Big]  
\,  |v-v_*| \,    m \dd v \dd v_* \dd \sigma
  \\
  +4\pi \int_{\R^3 \times \R^3}  \chi_{\delta^{-1}} \, \mu_* \, |h| \,
  |v-v_*| \,  m \dd v \dd v_*
  \\
  + 8\pi \, \int_{\R^3 \times \R^3} \chi_{\delta^{-1}} \, \mu \, |v-v_*| \,
  |h_*| \,   m \dd v \dd v_*, 
\end{multline}
and the  first term in the right hand side of \eqref{eq:AB2} is bounded thanks to 
Lemma~\ref{lem:Povznerbis} as 
 \begin{multline}  \label{eq:AB3}
  \int_{\R^3 \times \R^3 \times \mathbb S^2} \chi_{\delta^{-1}} \,
  \Big[\mu'_* \, |h'| + \mu' \, |h'_*| - \mu \, |h_*|  - \mu_* \,
  |h| \Big]  
\, |v-v_*| \,  m \dd v \dd v_* \dd \sigma
  \\
  = \int_{\R^3 \times \R^3} \chi_{\delta^{-1}} \, \mu_*\, |h| 
  \left( \int_{\R^3 \times \R^3 \times \mathbb S^2}  \Big[|v'_*|^k  +
    |v'|^k - |v_*|^k - |v|^k  \Big]   \dd \sigma \right)
  \, |v-v_*| \dd v \dd v_* 
  \\
  \le  
  \int_{\R^3 \times \R^3 } \chi_{\delta^{-1}} \, \mu_*\, |h| \, 
  C_k \, \left( |v|^{k-1} \, |v_*| + |v| \, |v_*|^{k-1}
  \right)
  \, |v-v_*|  \dd v \dd v_* 
  \\ 
 - (4\pi -\gamma_k) \int_{\R^3 \times \R^3} \chi_{\delta^{-1}} \, \mu(v_*) \, |h| \,  |v|^k  
  |v-v_*|   \dd v \dd v_*
\end{multline}
(observe that our characteristic function $\chi_{\delta^{-1}}$ is
invariant under the usual changes of variables as it only depends on
the kinetic energy and momentum).

Now using the elementary inequality 
\[
 \chi_{\delta^{-1}}(v,v_*) \le {\bf 1}_{|v|\ge \delta^{-1}/2} + {\bf
   1}_{|v_*|\ge \delta^{-1}/2} \le 2 \, \delta \, (|v| + |v_*|), 
\]
we easily and crudely bound from above the second and third terms of
the right hand side in \eqref{eq:AB2}, and the first term of the right hand side in
\eqref{eq:AB3}, in the following way
 \begin{multline}  \label{eq:AB4}
  4\pi \int_{\R^3 \times \R^3}  \chi_{\delta^{-1}} \,   |v-v_*|  \,
  \mu_* \, m \, |h| \,
  \dd v \dd v_*
  \\
  +
 8\pi \, \int_{\R^3 \times \R^3} \chi_{\delta^{-1}}  \, |v-v_*|  \,   m_* \, \mu_* \,
  |h| \dd v \dd v_* 
  \\
  +   C_k \, \int_{\R^3 \times \R^3 } \chi_{\delta^{-1}} \, |v-v_*|   \, 
 \left( |v|^{k-1} \, |v_*| + |v| \, |v_*|^{k-1}  \right)
  \, \mu_*\, |h|  \dd v \dd v_* 
  \\
  \le 4\pi \, \int_{\R^3 \times \R^3} \chi_{\delta^{-1}} \, \mu(v_*) \, |h| \,  |v|^k  
  |v-v_*|   \dd v \dd v_* \\ +
 8\pi \, \int_{\R^3 \times \R^3} \delta \, (|v| + |v_*|) \, |v-v_*|  \,   m_* \, \mu_* \,
  |h| \dd v \dd v_* 
  \\ +
  C_k \, \int_{\R^3 \times \R^3 } \delta \, (|v| + |v_*|) \, |v-v_*|    \, \langle v_* \rangle^{k-1}  \, \langle v \rangle^{k-1}  
  \, \mu_*\, |h|  \dd v \dd v_* 
  \\ 
 \le 4\pi \, \int_{\R^3 \times \R^3} \chi_{\delta^{-1}} \, \mu(v_*) \, |h| \,  |v|^k  
  |v-v_*|   \dd v \dd v_* +% \delta \, 
%   C_k \, \int_{\R^3  } \mu_*   \,     \langle v_* \rangle^{k+1}  \,  dv_* 
% \int_{\R^3  }  |h|  \, 
%     \, \langle v \rangle^{k+1}  \,dv =
    \OO(\delta) \, \| h
    \|_{L^1(\nu m)}. 
  \end{multline}
  Putting together the estimates \eqref{eq:AB1}, \eqref{eq:AB2},
  \eqref{eq:AB3} and \eqref{eq:AB4}, we get
\begin{multline*}
  \left\| \BB^2 _\delta h \right\|_{L^1(m)} \le \OO (\delta) \, \, \|
  h \|_{L^1(\nu m)} + \gamma_k \int_{\R^3 \times \R^3}
  \chi_{\delta^{-1}} \, \mu(v_*) \, |h| \, |v|^k
  |v-v_*|   \dd v \dd v_*\\
  \le \left( \OO (\delta) + \gamma_k \right) \, \, \| h \|_{L^1(\nu
    m)}
% \OO (\delta) \, \int_{\R^3  }  |h|  \, 
%     \, \langle v \rangle^{k+1}  \dd v +  \gamma_k \int_{\R^3}  (\mu *
%     |\cdot|)  \,  |h| \,  |v|^k  \dd v 
  \end{multline*}
which concludes the proof. 
 \end{proof}

 \subsection{Pointwise estimates on the remainder}
\label{sec:pointw-weight-estim}

The goal of the subsection is to establish estimates on $Q^+$ in
$L^\infty$ spaces with polynomial and exponential weights.  As a
preliminary step, we shall first establish a representation result for
the gain part of the collision operator $Q^+$ when applied to radially
symmetric functions.  The following result is adapted from \cite[Lemma
3.6]{MR1958975}, see also \cite{Semikov1995,Semikov1997}.  We give
however a full proof of the result for several reasons: the statement
as well as the step 1 of the proof are modified, and the final step 4
of the proof below was omitted in the quoted papers.
\begin{lem}\label{lem:Q+rad} 
 Let $F$ and $G \in L^1(\R^3)$ be some non-negative radially
  symmetric functions.  Then $Q^+(G,F)=Q^+(F,G)$ defined in
  \eqref{eq:decQ} is radially symmetric and, denoting $r=|v|$, we have
\begin{equation}\label{eq:Q+1}
  Q^+(G,F)(r) = \int_0 ^{+\infty} \!\! \int_0 ^{+\infty} {\bf 1}_{(r')^2 +
    (r'_*)^2 > r^2} \, B \, G(r'_*) \, F(r')  \dd r' \dd r'_*, 
\end{equation}
with 
$$
B := 64 \pi^2 \, \frac{ r' r'_*}{r} \, \min \{ r,r_*,r', r'_*\}, \quad
r_* := \sqrt{ (r')^2 + (r'_*)^2 - r^2}.
$$
\end{lem}

% \textcolor{red}{Je suis ok avec l'enonce, j'ai refait tous les
%   calculs. J'ai trouve deux coquilles de constantes dans la preuve,
%   mais qui se compensaient~!}

\begin{proof}[Proof of Lemma~\ref{lem:Q+rad}]
We proceed in several steps. 
\mk

\noindent
{\sl Step 1: Integral representation of the operator on the whole domain.} We claim that 
\begin{equation}\label{eq:Q+Physique}
  Q^+(F,G)(v) 
  = 8 \,  \int_{\R^3} \int_{\R^3} \int_{\R^3}  
  G( v'_* ) \, F ( v' ) \, \delta_{\mathsf C_m} \, \delta_{\mathsf C_e}  \dd v_* \dd v'
  \dd v'_*
\end{equation}
% \textcolor{red}{Attention changement d'un facteur 16 en 8 dans la
%   formule ci-dessus qui provient d'une coquille dans la formule qu'on
%   reprenait de Bobylev}
where
\[
\mathsf C_m := \left\{ (v,v_*,v',v'_*) \in (\R^3)^4, \ v + v_* = v' +v'_* \right\}
\]
and 
\[
\mathsf C_e := \left\{ (v,v_*,v',v'_*) \in (\R^3)^4, \ |v|^2 + |v_*|^2 = |v'|^2 +
|v'_*|^2 \right\}.
\]

In order to prove the claim, we use the identity (see~\cite[Lemma
1]{MR1728639} )
\begin{equation}\label{eq:IdentityDiracEnergy}
  \forall \, \Phi \in
  C(\R^3), \ \forall \, w \in \R^3, \quad \int_{\mathbb S^2} \Phi ( |w|
  \sigma - w) \dd \sigma = {1 \over |w|} \int_{\R^3}  \Phi (y) \,
  \delta_{y \cdot w + {1 \over 2} |y|^2 = 0 } \dd y. 
\end{equation}
% \textcolor{red}{Attention j'ai change le facteur $2$ en $1$ dans la formule
%   ci-dessus, en mettant le detail du calcul ci-dessous pour eviter les
%   erreurs} 
The proof is straightforward by completing the square in the Dirac
function
\begin{equation*}
  \int_{\R^3}  \Phi (y) \, \delta_{y \cdot w + {1 \over 2} |y|^2 = 0 } \dd y =  \int_{\R^3}  \Phi (y) \,
  \delta_{\frac{|y+\omega|^2-|\omega|^2}{2} = 0 } \dd y,
\end{equation*}
then changing variables to the spherical coordinates $y = - \omega
+ r \, \sigma$ 
\begin{equation*}
  \dots = \int_{r=0} ^{+\infty} \int_{\mathbb S^2}  \Phi (-\omega + r
  \sigma) \, \delta_{\frac{|r|^2-|\omega|^2}{2} = 0 } \, r^2 \dd \sigma \dd r,
\end{equation*}
and finally performing the change of variable $s =
(r^2-|\omega|^2)/2$ on the radial variable
\begin{equation*}
  \dots = \int_{s=-\frac{|\omega|^2}{2}} ^{+\infty} \int_{\mathbb S^2}  \Phi (-\omega + r
  \sigma) \, \delta_{s = 0 } \, r \dd \sigma \dd s 
  = |\omega|\, \int_{\mathbb S^2}  \Phi (|\omega|
  \sigma -\omega) \, \dd \sigma.
\end{equation*}
 
We start from the definition \eqref{eq:decQ}, \eqref{eq:defcoll} of
$Q^+$ and we write
\begin{multline*}
Q^+(G,F)(v) \\ 
= \int_{\R^3} \int_{\mathbb S^{2}} |v-v_*| \,
G( v_* -  (|w| \, \sigma - w)) \,  F ( v + (|w| \, \sigma - w)) \dd
v_* \dd \sigma
\\
= 2 \, \int_{\R^3} \int_{\R^3} 
F ( v + y) \, G( v_* -  y) \, \delta_{y \cdot w + {1 \over 2} |y|^2
  = 0 }  \dd v_* \dd y
\\
= 2 \,  \int_{\R^3} \int_{\R^3} \int_{\R^3}  
F ( v + y) \, G( v_* -  z) \, \delta_{y \cdot w + {1 \over 2} |y|^2
  = 0 } \, \delta_{y-z = 0} \dd v_* \dd y \dd z
\end{multline*}
where we have set $w := (v-v_*)/2$ and we have used
\eqref{eq:IdentityDiracEnergy}.  We conclude by performing the change
of variables $v' := v + y$, $v'_* := v_* - z$ and observing that
$$
\delta_{y \cdot w + {1 \over 2} |y|^2 = 0 } \, \delta_{y-z = 0}  
=  4\, \delta_{\mathsf C_m} \, \delta_{\mathsf C_e}  , 
$$
because $\delta_{y-z = 0} = \delta_{\mathsf C_m}$ and 
\begin{multline*}
  \fa (v,v'_*,v'_*,v'_*) \in \mathsf C_m, \quad {1 \over 4} \, \left(
    |v'|^2 + |v'_*|^2 - |v|^2  - |v_*|^2 \right) \\
  = {1 \over 4} \, \left( |v'-v-v_*|^2 + |v'|^2 - |v|^2  - |v_*|^2 \right)  \\
  = {1 \over 2} \{ (v'-v) \cdot (v-v_*) + |v'-v|^2 \} = y \cdot w - {1
    \over 2} \, |y|^2.
\end{multline*}

\smallskip\noindent {\sl Step 2. } The fact that $Q^+(G,F)$ is
radially symmetric when applied to two radial functions $F$ and $G$ is
straightforward by using rotational changes of variable in the
collision integral. The identity $Q^+(F,G) = Q^+(G,F)$ is obtained by
the change of variable $\sigma \to -\sigma$ in \eqref{eq:Q+Physique}.
We can then write for radially symmetric functions $F$ and
$G$ %and get from \eqref{eq:Q+Physique}
\[
Q^+(G,F)(r) = \int_0 ^{+\infty} \int_0 ^{+\infty} \int_0 ^{+\infty} K \,
\delta_{\mathsf C_e} \, G(r'_*) \, F(r') \dd r_* \dd r' \dd r'_*
\]
with 
\begin{equation}\label{eq:defK}
K := 8 \, (r_*)^{2} (r')^{2} (r'_*)^{2} \, 
\int_{\mathbb S^{2}} \int_{\mathbb S^{2}} \int_{\mathbb S^{2}}
\delta_{\mathsf C_m} \dd \sigma_* \dd \sigma' \,
\dd\sigma'_*
\end{equation}
with the transparent notation 
\[
\left\{ 
\begin{array}{l} \ds
r = |v|, \quad r_* = |v_*|, \quad r'=|v'|, 
\quad r'_*=|v'_*|, \vs \\ \ds 
\sigma_* = \frac{v_*}{|v_*|}, \quad \sigma' = \frac{v'}{|v'|}, \quad 
\sigma'_*= \frac{v'_*}{|v'_*|}.
\end{array}
\right.
\]
Using the distributional identity
\begin{equation*}
  \delta_{r_* ^2 = (r')^2 + (r'_*)^2 - r^2} \, {\bf 1}_{r_* \ge 0} =
  \frac{1}{2 r_*} \, \delta_{r_* = \sqrt{(r')^2 + (r'_*)^2 - r^2}} \,
  {\bf 1}_{(r')^2 + (r'_*)^2 - r^2 \ge 0}
\end{equation*}
we obtain
\begin{equation}\label{eq:Q+2}
Q^+(G,F)(r) = \int_0 ^{+\infty} \!\! \int_0 ^{+\infty} {\bf 1}_{(r')^2 +
  (r'_*)^2 > r^2}  {K \over 2  r_*} \, G(r'_*) \, F(r')  \dd r' \dd r'_*
\end{equation}
where now $r_*$ is defined by $r_* := \sqrt{(r')^2 + (r'_*)^2 - r^2}$. 

 \medskip

\noindent{\sl Step 3.} Let us prove that
  \begin{equation}\label{eq:intDeltaCm}
\int_{\mathbb S^{2}} \int_{\mathbb S^{2}} \int_{\mathbb S^{2}}
\delta_{\mathsf C_m}\dd \sigma_* \dd \sigma' \dd \sigma'_* =
{32 \pi  \over r r_* r' r'_*} A ,
\end{equation}
with 
$$
A := \int_0^{+\infty}  \sin (ru) \, \sin (r_* u) \, \sin (r' u) \, \sin (r'_* u) \,
\frac{{\rm d} u}{u^2} .
$$

We use the following representation of Dirac masses on $\R^3$: 
\[
\delta_{\mathsf C_m} = \frac{1}{(2\pi)^3} \int_{\R^3} e^{i (z, v+v_* - v' -
  v'_*)} \dd z
\]
which yields, thanks to a spherical change of variable on $z$ with $u
= |z|$ and $e = z/|z|$,
\begin{multline*}
  \int_{\mathbb S^{2}} \int_{\mathbb S^{2}} \int_{\mathbb S^{2}}
  \delta_{\mathsf C_m}\dd \sigma_* \dd \sigma' \dd \sigma'_* \\
  = \frac{1}{(2\pi)^3} \int_0 ^{+\infty} \int_{\mathbb S^{2}}
  \int_{\mathbb S^{2}} \int_{\mathbb S^{2}} \int_{\mathbb S^{2}} e^{i
    u (e, v+v_*-v' - v'_*)} \dd e \dd \sigma_* \dd \sigma' \dd
  \sigma'_* \, u^2 \dd u.
\end{multline*}
Observe that this formula is invariant under rotation of the variable
$v$: this can be proved by using appropriate rotations on the
integration variables $e$, $\sigma_*$, $\sigma'$, $\sigma'_*$. We can
therefore add an average over $\sigma = v/|v|$, and then remove the
spherical average over $e$, which is no more necessary:
\begin{multline*}
  \int_{\mathbb S^{2}} \int_{\mathbb S^{2}} \int_{\mathbb S^{2}}
  \delta_{\mathsf C_m}\dd \sigma_* \dd \sigma' \dd \sigma'_*
  \\
  = \frac{1}{(2\pi)^3} \int_0 ^{+\infty} \int_{\mathbb S^{2}}
  \int_{\mathbb S^{2}} \int_{\mathbb S^{2}} \int_{\mathbb S^{2}} e^{i
    u (e_0, v+v_*-v' - v'_*)} \dd \sigma \dd \sigma_* \dd \sigma' \dd
  \sigma'_* \, u^2 \dd u
\end{multline*}
for some fixed unit vector $e_0 \in \mathbb S^2$ (the volume of the
two spherical averages removed and added cancel). We then compute
\[
\int_{\mathbb S^2} e^{i u (e_0, w)} \dd \sigma = 2 \pi
\, \int_0 ^\pi e^{i u |w| \cos \theta} \, \sin \theta \dd \theta
= \frac{4\pi \sin (|w|u)}{ |w| u},
\]
and straightforwardly deduce \eqref{eq:intDeltaCm}. 

\medskip
\noindent{\sl Step 4. } We claim that for any $r,r_*,r',r'_* > 0$
satisfying the conservation of energy condition $r^2+r_*^2 = (r')^2 +
(r'_*)^2$, it holds
\begin{equation}\label{eq:ProdSinInt=}
A
= \frac \pi 2  \, \min \left\{ r,r_*,r',r'_* \right\}. 
\end{equation}
%\textcolor{red}{Ici chgt facteur 4 en 2 au denominateur, voir ci-dessous}
Indeed, by  Lebesgue dominated convergence theorem, we have 
\begin{equation*}
  A = \lim_{\eps \to 0} A_\eps  \ \mbox{ with } \ A_\eps 
  := \int_\eps^{+\infty}  \sin (ru) \, \sin (r_* u) \, \sin (r' u) \, \sin (r'_* u) \,
  \frac{{\rm d}u}{u^2}.
\end{equation*}
Using the identities $\sin z = (e^{iz} - e^{-iz})/(2i)$ and $\cos z = (e^{iz} + e^{-iz})/2$, we have 
\begin{eqnarray*}
&& 4 \, \sin (ru) \, \sin (r_* u) \, \sin (r' u) \, \sin (r'_* u) 
\\
&&\qquad= \cos (( r+r_*+r'+r'_*) u ) - \cos (( r+r_*+r'-r'_*) u ) 
\\
&&\qquad - \cos (( r+r_*-r'+r'_*) u ) + \cos (( r+r_*-r'-r'_*) u ) 
\\
&&\qquad- \cos (( r-r_*+r'+r'_*) u ) + \cos (( r-r_*+r'-r'_*) u ) 
\\
&&\qquad + \cos (( r-r_*-r'+r'_*) u ) - \cos (( r-r_*-r'-r'_*) u ).
\end{eqnarray*}
% \textcolor{red}{Attention il y avait une erreur d'un facteur 2 ici
%   j'ai remplace 8 par 4 a gauche ci-dessus : il faut 8 pour enlever
%   tous les denominateurs des sinus sous forme exponentielle, mais
%   ensuite il y a un facteur 2 pour recreer des cosinus a partir des
%   exponentielles} 
We observe that thanks to an integration by part, for any $a \in \R$, we have 
\begin{eqnarray*}
  \int_\eps^{+\infty} \cos (a \,  u ) \, \frac{{\rm d}u}{u^2} 
  &=& {\cos (a \,  \eps ) \over \eps} - a \, \int_\eps^\infty \sin (a
  \,  u ) \, \frac{{\rm d}u}{u} 
  \\
  &=& {1 \over \eps}- a \, \int_0^{+\infty} \sin (a \,  u ) \,
  \frac{{\rm d}u}{u}  + \OO(a^2 \, \eps)
  \\
  &=& {1 \over \eps}- {\pi \over 2} \, |a| + \OO(a^2 \, \eps).
\end{eqnarray*}
All together, we get 
\begin{eqnarray*}
- {8 \over \pi} \, A
&=& -  {8 \over \pi} \, \lim_{\eps\to0} A_\eps
\\
&=& \,\,\,\, | r+r_*+r'+r'_*| - | r+r_*+r'-r'_*| 
\\
&& \!\! - \, | r+r_*-r'+r'_*| + | r+r_*-r'-r'_*| 
\\
&&  \!\!- \, | r-r_*+r'+r'_*| + | r-r_*+r'-r'_*| 
\\
&& \!\!  +\,  | r-r_*-r'+r'_*| - | r-r_*-r'-r'_*|.
\end{eqnarray*}

Now assume first $r>r_*$, $r'>r'_*$ and $r>r'$, so that the energy
conservation condition implies that $ r>r'>r'_*>r_*$, and in
particular $r-r_* > r'-r'_* >0$. Hence any of the terms $r,r_*,r',
r'_*$ is smaller than the sum of the three other terms. Using all
these inequalities, the above expression then simplifies into
\begin{eqnarray*}
- {8 \over \pi} \, A
&=& \,\,\,\, ( r+r_*+r'+r'_*)  - ( r+r_*+r'-r'_*) 
\\
&&   \!\! - \, ( r+r_*-r'+r'_*) +| r+r_*-r'-r'_*|
\\
&&  \!\! - \, ( r-r_*+r'+r'_*) + ( r-r_*+r'-r'_*) 
\\
&&  \!\! + | r-r_*-r'+r'_*| + ( r-r_*-r'-r'_*)
\\
&=& - 2 \, r_* - 2 \, r'_* + | (r-r') - (r'_*-r_*)| + |(r-r') + (r'_*-r_*)|
\\
&=& - 2 \, r_* - 2 \, r'_* +2 \, \max \left\{ r-r',r'_*-r_* \right\}.
\end{eqnarray*}
Now, from the elementary inequality 
$$
\forall \, x,y \ge 1, \quad x^2 + y^2 - 1 \le (x+y-1)^2 ,
$$
we deduce that 
$$
r = r_* \, \sqrt{ \Bigl( \frac{r'}{r_*} \Bigr)^2 + \Bigl( \frac{r'_*}{r_*} \Bigr)^2 - 1} \le 
 {r_*} \, \left|  \frac{r'}{r_*}   +  \frac{r'_*}{r_*}- 1 \right| = r' + r'_* - r_*
$$
where we have removed the absolute value due to the inequalities
above. We thus obtain $\max\{ r-r', r'_*- r_* \} = r'_*-r_*$. As a
consequence, we get
\[
- {8 \over \pi} \, A =  - 2 \, r_* - 2 \, r'_* +2 \, (r'_*-r_*)  = - 4 r_*  
=  - 4 \,  \min \{ r,r_*,r',r'_*\}. 
\]
We then conclude \eqref{eq:ProdSinInt=} by using symmetries: the cases
$r < r_*$, $r' < r'_*$, and $r < r'$ are treated by using the three
swappings $v \leftrightarrow v_*$, $v' \leftrightarrow v'_*$ and
$(v,v_*) \leftrightarrow (v',v'_*)$ leaving invariant the energy
conservation identity.
 
\medskip
\noindent{\sl Step~5. Conclusion.} We conclude by gathering \eqref{eq:Q+2}
with \eqref{eq:defK}, \eqref{eq:intDeltaCm} and
\eqref{eq:ProdSinInt=}.
\end{proof}

We can now prove the pointwise estimates with polynomial weight on the
collision operator.

\begin{lem}\label{lem:decQ+rad} Assume $k >3$.
  Then we have the following \emph{bilinear estimate} on the $Q^+$
  operator defined in \eqref{eq:decQ}: 
\begin{multline}\label{eq:Q+rad-bil}
 \fa f,g \in L^\infty(\langle v \rangle^{k+1}), \quad \| Q^+(f,g) \|_{L^\infty(\langle v \rangle^k)} \\
 \le C(k) \,
  \left( \| f \|_{L^\infty(\langle v \rangle^{k+1})}  \| g
    \|_{L^\infty(\langle v \rangle^{k})} + \|
    g\|_{L^\infty(\langle v \rangle^{k+1})} \| f \|_{L^\infty(\langle
      v \rangle^k)} \right) 
\end{multline}
for some constant $C(k)>0$ depending on $k$.

Moreover, we have, for any $k>3$ and $\delta > 0$, the following more
precise linear estimate on the remainder operator $\BB^2_\delta$
(defined in \eqref{eq:def-remain}):
\begin{equation}\label{eq:Apointwise}
  \fa h \in L^\infty(\langle v \rangle^{k+1}),\quad 
  \left\| \BB^2_\delta  h \right\|_{L^\infty(\langle v \rangle^k)} \le
  \left( {4 \over k-1} + \eta(k,\delta) \right) \,   \| h
%  \|_{L^\infty(\langle v \rangle^{k+1})},
  \|_{L^\infty(\nu \, \langle v \rangle^{k})},
\end{equation}
for some constructive $\eta(k,\delta)$ such that
$\eta(k,\delta) \to 0$ as $\delta \to 0$ with $k$ fixed.
\end{lem}

\begin{rem} %\Blue
  Observe that a similar estimate is easily proved for the loss part
  of the collision operator $Q^-(g,f)$ as soon as $k >3$. These
  estimates for $Q^+$ recover, by another method, some estimates in
  \cite{MR900501}, in a more precise form and with the sharp constant
  (and weaker moment condition). They are different in nature from
  convolution-like estimates
  \begin{multline}\label{eq:Q+convol}
  \| Q^+(f,g) \|_{L^\infty(\langle v \rangle^k)} \\
 \le C \, \left(  \| g \|_{L^1(\langle v \rangle^{k+1})}  \| f
 \|_{L^\infty(\langle v \rangle^{k+1})}   
 +   \| f \|_{L^1(\langle v \rangle^{k+1})}  \| g \|_{L^\infty(\langle v \rangle^{k+1})}  \right)
\end{multline}
which hold for any $k \ge 2$ and any $f,g \in (L^1 \cap L^\infty)
(\langle v \rangle^{k+1})$, as proved for instance in \cite{MR711482}
or in \cite[Theorem 2.1 and Remark 3]{MR2081030}. % The decay of the
% constant in terms of $k$ in the estimate \eqref{eq:Apointwise} seems
% also to be a phenomenon which had been so far unoticed for such
% pointwise estimates on the Boltzmann collision operator.
\end{rem}

\begin{proof}[Proof of Lemma~\ref{lem:decQ+rad}]
% We first split the operator according to the decomposition 
% \[
% {\bf 1} = {\bf 1}_{|v-v_*| \le 1} + {\bf 1}_{|v-v_*| >1}.
% \]
% On the first part $|v-v_*| \le 1$, it is clear that 
% \[
% \| Q^+(g,f) \|_{L^\infty_k} \le C \, \| f \|_{L^\infty_{k}} \, \| g
% \|_{L^\infty_{k}}
% \]
% which easily allows to conclude the proof of
% \eqref{eq:Q+rad-bil}-\eqref{eq:Q+rad-lin}-\eqref{eq:Apointwise}. Let
% us now consider the remaining part and therefore assume in the sequel
% that $|v-v_*| >1$, and let us proceed in two steps.  \medskip
  We split the proof in two steps along the two parts of the
  statement.  \medskip

\noindent
{\it Step 1. The bilinear estimate \eqref{eq:Q+rad-bil}.} 
Define the  functions 
$$
\forall \, r > 0, \quad F(r) := \sup_{|v|=r} |f(v)|, \quad G(r) :=
\sup_{|v|=r} |g(v)|,
$$
so that 
$$
\left|Q^+(g,f)(v)\right| \le  Q^+(G,F)(|v|).
$$
Observing that now  $F$ and $G$ are radially symmetric functions, 
for $(r',r'_*) \in \R^2_+$ we get
$$
\{ (r')^2 + (r'_*)^2 \ge r^2 \} \subset \{ r' \ge r/\sqrt{2} \} \cup
\{ r'_* \ge r/\sqrt{2} \}.
$$
We can estimate $Q^+(G,F)$ by using the representation formula in
Lemma~\ref{lem:Q+rad} and the following splitting
 \begin{multline}\label{eq:Q+4}
   Q^+(G,F)(r) \le \frac{C_0}{r} \, \int_{r/\sqrt{2}} ^\infty \dd r'
   \int_0^\infty \dd r'_* G(r') \, F(r'_*) \, r' \, (r'_*)^2
   \\
   + \frac{C_0}{r} \,\int_{r/\sqrt{2}} ^\infty \dd r'_* \int_0^\infty
   \dd r' G(r') \, F(r'_*) \, (r')^2 \, r'_* =: I_1 + I_2
\end{multline}
where we have used $\min\{ r,r_*,r',r'_*\} \le r'_*$ in the first
term, $\min \{r,r_*,r',r'_* \} \le r'$ in the second term and we have
set $C_0 := 64\pi^2$.

For the first term, we set $m_k := (1+|v|^2)^{k/2}$ and we remark that
as soon as $k > 3$, we have, for $r \ge 1$, 
\begin{eqnarray*}
I_{1}  
&=&  
\frac{C_0}{r} \, \left( \int_{r/\sqrt 2} ^{+\infty}  r' \, G(r') \dd r'
\right) \, \left( \int_{0} ^{+\infty} F(r'_*) \, (r'_*)^2 \dd r'_* \right)
%{1 \over 4 \pi} \, \| F \|_{L^1(\R^3)}
\\
&\le& \frac{C_0}{r (k-3)}  \, \left[  \sup_{\R_+} (G m_{k+1}) \right] \, \left[
\sup_{\R_+} (F m_k) \right] \,  \int_{r/\sqrt 2} ^{+\infty}  
 {r' \dd r' \over (1 + (r')^2)^{k+1 \over 2}}  
%\Bigl(  \sup_{\R_+} F  \, m_k \Bigr) {1 \over 4 \pi} \, \| m_k^{-1} \|_{L^1(\R^3)}
 \\
 &\le& \frac{C_0 2^{(k-1)/2}}{(k-1)(k-3)} \, {1 \over m_k(r)}  \, \| g
 \|_{L^\infty(m_{k+1})} \, \| f \|_{L^\infty(m_k)},
\end{eqnarray*}
so that 
 $$
 \forall \, r > 0, \quad I_1(r) \, m_k(r) \le \frac{C_0
   2^{(k-1)/2}}{(k-1)(k-3)} \, \| g \|_{L^\infty(m_{k+1})} \, \| f
 \|_{L^\infty(m_k)}.
 $$
 Because the terms $I_1$ and $I_2$ are symmetric (the change of
 variable $(r',r'_*) \to (r'_*,r') $ exchanges the role played by $F$
 and $G$), we obtain the same estimate for $I_2$ where we exchange the
 role played by $f$ with $g$, and this concludes the proof of
 \eqref{eq:Q+rad-bil}.

 \medskip
 \noindent {\it Step 2.} Let us prove the following linearized
   estimate
\begin{equation}\label{eq:Q+rad-lin}
  \left\|\left[Q^+(\mu,f) + Q^+(f,\mu)\right] {\bf 1}_{|v| \ge \delta^{-1}}
  \right\|_{L^\infty(\langle v \rangle^k)} 
  \le  \left( {16 \pi  \over k-1} + \eta(k,\delta) \right) \,   \| f
  \|_{L^\infty(\langle v \rangle^{k+1})}
\end{equation}
for some constant $\eta(k,\delta) \to 0$ as $\delta \to 0$, for $k >3$
fixed. It implies the desired inequality \eqref{eq:Apointwise} since 
\begin{multline*}
4 \pi \, (1+|v|^2)^{1/2} \, {\bf 1}_{|v| \ge \delta^{-1}} \le 4 \pi \,
(1+|v|) \, {\bf 1}_{|v| \ge \delta^{-1}} \\ \le \nu(v) + 4 \pi \, {\bf
  1}_{|v| \ge \delta^{-1}} \le \nu(v) + \delta \, \nu(v). 
\end{multline*}

Setting $G := \mu$ and $F := m_{k+1}^{-1}$ we have 
$$
\begin{cases} 
|Q^+(\mu,f)| \le  \| f \|_{L^\infty(\langle v \rangle^{k+1})} \, Q^+
(G,F) \vs \\ 
|Q^+(f,\mu)| \le  \| f \|_{L^\infty(\langle v \rangle^{k+1})} \, Q^+
(F,G)
\end{cases}
$$
and since $Q^+ (G,F) = Q^+(F,G)$ (cf. Lemma~\ref{lem:Q+rad}), it is
enough to establish the estimate \eqref{eq:Q+rad-lin} for the term
$Q^+(G,F)$ only.

For any $\eps \in (0,1)$ and $(r',r'_*) \in
\R^2_+$, we have
$$
\{  (r')^2 + (r'_*)^2 \ge r^2 \} \subset \{  r' \ge \sqrt{\eps} \, r \} \cup \{  r'_* \ge (1-\eps) \, r \},
$$
so that we may   estimate  $Q^+(G,F)$ thanks to the following splitting 
 \begin{multline}\label{eq:Q+4}
   Q^+(G,F)(r) \le \frac{C_0}{r} \, \int_{ \sqrt{\eps} \, r
   }^{+\infty} \dd r' \int_0^\infty \dd r'_* \, G(r') \, F(r'_*) \, r' \,
   (r'_*)^2
   \\
   + \frac{C_0}{r} \,\int_{(1-\eps)\, r} ^\infty \dd r'_*
   \int_0^{+\infty} \dd r' \, G(r') \, F(r'_*) \, (r')^2 \, r'_* =: I_1 +
   I_2
\end{multline}
where we have used $\min \{ r,r_*,r',r'_* \} \le r'_*$ in the first
term, and $\min \{ r,r_*,r',r'_* \} \le r'$ in the second term. 

For the first term, we have
\begin{multline*}
I_{1}  =
\frac{C_0}{r} \left( \int_{\sqrt{\eps} r} ^{+\infty}  {r' \, e^{-(r')^2/2}
    \over (2\pi)^{3/2}} \dd r' \right) \, \left(  \int_0^{+\infty} {(r'_*)^2
    \over (1+(r'_*)^2)^{k+1 \over 2}} \dd r'_* \right) \\
\le \frac{C_0}{r} \, {e^{-\eps \, r^2 /2} \over (2\pi)^{3/2}} \, {\Theta \over k-2}, 
\end{multline*}
with $\Theta \in (0,1)$. On the other hand, for the second term, we
have for any $r \ge 1$
\begin{eqnarray*}
  I_{2}  
  &=&
  \frac{C_0}{r} \left( \int_0^{+\infty} (r')^2  \,  {e^{-(r')^2/2} \over
      (2\pi)^{3/2}} \dd r' \right) \, \left( \int_{(1-\eps) r}^{+\infty}
    {r'_* \over (1+(r'_*)^2)^{k+1 \over 2}}  \dd r'_* \right) 
  \\
  &=& \frac{C_0}{4 \pi r (k-1)} \,  {1 \over (1+ (1-\eps)^2 \, r^2)^{k-1 \over 2}} 
  \le {  16 \pi  \over k-1} \, {1 \over (1-\eps)^{k-1}} \, {1 \over r
    \, m_{k-1}(r)}%  &\le& {64 \pi  \over k-1} \, {1 \over (1-\eps)^{k-1}} \, {1 \over m_{k}(r)}.
\end{eqnarray*}
where we recall that $C_0 = 64 \pi^2$. 

By combining these two estimates together, we get for any $r \ge 1$ 
\[
Q^+(G,F)(r) \, m_k(r) \, {\bf 1}_{r \ge \delta^{-1}} \le {{16 \pi}
  \over k-1} + \phi (k,\delta,\eps)
\]
 with $\phi = \phi_1 + \phi_2$ and
 \begin{eqnarray*}
   \left\{ 
\begin{array}{ll} \ds 
  \phi_1 (k,\delta,\eps) :=&   {C_1 \over k-1} \left[ {1 \over
       (1-\eps)^{k-1}} \, \sup_{r \ge \delta^{-1}}  {m_1(r) \over r} - 1
   \right],
   \vs \\ \ds
   \phi_2 (k,\delta,\eps) :=& 
   {C_2  \over k-2} \, \left[ \sup_{r \ge \delta^{-1}} m_k( r)
     e^{-\eps \, r^2 /2}   \right], 
 \end{array}
\right.
\end{eqnarray*}
for some numerical constants $C_1, C_2 > 0$. We deduce that
\eqref{eq:Q+rad-lin} holds with $\eta(k,\delta) := \phi (k,\delta,
\delta)$ for instance.

 \medskip
 \noindent {\it Step 3.} Coming back to the definition of $\BB^2
 _\delta h$ we split it into three pieces
\begin{multline*} 
  |\BB^2 _\delta h (v)| \le \int_{\R^3 \times \mathbb S^2} {\bf
    1}_{|v| \ge R} \, ( \mu'_* \, |h'| + \mu'_* \, |h'| ) \, |v-v_*|
  \, \dd v_* \dd \sigma
  \\
  + \int_{\R^3 \times \mathbb S^2} {\bf 1}_{|v| \le R} \, (1 -
  \Theta_\delta) ( \mu'_* \, |h'| + \mu'_* \, |h'| ) \, |v-v_*| \, \dd
  v_* \dd \sigma
  \\
  + \int_{\R^3 \times \mathbb S^2} {\bf 1}_{|v| \le R} \, (1 -
  \Theta_\delta) \, \mu \, |h_*| \, |v-v_*| \dd v_* \dd \sigma =: I_1
  + I_2 + I_3.
\end{multline*}

For the first term $I_1$ we use \eqref{eq:Q+rad-lin} and we get
\[
\n{I_1}_{L^\infty (\langle v \rangle^k)} = \sup_{r \ge 0} \left(
  I_1(r) \, m_k(r) \right)  \le   {{16 \pi} \over k-1}+ \eta(k,R^{-1}).
 \]
 For the second term $I_2$ we use the sharp form of the convolution
 inequality \eqref{eq:Q+convol} as stated in \cite[Theorem
 2.1]{MR2081030} and we get for $k>3$ %for $\delta^{-1} \ge 2 R$
\begin{multline*} 
  I_2(r) \, m_k(r) \le m_k(R) \, \left\| \left(Q^+_\delta(\mu,|h|) +
      Q^+_\delta(|h|,\mu)\right) \, {\bf 1}_{|v| \le R}
  \right\|_{L^\infty}\\
 \le C \, m_k(R) \, \n{h}_{L^\infty(\langle v \rangle^k)} \sup_{|v|\le R} \int_{\R^3}
 \int_{\mathbb S^2} \left( 1- \Theta_\delta \right) \,
 \frac{1}{\langle v' \rangle^k \langle v'_*\rangle^k} |v-v_*| \dd v_*
 \dd \sigma \\ 
\le C \, m_k(R) \, \n{h}_{L^\infty(\langle v \rangle^k)} \sup_{|v|\le R} \int_{\R^3}
 \int_{\mathbb S^2} \left( 1- \Theta_\delta \right) \,
 \frac{1}{( 1+ |v'|^2 + |v'_*|^2)^{k/2}} |v-v_*| \dd v_*
 \dd \sigma \\ 
\le C \, m_k(R) \, \n{h}_{L^\infty(\langle v \rangle^k)} \sup_{|v|\le R} \int_{\R^3}
 \int_{\mathbb S^2} \left( 1- \Theta_\delta \right) \,
 \frac{1}{( 1+ |v|^2 + |v_*|^2)^{k/2}} |v-v_*| \dd v_*
 \dd \sigma
  % \\
  % \le m_k(R) \, \chi(k,\delta) \, \left( \| \mu \|_{L^1_{k-1 \over 2}}
  %   \, \| h \|_{L^\infty_{k-1 \over 2}} + \| h \|_{L^1_{k-1
  %       \over 2}} \, \| \mu \|_{L^\infty_{k-1 \over 2}} \right)
\end{multline*}
for some constant $C>0$. Observe that we can also write the same
control on the third term $I_3$ by a simpler argument: 
\begin{multline*} 
  I_3(r) \, m_k(r) \\
 \le C \, m_k(R) \, \n{h}_{L^\infty(\langle v \rangle^k)} \sup_{|v|\le R} \int_{\R^3}
 \int_{\mathbb S^2} \left( 1- \Theta_\delta \right) \,
 \frac{1}{\langle v \rangle^k \langle v_*\rangle^k} |v-v_*| \dd v_*
 \dd \sigma \\ 
\le C \, m_k(R) \, \n{h}_{L^\infty(\langle v \rangle^k)} \sup_{|v|\le R} \int_{\R^3}
 \int_{\mathbb S^2} \left( 1- \Theta_\delta \right) \,
 \frac{1}{( 1+ |v|^2 + |v_*|^2)^{k/2}} |v-v_*| \dd v_*
 \dd \sigma.
  % \\
  % \le m_k(R) \, \chi(k,\delta) \, \left( \| \mu \|_{L^1_{k-1 \over 2}}
  %   \, \| h \|_{L^\infty_{k-1 \over 2}} + \| h \|_{L^1_{k-1
  %       \over 2}} \, \| \mu \|_{L^\infty_{k-1 \over 2}} \right)
\end{multline*}

We then use 
\[
\left( 1- \Theta_\delta \right) \le \left( {\bf 1}_{|v-v_*| \ge \delta^{-1}} + {\bf 1}_{|v-v_*| \le 2
    \delta} + {\bf 1}_{ \cos \theta \ge 1 - 2\delta} \right)
\]
which gives rise to three terms to be controlled. The term associated
with the third part is $o(\delta)$ thanks to the $L^1$ integration on
the sphere, the second term is $O(\delta)$ thanks to the term
$|v-v_*|$ in the collision kernel, and for the first term, if we
assume $\delta$ small enough so that $\delta^{-1} \ge 2 R$, then we
deduce that $|v_*| \ge \delta^{-1}/2$ which gives a decay
$O(\delta^{k-2})$. We finally deduce that 
\begin{equation*}
  \n{I_2}_{L^\infty(\langle v \rangle^k)} + \n{I_3}_{L^\infty(\langle v \rangle^k)} \le o(\delta) \,
  \n{h}_{L^\infty(\langle v \rangle^k)}. 
\end{equation*}

Then the proof of \eqref{eq:Apointwise} follows by
% noting that
% $$
% \nu(v) := \int_{\R^3} \mu(v_*) \, 4\pi \, |v-v_*|\dd v_* \ge 4\pi \,
% \Bigl|\int_{\R^3} \mu(v_*) \, (v-v_*_\dd v_* \Bigr|= 4 \pi \, |v|,
%  $$
%  and
 gathering the preceding estimates on $I_1, I_2, I_3$. 
 \end{proof}

 \begin{rem}
   The reader can check that the above proof fails for Lebesgue spaces
   $L^q$, $q \in (1,+\infty)$: in fact the loss of weight in a
   bilinear inequality of the form $L^q \times L^q \to L^q$ seems
   strictly greater than what is allowed by $\nu$.
 \end{rem}

Let us now consider the case of a stretched exponential weight. 
\begin{lem}\label{lem:decQ+rad-exp} 
Consider the weight $m =  e^{\kappa |v|^\beta}$ with $\kappa >0$, $\beta \in (0,2)$.
Then we have the following \emph{bilinear estimate} on $Q^+$ 
  defined in \eqref{eq:decQ}:  
\begin{equation}\label{eq:Q+rad-bil-exp}
\| Q^+(g,f) \|_{L^\infty(\nu^{\beta} m)} 
   \le C \,
\left( \| f \|_{L^\infty(m)} \, \left\| g \right\|_{L^\infty(\nu  m)} + \|
  g\|_{L^\infty(m)} \, \left\| f  \right\|_{L^\infty(\nu m)} \right) , 
\end{equation}
for any $f,g \in L^\infty( \nu m)$  and for some constant $C>0$ depending on $m$.

Moreover, for any $\delta > 0$, we have the following \emph{linear
  estimate} on the remainder operator $\BB^2_\delta$: 
\begin{equation}\label{eq:Apointwise-exp}
 \fa h\in L^\infty( \nu m), \quad \left\| \BB^2_\delta  h \right\|_{L^\infty(m)} \le
\eta(\delta) \, \left\| h \right\|_{L^\infty(\nu m)},
\end{equation}
for some constructive constant $\eta(\delta)$ such that $\eta(\delta) \to 0$ as $\delta\to0$. 
\end{lem}

\begin{rem}
  Observe that by inspection $Q^-(h,\mu)$ is bounded in
  $L^\infty(m)$. However again such estimates are new for $Q^+$ to our
  knowledge.  They complement the $L^1$ integral estimates
  in~\cite{Mcmp}. These estimates show that the bilinear operator
  $Q^+$ is bounded for the norm $L^\infty(\nu m)$ for $\beta \in
  [1,2)$.
\end{rem}

\begin{proof}[Proof of Lemma~\ref{lem:decQ+rad-exp}]  We prove
  \eqref{eq:Q+rad-bil-exp} in step 1 and \eqref{eq:Apointwise-exp} in
  step 2. \mk

  \noindent {\sl Step 1. The bilinear estimate
    \eqref{eq:Q+rad-bil-exp}.}  We proceed as in step 1 of
  Lemma~\ref{lem:decQ+rad}. Consider $f,g \in L^\infty (\nu m)$ and
  introduce the associated radially symmetrized functions $F, G$ as
  before. We may estimate $Q^+(G,F)$ given by Lemma~\ref{lem:Q+rad}
  thanks to the following splitting
 \begin{multline}\label{eq:Q+4}
   Q^+(G,F)(r) \\ \le \frac{C_0}{r} \,\int_0^{+\infty} \!\!\int_0^{+\infty}
   {\bf 1}_{(r')^2 + (r'_*)^2 \ge r^2} \, {\bf 1}_{r' \ge r'_*} G(r')
   \, F(r'_*) \, r' \, (r'_*)^2 \dd r' \dd r'_*
   \\
   + \frac{C_0}{r} \,\int_0^{+\infty} \!\!\int_0^{+\infty} {\bf 1}_{(r')^2 +
     (r'_*)^2 \ge r^2} \, {\bf 1}_{r'_* \ge r'} G(r') \, F(r'_*) \,
   (r')^2 \, r'_* \dd r' \dd r'_* =: I_1 + I_2
\end{multline}
where we have used $\min\{ r,r_*,r',r'_* \} \le r'_*$ in the first
term, $\min \{ r,r_*,r',r'_* \} \le r'$ in the second term and we have
set again $C_0 := 64\pi^2$. 

We estimate the two terms in a symmetric way as
\begin{equation*}
\left\{ 
\begin{array}{l}\ds
I_1(r) \le   \| g \|_{L^\infty(m)} \, \| f
\|_{L^\infty(\langle v \rangle m)} \, J(r), \vs \\ \ds
I_2(r)  \le  \| g \|_{L^\infty(\langle v \rangle m)} \, \| f
\|_{L^\infty(m)} \, J(r),
\end{array}
\right.
\end{equation*}
with 
\begin{equation*}
  J(r)  =  \frac{C_0}{r}  \int_0^{+\infty} \!\! \int_0^{+\infty} {\bf 1}_{\rho \ge r} \,  {\bf 1}_{r'_* \ge r'}  
  \,  (m')^{-1} \, (m'_*)^{-1} \, (r')^2 \dd r' \dd r'_* 
\end{equation*}
where we denote $\rho^2 := (r')^2 + (r'_*)^2$. We introduce the
notations $x:= r'/\rho$, $y:= r'_*/\rho$, and we remark that by
inspection 
$$
\forall \, x \in [0, 1/\sqrt{2}], \quad x^\beta +
(1-x^2)^{\beta/2} -1 \ge \eta\, x^\beta
$$
for some explicit $\eta = \eta(\beta) \in (0,1)$.  As a consequence,
making the change of variables $(r',r'_*) \mapsto (r',\rho)$ and
noticing that the condition $r' \le r'_*$ is equivalent to the
condition $x \le 1/\sqrt{2}$, we get
\begin{multline*}
  J (r) = \frac{C_0}{r} \int_r^{+\infty} \dd \rho \int_0^{\rho/\sqrt{2}}
  \dd r' \, e^{-\kappa \, ((r')^\beta+(r'_*)^\beta)} \, (r')^2 \,
  {\rho \over r'_*} %1
  \\
  \le \frac{C_0 \sqrt{2}}{r} \int_r^{+\infty} e^{-\kappa \,
    \rho^\beta} \, \dd \rho \int_0^{+\infty} e^{-\kappa\eta \,
    (r')^\beta} \, (r')^2 \dd r' \le C \, \frac{e^{-\kappa \,
      r^\beta}}{r^\beta},
\end{multline*}
for some constant $C$ which depends on $C_0$, $\beta$, $\kappa$. 

Notice that in order to get the last inequality above we may proceed
as follows:
\begin{itemize} 
\item If $\beta \in (1,2)$ we use the inequality
  $1 \le \rho^{\beta-1}/r^{\beta-1}$ and we simply integrate exactly
  the resulting function by using its anti-derivative
  \begin{equation*}
    \int_r^{+\infty} e^{-\kappa \, \rho^\beta} \,  \dd \rho \le r^{1-\beta}
    \int_r^{+\infty} \rho^{\beta-1} e^{-\kappa \, \rho^\beta} \,  \dd \rho =
    r^{1-\beta} \frac{e^{-\kappa \, r^\beta}}{\beta-1}. 
  \end{equation*}
\item  If $\beta \in  (0,1)$, we write 
  \begin{multline*}
    I(r) := \int_r^{+\infty} e^{-\kappa \, \rho^\beta} \,  \dd \rho =
    \int_r^{+\infty} \rho^{1-\beta} \, \rho^{\beta-1} \, e^{-\kappa \,
      \rho^\beta} \,  \dd \rho \\
    = \left[ \rho^{1-\beta} \, \frac{e^{-\kappa \, \rho^\beta}}{-\kappa
        \beta} \right]_r ^{+\infty} + \frac{(1-\beta)}{\kappa \beta} \int_r^{+\infty}
    \rho^{-\beta} \, e^{-\kappa \, \rho^\beta} \,  \dd \rho \\
    \le \frac{r^{1-\beta} e^{-\kappa \, \rho^\beta}}{\kappa \beta} +
    \frac{r^{-\beta} (1-\beta)}{\kappa \beta}  \, I(r)
  \end{multline*}
which implies for $r \ge r_0$ with $r_0^{-\beta} (1-\beta) /(\kappa
\beta) \le 1/2$: 
\begin{equation*}
  \int_r^{+\infty} e^{-\kappa \, \rho^\beta} \,  \dd \rho \le
  \frac{2 r^{1-\beta} e^{-\kappa \, \rho^\beta}}{\kappa \beta}.
\end{equation*}
\end{itemize}

The estimate for small values of $r$, say $r \in [0,r_0]$, is a
consequence of \eqref{eq:Q+rad-bil}. This thus concludes the proof of
\eqref{eq:Q+rad-bil-exp}.

\medskip\noindent {\sl Step 2. The linearized estimate.} Estimate
\eqref{eq:Q+rad-bil-exp} implies the following \emph{linearized
  estimate}
\begin{equation}\label{eq:Q+rad-lin-exp}
  \left\|[Q^+(\mu,h) + Q^+(h,\mu)] {\bf 1}_{|v| \ge R} \right\|_{L^\infty(m)} \le
  O(\delta^\beta) \,   \left\| h  \right\|_{L^\infty(\nu m)}.
%+  C_k \, \| f \|_{L^\infty_{k}} 
\end{equation}
We then proceed as in the Step 3 of Lemma~\ref{lem:decQ+rad}:
\begin{multline*}
 |\BB^2 _\delta h (v)| \le \int_{\R^3 \times \mathbb S^2} {\bf
    1}_{|v| \ge R} \, ( \mu'_* \, |h'| + \mu'_* \, |h'| ) \, |v-v_*|
  \, \dd v_* \dd \sigma
  \\
  + \int_{\R^3 \times \mathbb S^2} {\bf 1}_{|v| \le R} \, (1 -
  \Theta_\delta) ( \mu'_* \, |h'| + \mu'_* \, |h'| ) \, |v-v_*| \, \dd
  v_* \dd \sigma
  \\
  + \int_{\R^3 \times \mathbb S^2} {\bf 1}_{|v| \le R} \, (1 -
  \Theta_\delta) \, \mu \, |h_*| \, |v-v_*| \dd v_* \dd \sigma =: I_1
  + I_2 + I_3.
\end{multline*}
The estimate \eqref{eq:Q+rad-lin-exp} implies 
\begin{equation*}
  \n{I_1}_{L^1(m)} \le O(\delta^\beta) \, \left\| h  \right\|_{L^\infty(\nu m)}.
\end{equation*}
Then the same estimates as  in the Step 3 of the proof of
Lemma~\ref{lem:decQ+rad} yield 
\begin{equation*}
  \n{I_2}_{L^\infty(m)} + \n{I_3}_{L^\infty(m)} \le o(\delta) \,
  \n{h}_{L^\infty(\langle v \rangle^k)}
\end{equation*}
(the truncation ${\bf 1}_{|v| \le R}$ means that any weight can be
chosen on the left hand side) which concludes the proof of
\eqref{eq:Apointwise-exp}.
\end{proof}

% \textcolor{red}{Ok sur cette preuve, j'ai verifie les calculs et
%   reecrit l'argument final pour le rendre plus simple et plus
%   explicite.}
\subsection{Dissipativity estimate on the coercive part}
\label{sec:diss-estim-b}

Let us summarize in the following lemma the estimates
available for $\BB^2 _\delta$.

\begin{lem} \label{lem:Asmall} Consider $p,q \in [1,\infty]$ and a
  weight function $m$ satifying one of the conditions {\bf (W1)}, {\bf
    (W2)}, {\bf (W3)} of Theorem~\ref{theo:LIBE2}.  Then the remainder
  collision operator $\BB^2_\delta$ (defined in
  \eqref{eq:def-remain}) satisfies
\begin{equation}\label{eq:Ac-exp-Lq}
  \forall \, h \in L^q _v(\nu m), \quad \left\|
   \BB^2_\delta  h \right\|_{L^q _v(m)} \le
  \Lambda_{m,q}(\delta)\, \left\| h \right\|_{L^q _v(\nu m)},
\end{equation}
and
\begin{equation}\label{eq:Ac-LpLq}
  \forall \, h \in L^q_vL^p_x(\nu \, m), \quad \left\|
    \BB^2_\delta h \right\|_{L^q_v L^p_x(m)} \le
  \Lambda_{m,q}(\delta) \, 
  \left\| h \right\|_{ L^q_vL^p_x(\nu m)},
\end{equation}
where $\Lambda_{m,q}(\delta)$ is some constructive constant (depending
on $m$ and $q$) such that
\begin{itemize}
\item $\Lambda_{m,q}(\delta) \to 0$ as $\delta \to 0$ for the conditions
{\bf (W1)} and {\bf (W2)};
\item $\Lambda_{m,q}(\delta) \to \phi_q(k) $ as $\delta \to 0$
for the condition {\bf (W3)} when $m := \langle v \rangle^k$, $k > 2$, where
$$
 \phi_q(k) := \Bigl( {4 \over k+2} \Bigr)^{1/q}  \Bigl( {4 \over k-1}  \Bigr)^{1-1/q} .
 $$
\end{itemize}
\end{lem}

\begin{rem}
  Remark that $\phi_q(k)$ goes to zero when $k$ goes to $+\infty$ and
  \[
  k > k^*_q := \frac{3 + \sqrt{49 - 48/q}}{2} \quad
  \Longrightarrow \quad \phi_q(k) <1,
  \]
  by the arithmetic-geometric inequality: we have 
\begin{equation*}
  \Bigl( {4 \over k+2} \Bigr)^{1/q}  \Bigl( {4 \over k-1}
  \Bigr)^{1-1/q} \le \frac1q \frac{4}{k+2} + \left( 1 - \frac1q
  \right) \frac{4}{k-1}
\end{equation*}
and 
\begin{equation*}
  \frac1q \frac{4}{k+2} + \left( 1 - \frac1q
  \right) \frac{4}{k-1} < 1 \Longleftrightarrow k > k_q ^*. 
\end{equation*}
\end{rem}

\begin{proof}[Proof of Lemma~\ref{lem:Asmall}] 
  We analyze separately the conditions {\bf (W1)}, {\bf (W2)} and {\bf
    (W3)} on the function $m$.

\medskip
\noindent
{\it Case~{\bf (W1)}: $p=q=2$ with Gaussian weight.} Arguing as in
\cite[Proposition~2.3]{Mcmp} one can prove the following
\begin{equation}\label{eq:Aremainder-classical}
\left\|  \AA^c _\delta h \right\|_{L^2(\mu^{-1/2})} \le o(\delta) \,
\| h \|_{L^2(\mu^{-1/2})}. 
\end{equation}
Let us recall the core of the proof, which relies on the careful
inspection of the explicit bound from above on the kernel of
$\AA^c_\delta$, inspired by the celebrated calculations of Hilbert and
Grad, as reported for instance in \cite[Chapter~7,
Section~2]{MR1307620}:
\[
\left| \AA^c _\delta h (v) \right| \le \int_{\R^3} k^c _\delta
(v,v') \, |h(v')| \dd v'
\]
with (when $\mu = (2\pi)^{-3/2} e^{-|v|^2/2}$) 
\begin{multline*}
  K^c _\delta (v,v') \le C \, (1-\Theta_\delta) \, \Bigg\{
  |v-v'|^{-1} \, \exp\left[ -\frac{|v-v'|^2}{8} - \frac{\left(
        |v|^2 - |v'|^2 \right)^2}{8 |v-v'|^2} \right] 
  \\ +
  |v-v'| \, \exp \left[ - \frac{\left(|v|^2 + |v'|^2\right)}{4}
  \right] \Bigg\}
\end{multline*}
from which \eqref{eq:Aremainder-classical} is easily deduced. 

\medskip\noindent {\it Cases~{\bf (W2)} and {\bf (W3)}.}  Recall that
\cite[Proposition~2.1]{Mcmp} establishes that for the stretch
exponential weight $m = e^{\kappa \, |v|^\beta}$ it holds
 \begin{equation}\label{eq:Ac-exp-L1}
   \forall \, h \in L^1(\nu m), \quad \left\|
     \BB_\delta ^2 h \right\|_{L^1(m)} \le
   \Lambda_{m,q}(\delta) \, \left\| h \right\|_{L^1(\nu m)}, \quad  
   \Lambda_{m,q}(\delta) \mathop{\longrightarrow}_{\delta\to0} 0,
\end{equation}
where however the definition of $\Theta_\delta$ is slightly different
from ours. But it is immediate to extend the proof to the
present situation.

Estimate \eqref{eq:Ac-exp-Lq} is then obtained by piling up
\eqref{eq:Aintegral}, \eqref{eq:Apointwise}, and
\eqref{eq:Apointwise-exp}, and using the Riesz-Thorin interpolation
theorem in order to obtain the $L^q$ estimate when $1 < q <\infty$.

\smallskip\noindent {\it Estimate \eqref{eq:Ac-LpLq}.}  Now observe
that all the estimates previously established on $\BB^2_\delta$ are
valid (with the same proofs) for $\tilde \BB^2_\delta$. Then, since
$\tilde \BB^2 _\delta$ is a nonnegative operator acting only in $v$, 
we have
\[
%\int_{\T^3} \left| \BB_\delta^2 h \right| \dd x \le 
\int_{\T^3} \left|\tilde \BB_\delta^2 h \right| \dd x \le 
\tilde \BB^2 _\delta \left( \int_{\T^3} |h| \dd x \right)
\]
and 
\[
\sup_{x \in \T^3} \left|\tilde \BB_\delta^2 h \right| \le \tilde
\BB^2 _\delta \left( \sup_{x \in \T^3} |h| \right)
\]
and therefore by interpolation
\begin{equation}\label{eq:interpA}
\left\|\tilde \BB_\delta^2 h \right\|_{L^p _x} \le \tilde
\BB^2 _\delta \left( \|h\|_{L^p_x} \right)
\end{equation}
for any $p \in [1,+\infty]$. We then conclude thanks to
\eqref{eq:Ac-exp-Lq} (used on $\tilde \BB^2 _\delta$): 
\begin{equation*}
  \left\|\BB^2_\delta h \right\|_{L^q_v L^p_x(m)} \le \left\|
    \tilde \BB^2_\delta h \right\|_{L^q_v L^p_x(m)} \le \left\|
    \tilde \BB^2_\delta \left( \n{h}_{L^p_x}\right)
  \right\|_{L^q_v(m)} \le \Lambda_{m,q}(\delta) \, \n{h}_{L^q_v
    L^p_x(m)}. 
\end{equation*}
\end{proof}

Let us now prove dissipativity estimates for the operator
$\BB_\delta$.
\begin{lem}\label{lem:decomp-new}
  Consider a weight $m$ and the space $\EE := W^{\sigma,q}_v W^{s,p}_x
  (m)$ with $p,q \in [1,+\infty]$ and $\sigma, s \in \N$, $\sigma \le
  s$. Then:
 \begin{itemize}
 \item[\textbf{(W1)}] When $m = \mu^{-1/2}$, $p=q=2$, there is
   $\lambda_0=\lambda_0(m,\delta) \in (0,\nu_0)$ such that
   $\lambda_0(m,\delta) \to \nu_0$ as $\delta \to 0$ and $(\BB_\delta
   + \lambda_0)$ is dissipative in $\EE$.
   \\
 \item[\textbf{(W2)}] When $m = e^{ \kappa \, |v|^\beta }$, $\kappa > 0$,
   $\beta \in (0,2)$ and $p,q \in [1,+\infty]$, there is
   $\lambda_0=\lambda_0(m,\delta) \in (0,\nu_0)$ such that
   $\lambda_0(m,\delta) \to \nu_0$ as $\delta \to 0$ and $(\BB_\delta
   +\lambda_0)$ is dissipative in $\EE$.
   \\
 \item[\textbf{(W3)}] When $m = \langle v \rangle^k$ with any $p,q \in
   [1,+\infty]$ and $k > k_q^*$, there is
   $\lambda_0=\lambda_0(k,q,\delta) \in (0,\nu_0)$ such that
\[
\left\{ 
\begin{array}{l}
  \lambda_0(k,q,\delta) \to
   \lambda_0 ^*(k,q) \in (0,\nu_0) \ \mbox{ when } \ \delta \to 0,
  \vspace{0.2cm} 
  \\
%  \bar \lambda_0(k,q) \to 0 \ \mbox{ when } \ k\to
%  \bar k_q, \vspace{0.2cm} 
%  \\ 
  \lambda_0^*(k,q) \to \nu_0 \ \mbox{ when } \ k\to +\infty,
\end{array}
\right.
\]
and $(\BB_\delta +\lambda_0)$ is dissipative in $\EE$.
 \end{itemize}
\end{lem}

\begin{rem}
  As in the previous statements, $k >k_q^*$ could be relaxed down to
  $k > k_q ^{**}$. 
\end{rem}

% \begin{rem} 
% In the case (iii) the choice of Lebesgue exponents $q \ge p$ is also
% possible, at the price of an additionnal condition $k \ge O(\nu_0
% ^{1-q})$. 
% \end{rem}

% \begin{rem} The lemma is mainly based on a smallness estimate on $\bar
%   \AA^c _\delta$ in appropriate spaces.  Such a smallness estimate on
%   $\bar\BB^2 _\delta$ was proved in $L^1(m)$ with stretched
%   exponential weight $m$ and in $L^2(\mu^{-1/2})$ in~\cite[Propositions
%   2.1, 2.3]{Mcmp}, see also \cite{Grad2} for the latter case.  We
%   extend this estimate to $L^\infty(m)$ with stretched exponential and
%   polynomial weight $m$ and in $L^1(m)$ with polynomial weight $m$.
%   That last result makes use of the so-called Povzner inequality.
% \end{rem}

\begin{proof}[Proof of Lemma~\ref{lem:decomp-new}] 
  We consider separately each case. Observe first that the
  $x$-derivatives commute with the operator $\BB_\delta$, therefore
  without restriction we do the proof for $s=0$. 
\medskip

\noindent {\it Case~{\bf (W1)}: $p=q=2$ with Gaussian weight.}  
We consider a solution $h_t$ to the linear equation
$$
\partial_t h_t = \BB_\delta \, h_t = \BB^2_\delta h_t - \nu \, h_t - v
\cdot \nabla_x h_t,
$$
with given initial datum $h_0$. We consider first $\sigma = 0$, and we
calculate
\[
\dt \| h_t \|^2_{L^2(\mu^{-1/2})} \le 2 \int_{\T^3 \times \R^3} 
\left| \BB^2 _\delta h\right| \, |h| \dd x \dd v - 2\int_{\T^3 \times
  \R^3} h^2 \, \nu \dd x \dd v 
\]
since the term involving $v \cdot \nabla_x$ cancels from its
divergence (in $x$) structure. This implies
\[
\dt \| h_t \|^2_{L^2(\mu^{-1/2})} \le - 2 \left(\nu_0-o(\delta)\right) \,
\| h_t \|^2_{L^2(\mu^{-1/2})}
\]
and concludes the proof of dissipativity. Since the $x$-derivatives commute with the
equation we have in the same manner
\[
\dt \left\| \nabla_x^{s} h_t
  \right\|^2_{L^2(\mu^{-1/2})} 
  \le - 2 \left(\nu_0-o(\delta)\right) \, \left\| \nabla_x^{s}h_t
  \right\|^2_{L^2(\mu^{-1/2})}.
\]

Then we consider the case of derivatives in $v$, say first $\sigma=1$
and $s \ge 1$. Note that we can reduce to the case $s=1$ by
differentiating in $x$ the equation (using that in the definition of
the norms \eqref{eq:def-norms} we sum over derivatives
$\partial^i_v \partial^j_x$ with $|i| \le \sigma$, $|j| \le s$,
$|i|+|j|\le \max\{\sigma; s\}$).  We compute the evolution of the
$v$-derivatives:
\begin{eqnarray*}
  \partial_t \partial_v h 
  &=& - v \cdot  \nabla_x \partial_v h - \partial_x h + \partial_v (\BB^2_\delta h - \nu h)   
  \\
%  &=& - v \cdot  \nabla_x \partial_v h - \partial_x h 
%  +  Q^+(\partial_v h, \mu) + Q^+(\mu,\partial_v h ) + Q^+(h,\partial_v\mu) 
%  \\
%  &&+Q^+(\partial_v\mu,h) - Q^-(\partial_v \mu,h) 
%  - Q^-( \mu,\partial_v h) -   \partial_v \left( \bar\AA_\delta (h)
%  \right) 
  \\
  &=& \BB_\delta ( \partial_v h) - \partial_x h + \RR h
\end{eqnarray*}
with 
\begin{equation}\label{eq:RRdef}
  \RR h := Q(h,\partial_v\mu) + Q(\partial_v\mu,h) - \left( \partial_v \AA_\delta\right)  (h)
  + \AA_\delta (\partial_v h) ,
\end{equation}
($\left( \partial_v \AA_\delta\right) (h)$ means that one
differentiates the kernel of the operator as opposed to its argument
$h$) where we have used twice the relation
\[
\BB^2_\delta h = Q^+(h,\mu) + Q^+( \mu,h) - Q^-(h,\mu) -
\AA_\delta(h),
\]
and the property 
\begin{equation}\label{eq:dvQ=Qdv}
\partial_v Q^\pm(f,g) = Q^\pm(\partial_vf, g) +
Q^\pm(f,\partial_v g)
\end{equation}
following from the translation invariance of the collision operator. 
We deduce that 
\begin{multline*}
  \dt \| \nabla_v h \|^2_{L^2(\mu^{-1/2})} \le - 2
  \left(\nu_0-o(\delta)\right) \, \| \nabla_v h \|^2_{L^2(\mu^{-1/2})}
  \\
  - \int_{\T^3 \times \R^3} \nabla_v h \cdot \nabla_x h \, \mu^{-1}
  \dd x \dd v + \|\RR h \|_{L^2(\mu^{-1/2})} \, \| \nabla_v h
  \|_{L^2(\mu^{-1/2})}.
\end{multline*}
Using one integration by parts and the regularizing property of the
operator $\AA_\delta$, we have
\[
\left\| \left(  \AA_\delta\right) (\partial_v h)
\right\|_{L^2(\mu^{-1/2})} ^2
+
\left\| \left( \partial_v \AA_\delta\right) (h)
\right\|_{L^2(\mu^{-1/2})} ^2 \le C \, \| h \|^2 _{L^2(\mu^{-1/2})} 
\] 
for some constant $C=C_\delta >0$ (depending on $\delta$). Moreover
using the computation of Hilbert and Grad (see above or again
\cite[Chapter~7, Section~2]{MR1307620}), we have
\[
\left\| Q^+(h,\partial_v\mu) + Q^+(\partial_v\mu,h) - Q^-(\partial_v
\mu,h) \right\|_{L^2(\mu^{-1/2})} ^2 \le C \, \| h \|_{L^2(\mu^{-1/2})} ^2
\]
for some constant $C>0$. Therefore the operator $\RR$ is bounded in
$L^2(\mu^{-1/2})$. Introducing the norm 
 $$
 \| h \|_{H^1_{x,v}(\mu^{-1/2})_\var} := \left( \| h \|^2_{L^2(\mu^{-1/2})} + \| \nabla_x h
 \|^2_{L^2(\mu^{-1/2})} + \var \, \| \nabla_v h \|^2_{L^2(\mu^{-1/2})} \right)^{1/2}
 $$
 for some given $\var >0$, we deduce
 \begin{multline*} 
\dt \| h \|^2 _{H^1_{x,v}(\mu^{-1/2})_\var} 
  \\ \le - 2  \left(\nu_0-o(\delta)\right) \, \left( \| h \|^2_{L^2(\mu^{-1/2})} +
     \left\| \nabla_x h \right\|^2_{L^2(\mu^{-1/2})} 
 + \var \, \left\|
       \nabla_v h \right\|^2_{L^2(\mu^{-1/2})}\right)
         \\ 
     + \var \, \left\| \nabla_v h \right\|_{L^2(\mu^{-1/2})} \,
     \left\| \nabla_x h \right\|_{L^2(\mu^{-1/2})} 
 + C \, \var \,
     \left\| \nabla_v h \right\|_{L^2(\mu^{-1/2})} \, \left\| h
     \right\|_{L^2(\mu^{-1/2})}
     \\
      \le - 2
   \left(\nu_0-o(\delta) - o(\eta) \right) \, \left( \| h \|^2_{L^2(\mu^{-1/2})} +
     \left\| \nabla_x h \right\|^2_{L^2(\mu^{-1/2})} 
     + \var \, \left\|
       \nabla_v h \right\|^2_{L^2(\mu^{-1/2})}\right) \\ 
   \le - 2
   \left(\nu_0-o(\delta) - o(\eta) \right) \, \| h \|^2
   _{H^1_{x,v}(\mu^{-1/2})_\var} 
\end{multline*} 
which concludes the proof by taking $\var$ small enough in terms of
$\delta$. The higher-order estimates can be performed with the norm
$$
\| h \|_{W^{\sigma,2}_v W^{2,s}_{x}(\mu^{-1/2})_\var} := \left(
  \sum_{|i| \le \sigma, \ |j| \le s, \ |i|+|j| \le \max\{\sigma; s\}}
  \var^{|i|} \, \left\| \partial ^i_v \partial ^j _x h \right\|^2_{L^2(\mu^{-1/2})}
\right)^{1/2}
 $$
for some $\var$ to be chosen small enough (in terms of $\delta$). 
\medskip

\noindent {\it Cases {\bf (W2)} and {\bf (W3)}: $p,q \in [1,+\infty]$ with
  stretched exponential and polynomial weights.} The proof of these
two cases are identical. We denote by $m$ either a polynomial weight or
a stretched exponential weight, using the respective estimates
established previously.

We consider again only the case $s=0$ since $x$-derivatives commute
with the equation, and we also look first at the case $\sigma=0$.

Consider first $1 \le p,q < +\infty$ and denote $\Phi'(z) := |z|^{p-1}
\, \hbox{sign} (z)$. We compute
\begin{multline*}
  \dt \| h_t \|_{L^q_v L^p_x(m)} = \| h \|_{L^q_v L^p_x(m)}
  ^{1-q} \, \times 
  \\
   \left( \int_{\R^3}
  \left( \int_{\T^3} (\BB_\delta (h) ) \, \Phi'( h) \dd x  \right) \,
  \left( \int_{\T^3} | h |^p \dd x \right)^{\frac{q}{p}-1} \, m^q \dd v \right).
\end{multline*}

Observing that
\begin{multline}\label{eq:mixB}
  \int_{\T^3} (\BB_\delta (h) ) \, \Phi'( h) \dd x = \int_{\T^3} \left[ 
  (\BB^2_\delta (h) ) \, \Phi'( h) - \nu \, |h|^p - \frac1p \, v
  \cdot \nabla_x \left( |h|^p \right) \right]  \dd x
  \\
  \le \left( \int_{\T^3} \left| \BB_\delta^2 (h) \right|^p \dd x \right)^{\frac{1}{p}}  \,
  \left( \int_{\T^3} | h |^p \dd x \right)^{1- \frac{1}{p}} - \nu \int_x |h|^p \dd x ,
\end{multline}
we deduce that 
\begin{multline}
\dt \| h_t \|_{L^q_v L^p_x(m)} 
  \\
   \le \| h \|_{L^q_v L^p_x(m)}
  ^{1-q} \, \left[ \left( \int_{\R^3} \left\|
     \BB_\delta ^2 (h) \right\|_{L^p_x} \, \| h \|_{L^p_x} ^{q-1}
  \, m^q \dd v \right) - \left(
  \int_{\R^3} \nu \, \| h \|_{L^p_x} ^q \, m^q \dd v \right) \right].
\end{multline}

Denoting $H = \|h\|_{L^p_x}$, we obtain thanks to \eqref{eq:interpA}
\begin{multline*}
\dt \| h_t \|_{L^q_v L^p_x(m)} 
  \\
   \le \| h \|_{L^q_v L^p_x(m)}
^{1-q} \, \left[ \left( \int_{\R^3} \tilde \BB_\delta ^2(H) \, \nu^{-1/q'} \, m \, H ^{q-1}
    \, m^{q-1} \, \nu^{1/q'}\dd v \right) - \int_{\R^3} \nu \, H ^q \, m^q \dd v 
\right] 
  \\
\le
\| h \|_{L^q_v L^p_x(m)}
^{1-q} \, \left[
 \left\| \tilde \BB_\delta ^2 (H) \right\|_{L^q _v(m \nu^{-1/q'})} \,  \| H\|_{L^q_v (m \nu^{1/q})}
^{q-1} -  \int_{\R^3} \nu \, H ^q \, m^q \dd v  \right]
.
% \left( \int_{\R^3} \tilde \AA_\delta ^c H \dd v \right) - \left(
%   \int_{\R^3} \nu \, H \dd v \right).
\end{multline*}
Using then \eqref{eq:Ac-LpLq}  and   
\[
\| h \|_{L^q_v L^p_x(m)}   \le \nu_0^{-1/q} 
\| h \|_{L^q_v L^p_x(m \nu^{1/q} )}, 
\]
we  finally deduce that
\begin{eqnarray}\nonumber
  \dt \| h_t \|_{L^q_v L^p_x(m)} 
 & \le &
  \| h \|_{L^q_v L^p_x(m)}^{1-q} \, 
\left[ \Lambda_{m\nu^{-1/q'},q} (\delta) - 1 \right] \,   
  \| h \|_{L^q_v L^p_x(m\nu^{1/q})}^q
\\ \label{eq:BL:gaindissip}
  &\le&  \nu_0^{1/q-1} \left[ \Lambda_{m\nu^{-1/q'},q} (\delta) - 1 \right] \,     \| h \|_{L^q_v L^p_x(m\nu^{1/q})}
\\ \nonumber
  &\le& - \nu_0^{-1} \,[ 1 - \Lambda_{m \nu^{-1/q'},q} (\delta)]  \,  \| h_t \|_{L^q_v L^p_x(m)}, 
\end{eqnarray}
which concludes the proof of dissipativity in this case.

The cases $p = +\infty$ and $q = +\infty$ are then obtained by taking
the corresponding limits in the above estimate. The $v$-derivatives
can be treated with the same line of arguments as in the case {\bf
  (W1)}. Arguing as before we obtain
\begin{multline*} 
\dt \left( \| h \| _{L^q_v L^p_x(m)} + \|   \nabla_x h \| _{L^q_v L^p_x(m)} \right)
 \\
    \le
  - \nu_0^{1/q-1} \, [ 1 - \Lambda_{m \nu^{-1/q'},q} (\delta)]   \, 
  \Bigl( \| h_t \|_{L^q_v L^p_x(m\nu^{1/q})}  +  \| \nabla_x h_t \|_{L^q_v L^p_x(m\nu^{1/q})} \Bigr)
%  \\
%  \left[ \left( \int_{\R^3} \|h \|_{L^p_x} ^q \, \nu \, m^q \dd v
%    \right)^{\frac{1}{q}} + \left( \int_{\R^3} \| \nabla_x h
%      \|_{L^p_x} ^q \, \nu \, m^q \dd v \right)^{\frac{1}{q}} \right]
\end{multline*}
and 
\begin{multline*}
\dt \| \nabla_v h \| _{L^q_v L^p_x(m)} \\ \le
- \nu_0^{1/q-1} \,[ 1 - \Lambda_{m \nu^{-1/q'},q} (\delta)]  \,   \| \nabla_v h \| _{L^q_v L^p_x(m)}  
   + \|\nabla_x h \|_{L^q_v L^p_x(m)}
  + \|\RR h \|_{L^q_v L^p_x(m)}, 
\end{multline*}
where $\RR$ is defined in \eqref{eq:RRdef}. Using the
Lemmas~\ref{lem:Aintegral} and \ref{lem:decQ+rad} when $m$ is a
polynomial weight, and \eqref{eq:Ac-exp-L1} and
Lemma~\ref{lem:decQ+rad-exp} when $m$ is an exponential weight, and
the regularization property of the operator $\AA_\delta$, we prove
that
\[
\|\RR h \|_{L^q_v L^p_x(m)} \le C \, \left( \int_{\R^3} \| h_t
  \|_{L^p_x} ^q \, \nu \, m^q
  \dd v \right)^{\frac{1}{q}} , 
\] 
for some constant $C=C_\delta >0$  (depending on $\delta$). We then
introduce the norm
 $$
 \| h \|_{W^{1,q}_v W^{1,p}_x(m)_\var} := \| h \|_{L^q_v L^p_x(m)} + \| \nabla_x h
 \|_{L^q_v L^p_x(m)} + \var \, \| \nabla_v h \|_{L^q_v L^p_x(m)},
 $$
 for some $\var>0$ to be fixed later, and we deduce 
 \begin{multline*} 
\dt  \| h \|_{W^{1,q}_v W^{1,p}_x(m)_\var} \le
   - \nu_0^{1/q-1} \,[ 1 - \Lambda_{m \nu^{-1/q'},q} (\delta)] \, 
   \Bigg[ \left( \int_{\R^3} \|h \|_{L^p_x} ^q \, \nu \, m^q \dd v
     \right)^{\frac{1}{q}} 
  \\
   + \left( \int_{\R^3} \| \nabla_x h
       \|_{L^p_x} ^q \, \nu \, m^q \dd v \right)^{\frac{1}{q}} + \var \, \left(
       \int_{\R^3} \| \nabla_v h \|_{L^p_x} ^q \, \nu \, m^q
       \dd v \right)^{\frac{1}{q}} \Bigg]
  \\
   + C \, \var \, \left( \int_{\R^3} \|h \|_{L^p_x} ^q \, \nu \, m^q \,
     dv \right)^{\frac{1}{q}} + \var \, \| \nabla_x h \|_{L^q_v
     L^p_x(m)} 
  \\
   \le
   - \left(\nu_0^{1/q-1} \,[ 1 - \Lambda_{m \nu^{-1/q'},q} (\delta)]  - o(\var) \right) \, 
   \Bigg[ \left( \int_{\R^3} \|h \|_{L^p_x} ^q \, \nu \, m^q \dd v
     \right)^{\frac{1}{q}} 
  \\
   + \left( \int_{\R^3} \| \nabla_x h
       \|_{L^p_x} ^q \, \nu \, m^q \dd v \right)^{\frac{1}{q}} + \var \, \left(
       \int_{\R^3} \| \nabla_v h \|_{L^p_x} ^q \, \nu \, m^q \dd v
     \right)^{\frac{1}{q}} \Bigg] \\ 
     \le
   - \left(\nu_0^{1/q-1} \,[ 1 - \Lambda_{m \nu^{-1/q'},q} (\delta)]  - o(\var)
   \right) \, \| h \|_{W^{1,q}_v W^{1,p}_x(m)_\var}
\end{multline*} 
which concludes the proof by taking $\var$ small enough in terms of
$\delta$. The higher-order estimates are performed with the norm
$$
 \| h \|_{W^{s,q}_v W^{s,p}_x(m)_\var} := \sum_{|i|\le \sigma, \ |j|
   \le s, \ |i|+|j|\le \max\{\sigma; s\}}
 \var^{|i|} \, \left\| \partial_v  ^i \partial ^j _x h \right\|_{L^q_v L^p_x(m)}
 $$
for some $\var >0$ to be chosen small enough (in terms of $\delta$). 
\end{proof}

% \textcolor{red}{J'ai verifie la preuve en details, tout me semble
%   juste, mais la fin est redigee de maniere peut-etre un peu trop
%   compacte, a voir eventuellement plus tard en fonction des retours.}

\subsection{Regularization estimates in the velocity variable}
\label{sec:veloc-regul-estim}

In this subsection we prove a regularity estimate on the truncated
operator $\AA_\delta$, which improves the result
~\cite[Proposition~2.4]{Mcmp}. In the latter paper, it was established
in~\cite[Proposition 2.4 (iii)]{Mcmp}, for a slightly weaker
truncation function $\Theta_\delta$ (and the same proof would apply
here), the boundedness of the operator $\AA_\delta$ from $L^1(\langle
v \rangle^\gamma)$ into the space of $W^{1,1}_v $ functions with
compact support. We prove here:
\begin{lem}
  \label{lem:regA}
  For any $s \in \N$ the operator $\AA_\delta$ maps $L^1_v
  (\langle v \rangle)$ into $H^s_v$ functions with compact support,
  with explicit bounds (depending on $\delta$) on the $L^1 _v (\langle
  v \rangle) \to H^s _v$ norm and on the size of the support.

  More precisely, there are two constants $C_{s,\delta}, R_\delta>$ so that 
  $$
  \fa h \in L^1_v(\langle v \rangle), \quad \hbox{{\em supp}} \,
  \AA_\delta h \subset B(0,R_\delta), \quad \| \AA_\delta h
  \|_{H^s_v} \le C_{s,\delta} \, \| h \|_{ L^1_v(\langle v \rangle)
  }.
  $$
\end{lem}

\begin{proof}[Proof of Lemma~\ref{lem:regA}] 
  On the one hand, it is clear that the range of the operator
  $\AA_\delta$ is included into compactly supported functions thanks
  to the truncation, with a bound on the size of the support related
  to $\delta$.

  On the other hand, the proof of the smoothing estimate is a
  straightforward consequence of the regularization property of the
  gain part $Q^+$ of the collision operator discovered by P.-L. Lions
  \cite{MR1284432,MR1295942}, and we only sketch it. Let us recall
  that
\[
\AA_\delta h = Q^+_{B_\delta}(\mu,h) + Q^+_{B_\delta}(h,\mu) 
- Q^-_{B_\delta} (\mu,h)
\]
where $Q^+_{B_\delta}$ (resp. $Q^-_{B_\delta}$) is the gain
  (resp. loss) part of the collision operator associated to the
mollified collision kernel $B_\delta = \Theta_\delta \, B$. More
precisely, we have
\[
Q^+ _{B_\delta} (f,g) := \int_{\R^3} \int_{\mathbb{S}^2} \Theta_\delta
\, f(v') \, g(v'_*) \, |v-v_*|^\gamma \, b(\cos \theta) \dd v_* \dd \sigma
\]
and, since we can decompose the truncation as $\Theta_\delta =
\Theta^1_\delta(v) \, \Theta^2_\delta(v-v_*) \,
\Theta^3_\delta(\cos\theta)$, we have the formula
\begin{eqnarray*} Q^-_{B_\delta} (\mu,h) &:=&
\int_{\R^3} \int_{\mathbb{S}^2} \Theta_\delta
\,\mu(v)   \, h(v_*) \,  |v-v_*|^\gamma \, b(\cos \theta) \dd v_* \dd \sigma 
  \\
&=& \mu(v) \, \Theta^1_\delta(v) \, (f * \nu_\delta) (v), \quad
\nu_\delta \in C_c(\R^3).
\end{eqnarray*}
The regularity estimate is trivial for $Q^-_{B_\delta} (\mu,h)$ thanks
to the truncation and convolution structure, and the regularity
estimate for $Q^+_{B_\delta}$ follows immediately from the result
discovered in \cite{MR1284432,MR1295942}  in the form proven in
\cite[Theorem~3.1]{MR2081030}.
\end{proof}

% \textcolor{red}{A nouveau la preuve est tres compacte, mais on verra
%   plus tard si ca pose soucis.}

\subsection{Iterated averaging lemma}
\label{sec:new-class-iterated}

In this subsection we prove the key regularity results for our
factorization and enlargement theory.  We begin with an ``averaging
lemma'' (in the spirit of \cite{MR923047,MR1949176}) for the free
transport equation. This first result requires regularity in the
velocity variable. We shall then show how to get rid of the assumption
by a new \textit{iterated} averaging lemma.

\begin{lem}\label{lem:BLaveraginglemma}
  Consider $f \in L^1([0,T]; L^1(\T^d \times \R^d))$ and $\initem{f} \in
  L^1(\T^d \times \R^d))$ such that $\nabla_v \initem{f} \in L^1(\T^d
  \times \R^d))$ and (in the weak sense)
$$
\partial_t f + v \cdot \nabla_x f = 0 \quad \hbox{on}\quad [0,T)
\times \T^d \times \R^d, \quad f_{|t=0} = \initem{f} \quad
\hbox{on}\quad \T^d \times \R^d.
$$
For any fixed $\varphi \in \DD(\R^d)$, let us define 
$$
\rho_\varphi (t,x) := \int_{\R^d} f_t(x,v) \, \varphi(v) \dd v.
$$
Then $\rho_\varphi$ satisfies
\begin{equation}\label{eq:transfert}
\left\|\rho_\varphi (t,\cdot) \right\|_{W^{1,1}_x} \le \left( 1 + {1 \over t}
\right) \, \|\varphi \|_{W^{1,\infty}} \, \left( \|\initem{f}
\|_{L^1_{x,v}}+ \|\nabla_v \initem{f} \|_{L^1_{x,v}}\right).
\end{equation}
\end{lem}

\begin{rem} It is worth mentioning that a similar result holds in
  $L^2$.  It may be compared with the classical averaging lemma for
  the free transport equation: a typical statement (see
  \cite{MR1798557,MR1669221} as well as
  \cite{MR923047,MR1127927,MR1634024,MR2597793} and the references
  therein for more details) is
\begin{equation}\label{eq:ClassicAverageLemma}
  \left\|\rho_\varphi (t,\cdot) \right\|_{H^{1/2}_x} \le  (1 + t) \, \|\varphi \|_{W^{1,\infty}}  \, \|f_0 \|_{L^2_{x,v}}. 
\end{equation}
Hence the gain of derivability in the $x$ variable is weaker compared
to~\eqref{eq:transfert}, but there is no regularity assumption on the
initial datum. However, it is well known that
\eqref{eq:ClassicAverageLemma} is false for $p=1$ (see the discussion
in \cite{MR923047} and the related work \cite{MR1903763}). In the
estimate~\eqref{eq:transfert} we can cover the critical $L^1$ case at
the price of assuming more initial regularity on the velocity
variable. It shares some similarity with the results in
\cite{MR1949176}. The proof makes use of the ``gliding norms''
introduced in \cite{MR2863910}.
\end{rem}

\begin{proof}[Proof of Lemma~\ref{lem:BLaveraginglemma}]
Introducing the differential operator 
\begin{equation}\label{def:Dt}
D_t := t \nabla_x + \nabla_v, 
\end{equation}
we observe that $D_t$ commutes with the free transport operator
$\partial_t + v \cdot \nabla_x$, so that
\[
\partial_t (D_t f) + v \cdot \nabla_x (D_t f) = 0.
\]
From the mass preservation for the free transport flow on $f_t$ and $D_t
f_t$, we deduce
$$
\fa t \ge 0, \quad \|f _t \|_{L^1} = \|f_0 \|_{L^1}, \quad \left\|D_t f _t \right\|_{L^1} =
\left\|D_0 f_0 \right\|_{L^1} = \left\|\nabla_v f _0 \right\|_{L^1}.
$$
Finally we calculate
\begin{eqnarray*}
  \nabla_x \rho_\varphi(t,x)  &= & \int_{\R^d} \left( {D_t \over t}- \nabla_v \right) f_t(x,v) \, \varphi(v) \dd v
  \\
  &= & {1 \over t}\int_{\R^d}   \left( D_t   f\right) (t,x,v) \, \varphi(v) \dd v + \int_{\R^d}   f (t,x,v) \, \nabla_v  \varphi(v) \dd v ,
 \end{eqnarray*}
and we conclude the proof thanks to the previous estimates. 
\end{proof}

Let us recall the notation $T_n (t) := (\AA_\delta
S_{\BB_\delta})^{(*n)}$ for $n \ge 1$, where $S_{\BB_\delta}(t)$ is
the semigroup generated by the operator $\BB_\delta$. We remind the
reader that the $T_n(t)$ operators are merely time-indexed family of
operators which do not have the semigroup property in general.

%-------------------  lemme T -------------------------------------------- 

\begin{lem}
  \label{lem:regularization-lbe} 
  Consider $s \in \R_+$, % $p,q \in [1,+\infty]$
  and a weight $m$ so that the assumptions of Lemma
  \ref{lem:decomp-new} are satisfied (hence $\BB_\delta$ is
  dissipative in $W^{s',1}_{x,v}(m)$ for $s' \in [0,s+4] \cap \N$). 

  Then the time indexed family $T_n$ of operators satisfies the following: for
  any $\lambda_0 ' \in (0,\lambda_0)$ where $\lambda_0$ is provided by
  Lemma~\ref{lem:decomp-new}, there is some constructive constants
  $C=C(\lambda_0',\delta)>0$ and $R=R(\delta)$ such that for any $t \ge 0$
 % that for any there is
 %  $C= C(s,\sigma,p,q,m) > 0$ and $a = a(s,\sigma,p,q,m) \in
 %  (-\nu_0,0)$ such that
 \[
\hbox{{\em supp}} \,  T_n(t) h  \subset K := B(0,R), 
\]
and 
 \begin{eqnarray} \label{estim:T1gainW} 
& \fa t \ge 0, &  
\left\| T_1(t) h \right\|_{W^{s+1,1}_{x,v} (K)} \le C \,
    { e^{-\lambda_0' \, t}  \over t}\, \| h \|_{W^{s,1}_{x,v}(m)}, \quad \hbox{if} \,\, s \ge 1; 
\\
\label{estim:T2gainW} 
& \fa t \ge0, &
\left\| T_2(t) h \right\|_{W^{s+1/2,1}_{x,v} (K)} \le C \,
    e^{-\lambda_0' \, t} \, \| h \|_{W^{s,1}_{x,v}(m)},  \quad \hbox{if} \,\, s \ge 0. 
\eear
\end{lem}

\begin{rem}
  Our proof extends verbatim to the case of $W^{s,p}_{x,v}$ spaces in
  \eqref{estim:T2gainW}, with $p \in [1,+\infty)$. The important
  aspect of our estimates is the optimal time
  decay. % This iterated kind of averaging lemma does not seem known
  % in the literature in the $L^1$ case which is notoriously degenerate,
  % see for instance \cite{MR923047,MR1903763}.
  The core idea is to exploit correctly the combination of a
  $v$-regularizing operator $\AA_\delta$ and a transport semigroup
  $\SS_{\BB_\delta}$. However the usual averaging lemma degenerate in
  $L^1$, where only a mere compactness property in space is
  retained. We here show that by using the propagation of a
  \emph{time-dependent phase space} regularity (thanks to the
  introduction of the operator $D_t$), one can still keep track of
  some velocity regularity, and \emph{transfer} it to the space
  variable, while preserving at the same time the correct time decay
  asymptotics.  
\end{rem}

\begin{proof}[Proof of Lemma~\ref{lem:regularization-lbe}]
  Let us consider $h \in W^{s,1} _{x,v}(m)$, $s \in \N$. We have
  from Lemma~\ref{lem:regA} and the fact that the $x$-derivatives
  commute with $T_1(t)$:
 \[
   \left\| T_1(t) h \right\|_{W^{s,1} _x
    W^{s+1,1}_v (K)} = \left\| \AA_\delta \, S_{\BB_\delta}(t) \, h_0
  \right\|_{W^{s,1} _x W^{s+1,1}_v (K)} \le C \, \left\| S_{\BB_\delta}(t) \,
  h \right\|_{W^{s,1} _{x,v} (m)}.
\]
Using that $\BB+\lambda_0$ is dissipative in $W^{s,1} _{x,v}(m)$, with
$\lambda_0 >0$, from Lemma~\ref{lem:decomp-new}, we get
\begin{equation}\label{TL1xW1v}
\|T_1(t) h \|_{W^{s,1} _x W^{s+1,1}_v (K)} \le C \, e^{-\lambda_0 \, t}
\, \| h \|_{W^{s,1} _{x,v} (m)}.  
\end{equation}
  
Assume now $h \in W^{s,1} _x W^{s+1,1}_v(m)$ and consider the function
$g_t = S_{\BB_\delta}(t) (\partial_x ^\alpha h)$, for any $|\alpha|\le
s$. Such function satisfies
\[
\partial_t g_t + v \cdot \nabla_x g_t = Q(\mu,g_t) + Q(g_t,\mu) -
\AA_\delta g_t.
\]

Using (1) that the operator $D_t$ defined in \eqref{def:Dt}
commutes with the free transport equation, and (2) the translation
invariance property \eqref{eq:dvQ=Qdv} of the collision operator, we
have
\begin{eqnarray*}
\partial_t (D_t g_t) + v \cdot \nabla_x (D_t g_t) 
&=& Q(\nabla_v \mu,g_t) + Q(g_t,\nabla_v \mu) 
\\
&&+ \,\, Q(\mu, D_t g_t) + Q(D_t g_t, \mu) - D_t \left(\AA_\delta g_t\right).
\end{eqnarray*}
With the notation of \eqref{eq:def-barAAd}, we rewrite the last term
as
\begin{eqnarray*}
  D_t \left( \AA_\delta g_t \right)
  &=& D_t \int_{\R^3} k_\delta(v,v_*) \, g_t(v_*) \dd v_* 
  \\
  &=& \int_{\R^3} \nabla_v k_\delta(v,v_*) \, g_t(v_*) \dd v_* 
  - \int_{\R^3} k_\delta(v,v_*) \, \nabla_{v_*} g_t(v_*) \dd v_*  
  \\
  &&\qquad
  + \int_{\R^3} k_\delta(v,v_*)  \, (D_t g_t)(v_*) \dd v_*  
  \\
  \\
  &=& \AA^1_\delta g_t + \AA^2 _\delta g_t + \AA_\delta (D_t g_t), 
\end{eqnarray*}
where we have performed one integration by part in the term of the
middle and where $\AA^1_\delta$ stands for the integral operator
associated with the kernel $\nabla_v k_\delta$ and $\AA^2_\delta$
stands for the integral operator associated with the kernel $
\nabla_{v_*} k_\delta$.  All together, we may write
\begin{equation}\label{eq:expBtDt}
  \partial_t (D_t g_t) = \BB_\delta (D_t g_t) + \JJ_\delta (g_t)
\end{equation}
with 
$$
\JJ_\delta f := Q(\nabla_v \mu,f) + Q(f,\nabla_v \mu) + \AA^1 _\delta f
+\AA^2 _\delta f.
$$
On this last term we have the following $\delta$-dependent estimate
obtained by gathering Lemmas~\ref{lem:Aintegral} and \ref{lem:regA}:
\[
\| \JJ_\delta f \|_{L^1(m)} \le C_\delta \, \| f \|_{L^1( \nu m )}.
\]
% Using the decomposition $B = \Theta_\delta \, B + (1-\Theta_\delta ) \, B$ in order to estimate the two first terms, we 
% have 
% $$
% \| \JJ f \|_{L^1(m)} \le \eta_\delta \, \| f \|_{L^1(m \, \nu)} +
% C_\delta \, \| f \|_{L^1(m )},
% $$
% with $\eta_\delta \to 0$ (and $C_\delta \to 0$) when $\delta\to0$. 

Then arguing as in Lemma~\ref{lem:decomp-new}, we have 
$$
{{\rm d} \over {\rm d}t} \int_{\T^3 \times \R^3} |D_t g_t| \, m \dd x
\dd v \le - \frac{\lambda_0}{\nu_0} \, \int_{\T^3 \times \R^3} |D_t
g_t| \, \nu \, m \dd x \dd v + C \, \| g_t \|_{L^1( \nu m )}
$$
and 
$$
{{\rm d} \over {\rm d}t} \int_{\T^3 \times \R^3} |g_t| \, m \dd x \dd v \le
- \frac{\lambda_0}{\nu_0} \, 
\int_{\T^3 \times \R^3} |g_t| \, \nu \, m \dd x \dd v.
$$

Combining that last two differential inequalities we obtain, for any
$\lambda_0 ' \in (0,\lambda_0)$ and for $\var$ small enough
$$
{{\rm d} \over {\rm d}t} \left( e^{\lambda_0 ' t} \, \int_{\T^3
    \times \R^3} \left(\var \, |D_t g_t| + |g_t|\right) \, m \dd x \dd
  v \right) \le 0,
$$
which implies
\begin{equation}\label{eq:DthL1}
\forall \, t \ge 0, \quad \left\| D_t g_t \right\|_{L^1(m)} + \left\| g_t \right\|_{L^1(m)}  
\le \var^{-1} \, e^{-\lambda_0 ' t} \, \left\| h \right\|_{W^{s,1} _x W^{1,1}_v(m)}.
\end{equation}

Then we write 
\begin{eqnarray*}
  t \, \nabla_x T_1(t) (\partial_x ^\alpha h) 
  &=& \int_{\R^3}  k_\delta (v,v_*) \, \left[ (D_t g_t) - \nabla_{v_*} g_t \right] (x,v_*)  \dd v_*
  \\
  &=& \AA_\delta \left( D_t g_t \right) + \AA^2 _\delta g_t,
\end{eqnarray*}
so that thanks to \eqref{eq:DthL1}
\begin{eqnarray*}
t \, \left\|\nabla_x T_1(t) (\partial_x ^\alpha h) \right\|_{L^1(K)} 
&\le& C \, \left[ \left\|D_t g_t \right\|_{L^1(m)} + \left\| g_t \right\|_{L^1(m)} \right]
\\
&\le& C \, \var^{-1} \, e^{-\lambda_0 ' t} \, \left\| h \right\|_{W^{s,1} _x W^{1,1}_v(m)}. 
\end{eqnarray*}
Together with estimate \eqref{TL1xW1v} and Lemma~\ref{lem:regA}, for
$s \ge 0$, we conclude that
\[
  \| T_1(t) (\partial_x ^\alpha h) \|_{W^{1,1}_x W^{s+1,1}_v (K)} \le {C \, e^{-\lambda_0'
      t} \over t}  \, \left\| h \right\|_{W^{s,1} _x W^{1,1}_v(m)}, 
\]
which in turns implies   \eqref{estim:T1gainW}. 

We now interpolate between the last inequality for a given $s \in
[0,1]$, i.e. 
\[
  \| T_1(t) (h) \|_{W^{s+1,1}_x W^{s+1,1}_v (K)} \le {C \, e^{-\lambda_0'
      t} \over t}  \, \left\| h \right\|_{W^{s,1} _x W^{1,1}_v(m)}
\]
and 
\[
\|T_1(t) h \|_{W^{s,1} _x W^{s+1,1}_v (K)} \le C \, e^{-\lambda_0 \, t}
\, \| h \|_{W^{s,1} _x W^{1,1}_v(m)}  
\]
obtained from \eqref{TL1xW1v} written for the same $s$, which gives
\begin{multline}\label{estim:htW1/2x1v} 
  \| T_1(t) h \|_{W^{s+1/2,1}_{x,v}
    (K)} \\ \le C \left( \frac{e^{-\lambda_0'
      t}}{t} \right)^{1/2} \left( e^{-\lambda_0 \, t} \right)^{1/2} \,
\| h \|_{W^{s,1} _x W^{1,1}_v(m)} \le 
{C \, e^{-\lambda_0 ' t} \over \sqrt t} \, \| h \|_{W^{s,1} _x W^{1,1}_v(m)}.  
\end{multline} 
Putting together \eqref{estim:htW1/2x1v} and \eqref{TL1xW1v}, for $s
\in [0,1]$, we get
\begin{eqnarray*} 
  \left\| T_2(t) h \right\|_{W^{s+1/2,1}_{x,v}(K)} &\le& \int_0^t \left\|
    T_1(t-\tau) \, T_1(\tau) h \right\|_{W^{s+1/2,1}_{x,v}(K)} \dd \tau
  \\
  &\le& C \int_0^t \frac{e^{-\lambda_0'(t-\tau)}}{(t-\tau)^{1/2}} \, 
  \| T_1(\tau) h \|_{W^{s,1} _x W^{1,1}_v(m)}  \dd \tau
  \\
  &\le& C \, \left( \int_0^t \frac{e^{-\lambda_0'(t-\tau)}}{(t-\tau)^{1/2}} \, 
    e^{-\lambda_0 \tau} \dd \tau \right) \,  \left\| h \right \|_{W^{s,1}_{x,v}(m)}
  \\
  &\le& C \, e^{-\lambda_0 't} \, \left( \int_0^t
    \frac{e^{-(\lambda_0-\lambda_0')\, \tau}}
    {(t-\tau)^{1/2}} \dd \tau \right) \,  \left\| h \right \|_{W^{s,1}_{x,v}(m)}
  \\
  &\le& C' \, e^{-\lambda_0 't} \, \| h \|_{W^{s,1}_{x,v}(m)}, 
\end{eqnarray*} 
for some other constant $C'>0$, which concludes the proof.
\end{proof}

\begin{rem} The case when the Lebesgue integrability exponent $p \in
  (1,+\infty)$ is different from $p=1$ is less degenerate, and the
  regularization result in finite time can also be obtained thanks to
  classical averaging lemmas \cite{MR923047}. However we both need the
  precise asymptotic estimates and the case $p=1$ in the sequel of
  this paper. 

  Let us explain briefly the alternative argument for the regularity
  in the simplest case, namely when $p=2$ and $s = 0$. The classical
  averaging lemma (see \cite[Lemma 1]{MR1798557} and the proof of
  \cite[Theorem 2.1]{MR1669221}) can be stated as follows in its
  simplest form: any solution $f \in C([0,T]; L^2(\T^3 \times \R^3))$
  to the kinetic equation
$$
\partial_t f_t + v \cdot \nabla_x f_t = g_t, \quad f_{|t=0} = h, 
$$
satisfies for any $\psi \in \DD(\R^3)$ the estimate
$$
\left\|\int_{\R^3} f_t(x,v_*) \, \psi(v_*) \dd v_*
\right\|_{L^2_t\left(H^{1/2}_x\right)} \le C \, \left( \| h \|_{L^2_{x,v}}+
\| g \|_{L^2_{t,x,v}}\right)
$$
where $L^2 _t$ means the $L^2$ norm on the whole real line of times.
Observing that $f_t = S_{\BB_\delta}(t) h$ satisfies the above
kinetic equation with 
\[
g_t:= \BB_\delta f_t = - \nu \, f_t - \BB^2_\delta f_t
\]
and that
\[
\| g_t \| _{L^2(m)} \le C \, \| f_t\| _{L^2(\nu^2 m)} \le C \,
e^{-\lambda_0 t}
\, \| h \|_{L^2(\nu^2 m)},
\]
we deduce that 
$$
\| T_1 (t) h \|_{L^2_t(H^{1/2}_{x,v}(K))} \le C \,  \| h
\|_{L^2(\nu^2 m)}.  
$$
Now, using the Cauchy-Schwarz inequality, we have 
\begin{multline*}
  \| T_2 (t) h \|_{H^{1/2}(K)} \\
  \le \|  h \|_{L^2(\nu^2m)} \, \int_0^t \| T_1(t-s) \|_{L^2(\nu^2 m) \to
    H^{1/2}(K)}  
  \| T_1(s) \|_{L^2(\nu^2m)} \dd s 
  \\
  \le \|  h \|_{L^2(\nu^2m)} \,  \left(  \int_0^t \| T_1(s)
    \|_{L^2(\nu^2 m) \to H^{1/2}(K)}^2 \dd s\right)^{1/2} \times
  \left(  \int_0^t \|T_1(s) \|_{L^2(\nu^2 m)}^2 \dd s\right)^{1/2}  
  \\
  \le C \, \|  h \|_{L^2(\nu^2 m)}
\end{multline*}
which allows to recover pointwise in time estimates. 
%and we conclude similarly as in the proof of Lemma~\ref{lem:regularization-lbe}. 
\end{rem}

\subsection{Proof of the main hypodissipativity result}
% for the  linearized Boltzmann equationDecomposition of the linearized Boltzmann operator and estimates}
\label{sec:decomp-line-boltzm}

We may now conclude the proof of Theorem~\ref{theo:LIBE2}. We consider
$p,q,s,\sigma$ and $m$ that satisfy the assumptions of the theorem. We
set $\EE = W^{\sigma,q} _vW^{s,p}_x(m)$ and $E :=
H^{s'}_{x,v}(\mu^{-1/2})$ with $s' \in \N^*$ large enough.

We apply Theorem~\ref{theo:EnlargingSGdecay}. On the one hand, for
$s'$ large enough, we have $E \subset \EE$. Then we see that {\bf
  (A3)} is fulfilled and {\bf (A1)} is nothing but
\cite[Theorem~3.1]{MNeu}.  On the other hand, assumption {\bf (A2)} is
a direct consequence of Lemma~\ref{lem:regA},
Lemma~\ref{lem:decomp-new} and Lemma~\ref{lem:regularization-lbe},
together with Lemma~\ref{lem:Tn}. Indeed, from
Lemma~\ref{lem:regularization-lbe} and Lemma~\ref{lem:Tn} we have for
instance
$$
 \left\| T_n(t) h \right\|_{H^{s'} _{x,v}(\mu^{-1/2})} \le C \,
 e^{-\lambda_0 ' t}
 \| h \|_{L^1_{x,v}(\langle v \rangle^3)},
$$
so that 
$$
\left\| T_{n+1}(t) h \right\|_{E}  \le C \, e^{-\lambda_0 ' t} \|  h \|_{\EE}.
$$

This proves the exponential decay on the semigroup in $\EE$. Then one
obtains a rate of decay in $\EE$ equal to the one in $E$ as soon as
$\lambda_0$ (provided by Lemma~\ref{lem:decomp-new}) is strictly
greater than the spectral gap $\lambda \in (0,\nu_0)$ in $E$ (which
required the condition $k$ is large enough on the exponent of the
weight in case of a polynomial weight), which also then allows to take
$\lambda_0'$ strictly greater than the sepctral gap in $E$ in
Lemma~\ref{lem:regularization-lbe} and Lemma~\ref{lem:Tn}. This proves
the last claim in the statement of Theorem~\ref{theo:LIBE2}.

\subsection{Structure of singularities for the linearized flow}
\label{sec:study-sing-struct}

From the previous study of the decay rate of the linearized flow, we
have obviously the following decomposition of the solution $h_t :=
S_\LL(t) \init{h}$:
\[
h_t = \Pi \init{h} + \left( h_t - \Pi \init{h}\right).
\]
In this decomposition the first part is infinitely regular, say in
$H^\infty(\mu^{-1/2})$, and the second part decays like $O(e^{-\lambda
  \, t})$, where $\lambda>0$ denotes the optimal spectral gap (for
polynomial moments this requires the condition $k > k^* _q$). We shall
now make more precise the singularity structure of the second part,
showing on the one hand that its dominant part in this asymptotic
behavior is as regular as wanted, and on the other hand that its worst
singularities are supported by the free motion characteristics. One
way to understand these statements is through a spectral decomposition
of the semigroup, and the method we expose here can be considered as a
quantitative spectral decomposition in this context.

\subsubsection{Asymptotic amplitude of the singularities}
\label{sec:asympt-ampl-sing}

Let us consider for instance the space $L^1_{x,v}(m)$ where the weight
$m$ satisfies the assumptions of Theorem~\ref{theo:LIBE2}.  Other
spaces can be considered, provided that they fall within the scope of
Theorem~\ref{theo:LIBE2}. We % assume without loss of generality that
% $\Pi \init{h} =0$ (which implies $\Pi h_t =0$ for any later time), and
% we
start from the following decomposition formula of the semigroup 
 \begin{multline*}
   S_\LL(t) = \Pi_{\LL,0} + \sum_{\ell=0}
   ^{n-1} (-1)^\ell \, (\mbox{Id} - \Pi_{\LL,0}) \, S_\BB \ast \left( \AA S_\BB \right)^{\ast \ell}(t) \\
   + (-1)^n \left[ (\mbox{Id} - \Pi_{L,0}) S_L\right] \ast
   \left( \AA S_\BB \right)^{\ast n}(t)
 \end{multline*}
 that has been proved. We then use on the one hand that, given any $s
 \in \N$ and $\var >0$, there is $n$ large enough so that
% \begin{equation*}
%   \supp \, \left( \AA S_\BB \right)^{\ast n}(t) \, h \subset \T^3 \times
%   B(0,R)
% \end{equation*}
% and 
\begin{equation*}
\n{\left( \AA S_\BB \right)^{\ast n}(t)
    h}_{H^s_{x,v}(\mu^{-1/2})} \le C
  \, e^{-(\nu_0-  \var) \, t} \, \n{h}_{L^1_{x,v}(m)}
\end{equation*}
thanks to the previous study, and 
\begin{equation*}
  \n{\left[ (\mbox{Id} - \Pi_{L,0}) S_L\right] h}_{H^s_{x,v}(\mu^{-1/2})}
  \le C \, e^{-\lambda \, t} \, \n{h}_{H^s_{x,v}(\mu^{-1/2})}
\end{equation*}
with the optimal rate $\lambda$. Since $\nu_0 > \lambda$, by choosing
$\var>0$ small enough we deduce that 
\begin{equation*}
  \n{\left[ (\mbox{Id} - \Pi_{L,0}) S_L\right] \ast
   \left( \AA S_\BB \right)^{\ast n}(t) h}_{H^s_{x,v}(\mu^{-1/2})}
  \le C \, e^{-\lambda \, t} \, \n{h}_{L^1_{x,v}(m)}
\end{equation*}
with the optimal rate $\lambda$. On the other hand, for all the other
terms in the decomposition we use the decay of $S_\BB(t)$ with
exponential rate as close as wanted to $-\nu_0$ to deduce that, for
any $\var>0$ 
\begin{equation*}
  \n{\sum_{\ell=0}
   ^{n-1} (-1)^\ell \, (\mbox{Id} - \Pi_{\LL,0}) \, S_\BB \ast \left(
     \AA S_\BB \right)^{\ast \ell}(t) h}_{L^1_{x,v}(m)} \le C \,
 e^{-(\nu_0-\var) \, t} \, \n{h}_{L^1_{x,v}(m)}. 
\end{equation*}
This thus shows that for any $s \in \N$ and $\var >0$ there is a
decomposition of  the linearized flow as
\begin{equation*}
  S_\LL(t) = \Pi_{\LL,0} + S^s_\LL(t) + S^r_\LL(t)
\end{equation*}
where $S^s_\LL(t)$ satisfies 
\[
%\supp \, S^s_\LL(t) \, h \subset B(0,R), \quad 
\n{ S^s_\LL(t) \, h}_{H^s_{x,v}(\mu^{-1/2})} \le C \, \n{h}_{L^1_{x,v}(m)} \,
e^{-\lambda \, t}
\]
with the sharp rate $\lambda>0$ and where $S^r_\LL(t)$ satisfies
\[
\n{ S^r_\LL(t) \, h}_{L^1_{x,v}(m)} \le C \, \n{h}_{L^1_{x,v}(m)} \, e^{-(\nu_0-\var) \, t}. 
\]
In words, the part $S^s$ is as smooth as wanted, with Gaussian
localization as in the small linearization space, and decays in time
with the sharp rate $\lambda$, and the part $S^r$ decays in time
exponentially fast in the original space $L^1_{x,v}(m)$ with a rate as
close as wanted to $\nu_0$, which corresponds to the onset of the
continuous spectrum. The latter part $S^r$ carries all the
singularities of the flow. 

\subsubsection{Localization of the $L^2$ singularities}
\label{sec:local-l2-sing}

We consider now the space $L^2_{x,v}(m)$ with a weight $m$ so that the
assumptions of Theorem~\ref{theo:LIBE2} are satisfied. (Again other
spaces could be considered). We know that the solution $h_t$ to the
linearized problem remains uniformly bounded in this space along
time. We now consider the decomposition 
\begin{equation*}
  \LL = \KK - v \cdot \nabla_x - \nu := \KK + \BB_0
\end{equation*}
and apply our decomposition at order one: 
\begin{equation*}
   S_\LL(t) = \Pi_{\LL,0} + (\mbox{Id} - \Pi_{\LL,0}) \, S_{\BB_0} (t) 
   - \left[ (\mbox{Id} - \Pi_{L,0}) S_L\right] \ast
   \left( \KK S_{\BB_0} \right)(t). 
 \end{equation*}
Then one checks with the help of the explicit formula
\[
S_{\BB_0}(t) h(x,v) = e^{-\nu(v) \, t} \, h(x-vt,v)
\]
that the second term in the right hand side propagates the singularity along the
characteristic lines of the transport flow while damping their
amplitude like $e^{-\nu(v) \, t}$. Finally for the third term we use
that by interpolation and averaging lemma (as in \cite{MR2081030} and
\cite{MR1798557})
\[
\n{\left( \KK S_{\BB_0} \right)(t) h}_{H^\alpha_{x,v,loc}} \le \frac{C}{\min\{t^\theta ; \, 1 \}} \, 
\n{h}_{L^2_{x,v}(m)}
\]
for some small but non-zero $\alpha>0$ and some $\theta>0$. This
proves the decomposition
\[
S_\LL(t) h \in \left[ \Pi_{\LL,0} + (\mbox{Id} - \Pi_{\LL,0}) \, \left(
  e^{-\nu(v) \, t} \, h(x-vt,v) \right) \right] + O(t^{-\theta}) \, H^\alpha_{x,v,loc}
\]
where $H^\alpha_{x,v,loc}$ denotes some function which belongs to the
fractional Sobolev space $H^\alpha_{x,v}$ when restricted to any
compact set. This captures the localization of $L^2$ singularities. 

%%%%%%%%%%%%%%%%%%%% Boltzmann %%%%%%%%%%%%%%%%%%%%%%%%%%%%%%
\section{The nonlinear Boltzmann equation}
\label{sec:boltzmann}
\setcounter{equation}{0}
\setcounter{theo}{0}

In this section, we are concerned with the proof of the main outcome
of our theory: two new Cauchy results for the nonlinear Boltzmann equation
with optimal decay rates, and the proof of the exponential $H$-theorem
under a priori assumptions.

% \Red SM: I am agree with the statements and proofs below when $q =
% 1$, $s \in \N$.  I believe that the results can also be easily
% proved in the case $q=\infty$ and $m = \exp (\alpha \, |v|^\beta)$,
% $\beta \in (1,2)$. As I have pointed out to Cl\'ement, there is
% still some gaps in the proofs in order to obtain similar results
% when $q \in (1,\infty)$ and when $q=\infty$, $m$ grows moderately.
% Do we try to fill the gap?

%\Black

%----------------------------------------------------------------------------------------------------
\subsection{The main results} 

We consider the fully non-linear problem \eqref{eq:NLBE}, first in the
close-to-equilibrium regime, then in the weakly inhomogeneous regime,
and finally the far-from-equilibrium regime with a priori bounds. Here
and below we call \emph{normalized distribution} a distribution with
zero momentum, and mass and temperature normalized to one (remember
that the volume of the torus is normalized to one, and therefore this
definition is unchanged for spatially homogeneous distributions). This
normalization induces no loss of generality thanks to the conservation
laws of the nonlinear flow. Let us first define the notion of
solutions we shall use
\begin{defin}[Conservative solution]\label{def:sol}
  For some non-negative inital data $\init{f} \in L^1_v 
  L^\infty_x(1+|v|^2)$, we say that for $T \in (0,+\infty]$, 
  \[
  0 \le f \in L^1_{t,loc} \left([0,T),L^1_v L^\infty_x(1+|v|^2)\right) \cap
  C^0 _t\left([0,T),L^1_v L^\infty_x(1+|v|^1)\right)
  \]
  is a \emph{conservative (distributional) solution} on $[0,T)$ if it
  satisfies
  \begin{equation*}
    \left\{ 
      \begin{array}{l} \ds 
        \partial_t f + v \cdot \nabla_x f = Q(f,f) \quad \mbox{in the
          sense of distributions,} \vs \\ \ds 
        f_{|t=0} = \init{f} \quad \mbox{almost everywhere,}
      \end{array}
\right.    
\end{equation*}
and satisfies the conservation law
  \begin{equation*}
    \fa t \ge 0, \quad \int_{\T^3 \times \R^3} f_t(x,v) (1+|v|^2) \dd
    x \dd v = \int_{\T^3 \times \R^3} \init{f}(x,v) (1+|v|^2) \dd
    x \dd v.
  \end{equation*}
\end{defin}

\begin{rem}
  The solutions can also understood in the renormalized sense and in
  the mild sense, that is in the sense of the almost everywhere
  equality
$$
f_t(x,v) = \init{f}(x-vt,v) + \int_0^t Q(f_\tau,f_\tau)(x-v(t-\tau),v)
\dd \tau.
$$
Observe that thanks to the bilinear estimates available on $Q$, for
solutions in $L^1_{t,loc}([0,T),L^1_vL^\infty_x(1+|v|^2))$, the last term of the
right hand side is always well-defined as a measurable function.
\end{rem}

\begin{theo}[Nonlinear stability]
  \label{theo:NLBE} 
We divide our main result into:

\mk

\noindent {\bf (I) A priori properties of conservative solutions.}
Consider a conservative solution as defined above on $[0,T)$, $T \in
(0,+\infty]$, with a uniform bound from below on the initial
distribution
\begin{equation}\label{eq:hyp-lower}
  \fa x \in \T^3, \ v \in \R^3, \quad \initem{f}(x,v) \ge \varphi(v)
  \ge 0,
  \quad \int_{\R^3} \varphi(v) \dd v \in (0,+\infty). 
\end{equation}
Then this solution satisfies for any positive time $t>0$: 
\[
\left\{ 
  \begin{array}{l}\ds
    \fa k >0, \quad \n{f_t}_{L^1_{x,v}(1+|v|^k)} <+\infty \vs \\ \ds
    \fa x \in \T^3, \ v \in \R^3, \ f_t(x,v) \ge K_1 \, e^{-K_2 |v|^2} 
  \end{array}
\right.
\]
for some $K_1, K_2>0$. In the case of a global solution
($T=+\infty$), these estimates are uniform as time goes to infinity.

Moreover when the initial data belongs to $L^1_v W^{3,1}_x(1+|v|^2)$
the moment estimate can be (strongly) improved into $\n{f_t}_{L^1_v
  W^{3,1}_x (e^{\kappa |v|})} <+\infty$ for some $\kappa>0$. However
for higher-order exponential moments $L^1_v W^{3,1}_x (e^{\kappa
  |v|^\beta})$, $\beta \in (1,2]$, $\kappa >0$, if they are not finite
initially they remain infinite for all times.

Finally these conservative solutions are a priori unique (without
perturbative assumptions) at least when restricted to $L^1_{t,loc}
L^1_v L^\infty_x(1+|v|^k) \cap C^0_t L^1_v L^\infty_x (1+|v|^{k-1})$,
$k>2$, or, in the critical case $k=2$, when restricted to $L^1_{t,loc} 
L^1_v W^{3,1}_x(1+|v|^2) \cap C^0_t L^1_v W^{3,1}_x (1+|v|)$.
\mk

\noindent {\bf (II) Nonlinear stability.} For any $k>2$, there is some
constructive constant $\epsilon = \epsilon(k) >0$ such that for any
normalized non-negative initial data % $\initem{f} \in
% L^1_vL^\infty_x(1+|v|^2)$
satisfying
\[
\left\| \initem{f} - \mu \right\|_{L^1_vL^\infty_x(1+|v|^k)} \le \epsilon(k), 
\]
where $\mu$ is the Maxwellian equilibrium defined in
\eqref{def:secBL:mu}, there exists a unique global conservative
solution in $L^\infty_t L^1_v L^\infty_x (1+|v|^k) \cap C_t^0 L^1_v
L^\infty_x$ to \eqref{eq:NLBE} with initial $\initem{f}$, which
satisfies
\begin{equation}\label{eq:secBNL:th1}
  \forall \, t \ge 0, \quad 
  \left\| f_t -  \mu \right\|_{L^1_vL^\infty_x(1+|v|^k)} \le C_1 \, 
  e^{ -\lambda \, t} \, \left\|
    \initem{f} -\mu \right\|_{L^1_vL^\infty_x(1+|v|^k)}
\end{equation}
where $\lambda$ is the optimal linearized rate in
Theorem~\ref{theo:LIBE2} and for some explicit constant $C_1\ge
1$.

\mk

\noindent {\bf (III) Stability in stronger norms.} Consider for $p,q
\in [1,+\infty)$ any functional space
\[
\EE = W^{\sigma,1}_v W^{s,p}_x (m) \cap W^{\sigma,q}_v W^{s,p}_x (m)
\subset L^1_vL^\infty_x(1+|v|^2)
\]
with $s,\sigma \in \N$, $\sigma \le s$, $s > 6/p$ and $m$ satisfying
one of the assumptions {\bf (W1)}, {\bf (W2)}, {\bf (W3)} in
Theorem~\ref{theo:LIBE2}. In the case $p=+\infty$ one can consider the
same spaces but including additionally the case $s \ge 0$. Finally in
the case $q=+\infty$ of {\bf (W2)} or {\bf (W3)}) then consider the
simpler functional spaces
\[
\EE = W^{\sigma,\infty}_v W^{s,p}_x (m) \subset
L^1_vL^\infty_x(1+|v|^2). 
\]

Then there is some constructive constant $\epsilon = \epsilon(\EE) >0$
such that if the previous initial data satisfies furthermore $\|
\initem{f} -\mu \|_\EE \le \epsilon(\EE)$, we have the estimate
\begin{equation}\label{eq:secBNL:th1}
  \forall \, t \ge 0, \quad 
  \left\| f_t -  \mu \right\|_{\EE} \le 
  C_2 \, e^{ -\lambda \, t} \, \left\| \initem{f} -\mu \right\|_{\EE}.
\end{equation}
with the optimal rate $\lambda$ and for some constructive constant
$C_2 \ge 1$. 
\end{theo}

\begin{rems}
  \begin{enumerate}
  \item The rate $\lambda$ and constants in Theorem~\ref{theo:NLBE} on
    the nonlinear flow are obtained in a constructive way and the rate
    is the same as for the linearized flow. In turn we have given
    sufficient conditions in Theorem~\ref{theo:LIBE2} for this rate to
    be the same as the sharp rate in the space
    $L^2(\mu^{-1/2})$. Finally in the latter space, the decay rate and
    constants were proved in~\cite{MR0445138} by non-constructive
    argument based on Weyl's theorem, and then the series of papers
    \cite{Baranger-Mouhot,MR2254617,MR2301289,MNeu} provided
    constructive proof with explicit constants and estimates on the
    rate $\lambda$. 
  \item Some refinements of these theorems could be considered: (1)
    extend these results to variable hard potentials ($\gamma \in
    (0,1]$); (2) extend these results to solutions $M^1_{v}
    W^{s,p}_x(m)$ that are merely measures in the velocity variable,
    by using the recent works \cite{lm1,lm2} at the spatially
    homogeneous level\footnote{Note that in this case the lower bound
      assumption~\eqref{eq:hyp-lower} should be changed into:
      $\varphi$ non-negative measure with positive mass and different
      from a single Dirac mass.}. We did not include these natural
    extensions in the statement as it is already long enough.
\item It seems also that in the spatially homogeneous setting the
  optimal rate in $(W_v ^{\sigma,1} \cap W_v^{\sigma,q})(m)$, $\sigma
  \ge 0$, $q \in [1,+\infty]$, with $m$ satisfying {\bf (W3)},
  provided by Theorem~\ref{theo:NLBE+} is new (whereas it was proved
  in the case {\bf (W2)} in \cite{Mcmp}).
\item The fact that Gaussian moments do not appear in part {\bf (I)}
  justifies the need for enlarging the functional space of the decay
  estimates on the linearized flow. An interesting open question is to
  clarify whether the nonlinear Boltzmann equation (starting with the
  spatially homogeneous case) is indeed \emph{ill-posed} in
  $L^2(\mu^{-1/2})$ in the non-perturbative regime.
  \end{enumerate}
\end{rems}

\begin{theo}[Weakly inhomogeneous solutions]
 \label{theo:NLBE+}
 Consider a normalized non-negative spatially homogeneous distribution
 $\initem{g} = \initem{g}(v) \in L^1_v (1+|v|^k)$, $k >2$. Then there
 is some constructive constant $\epsilon >0$ depending on the mass,
 energy and $k$-moment of $\initem{g}$, such that for any normalized
 non-negative initial data $\initem{f} \in L^1_v L^\infty_x(1+|v|^k)$
 satisfying
\[
\left\| \initem{f} - \initem{g} \right\|_{L^1_v L^\infty _x(1+|v|^k)}
\le \epsilon,
\]
there exists a unique global conservative solution in $L^\infty_t
L^1_v L^\infty_x (1+|v|^2) \cap C_t^0 L^1_v L^\infty_x(1+|v|)$ to
\eqref{eq:NLBE} with initial data $\initem{f}$, which satisfies
\begin{equation}\label{eq:secBNL:th2}
\forall \, t \ge 0, \quad \left\| f_t - g_t \right\|_{L^1_v L^\infty_x(1+|v|^2)} \le
C \, \epsilon, 
\end{equation}
where $g_t$ is the unique conservative solution to the spatially
homogeneous Boltzmann equation starting from $\initem{g}$, as well as
the properties {\bf (I)} above and 
\begin{equation*}
  \forall \, t \ge 0, \quad 
  \left\| f_t -  \mu \right\|_{L^1_v L^\infty_x(1+|v|^2)} \le C \, 
  e^{ -\lambda \, t}
\end{equation*}
where $\lambda>0$ is the optimal linearized rate in
Theorem~\ref{theo:LIBE2} and for some constant $C>0$.
\end{theo}

\begin{rems}
\begin{enumerate}
% \item The smallness parameter $\epsilon$ depends more precisely on a
%   positive lower bound on $\int_{\R^3} \init{g}(v) \, |v| \dd v$ and an upper
%   bound on $\int_{\R^3} \init{g}(v) \, |v|^2 \dd v$. In fact any
%   probability measure $\init{g}$ with finite second moment and different from
%   $\delta_0$ would satisfy such bounds. 

\item It is possible to prove a posteriori estimates on $f_t$ in
  spaces of the form
\[
W^{\sigma,1}_v W^{s,p}_x (1+|v|^k) \cap W^{\sigma,q}_v W^{s,p}_x
  (1+|v|^k) \subset L^1_v L^\infty_x(1+|v|^k)
\] 
(with the conditions {\bf (W3)} on $s,\sigma,p,q$ and $k$), by using
some refined technical convolution inequalities on the collision
operator from \cite{MR2081030}. We leave this question, as well as
that of a general a posteriori regularity theory, to further studies.

\item Theorems~\ref{theo:NLBE} and \ref{theo:NLBE+} provide the
  largest class of unique solutions constructed so far to our
  knowledge (in $L^1_v L^\infty_x (1+ |v|^{2+0})$ close to equilibrium
  or close to spatially inhomogeneous solutions). It is an interesting
  open question whether existence and uniqueness can be obtained in
  the space $L^1_v L^\infty_x (1+ |v|^2)$ (or $L^1_v
  W^{3,1}_x(1+|v|^2)$ where we have proved above that a priori
  uniqueness holds for conservative solutions) with a perturbation
  condition.
\end{enumerate}
\end{rems}

\begin{theo}[Exponential $H$-theorem with a priori bounds]\label{theo:expEB}
  Let $(f_t)_{t \ge 0}$ be a normalized non-negative smooth solution
  of \eqref{eq:NLBE} such that for $k,s$ large enough
\[
\sup_{t \ge 0} \left( \| f_t \|_{H^s(\T^d \times \R^3)} + \| f_t
  \|_{L^1(1+|v|^k)} \right) < +\infty,
\]
and such that its spatial density 
\[
\fa x \in \T^3, \quad \initem{\rho}(x) = \int_{\R^d} \initem{f}(x,v)
\dd v \ge \alpha >0
\]
is uniformly positive on the torus.
\smallskip

Then this solution satisfies 
\begin{equation*}
  \forall \, t \ge 0, \quad \n{f_t}_{L^1_v L^\infty_x (1+|v|^2)} \le C \,
  e^{-\lambda \, t}
\end{equation*}
and 
\begin{equation*}
\forall \, t \ge 0, \quad \int_{\T^d \times \R^3} f_t \, \log
\frac{f_t}{\mu} \dd x \dd v \le C \,
e^{-\lambda \, t}
\end{equation*}
for some constructive constant $C >0$, and where $\lambda >0$ is the
optimal linearized rate in Theorem~\ref{theo:LIBE2}. 
\end{theo}

\begin{rem}
  Our relaxation rate in $L^1_v L^\infty_x (1+|v|^2)$ norm is
  optimal. However the linearization of the relative entropy would
  suggest the relaxation rate $O(e^{-2 \, \lambda \, t})$ for the
  relative entropy since
  \begin{equation*}
    \int_{\T^d \times \R^3} f_t \, \log
\frac{f_t}{\mu} \dd x \dd v = \int_{\T^d \times \R^3} \left( \frac{f_t}{\mu} \, \log
\frac{f_t}{\mu} - \frac{f_t}{\mu} +1 \right) \dd x \dd v
  \end{equation*}
  and $z \log z - z +1 \sim z^2/2$ at $z =1$. This statement needs
  however proper justification; first of all in order to be true it
  would require for the solution $f_r$ to have tails decaying as
  $\mu$, which is expected to be wrong outside specific perturbative
  regimes. Therefore it is an interesting open question to know
  whether the relaxation rate of the relative entropy for perturbative
  solutions with {\em{polynomial}} tail lies between $e^{-\lambda \,
    t}$ and $e^{- 2 \, \lambda \, t}$.  The importance of tail's decay
  was already outlined by Cercignani in his conjecture~\cite{Ce82}.
\end{rem}

\subsection{Strategy of the Proof of Theorem~\ref{theo:NLBE}}
\label{sec:strategy}
\mbox{ } 
\sk

\noindent
Part {\bf (I)}: The moment bounds are inspired by the arguments in the
spatially homogeneous case \cite{MR2264623,Mcmp,lm1,AlCaGaMo} and more
precisely by the techniques developed in \cite{AlCaGaMo}. The lower
bounds is obtained from the results in
\cite{MR1461954,Mlower,MR1697562,MR1697495,briant}. The a priori
uniqueness is inspired by the proof of uniqueness in the spatially
homogeneous case \cite{MR1697562,MR1716814,lm1}: more precisely it
extends to the spatially inhomogenenous case the method presented by
Lu in \cite{MR1716814} (see also \cite{lm1}). \mk

\noindent
Part {\bf (II)} and {\bf (III)}: The study of the nonlinear stability
is based on \emph{energy methods}. Such methods are often used in
nonlinear PDE's, and use the coercivity properties of the linearized
operator. However in the present situation the time decay estimates
obtained on the linearized semigroup do \emph{not} imply coercivity
inequalities on some Dirichlet form due to the absence of symmetry
structure. To resolve this issue we introduce a new
\emph{non-symmetric energy method}. We introduce in the next
subsection a \emph{dissipative Banach norm}, for which some suitable
coercivity is recovered. This norm involves the linearized evolution
flow for all times. More precisely we prove:
\begin{itemize}
\item[(1)] Bilinear estimates to control the nonlinear remainder in
  the equation for any given initial datum $\init{g}\in
  W_v^{\sigma,q}(m)$.
\item[(2)] The key a priori estimate for $k>2$ moments which provides
  the ``linearisation trap''.
% \item[(3)] the key a priori estimate for the critical case $k=2$ used
%   only for initial times,
\item[(3)] A local-in-time existence result. 
%\item[(4)] some moments bounds in a spatially inhomogeneous setting.
\end{itemize}
We then conclude the proof by standard continuation method.

The proof of Theorem~\ref{theo:NLBE+} is based on the previous
linearized stability estimates in functional spaces large enough to be
compatible with the Cauchy theory of the spatially homogeneous
equation in the large, and a classical argument on the dynamics,
inspired from \cite{MR900501}. It is sketched in Figure~1: the
spatially homogeneous solutions are represented as a subset a general
solutions. By proving local-in-time stability in $L^1_v L^\infty_x$
spaces, we can capture a general solution around this subset. If this
time is large enough, which is granted if the perturbation between
$\init{f}$ and $\init{g}$ is small enough, then $f_t$ is driven
towards equilibrium thanks to the relaxation estimates known for
$g_t$. Finally we use the linearized stability estimates once the
stability neighborhood is entered by $f_t$.
\begin{figure}[h]\label{fig:weakly-inhom}
 \hspace{-1.5cm} \includegraphics[width=14cm]{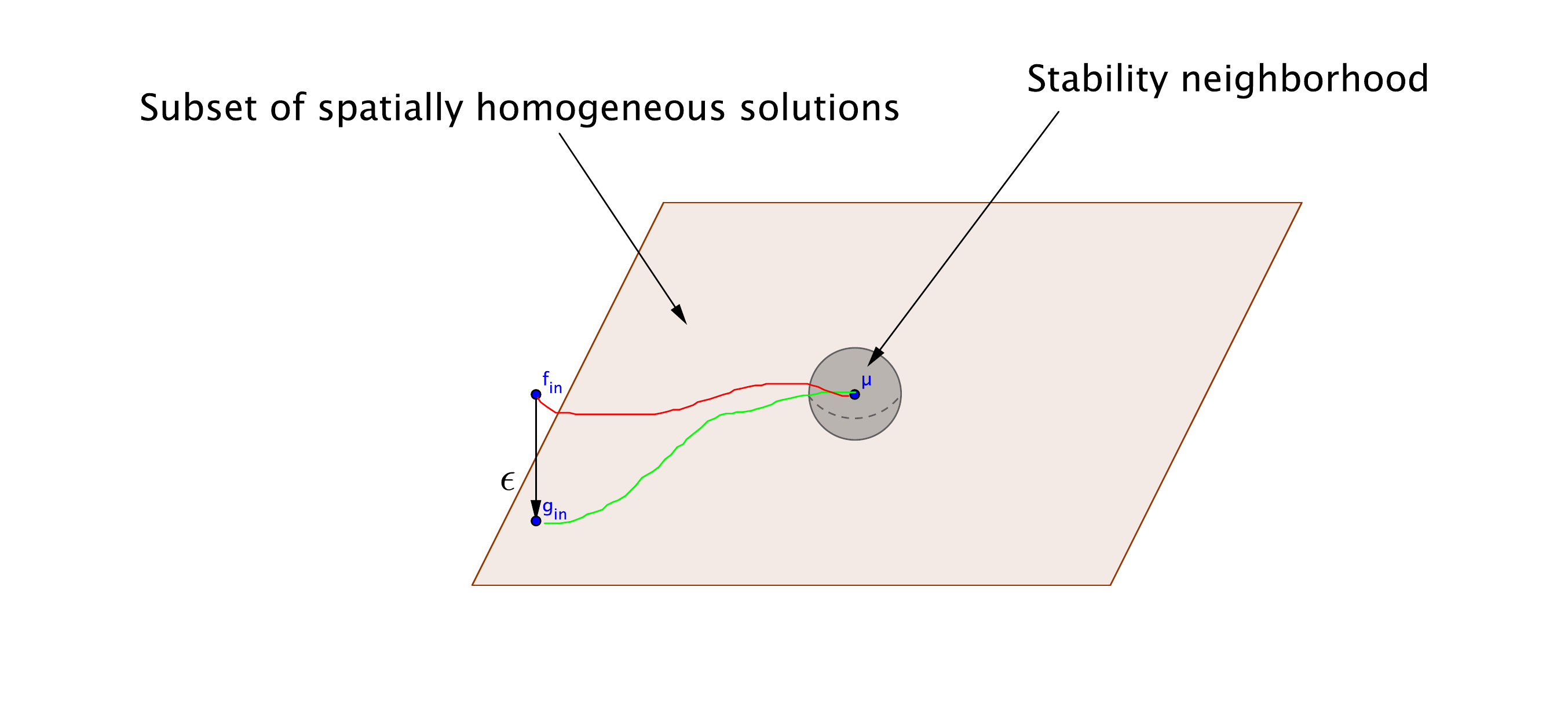}
  \caption{Sketch of the construction of weakly inhomogeneous
    solutions.}
\end{figure}
%----------------------------------------------------------------------------------------------------

\subsection{Proof of Theorem~\ref{theo:NLBE}, part {\bf (I)}}
\label{sec:priori-moment-bounds}

\subsubsection{A priori moment bounds}
\label{sec:moment-bounds}

Polynomial moments estimates are now a classical tool in the theory of
the spatially homogeneous Boltzmann equation. Exponential moments
estimates for the spatially homogeneous Boltzmann equation are more
recent, see \cite{MR1478067,MR2096050,MR2533928} and the references
therein. In the latter references exponential moments (in integral or
pointwise forms) are shown to be {\em propagated}. In the papers
\cite{MMR,Mcmp,lm1,AlCaGaMo} a theory of {\em appearance} of
exponential moments was developed, still in the spatially homogeneous
case. We shall extend this theory to the inhomogeneous framework,
taking advantage of the a priori bounds on the solutions.

 \begin{lem}\label{lem:mts}   
   Consider for $T \in (0,+\infty]$ a conservative solution 
   \[
   0 \le f \in L^1 _{t,loc} \left([0,T),L^1_v L^\infty_x(1+|v|^2)\right) \cap
   C^0 _t \left([0,T),L^1_v L^\infty_x(1+|v|)\right)
  \]
  with initial datum bounded uniformly from below as
  in~\eqref{eq:hyp-lower}.

  Then the solution $f$ has the following properties: for any $k>2$
  and $T' \in (0,T)$ there is an explicit constant $C(k,T')>0$
  depending on $k>2$, on the $L^\infty _t ([0,T'],L^1_v
  L^\infty_x(1+|v|))$ norm of the solution, on the lower bound
  \eqref{eq:hyp-lower} on the initial datum, and on $T'$, so that
\begin{equation}
  \label{eq:poly-moments}
  \forall \, t \in (0,T'], \quad \int_{\T^3 \times \R^3} f_t(x,v) \,
  |v|^k \dd x \dd v \le C(k,T') \, \max \left\{
    \frac{1}{t^{k-2}}, 1 \right\}. 
\end{equation}
\end{lem}
\begin{rem}
  Observe that our moment estimate is not uniform in time. This is due
  to the lack of known uniform-in-time estimates from below on
  solutions to the nonlinear Boltzmann equation with such a low
  regularity. This will however not cause any problem for our
  uniform-in-time stability results since the ``trapping mechanism''
  around the linearized regime takes over in finite time for the
  solutions we considered.
\end{rem}

\begin{proof}[Proof of Lemma~\ref{lem:mts}]
  Using the Duhamel formulation and the above bounds on the solution
  we have for $T' \in (0,T)$: 
  \begin{multline*}
    \fa t \ge \in [0,T'], \ x \in \T^3, \ v \in \R^3, \\ f_t(x,v) =
    e^{-\int_0 ^t Q^-(f_\tau,f_\tau)(x-v(t-\tau),v) \dd \tau} \,
    \init{f}(x-vt,v) \\ + \int_0 ^t e^{-\int_0 ^\tau
      Q^-(f_{\tau'},f_{\tau'})(x-v(\tau-\tau'),v) \dd
      \tau'} \, Q^+(f_\tau,f_\tau)(x-v(t-\tau),v) \dd \tau \\
    \ge e^{-c(T') t (1+|v|)} \, \varphi(v)
  \end{multline*}
for some constant $c(T')>0$ depending on $T'$ (through the
$L^\infty_t([0,T'],L^1_v L^\infty_x(1+|v|^2))$ norm of the solution). 

  We deduce that there is a constant $K(T') >0$ so that
\begin{equation*}
  \fa t \in [0,T'], \ x \in \T^3, \ v \in \R^3, \quad
          \int_{\R^3} f_t(x,v_*) \, |v-v_*| \dd v_* \ge K(T') \,
          (1+|v|). 
  \end{equation*}

Consider now the moments of the solutions 
\begin{equation*}
  M_k[f_t] := \int_{\T^3 \times \R^3} f_t(x,v) \, (1+|v|^k) \, \dd x \dd v,
  \quad k \ge 0,
\end{equation*}
and apply the Lemma~\ref{lem:Povznerbis} to get for $k>2$ the
following inequality in the sense of distribution
\begin{multline*}
  \dt M_k[f_t] = \int_{\T^3 \times \R^3 \times \R^3 \times \mathbb
    S^2} f_t(x,v) \, f_t(x,v_*) \\ \times 
  \left[ |v'|^k + |v'_*|^k - |v|^k - |v_*|^k \right] \, |v-v_*| \dd v \dd v_* \dd
  \sigma \dd x \\
\le C_k \, \int_{\T^3 \times \R^3 \times \R^3} f_t(x,v) \, f_t(x,v_*) \,
(1+|v|^k) \, (1+|v_*|^2) \dd x \dd v \\ - 2 \, \int_{\T^3 \times \R^3 \times \R^3}
f_t(x,v) \, f_t(x,v_*) \, |v|^k \, |v-v_*| \dd x \dd v \dd v_*
\\ \le C'_k \, M_k[f_t] - K_k \, M_{k+1}[f_t] 
 \le C' _k \, M_k[f_t] - K_k ' \, M_{k} ^{\frac{k-1}{k-2}} [f_t]
\end{multline*}
for some constants $C_k, C_k', K_k, K_k' >0$ depending on the
$L^1_vL^\infty_x(1+|v|^2)$ upper bound on the solution and the
previous lower bound. By standard interpolation and Gronwall
inequality argument this leads to the bound
\begin{equation*}
  \fa t \in (0,T'], \quad M_k[f_t] \le \frac{C(k,T')}{t^{k-2}}
\end{equation*}
for some constant $C(k,T')>0$ which depends on $k>2$, $T'>0$ and on
the bounds on the solution.
\end{proof}

\begin{lem}\label{lem:exp-moments}
  Consider for $T \in (0,+\infty]$ a conservative solution 
   \[
   0 \le f \in L^1_{t,loc} \left([0,T),L^1_v W^{3,1}_x(1+|v|^2)\right) \cap
   C^0_t\left([0,T),L^1_v W^{3,1}_x(1+|v|)\right)
  \]
  with initial datum bounded uniformly from below as in
  ~\eqref{eq:hyp-lower}. 

  Then for any $T' \in (0,T)$, there exist explicit constants $\kappa,
  C >0$ (depending on the bounds assumed on the solution, on the lower
  bound \eqref{eq:hyp-lower} on the initial datum, and on $T'>0$) such
  that
     \begin{equation}\label{eq:mts}
       \forall \, t \in [0,T'], \quad \n{f_t}_{L^1_v W^{3,1}_x(e^{\kappa 
           \min\{t,1\} |v|})} \le C.
     \end{equation}
\end{lem}

\begin{proof}[Proof of Lemma~\ref{lem:exp-moments}]
  As a first step let us extend the polynomial moment bounds to the
  derivatives of the solution. Let us define 
  \begin{equation*}
    \tilde M_k(t) := \sum_{|\alpha|\le 3} c_\alpha \, M_k[\partial_x ^\alpha
    f_t]
  \end{equation*}
  for some constants $c_\alpha >0$ to be fixed later. Arguing as in
  the previous lemma and using the Sobolev embedding $W^{3,1}_x
  \hookrightarrow L^\infty_x$, we get
\begin{equation*}
  \dt M_k[f_t] \le C' _k \, M_k[f_t] - K_k ' \, M_{k}
  ^{\frac{k-1}{k-2}} [f_t]
\end{equation*}
for some constants depending on time. For the first derivatives we
write (with the notation $s = \mbox{sign}(\partial_x f)$)
\begin{multline*}
  \dt M_k[\partial_x f_t] = \int_{\T^3 \times \R^3 \times \R^3 \times \mathbb
    S^2} \partial_x f_t(x,v) \, f_t(x,v_*) \\ \times 
  \left[ |v'|^k s' + |v'_*|^k s'_* - |v|^k s - |v_*|^k s_* \right] \, |v-v_*| \dd v \dd v_* \dd
  \sigma \dd x \\
\le C_k \, \int_{\T^3 \times \R^3 \times \R^3} |\partial_x f_t(x,v)| \, f_t(x,v_*) \,
(1+|v|^k) \, (1+|v_*|^2) \dd x \dd v \\ - 2 \, \int_{\T^3 \times \R^3 \times \R^3}
|\partial_x f_t(x,v)| \, f_t(x,v_*) \, |v|^k \, |v-v_*| \dd x \dd v
\dd v_* \\ 
+ 2 \, \int_{\T^3 \times \R^3 \times \R^3}
|\partial_x f_t(x,v)| \, f_t(x,v_*) \, (1+|v|) \, (1+|v_*|^{k+1}) \dd x \dd v
\dd v_* 
\\ \le C_k \, M_k[\partial_x f_t] - K_k \, M_{k+1}[\partial_x f_t] +
C \, M_{k+1}[f_t].
 % \le C' _k \, M_k[\partial_x f_t] - K_k ' \, M_{k} ^{\frac{k-1}{k-2}}
%  [\partial_x f_t] +
% C'' \, M_{k+1}[f_t].  
\end{multline*}
We calculate similarly for any $|\alpha|\le 3$: 
\begin{equation*}
  \dt M_k[\partial_x ^\alpha f_t] \le C _k \, M_k[\partial_x ^\alpha f_t] - K_k  \, M_{k+1} 
 [\partial_x ^\alpha f_t] +
 C \, \sum_{\beta < \alpha} M_{k+1}[\partial_x ^\beta f_t].  
\end{equation*}
Finally choosing suitable constants $c_\alpha>0$, we deduce 
\begin{equation*}
  \dt \tilde M_k(t) \le C_k ' \, \tilde M_k(t) - K_k ' \, \tilde M_{k+1}(t)
  \le C_k ' \, \tilde M_k(t) - K_k '' \, \tilde M_{k} (t) ^{\frac{k-1}{k-2}}
\end{equation*}
which shows that 
\begin{equation*}
  \fa t \in (0,T'], \quad \tilde M_k(t) \le C_k \, \max\left\{
    \frac{1}{t^{k-2}}, 1 \right\}. 
\end{equation*}

We now consider exponential moments and extend the argument in
\cite{AlCaGaMo} to spatially inhomogeneous solutions in the torus.
Our goal is to estimate the quantity
\begin{multline*}
  E(t,z)
  := \sum_{|\alpha|\le 3} c_\alpha \int_{\T^3 \times \R^3} \left| \partial_x
    ^\alpha f_t(x,v)\right| \exp\big( z |v| \big) \dd x \dd v \\
  = \sum_{|\alpha|\le 3} c_\alpha \sum_{k=0}^\infty M_k[\partial_x ^\alpha f_t] \frac{z^k}{k!}
\end{multline*}
where $z$ will depend on time. For use below let us define the
truncated sum as
\begin{equation*}
  E^n (t,z)  :=   \sum_{|\alpha|\le 3} c_\alpha \sum_{k=0}^n M_k[\partial_x ^\alpha f_t] \frac{z^k}{k!}
\end{equation*}
for $n \in \mathbb N$, $z \ge 0$, and $t \geq 0$. We also define
\begin{equation*}
  I^n (t,z) :=  \sum_{|\alpha|\le 3} c_\alpha \sum_{k=0}^n M_{k+1}[\partial_x ^\alpha f_t] \frac{z^k}{k!}
\end{equation*}
and 
\begin{equation*}
    S_{\ell}(t)
    := \sum_{|\alpha|\le 3} c_\alpha \sum_{k=1}^{k_\ell} {\ell \choose k}
    \left( M_{k+1}[\partial_x ^\alpha f_t] M_{\ell-k}[\partial_x ^\alpha
      f_t]  + M_{k}[\partial_x ^\alpha f_t]  M_{\ell-k+1}[\partial_x
      ^\alpha f_t] \right),
  \end{equation*}
where $k_\ell$ is the integer part of $(\ell+1)/2$. 

Let us prove the following inequality: there exists some constant
$C>0$ independent of $n$ such that for any $\ell_0 \ge 2$ the
following holds:
  \begin{equation}
    \label{eq:convolution}
    \sum_{\ell=\ell_0} ^n \frac{z^\ell}{\ell!} S_{\ell}(t) \le C \,
    E^n (t,z) \, I^n (t,z).
  \end{equation}
  The first part of the sum in the left hand side of
  \eqref{eq:convolution} can be bounded as:
\begin{multline*}
    \sum_{|\alpha|\le 3} c_\alpha \sum_{\ell=\ell_0}^n \frac{z^\ell}{\ell!} \sum_{k=1}^{k_\ell} {\ell \choose k}
  M_{k+1}[\partial_x ^\alpha f_t] M_{\ell-k}[\partial_x ^\alpha
      f_t] \\
  = \sum_{|\alpha|\le 3} c_\alpha \sum_{\ell=\ell_0}^n \sum_{k=1}^{k_\ell}
  M_{k+1}[\partial_x ^\alpha f_t]  \, 
  \frac{z^{k}}{k!}  \, M_{\ell-k}[\partial_x ^\alpha f_t] \, \frac{z^{\ell-k}}{(\ell-k)!}
  \\
  \leq \sum_{|\alpha|\le 3} c_\alpha \sum_{k=1}^n M_{k+1}[\partial_x ^\alpha
  f_t]  \, 
  \frac{z^{k}}{k!}  \sum_{\ell=\max\{\ell_0, 2k-1\}}^n M_{\ell-k}[\partial_x ^\alpha f_t]   
  \frac{z^{\ell-k}}{(\ell-k)!} \\ \leq C \, I^n (t,z) \, E^n (t,z).
\end{multline*}
%where we have used $\beta \le s$ in the last line. 
We carry out a similar estimate for the other part:
\begin{multline*}
   \sum_{|\alpha|\le 3} c_\alpha \sum_{\ell=\ell_0}^n \frac{z^\ell}{\ell!} \sum_{k=1}^{k_\ell} {\ell \choose k}
  M_{k}[\partial_x ^\alpha f_t] M_{\ell-k+1}[\partial_x ^\alpha
      f_t] \\
  = \sum_{|\alpha|\le 3} c_\alpha \sum_{\ell=\ell_0}^n \sum_{k=1}^{k_\ell}
  M_{k}[\partial_x ^\alpha f_t]  \, 
  \frac{z^{k}}{k!}  \, M_{\ell-k+1}[\partial_x ^\alpha f_t] \, \frac{z^{\ell-k}}{(\ell-k)!}
  \\
  \leq \sum_{|\alpha|\le 3} c_\alpha \sum_{k=1}^n M_{k}[\partial_x ^\alpha
  f_t]  \,  \frac{z^{k}}{k!}  \sum_{\ell=\max\{\ell_0, 2k-1\}}^n M_{\ell-k+1}[\partial_x ^\alpha f_t]   
  \frac{z^{\ell-k}}{(\ell-k)!} \\ \leq C \, E^n (t,z) \, I^n (t,z).
\end{multline*}
This concludes the proof of \eqref{eq:convolution}. 

First we notice that in order to prove \eqref{eq:mts} it is enough to
prove the following: there are some constants $T_0\in(0,T)$ and
$\kappa, C > 0$ (which depend only on $b$ and the initial mass and
energy) such that
  \begin{equation}
    \label{eq:exp-generation}
    \n{f_t}_{L^1_v W^{3,1}_x(e^{\kappa t |v|})} \le C
    \quad
    \text{ for } t \in [0, T_0].
  \end{equation}
  Indeed, since the assumptions are satisfied on the whole time
  interval $[0,T)$, for $t \ge T_0$ it is then possible to apply
  \eqref{eq:exp-generation} starting at time $(t-T_0)$.

  Hence, we aim at proving the estimate \eqref{eq:exp-generation}. Let
  us denote 
\[
E_0 = E^n(0,0) = E(0,0) = \n{\init{f}}_{L^1_v W^{3,1}_x}.
\]
Consider $\kappa >0$ to be fixed later, $n \in \mathbb N$ and define $T_0 >
0$ as
\begin{equation*}
  \label{eq:def-T}
  T_0 := \min\Big\{ 1 \ ; \ \sup\big\{t > 0 \ \text{ s.t. } \ E^n
  (t,\kappa t) < 4 \, E_0 \big\}\Big\}.
\end{equation*}
The definition is consistent and the previous polynomial moment
estimates ensure that $T_0>0$ for each given $n$. The
bound of $1$ is not essential, and is included just to ensure that
$T_0$ is always finite. We note that a priori such $T_0$ depends on
the index $n$ in the sum $E^n$ but we will prove a uniform bound on
it.

Choose an integer $\ell_0 \ge 3$, to be fixed later. Arguing as in
\cite{AlCaGaMo}, by classical functional inequalities we have
\begin{equation*}
 \fa t \in [0,T), \ \ell \ge \ell_0, \quad \dt M_\ell[f_t] 
  \leq A_\ell \, S_\ell(t)
  - K \, M_{\ell+1}[f_t] 
%+ K_2  m_{\beta p}
\end{equation*}
with $S_\ell$ defined as before, $K>0$ uniform, and $A_\ell$ positive
decreasing and going to zero as $\ell \to \infty$. We can extend this
argument to higher derivatives at the price of an additional error
term as before:
\begin{multline*}
  \fa t \in [0,T), \ \ell \ge \ell_0, \quad \dt
  M_\ell[\partial^\alpha_x f_t] \\
  \leq A_\ell \, S_\ell(t)
  - K \, M_{\ell+1}[\partial_x ^\alpha f_t] 
+ C \, \sum_{\beta < \alpha} M_{\ell+1}[\partial^\beta _x f_t]. 
\end{multline*}
By linear combination with careful choice of the constants $c_\alpha$
we deduce that 
\begin{equation*}
\fa t \in [0,T), \ \ell \ge \ell_0, \quad  \dt \tilde M_\ell(t)
  \leq A_\ell \, S_\ell(t)
  - K \, \tilde M_{\ell+1}(t) 
\end{equation*}
for some uniform $K>0$ and $A_\ell$ positive decreasing going to zero
as $\ell \to 0$. 

In addition, we know from the previous polynomial estimates that 
\begin{equation}
  \label{eq:small-moments}
  \fa t \in [0,T), \quad \sum_{\ell=0}^{\ell_0} \tilde M_\ell(t)
  \,t^\ell \leq C_{\ell_0}. 
\end{equation}
% Here and in the rest of this proof $C$ will denote any constant that
% depends on $p_0$.
%
Taking any $\kappa \in (0,1)$ and using the product rule we get:
\begin{multline*}
  %\label{eq:di5}
  \td{}{t} \sum_{\ell=\ell_0}^n \tilde M_\ell(t) \, \frac{(\kappa
    t)^\ell}{\ell!}
  \\
  \leq \sum_{\ell=\ell_0}^n \frac{(\kappa t)^\ell}{\ell!}  \left(
    A_\ell \, S_\ell(t) - K \, \tilde M_{\ell+1}(t) \right) + \kappa
  \, \sum_{\ell=\ell_0}^n \tilde M_\ell (t) \, \frac{(\kappa
    t)^{\ell-1}}{(\ell-1)!}
  \\
  \leq \sum_{\ell=\ell_0}^n \frac{(at)^\ell}{\ell !} A_\ell \,
  S_\ell(t) + (\kappa - K) \, I^n (t,\kappa t) + (K +\kappa) \, \sum_{\ell=1}^{\ell_0}
  \tilde M_\ell(t)\,\frac{(\kappa t)^{\ell-1}}{(\ell-1)!}
  \\
  \leq \sum_{\ell=\ell_0}^n \frac{(\kappa t)^\ell}{\ell !} A_\ell \, S_\ell
  + (\kappa- K)\, I^n(t,\kappa t) +\frac{(K+\kappa )}{t} \, C_{\ell_0},
\end{multline*}
where we have used that $\kappa < 1$ and
inequality~\eqref{eq:small-moments} in the last step. Hence, from the
inequality~\eqref{eq:convolution} we obtain
\begin{equation*}
  \td{}{t} \sum_{\ell=\ell_0}^n \tilde M_\ell(t) \frac{(\kappa t)^\ell}{\ell!}
   \leq I^n(t,\kappa t) \Big[
    C \, A_{\ell_0} \, E^n(t,\kappa t) + (\kappa - K) \Big] 
+ \frac{(K+\kappa )}{t} \,  C_{\ell_0}.
 \end{equation*}
 Next, choose $\kappa \le\min\{1, K/2\}$ and $\ell_0$ large enough so that 
 \begin{equation*}
   \fa t \in [0,T_0], \quad C \, A_{\ell_0} \, E^n(t,\kappa t) \le C \,
   A_{\ell_0} \, 4 \, E_0 \le \frac{K}{4}.
 \end{equation*}
Hence
\begin{multline}
  \label{eq:di7}
  \td{}{t} \sum_{\ell=\ell_0}^n \tilde M_\ell(t) \frac{(\kappa t)^\ell}{\ell!}
    \leq - \frac{K}{4} \, I^n (t,\kappa t)
  + \frac{(K+\kappa )}{t} \, C_{\ell_0}
  \\
  \leq 
  - \frac{1}{t}  \left[
    \frac{K}{4\kappa } \left(E^n  (t,\kappa t) - E_0\right)
    - (K+\kappa ) \, C_{\ell_0} \right]
\end{multline}
where for the last inequality we have used that (thanks to the
conservation of the total mass) 
\[
I^n (t,\kappa t) \geq \frac{(E^n (t,\kappa t) - E_0)}{\kappa t}.
\]
We make the additional restriction that $\kappa < E_0 /
(6C_{\ell_0})$, which together with $\kappa \le\min\{1, K/2\}$ implies
that
\begin{equation*}
  \frac{K}{4 \kappa} \, E_0 > (K+\kappa) \, C_{\ell_0}.
\end{equation*}
Then we have 
% and using again that
% $E^n _\beta (t,at) \leq 4m_0$ on $[0,T]$, \eqref{eq:di7} implies that
\begin{equation}
  \label{eq:di8}
  \td{}{t} \sum_{\ell=\ell_0}^n \tilde M_\ell (t) \frac{(\kappa t)^\ell}{\ell!} 
  \leq 0
%   K_2 E^n _\beta (t,at)
% %  + (K_2+a) m_1
% %  \\
%   \leq
%   4 K_{2} m_0
%   % + (K_2+a) m_1
%   =: K_{3}
\end{equation}
for any time $t \in [0,T_0]$ for which $E^n (t,\kappa t) \geq 2 \,
E_0$ holds. This is true in particular when $\sum_{\ell=\ell_0}^n
\tilde M_\ell(t) \, \frac{(\kappa t)^\ell}{\ell!} \geq 2 \, E_0$. We deduce
that
\begin{equation}
  \label{eq:large-moments-bound}
  \fa t \in [0,T_0], \quad \sum_{\ell=\ell_0}^n \tilde M_\ell
  \frac{(\kappa t)^\ell}{\ell!} 
  \leq 2\, E_0.
\end{equation}
In order to finish the argument we need to bound the initial part of
the full sum (from $\ell=0$ to $\ell_0-1$.) Indeed, we note that from
\eqref{eq:small-moments},
\begin{equation}
  \label{eq:small-moments-2}
  \fa t \in [0,T_0], \quad \sum_{\ell=0}^{\ell_0-1} \tilde M_\ell(t)
  \, \frac{(\kappa t)^\ell}{\ell!}
  \leq E_0 + \kappa \, C_{\ell_0}, 
\end{equation}
so, recalling that $6 \kappa C_{\ell_0} < E_0$ and using
\eqref{eq:large-moments-bound} and \eqref{eq:small-moments-2} we get
\begin{equation*}
  E^n (t,\kappa t) =
  \sum_{\ell=0}^{\ell_0-1} \tilde M_\ell(t) \, \frac{(\kappa t)^\ell}{\ell !}
  + \sum_{\ell=\ell_0}^n \tilde M_\ell(f) \, \frac{(\kappa t)^\ell}{\ell!}
  \leq 3 \, E_0 + \kappa \, C_{\ell_0} \leq \frac{19}{6} \, E_0 
\end{equation*}
for $t \in [0, T_0]$, uniformly in $n$. %This shows that $T \geq 2 m_0 / K_3$, where $K_3$ is a constant which
%depends only on $b$ and the initial mass and energy. 
%
Finally, gathering all conditions imposed along the proof on the
parameter $\kappa$, we choose
\begin{equation}\label{choose-a}
  \kappa :=
  \min \left\{1, \frac{K}{2}, \frac{E_0}{6C_{\ell_0}} \right\}
\end{equation}
independently of $n$.  We conclude, from the definition of $T_0$, that
$T_0=1$ for all $n$.  Sending $n\rightarrow\infty$, we deduce the
result. 
\end{proof}

\subsubsection{Non appearance of ``superlinear'' exponential moments}
\label{sec:non-appe-superl}

\begin{lem}\label{lem:exp-moments-bis}
  Consider for $T \in (0,+\infty]$ a conservative solution 
   \[
   0 \le f \in L^1_{t,loc}\left([0,T),L^1_v W^{3,1}_x(1+|v|^2)\right) \cap
   C^0_t\left([0,T),L^1_v W^{3,1}_x(1+|v|)\right)
  \]
  with initial datum bounded uniformly from below as in
  ~\eqref{eq:hyp-lower}.  Assume that for $\beta \in (1,2]$ the
  initial condition satisfies
  \begin{equation*}
    \fa \kappa >0, \quad \n{\initem{f}}_{L^1_v W^{3,1}_x(e^{\kappa \, |v|^\beta})} = +\infty.
  \end{equation*}
  Then we have 
\begin{equation*}
  \fa t \ge 0, \ \fa \kappa >0, \quad \n{f_t}_{L^1_v W^{3,1}_x(e^{\kappa \, |v|^\beta})} = +\infty.
  \end{equation*}
\end{lem}

\begin{proof}[Proof of Lemma~\ref{lem:exp-moments-bis}]
  We only sketch the proof in the case $\beta =2$ and leave to the
  reader the general case. The key idea is to define
  \begin{equation*}
    E^n _R(t,z)  :=   \sum_{|\alpha|\le 3} c_\alpha \sum_{k=0}^n 
    \left( \int_{\T^3 \times \R^3} \left|\partial_x ^\alpha f_t\right|
      \, (1+|v|)^{2k} \,
    {\bf 1}_{|v| \le R} \, \dd x \dd v \right) \frac{z^k}{k!}
  \end{equation*}
  for some parameter $R >0$, and then consider
  $E^n_R(t,\kappa(1+\kappa' t))$ with $\kappa$ arbitrary and $\kappa'$
  to be fixed later. Observe that $E^n_R(t,z)$ is always well-defined
  and finite for all time and value of $z$ due to the truncations. We
  calculate (dropping out the positive terms)
\begin{multline*}
  \dt E^n_R\left(t,\kappa(1+ \kappa' t)\right) \\ \ge - K \, \sum_{|\alpha|\le 3}
  c_\alpha \sum_{k=0}^n \left( \int_{\T^3 \times \R^3}
    \left|\partial_x ^\alpha f_t\right| \, {\bf 1}_{|v| \le R} \,
    (1+|v|)^{2k+1} \, \dd x \dd v \right) \,
  \frac{\left(\kappa(1+\kappa' t)\right)^k}{k!} \\
  + \kappa \kappa' \, \sum_{|\alpha|\le 3} c_\alpha \sum_{k=1}^n \left(
    \int_{\T^3 \times \R^3} \left|\partial_x ^\alpha f_t\right| \,
    {\bf 1}_{|v| \le R} \, (1+|v|)^{2k} \, \dd x \dd v \right) \,
  \frac{\left(\kappa(1+\kappa' t)\right)^{k-1}}{(k-1)!}.
\end{multline*}
We deduce that for $\kappa'$ large enough
\begin{multline*}
  \dt E^n_R\left(t,\kappa(1+ \kappa't)\right) \\ \ge - K \, \sum_{|\alpha|\le 3} c_\alpha
  \left( \int_{\T^3 \times \R^3} \left|\partial_x ^\alpha f_t\right|
    \, {\bf 1}_{|v| \le R} \, (1+|v|)^{2n+1} \, \dd x \dd v \right) \,
  \frac{\left(\kappa(1+\kappa' t)\right)^n}{n!}.
\end{multline*}
Since the right hand side goes to zero as $n \to +\infty$ we deduce the a priori
estimate
\begin{equation*}
  \dt E_R ^\infty\left(t,\kappa(1+ \kappa' t)\right) \ge 0. 
\end{equation*}
We hence deduce by passing to the limit $R \to \infty$ that $E^\infty
_\infty(t,\kappa(1+\kappa' t)) = +\infty$ for $t \ge 0$ which concludes
the proof.
  \end{proof}

\subsubsection{A priori lower bounds}
\label{sec:priori-lower-bounds}

The proof of the Maxwellian lower bound in part {\bf (I)} of
Theorem~\ref{theo:NLBE} is a straightforward application of
\cite{Mlower} and we shall therefore skip the proof. In the
paper~\cite{Mlower} an a priori bound was assumed on the entropy but
it can be removed using the non-concentration estimates on the
iterated gain term first discovered in \cite{MR1697562} and then
developed in \cite{MR1697495}. We refer to the more recent preprint
\cite{briant} where these issues are discussed.

\subsubsection{A priori uniqueness for conservative solutions}
\label{sec:priori-uniq-crit}

This subsection is related to the Cauchy theory for unique solution to
the spatially homogeneous Boltzmann for hard spheres in
$L^1_v(1+|v|^2)$. Let us refer first to \cite{MR0356798} for the idea
of the key a priori estimate on moment of the difference of two
solutions and \cite{MR0339665,MR0339666} for the first uniqueness
result in a space of the form $L^1_v(1+|v|^k)$ (with $k>2$). Then we
refer to \cite{MR1697562,MR1716814} (and later \cite{lm1} following
the same approach) for the more recent optimal results. In these
papers, there are mainly two approaches. The first one
\cite{MR1697562} relies on a subtle variants of the Povzner
inequality, and the second one \cite{MR1716814} (see also \cite{lm1})
is more direct and relies on the estimate of the tail of the
distribution at initial times. We shall elaborate upon this second
approach in this subsection.

\begin{lem}[A priori uniqueness in $L^1_vL^\infty_x (1+|v|^k)$, $k>2$]
\label{lem:presque-energy}
Consider for $T \in (0,+\infty]$ and $k>2$ two conservative
distributional solutions
\[
f_t, g_t \in L^1_{t,loc}\left([0,T),L^1_v L^\infty_x(1+|v|^k)\right) \cap
C^0_t\left([0,T),L^1_v L^\infty(1+|v|^{k-1})\right)
\]
with initial data $\initem{f}, \initem{g}$ satisfying the lower bound
assumption~\eqref{eq:hyp-lower}. Then for any $T' \in [0,T)$ there is
some constant $C(T')>0$ (depending on the bounds assumed on the
solutions, on the lower bound \eqref{eq:hyp-lower} on the initial
datum, and on $T'>0$) such that
  \begin{equation}
    \label{eq:estimate-presque-energy}
    \fa t \in [0,T'], \quad \n{f_t -g_t}_{L^1_{x,v}(1+|v|^2)}
    \le C(T') \, \n{\initem{f} - \initem{g}}_{L^1_{x,v}(1+|v|^2)}.
  \end{equation}
\end{lem}
\begin{proof}[Proof of Lemma~\ref{lem:presque-energy}]
%   The following lower bound estimate on $f_t$ holds and is obtained
%   through the Duhamel principle with the bound on the energy:
% \begin{equation*}
%   \fa t \ge 0, \quad f_t(x,v) \ge \init{f}(x,v) \, e^{-C \, t}.
% \end{equation*}
Arguing as before we get 
\begin{equation*}
  \fa t \in [0,T'], \quad \int_{\R^3} f_t(x,v_*) \, |v-v_*| \dd v_* % \ge
  % e^{-C \, t} \int_{\R^3} \init{f}(x,v_*) \, |v-v_*| \dd v_*  \\
  % \ge e^{-C \, t} \, \int_{\R^3} \varphi(v_*) \, |v-v_*| \dd v_* 
  \ge
  K(T') \, (1+|v|)
\end{equation*}
for some constant depending on the $L^\infty_t([0,T'],L^1_v
  L^\infty(1+|v|^{k-1}))$ norm of $f$ and the lower bound
  \eqref{eq:hyp-lower} on $\init{f}$. 

We then write the estimate (arguing as in the previous section) 
\begin{multline*}
  \fa t \in [0,T'], \quad \dt \n{f_t}_{L^1_v L^\infty_x (1+|v|^k)} \\ \le C \, \n{f_t}_{L^1_v
    L^\infty_x (1+|v|^2)} \, \n{f_t}_{L^1_v L^\infty_x (1+|v|^k)} - K
  \, \n{f_t}_{L^1_v L^\infty_x (1+|v|^{k+1})}
\end{multline*}
which shows that $\n{f_t}_{L^1_v L^\infty_x (1+|v|^{3})}$ is
time-integrable. Similarly we deduce that $\n{g_t}_{L^1_v L^\infty_x
  (1+|v|^{3})}$ is time-integrable on $[0,T']$. Finally we obtain the
continuity of the flow in the topology $L^1_v L^\infty_x(1+|v|^2)$:
\begin{multline*}
  \dt \n{f_t - g_t}_{L^1_v L^\infty_x(1+|v|^2)} \\ \le C \, \left( \n{f_t}_{L^1_v
      L^\infty_x(1+|v|^3)} + \n{g_t}_{L^1_v L^\infty_x(1+|v|^3)}
  \right) \, \n{f_t - g_t}_{L^1_v L^\infty_x(1+|v|^2)} 
\end{multline*}
and thus 
\begin{multline*}
  \fa t \in [0,T'], \quad \n{f_t - g_t}_{L^1_v L^\infty_x(1+|v|^2)}  \\
  \le C \, e^{\int_0 ^t \left( \n{f_\tau}_{L^1_v L^\infty_x(1+|v|^3)}
      + \n{g_\tau}_{L^1_v L^\infty_x(1+|v|^3)} \right) \dd \tau} \,
  \n{\init{f} - \init{g}}_{L^1_v L^\infty_x(1+|v|^2)}
\end{multline*}
and the claimed uniqueness property follows.
\end{proof}

The next lemma follows an idea first introduced for the spatially
homogeneous Boltzmann equation in \cite{MR1697562,MR1716814}, using
the reformulation in \cite{lm1}.

\begin{lem}[A priori uniqueness in the critical case  $k=2$]\label{lem:energy}
  Consider for $T \in (0,+\infty]$ two conservative distributional solutions 
\[
f_t, g_t \in L^1_{t,loc} \left([0,T),L^1_v W^{3,1}_x(1+|v|^2)\right) \cap
C^0_t \left([0,T),L^1_v W^{3,1}_x(1+|v|)\right)
\]
with initial data $\initem{f}, \initem{g}$ satisfying the lower bound
assumption~\eqref{eq:hyp-lower}. 

Then for any $T' \in (0,T)$, there is an explicit function $\Psi :
\R_+ \to \R_+$ which depends on $T'>0$, $\initem{f}$ and $\initem{g}$,
which is continuous and satisfies $\Psi(0)=0$ and $\Psi(r)>0$ for
$r>0$, such that
  \begin{equation}
    \label{eq:estimate-energy}
    \fa t \in [0,T'], \quad \n{f_t -g_t}_{L^1_{x,v}(1+|v|^2)}
    \le \Psi\left( \n{\initem{f} - \initem{g}}_{L^1_{x,v}(1+|v|^2)} \right).
  \end{equation}
\end{lem}

\begin{proof}[Proof of Lemma~\ref{lem:energy}]
  We fix $T'\in (0,T)$ for the whole proof. Arguing exactly as in the
  first part of Lemma~\ref{lem:mts}, we deduce that there is a
  constant $K(T') >0$ so that
\begin{equation*}
  \fa t \in [0,T'], \ x \in \T^3, \ v \in \R^3, \quad \left\{ 
    \begin{array}{l} \ds
      \int_{\R^3} f_t(x,v_*) \, |v-v_*| \dd v_* \ge K(T') \, (1+|v|), \vs \\ \ds 
      \int_{\R^3} g_t(x,v_*) \, |v-v_*| \dd v_* \ge K(T') \, (1+|v|),
    \end{array}
  \right.
  \end{equation*}
  and
\begin{equation*}
  \fa t \in (0,T'], \quad 
  \left\{ 
    \begin{array}{l}\ds 
      \tilde M_k(t) \le C_k(T') \, \min\left\{ \frac{1}{t^{k-2}}, 1 \right\} \vs \\ \ds
      \tilde M_k(t) \le C_k(T') \, \min\left\{ \frac{1}{t^{k-2}}, 1 \right\}
    \end{array}
  \right.  
\end{equation*}
for some constant $C_k(T')$ depending on $T'>0$ and $k>2$, and where
$\tilde M_k$ was defined in the proof of Lemma~\ref{lem:mts} (recall
that it involves the derivatives $\partial_x ^\alpha$, $|\alpha| \le
3$).

Let us denote $d_t := f_t - g_t$ and $s_t := f_t + g_t$. 

We have by usual calculations
\begin{multline*}
  \dt \int_{\T^3 \times \R^3} |d_t| \, (1+|v|^2) \dd x \dd v \\ \le C
  \, \left( \int_{\T^3 \times \R^3} |d_t| \dd x \dd v \right) \,
  \left( \sup_{x \in \T^3} \int_{\R^3} |s_t| \, (1+|v|^3) \dd v
  \right) \\ + C \, \, \left( \int_{\T^3 \times \R^3} |d_t| \,
    (1+|v|) \dd x \dd v \right) \, \left( \sup_{x \in \T^3}
    \int_{\R^3} |s_t| \, (1+|v|^2) \dd v
  \right)\\
  \le C_1 \, \min\left\{ \frac{1}{t}, 1 \right\} \, \left( \int_{\T^3
      \times \R^3} |d_t| \dd x \dd v \right) + C_2 \, \left(
    \int_{\T^3 \times \R^3} |d_t| \, (1+|v|) \dd x \dd v \right)
\end{multline*}
which provides a simple Gronwall-like estimates for times bounded away
from zero. 

Let us now consider small times. Define 
\[
r := \min \left\{ \n{\init{d}}_{L^1_{x,v}(1+|v|^2)} \, ; \, T' \right\}
\]
and let us estimate the $L^1_{x,v}(1+|v|^2)$ norm of the difference
for the times $t \in [0,r]$. Then calculate
\begin{multline*}
  \fa t \in [0,r], \quad \n{d_t}_{L^1_{x,v}(1+|v|^2)} \\ \le
  \int_{\T^3 \times \R^3} d_t \, (1+|v|^2) \dd x \dd v + 2 \,
  \int_{\T^3
    \times \R^3}(d_t)_+ \, (1+|v|^2) \dd x \dd v \\
  \le \int_{\T^3 \times \R^3} \init{d} \, (1+|v|^2) \dd x \dd v + 2 \,
  \int_{\T^3
    \times \R^3}f_t \, (1+|v|^2) \dd x \dd v \\
  \le \int_{\T^3 \times \R^3} |\init{d}| \, (1+|v|^2) \dd x \dd v + 2
  \, \int_{\T^3
    \times \R^3}f_t \, (1+|v|^2) \dd x \dd v \\
  \le r + 2 \, \int_{|v| \le R} f_t \, (1+|v|^2) \dd x \dd v + 2 \,
  \int_{|v| > R} f_t \, (1+|v|^2) \dd x \dd v \\ 
  \le r + 2 \, (1+R^2) \, \n{d_t}_{L^1_{x,v}} + 2 \,
  \int_{|v| > R} f_t \, (1+|v|^2) \dd x \dd v
\end{multline*}
for some parameter $R>0$ to be chosen later, where we have used the
conservation of the energy of our solutions and the inequality
$(d_t)_+ \le f_t$.

The second term in the right hand side above can be estimated as
\begin{multline*}
  \dt \int_{\T^3 \times \R^3} |d_t| \dd x \dd v \le C\, \int_{\T^3
    \times \R^3} d_t (s_t)_* \, |v-v_*| \dd x \dd v \dd v_* % \\ \le C
  % \, \left( \int_{\T^3 \times \R^3} |d_t| \, (1+|v|) \dd x \dd v
  % \right) \, \left( \int_{\T^3 \times \R^3} |s_t| \, (1+|v|) \dd x
  %   \dd v
  % \right)
  \\
  \le C' \, \left( \int_{\T^3 \times \R^3} |d_t| \, (1+|v|) \dd x \dd
    v \right).
\end{multline*}
Hence
\begin{multline*}
  \fa t \in [0,r], \quad \n{d_t}_{L^1_{x,v}} \\ \le \n{\init{d}}_{L^1_{x,v}(1+|v|^2)} + C' \,
  \int_0 ^t \left( \n{f_\tau}_{L^1_{x,v}(1+|v|^2)}
    \n{g_\tau}_{L^1_{x,v}(1+|v|^2)} \right) \dd \tau 
  \le C'' \, r.
\end{multline*}
Finally the third term of the right hand side can be estimated as
\begin{multline*}
  \int_{|v| > R} f_t \, (1+|v|^2) \dd x \dd v = \int_{\T^3 \times
    \R^3} f_t \, (1+|v|^2) \dd x \dd v - \int_{|v| \le R} f_t \,
  (1+|v|^2) \dd x \dd v
  \\
  = \int_{\T^3 \times \R^3} \init{f} \, (1+|v|^2) \dd x \dd v -
  \int_{|v| \le R} f_t \, (1+|v|^2) \dd x \dd
  v \\
  = \int_{\T^3 \times \R^3} \init{f} \, (1+|v|^2) \dd x \dd v -
  \int_{|v| \le R} \init{f} \, (1+|v|^2) \dd x \dd
  v  \\
  - \int_0 ^t \int_{|v| \le R} Q(f_\tau,f_\tau) \, (1+|v|^2) \dd x \dd
  v \dd \tau  \\
  \le \int_{|v| > R} \init{f} \, (1+|v|^2) \dd x \dd v + \int_0 ^t
  \int_{|v| \le R} Q^-(f_\tau,f_\tau) \, (1+|v|^2) \dd x \dd v \dd
  \tau \\
  \le \int_{|v| > R} \init{f} \, (1+|v|^2) \dd x \dd v + C''' \, r \, (1+R^2) 
\end{multline*}
where we have used again the conservation of energy and the evolution
equation integrated against $(1+|v|^2) \, {\bf 1}_{|v| \le R}$. 

Combining the three estimates we deduce that
\begin{multline*}
  \fa t \in [0,r], \quad \n{d_t}_{L^1_{x,v}(1+|v|^2)} \le \\
r + 2 \, C'' \, r \, (1+R^2) + 2 \, \int_{|v|>R} \init{f} \, (1+|v|^2)
\dd x \dd v + C''' \, r \, (1+R^2).
\end{multline*}
We finally choose for instance $R = r^{-1/3}$ and define
\begin{multline*}
  \Psi_0(r) := r + 2 \, C'' \, r \,
  \left(1+r^{-2/3}\right)  \\
  + 2 \, \int_{|v|>r^{-1/3}} \init{f} \, (1+|v|^2)
\dd x \dd v + C''' \, r \, \left(1+r^{-2/3}\right)
\end{multline*}
which depends on the profiles $\init{f}$ and $\init{g}$ via the tail
estimate in the right hand side and also via the constants depending
on the mass and energy.

We have therefore 
\begin{equation*}
  \fa t \in [0,r], \quad \n{d_t}_{L^1_{x,v}(1+|v|^2)} \le \Psi_0\left(
    r \right). 
\end{equation*}

To conclude with the final stability estimate in the case $r < T'$, we
write
\begin{multline*}
  \fa t \in [0,T'], \quad \n{d_t}_{L^1_{x,v}(1+|v|^2)} \\ \le
  \n{d_r}_{L^1_{x,v}(1+|v|^2)} + \int_r ^t \left( 
  \frac{{\rm d}}{{\rm d}\tau} \n{d_\tau}_{L^1_{x,v}(1+|v|^2)} \right) \dd \tau \\ 
\le \Psi_0\left( r \right) + \int_r ^t
  \left( C_1 \, \min\left\{ \frac{1}{\tau}, 1 \right\} \,
  \n{d_\tau}_{L^1_{x,v}} + C_2 \, 
  \n{d_\tau}_{L^1_{x,v}(1+|v|)} \right) \dd \tau. 
\end{multline*}
If $T' \ge r \ge 1$ the proof is clear by a Gronwall estimate, for $r
<1$ we write first (assuming $T' \ge 1$ for notational simplicity, the
case $T'<1$ is similar)
\begin{multline*}
  \fa t \in [0,T'], \quad \n{d_t}_{L^1_{x,v}(1+|v|^2)} 
\le \Psi_0\left(r \right) \\ + C_1 \, \int_r
  ^{1} 
    \n{d_\tau}_{L^1_{x,v}} \frac{\dd \tau}{\tau} + (C_1 + 2 C_2)\, \int_{1} ^{T'}
    \n{d_\tau}_{L^1_{x,v}(1+|v|^2)} \dd \tau
\end{multline*}
and for the second term of the right hand side we use the estimate on
$\n{d_\tau}_{L^1_{x,v}}$: 
\begin{multline*}
  \int_r ^{1} \n{d_\tau}_{L^1_{x,v}} \frac{\dd \tau}{\tau}
\le \int_r ^{1} \left( \n{\init{d}}_{L^1_{x,v}} + C \,
  \int_0 ^\tau \n{d_{\tau'}}_{L^1_{x,v}} \dd \tau' \right) \frac{\dd
  \tau}{\tau} \\ 
\le r |\ln r| + C \, \int_0 ^1 \n{d_{\tau'}}_{L^1_{x,v}} \, |\ln
\tau'| \dd \tau'.
\end{multline*}

We thus deduce 
\begin{multline*}
  \fa t \in [0,T'], \quad \n{d_t}_{L^1_{x,v}(1+|v|^2)} 
\le \Psi_0\left(r \right) + r |\ln r| \\ + C_1' \, \int_0
  ^{1} 
    \n{d_\tau}_{L^1_{x,v}} \, |\ln \tau| \dd \tau + C_2 ' \, \int_{1} ^{T'}
    \n{d_\tau}_{L^1_{x,v}(1+|v|^2)} \dd \tau
\end{multline*}
which yields the result for some nonlinear function $\Psi=\Psi(r)$ by
the Gronwall lemma. 
\end{proof}

\subsection{Proof of Theorem~\ref{theo:NLBE}, parts {\bf (II)} and
  {\bf (III)}}

\subsubsection{A dissipative Banach norm}
\label{sec:diss-banach-norm}

In this subsection we construct a Banach norm for which the semigroup
is not only dissipative, but also has a stronger dissipativity
property: the damping term in the energy estimate controls the
\emph{norm of the graph} of the collision operator.

Observe that in this theorem, the rate of decay is possibly worse than
in Theorem~\ref{theo:LIBE2}. It shall not however cause any problem
when searching for the rate of decay of the nonlinear equation, as the
latter can be recovered by a bootstrap argument once the stability is
proved.

\begin{prop}\label{theo:LIBE3} 
  Consider the space $\EE = W^{\sigma,q}_v W^{s,p}_x (m)$ with the
  same assumptions as in Theorem~\ref{theo:LIBE2}, with a norm denoted
  by $\| \cdot \|_\EE$, and define the equivalent norm
\begin{equation}\label{def:NtEE}
\Nt h \Nt_{\EE} :=   \eta \, \| h\|_{\EE}   + \int_0 ^{+\infty} \left\|
  S_{\LL} (\tau) h \right\| _{\EE} \dd \tau, \quad \eta >0. 
\end{equation}

Then there exists $\eta > 0$ (small enough) and $\lambda_1 \in
(0,\lambda)$ such that for any $\initem{h} \in \EE$, $\Pi \initem{h} =
0$ (let us recall that $\Pi$ is the projection on the eigenspace
associated to the eigenvalue $0$ thanks to the formulas
\eqref{def:SpectralProjection} and \eqref{eq:project-explicit}), the
solution $h(t) := S_\LL(t) \initem{h}$ to the linearized flow
\eqref{eq:LIBE} satisfies:
\begin{equation*}
  \forall \, t \ge 0, \quad 
  {{\rm d} \over {\rm d}t} \Nt h_t \Nt_\EE \le -\lambda_1 \,  \Nt h_t \Nt_{ \EE_\nu},
\end{equation*}
where 
\[
\EE_\nu := W^{\sigma,q}_v W^{s,p}_x (\nu^{1/q} \, m) \subset \EE
\]
and $\Nt \cdot \Nt_{\EE_\nu}$ is defined from $\| \cdot \|_{\EE_\nu}$
as in \eqref{def:NtEE}:
\[
\Nt h \Nt_{\EE_\nu} :=   \eta \, \| h\|_{\EE_\nu}   + \int_0 ^{+\infty} \left\|
  S_{\LL} (\tau) h \right\| _{\EE_\nu} \dd \tau.
\]
\end{prop}

\begin{proof}[Proof of Theorem~\ref{theo:LIBE3}] 
  From the decay property of $\LL$ provided by
  Theorem~\ref{theo:LIBE2} we have
\[
\left\| S_{\LL} (\tau) h \right\|_{\EE} \le C \, e^{- \lambda \, t} \,
\| h \|_{\EE}.
\]
Therefore we deduce that
\[
C_1(\eta) \, \| h \|_{\EE} \le \Nt h \Nt_{\EE}   \le C_2(\eta) \, \|
h\|_{\EE}
\]
for some constants $C_1(\eta), C_2(\eta) >0$ depending on $\eta$,
i.e. the norms $\| \cdot \|_{\EE}$ and $\Nt \cdot \Nt_{\EE} $ are
equivalent for any $\eta > 0$. 

Let us now compute the time derivative of the norm $\Nt \cdot \Nt_\EE$
along $h_t$ which solves the linear evolution problem
\eqref{eq:LIBE}. Observe that $\Pi h_t =0$ for any time $t \ge 0$ due
to the mass, momentum and energy conservation of the linearized
Boltzmann equation.

Since the $x$-derivatives commute with the linearized operator, we can
set $s=0$ without loss of generality. We consider first $\sigma =0$
and $p,q \in [1,+\infty)$.  We denote again $\Phi'(z) := |z|^{p-1} \,
\hbox{sign} (z)$ and we have
\begin{multline*}
  \frac{{\rm d}}{{\rm d}t} \Nt h_t \Nt_\EE = \eta \, \| h_t
  \|_\EE^{1-q} \, \int_{\R^3} \left( \int_{\T^3} \LL (h_t) \,
    \Phi'(h_t) \dd x \right) \, \| h_t \|_{L^p_x} ^{q-p} \, m^q \dd v
  \\
  + \int_0 ^{+\infty} \frac{\partial}{\partial t} \left\| h_{t+\tau}
  \right\|_\EE \dd \tau =: I_1 + I_2 .
\end{multline*}

Concerning the first term $I_1$ we have, arguing as in the proof of
Lemma~\ref{lem:decomp-new} (cases {\bf (W2)}-{\bf (W3)}):
\[
I_1 = \eta  \, \| h_t \|_\EE^{1-q} \,
  \int_{\R^3} \left( \int_{\T^3} \left(  \AA_\delta + 
      \BB_\delta \right) (h_t) \, \Phi'(h_t) \dd x \right)
  \, \| h_t \|_{L^p_x} ^{q-p} \, m^q \dd v 
\]
where we have dropped the transport term thanks to its divergence
structure. Thanks to the dissipativity of $B_\delta$ proved in
Lemma~\ref{lem:decomp-new} and the bounds on $\AA_\delta$ in
Lemma~\ref{lem:regA} we get
\[
I_1  \le \eta \, \left( C \, \| h \|_\EE - K \, \| h \|_{\EE_\nu}
\right) 
\]
for some constants $C, K>0$. 

The second term is computed exactly: 
\[
I_2 = \int_0 ^{+\infty} \frac{\partial}{\partial t} \left\|
      h_{t+\tau}  \right\|_\EE \dd \tau = \int_0 ^{+\infty} \frac{\partial}{\partial \tau} \left\|
      h_{t+\tau}  \right\|_\EE \dd \tau = - \| h \|_\EE.
\]

The combination of the two last equations yields the desired result 
\[
\frac{{\rm d}}{{\rm d}t} \Nt h_t \Nt_\EE \le - K \, \Nt h_t
\Nt_{\EE_\nu}
\]
with $K>0$, by choosing $\eta$ small enough.

Then the cases $p=+\infty$ and $q = +\infty$ are obtained by passing
to the limit. 

Finally the case of a higher-order $v$-derivative is treated by an
argument close to the one in Lemma~\ref{lem:decomp-new}. For instance
the case $\sigma=s=1$ is proved by introducing the norms
$$
\left\{ 
\begin{array}{lcl} \ds
 \Nt h \Nt_{\EE_\var} & := & \Nt h \Nt_\EE + \Nt \nabla_x h
 \Nt_\EE + \var \, \Nt \nabla_v h \Nt_\EE, \vs \\ \ds 
\Nt h \Nt_{\EE_{\nu,\var}} & := & \Nt h \Nt_{\EE_\nu} + \Nt \nabla_x h
 \Nt_{\EE_\nu} + \var \, \Nt \nabla_v h \Nt_{\EE_\nu}, 
\end{array}
\right.
$$
for some second parameter $\var>0$ small enough. Arguing as before we obtain
\begin{multline*} 
  \dt \left( \Nt h_t \Nt_{L^q_v L^p_x(m)} + \Nt \nabla_x h_t \Nt _{L^q_v
      L^p_x(m)} \right)
  \\
  \le - K_1 \, \Bigl( \Nt h_t \Nt_{L^q_v L^p_x(m\nu^{1/q})} + \Nt
  \nabla_x h_t \Nt_{L^q_v L^p_x(m\nu^{1/q})} \Bigr)
%  \\
%  \left[ \left( \int_{\R^3} \|h \|_{L^p_x} ^q \, \nu \, m^q \dd v
%    \right)^{\frac{1}{q}} + \left( \int_{\R^3} \| \nabla_x h
%      \|_{L^p_x} ^q \, \nu \, m^q \dd v \right)^{\frac{1}{q}} \right]
\end{multline*}
and 
\begin{multline*}
\dt \Nt \nabla_v h_t \Nt _{L^q_v L^p_x(m)} \\ \le
- K_2  \,   \Nt \nabla_v h_t \Nt _{L^q_v L^p_x(m\nu^{1/q})}  
   + \Nt \nabla_x h_t \Nt_{L^q_v L^p_x(m)}
  + \Nt \RR h_t \Nt_{L^q_v L^p_x(m)}, 
\end{multline*}
where $\RR$ is defined in \eqref{eq:RRdef}. Using (a) the
Lemmas~\ref{lem:Aintegral} and \ref{lem:decQ+rad} when $m$ is a
polynomial weight, (b) \eqref{eq:Ac-exp-L1} and
Lemma~\ref{lem:decQ+rad-exp} when $m$ is an exponential weight, (c)
the regularization property of the operator $\AA_\delta$, (d) the
equivalence of the norms $\Nt \cdot \Nt$ and $\n{\cdot}$, we prove
that
\[
\Nt\RR h_t \Nt_{L^q_v L^p_x(m)} \le C \, \Nt h_t \Nt_{L^q_v L^p_x(m\nu^{1/q}))}
\] 
for some constant $C >0$. We deduce that for $\var$ small enough
 \begin{equation*} 
\dt  \Nt h_t \Nt_{\EE_\var} \le
   -K_3 \, \Nt h_t \Nt_{\EE_{\nu,\var}}
\end{equation*} 
for some $K_3 >0$. The higher-order estimates are performed with the
norm
$$
 \left\{ 
   \begin{array}{l} \ds \Nt h \Nt_{\EE_\var} := \sum_{|i|\le \sigma, \
       |j|\le s, \ |i|+|j| \le \max\{\sigma; s\}}
     \var^{|i|} \, \Nt \partial ^i _v \partial ^j _x h \Nt_{L^q_v L^p_x(m)} \vs \\
     \ds \Nt h \Nt_{\EE_{\nu,\var}} := \sum_{|i|\le \sigma, \ |j|\le
       s, \ |i|+|j| \le \max\{\sigma; s\}} \var^{|i|} \, \Nt \partial ^i
     _v \partial ^j h \Nt_{L^q_v L^p_x(m \nu^{1/q})}
\end{array}
\right.
$$
for some $\var >0$ to be chosen small enough. 
\end{proof}

%%%%%%%%%%%%%%%%%%%%%%%%%%%%%%%%%%%%%%

\subsubsection{The bilinear estimates}
\label{sec:bilinear-estimate}

Let us summarize the bilinear estimate available on the nonlinear term
in the equation \eqref{eq:NLBE}.

\begin{lem}
  \label{lem:bilinear}
  Consider the space $W^{\sigma,q}_v W^{s,p}_x (m)$ with $s,\sigma \in
  \N$, $\sigma \le s$, $s > 6/p$, $s \ge 0$ when $p=+\infty$, with $m$
  satisfying one of the assumptions {\bf (W1)}, {\bf (W2)}, {\bf (W3)} of
  Theorem~\ref{theo:LIBE2}. Then in the case $q<+\infty$ we have
    \begin{multline*}
      \left\| Q(g,f) \right\|_{W^{\sigma,q}_v W^{s,p}_x (m\nu^{1/q-1})} \\
      \le C \, \Big( \| g \|_{W^{\sigma,1}_v W^{s,p}_x (m)} \, \| f
      \|_{W^{\sigma,q}_v W^{s,p}_x (m \nu^{1/q})} + \| g
      \|_{W^{\sigma,q}_v W^{s,p}_x (m \nu^{1/q})}
      \, \| f \|_{W^{\sigma,1}_v W^{s,p}_x (m)} \\
      \| g \|_{W^{\sigma,1}_v W^{s,p}_x (m\nu)} \, \| f
      \|_{W^{\sigma,q}_v W^{s,p}_x (m)} + \| g \|_{W^{\sigma,q}_v
        W^{s,p}_x (m)} \, \| f \|_{W^{\sigma,1}_v W^{s,p}_x (m
        \nu)}\Big)
     \end{multline*}
    for some constant $C>0$, which implies 
\begin{multline*}
      \left\| Q(g,f) \right\|_{W^{\sigma,1} _v W^{s,p}_x(m) \cap W^{\sigma,q}_v W^{s,p}_x (m\nu^{1/q-1})} \\
      \le C \, \Big( \| g \|_{(W^{\sigma,1}_v \cap W^{\sigma,q}_v) W^{s,p}_x (m)} \, \| f
      \|_{W^{\sigma,1} _v W^{s,p}_x (m \nu) \cap W^{\sigma,q} _v W^{s,p}_x (m \nu^{1/q})}
      \\ + \| g \|_{W^{\sigma,1}_v W^{s,p}_x (m \nu) \cap W^{\sigma,q}_v W^{s,p}_x (m \nu^{1/q}) }
      \, \| f \|_{(W^{\sigma,1}_v \cap W^{\sigma,q}_v) W^{s,p}_x (m)}
      \Big)
     \end{multline*}
and in the case $q=+\infty$ we have simply
\begin{equation*}
      \left\| Q(g,f) \right\|_{W^{\sigma,\infty}_v W^{s,p}_x (m\nu^{-1})} 
      \le C \, % \Big(
      \| g \|_{W^{\sigma,\infty}_v W^{s,p}_x (m)} \, \| f
      \|_{W^{\sigma,\infty} _v W^{s,p}_x (m)}.
      % \\ + \| g \|_{W^{\sigma,\infty}_v W^{s,p}_x (m \nu) }
      % \, \| f \|_{W^{\sigma,\infty}_v W^{s,p}_x (m)}
      % \Big). 
     \end{equation*}
  %\end{itemize}
% On the gain term we have the more precise estimate 
% \begin{multline*}
% \left\| Q^+(g,f) \right\|_{W^{\sigma,q}_v W^{s,p}_x (m)} \\ 
% \le C_\delta \, \Big(
%   \| g \|_{W^{\sigma,1}_v W^{s,p}_x (m)} \, \| f \|_{W^{\sigma,q}_v
%     W^{s,p}_x (m)} + \| g \|_{W^{\sigma,q}_v W^{s,p}_x (m)} 
%   \, \| f \|_{W^{\sigma,1}_v W^{s,p}_x (m)} \Big)\\ 
% + \delta \Big( \| g \|_{W^{\sigma,1}_v W^{s,p}_x (m)} \, \| f \|_{W^{\sigma,q}_v
%     W^{s,p}_x (m \nu^{1/q})} + \| g \|_{W^{\sigma,q}_v W^{s,p}_x (m \nu^{1/q})} 
%   \, \| f \|_{W^{\sigma,1}_v W^{s,p}_x (m)} \\  
% + \| g \|_{W^{\sigma,1}_v W^{s,p}_x (m\nu)} \, \| f \|_{W^{\sigma,q}_v
%     W^{s,p}_x (m)} + \| g \|_{W^{\sigma,q}_v W^{s,p}_x (m)} 
%   \, \| f \|_{W^{\sigma,1}_v W^{s,p}_x (m \nu)}\Big)
% \end{multline*}
% for any $\delta>0$ and some corresponding constant $C_\delta>0$. 
% More precisely,  for any $\delta  > 0$, there exists  $C_\delta$ such that  for any 
% $f,g \in \EE_{\ell+1}$, $\ell \in [-1,0]$, we may split 
% \begin{eqnarray*}
% Q^{+,*}(g,f) &:=&Q(f,g) + Q(g,f) + f \, (g * 4\pi |\cdot |) 
% \\
% &=& Q_{S,\delta}^{+,*}(g,f) + Q_{R,\delta}^{+,*}(g,f)
% \end{eqnarray*}
% with 
% \bear\label{eq:QintegralS}
% \qquad \left\|
%   Q_{S,\delta}^{+,*}(g,f)  \right\|_{\EE} 
%   &\le& C_\delta \,   \left\| f \right\|_{\EE}  \, \left\| g  \right\|_{\EE}  
%   \\  \label{eq:QintegralR}
%   \qquad \left\|
%   Q_{R,\delta}^{+,*}(g,f)  \right\|_{\EE_\ell} 
%   &\le& \delta  \left\| f \right\|_{\EE_{\ell+1}}  \, \left\| g  \right\|_{\EE_{\ell+1}}.  
%   \eear
\end{lem}

\begin{proof}[Proof of Lemma~\ref{lem:bilinear}]
  For $\sigma=s=0$ and $q<\infty$ this estimate is an immediate
  consequence of the convolution inequalities on $Q$ established in
  \cite{MR946973}, together with the inequality $m (m' m'_*)^{-1} \le
  C \, m_*$.  (For the specific case of stretch exponential weight $m
  = e^{\kappa \, |v|^\beta}$, $\kappa >0$ and $\beta \in (0,2)$, we
  also refer to \cite{Mcmp} where the proof is explicitely
  written). In the case $q=+\infty$ we use Lemmas~\ref{lem:decQ+rad}
  and~\ref{lem:decQ+rad-exp}. 

  Finally the $x$ and $v$ derivatives are treated thanks to the
  distributive properties
\begin{equation*}
  \left\{ 
\begin{array}{l} \ds
  \nabla_x Q(g,f) = Q(\nabla_x g, f) + Q(g,\nabla_x f) \vs \\ \ds 
  \nabla_v Q(g,f) = Q(\nabla_v g, f) + Q(g,\nabla_v f)
\end{array}
\right.
\end{equation*}
and Sobolev embeddings. 
\end{proof}
 
\subsubsection{The a priori stability estimate}
\label{sec:key-priori-estimate-1}

\begin{lem}[A priori stability estimate]
\label{lem:apriori} 
Consider $s,\sigma \in \N$, $p,q \in [1,+\infty]$ with $\sigma \le s$,
$s > 6/p$, or $s \ge 0$ when $p=+\infty$, with $m$ satisfying one of
the assumptions {\bf (W1)}, {\bf (W2)}, {\bf (W3)} of
Theorem~\ref{theo:LIBE2}. Then consider the spaces
\[
\left\{ 
\begin{array}{l} \ds
\EE^q  := W^{\sigma,q}_v W_x^{s,p}
(m)\vs \\ \ds
\EE_\nu ^q := W^{\sigma,q}_v
W^{s,p}_x (m\nu^{1/q})
\end{array}
\right.
\]
if $q <+\infty$, or simply 
\[
\EE^\infty := \EE^\infty _\nu = W^{\sigma,\infty}_v W_x^{s,p}
(m)
% \left\{ 
% \begin{array}{l} \ds
% \EE := \EE_\nu = W^{\sigma,\infty}_v W_x^{s,p}
% (m)\vs \\ \ds
% \EE_\nu := W_v^{\sigma,\infty} W_x^{s,p}(m\nu) 
% \end{array}
% \right.
\]
if $q=+\infty$. Consider a solution
\[
f_t = \mu + h_t \in \EE
\]
to the nonlinear Boltzmann equation, with $\Pi h_t = 0$. 

Then for $q <+\infty$, $h_t$ satisfies the estimate
\begin{equation}\label{ineq:NLBEaprioriBd}
  \frac{{\rm d}}{{\rm d}t} \Nt h_t \Nt_{\EE^q} \le 
  \left( C \,  \Nt h_t \Nt_{\EE^q \cap \EE^1} - K \right) \, \Nt h_t \Nt_{\EE^q} ^{1-q}
  \, \Nt h_t \Nt_{\EE^q _\nu} ^q \le \left( C \,  \Nt h_t \Nt_{\EE^q
      \cap \EE^1} - K \right) \, \Nt h_t \Nt_{\EE^q} 
\end{equation}
for some constants $C, K>0$, which also writes 
\[
\frac{{\rm d}}{{\rm d}t} \left( \frac1q \, \Nt h_t \Nt_{\EE^q} ^q \right) \le 
  \left( C \,  \Nt h_t \Nt_{\EE^q \cap \EE^1} - K \right) \, \Nt h_t
  \Nt_{\EE^q _\nu} ^q{\color{red}{.}}
\]

When $q=+\infty$ we have the cleaner estimate
\begin{equation}\label{ineq:NLBEaprioriBd-infty}
  \frac{{\rm d}}{{\rm d}t} \Nt h_t \Nt_{\EE^\infty} \le 
  \left( C \,  \Nt h_t \Nt_{\EE^\infty} - K \right) \, \Nt h_t \Nt_{\EE^\infty}
\end{equation}
for some constants $C, K>0$. 
\end{lem}

\begin{proof}[Proof of Lemma~\ref{lem:apriori}] Assume first
  $\sigma=s=0$ and consider $q\in[1,+\infty)$ and $p \in [1,+\infty)$,
  and denote $\Phi(z) = |z|^p/p$. We calculate
\begin{eqnarray*}
\frac{{\rm d}}{{\rm d}t} \Nt h_t \Nt_{L^q_v L^p _x(m)}   = I_1 + I_2
\end{eqnarray*}
with 
\begin{multline*}
  I_1 := \eta \, \n{h_t}_{L^q_v L^p_x(m)}^{1-q} \, \int_{\R^3} \left(
    \int_{\T^3} \LL h_t \, \Phi'(h_t) \dd x \right) \, \| h_t
  \|_{L^p_x} ^{q-p} \, m^q \dd v
  \\
  + \n{h_t}_{L^q_v L^p_x(m)}^{1-q} \, \int_0^{+\infty} \int_{\R^3}
  \left( \int_{\T^3} ( \SS_\LL(\tau) \,(\LL h_t ) \, \Phi'(e^{\tau\LL}
    \, h_t) \dd x \right) \, \n{\SS_\LL(\tau) \, h_t}_{L^p_x} ^{q-p}
  \, m^q \dd v \dd \tau
\end{multline*}
 and
\begin{multline*}
  I_2 := \eta \, \n{h_t}_{L^q_v L^p_x(m)}^{1-q} \, \int_{\R^3} \left(
    \int_{\T^3} Q(h_t,h_t) \, \Phi'(h_t) \dd x \right) \, \| h_t
  \|_{L^p_x} ^{q-p} \, m^q \dd v
  \\
  + \n{h_t}_{L^q_v L^p_x(m)}^{1-q} \,\int_0^{+\infty} \int_{\R^3}
  \left( \int_{\T^3} (\SS_\LL(\tau) \, Q(h_t,h_t)) \, \Phi'(e^{\tau\LL}
    \, h_t) \dd x \right) \, \n{\SS_\LL(\tau) \, h_t}_{L^p_x} ^{q-p}
  \, m^q \dd v \dd\tau.
\end{multline*}
In Proposition ~\ref{theo:LIBE3} we proved that choosing $\eta > 0$ it
holds
$$
I_1 \le - K \, \n{h_t}_{L^q_v L^p_x(m)}^{1-q} \, \Nt h_t
\Nt_{L^q_vL^p_x(m \nu^{1/q})} ^q \quad \mbox{ for some
} K>0. 
$$
For the second term, the H\"older inequality implies 
\begin{multline*}
  \int_{\R^3} \left( \int_{\T^3} Q(h_t,h_t) \, \Phi'(h_t)
    \dd x \right) \, \| h_t \|_{L^p_x} ^{q-p} \, m^q \dd v \\ \le 
  \int_{\R^3} \n{Q(h_t,h_t)}_{L^p_x}  \n{h_t}_{L^p_x} ^{q-1} \, m^q \dd v
  \\
  \le \n{Q(h_t,h_t)}_{L^q_v L^p_x(m \nu^{1/q-1})}  \n{h_t}_{L^q _v
    L^p_x(m\nu^{1/q})} ^{q-1}
\end{multline*}
and similarly 
\begin{multline*}
  \int_{\R^3} \left( \int_{\T^3}  (e^{\tau\LL} \,
    Q(h_t,h_t)) \, \Phi'(e^{\tau\LL} \, h_t) \dd x \right) \, \|
  e^{\tau\LL} \, h_t \|_{L^p_x} ^{q-p} \, m^q \dd v \\ \le 
  \int_{\R^3} \n{e^{\tau\LL} \, Q(h_t,h_t)}_{L^p_x}  \n{e^{\tau\LL} \, h_t}_{L^p_x} ^{q-1} \, m^q \dd v
  \\
  \le \n{e^{\tau\LL} \, Q(h_t,h_t)}_{L^q_v L^p_x(m \nu^{1/q-1})}
  \n{e^{\tau\LL} \, h_t}_{L^q_v L^p_x(m\nu^{1/q})} ^{q-1}.
\end{multline*}

Using the bilinear estimate in Lemma~\ref{lem:bilinear} and the
semigroup decay in Theorem~\ref{theo:LIBE2} (noticing that $\Pi
Q(h_t,h_t) = 0$) we get the following estimates
\begin{equation*}
  \n{Q(h_t,h_t)}_{L^q_vL^p_x(m \nu^{1/q-1})} \le C \, \|h_t \|_{(L^1_v \cap L^q_v)L^\infty_x(m)} \, \|h_t
  \|_{L^1_v L^p_x(m\nu)\cap L^q_v L^p_x(m\nu^{1/q})}
\end{equation*}
and 
\begin{multline*}
  % \Nt Q(h_t,h_t) \Nt_{L^q_vL^p_x(m \nu^{1/q-1})} \\ \le \eta \, \left\| Q(h_t,h_t)
  % \right\|_{L^q_vL^p_x(m)} +
  \int_0 ^{+\infty} \left\| S_\LL(\tau)
    Q(h_t,h_t) \right\|_{L^q_vL^p_x(m \nu^{1/q-1})} \dd \tau \\
  \le % C \, \eta \, \|h_t \|_{(L^1_v \cap L^q_v)L^\infty_x(m)} \, \|h_t
  % \|_{L^1_v L^p_x(m\nu)\cap L^q_v L^p_x(m\nu^{1/q})} \\
  C' \, \left( \int_0 ^{+\infty} e^{-\lambda \, \tau} \dd \tau
  \right) \, \|h_t \|_{(L^1_v \cap L^q_v) L^\infty _x(m)} \, \|h_t \|_{L^1_v
    L^p_x(m\nu)\cap L^q_v L^p_x(m\nu^{1/q})}
  \\
  \le % C' \, \|h_t \|_{(L^1_v \cap L^q_v)L^\infty_x(m)} \, \|h_t
  % \|_{L^1_v L^p_x(m\nu)\cap L^q_v L^p_x(m\nu^{1/q})} \\
  % \le 
C'' \, \Nt h_t \Nt_{(L^1_v \cap L^q_v)L^\infty_x(m)} \, \Nt h_t
  \Nt_{L^1_v L^p_x(m\nu)\cap L^q_v L^p_x(m\nu^{1/q})}
\end{multline*}
for some constant $C,C', C''>0$. We deduce that 
\begin{equation*}
  I_2  \le C''' \, \n{h_t}_{L^q_v L^p_x(m)}^{1-q} \, \Nt h_t \Nt_{(L^1_v \cap L^q_v)L^\infty_x(m)} \, \Nt h_t
  \Nt_{L^1_v L^p_x(m\nu)\cap L^q_v L^p_x(m\nu^{1/q})} ^q
  % \|Q(h_t,h_t) \|_{L^q_v L^p_x (m)} 
  % + \int_0^{+\infty}  \|  e^{\tau\LL} \,  Q(h_t,h_t) \|_{L^q_vL^p_x(m)}\dd\tau 
\end{equation*}
and thus (using Sobolev embeddings or passing to the limit $p \to \infty$)
\begin{eqnarray*}
  \frac{{\rm d}}{{\rm d}t} \Nt h_t \Nt_{\EE}  \le  \left( C \, \Nt h_t \Nt_{\EE} - K
\right) \, \n{h_t}_{\EE}^{1-q} \, \Nt h_t \Nt_{\EE_\nu} ^q. 
\end{eqnarray*}

This concludes the proof in the case $\sigma=s=0$, $q<+\infty$ and
$p=+\infty$. In the case $p <+\infty$ and $0 < \sigma \le s$, one uses
the distributive property of the derivatives and Sobolev embeddings.

The case $q=+\infty$ is handled similarly by using the final estimates
in Lemma~\ref{lem:bilinear}. We use the previous argument with
$q<+\infty$ unchanged and take the limit $q \to \infty$ in the
bilinear estimates to get
\begin{equation*}
\frac{{\rm d}}{{\rm d}t} \Nt h_t \Nt_{L^\infty_v L^p _x(m)} \le - K \,
\Nt h_t \Nt_{L^\infty_v L^p _x(m)}  + \Nt Q(h_t,h_t) \Nt_{L^\infty_v L^p_x(m\nu^{-1})}.
\end{equation*}
The bilinear estimate in Lemma~\ref{lem:bilinear}
for $q=+\infty$ and the semigroup decay in Theorem~\ref{theo:LIBE2}
(noticing that $\Pi Q(h_t,h_t) = 0$) yield 
\begin{equation*}
  \n{Q(h_t,h_t)}_{L^\infty_vL^p_x(m \nu^{-1})} \le C \, \|h_t \|_{L^\infty_vL^\infty_x(m)} \, \|h_t
  \|_{L^\infty_v L^p_x(m)}
\end{equation*}
and 
\begin{multline*}
  \int_0 ^{+\infty} \left\| S_\LL(\tau)
    Q(h_t,h_t) \right\|_{L^\infty_vL^p_x(m \nu^{-1})} \dd \tau \\
  \le   C' \, \left( \int_0 ^{+\infty} e^{-\lambda \, \tau} \dd \tau
  \right) \, \|h_t \|_{L^\infty_v L^\infty _x(m)} \, \|h_t \|_{L^\infty_v L^p_x(m)}
  \\
  \le C'' \, \|h_t \|_{L^\infty_v L^\infty _x(m)} \, \|h_t \|_{L^\infty_v L^p_x(m)}
\end{multline*}
for some constant $C,C', C''>0$. We deduce that 
\begin{equation*}
  \Nt Q(h_t,h_t) \Nt_{L^\infty_v L^p_x(m\nu^{-1})} \le C''' \, \|h_t
  \|_{L^\infty_v L^p_x(m)} ^2
\end{equation*}
and thus (using Sobolev embeddings or passing to the limit $p \to \infty$)
\begin{eqnarray*}
  \frac{{\rm d}}{{\rm d}t} \Nt h_t \Nt_{\EE^\infty}  \le  \left( C \, \Nt h_t \Nt_{\EE^\infty} - K
\right) \, \n{h_t}_{\EE^\infty}. 
\end{eqnarray*} 
\end{proof}

\subsubsection{Final proof}
\label{sec:proof-theor-refth}

We consider the close-to-equilibrium regime and the spaces $\EE$ and
$\EE_\nu$ as before.  We will construct solutions through the
following iterative scheme
$$
\partial_t h^{n+1} = \LL \, h^{n+1} + Q(h^{n+1},h^n), \quad n \ge 1, 
$$
with the initialization 
\[
\partial_t h ^0= \LL \, h^0, \qquad h^0 _0 = \init{h}^0 =\init{h}, \quad \Nt
\init{h} \Nt_{\EE^q} \le \eps/2.
\]
The functions $h^n$, $n \ge 0$ are well-defined in $\EE$ for all times
$t \ge 0$ thanks to the study of the semigroup in
Theorem~\ref{theo:LIBE2} and the stability estimates proven below.

\smallskip We split the proof into four steps. The first two steps of
the proof establish the stability and convergence of the iterative
scheme, and they are mainly an elaboration upon the key a priori
estimate of the previous subsection. The third step consists of a
bootstrap argument in order to recover the optimal decay rate of the
linearized semigroup. The fourth step details the modifications to the
argument for $q=+\infty$.

% \medskip
% \noindent
% {\em Step~1. Local-in-time solutions.} 
% We consider the equation 
% \begin{equation*}
%   \partial_t f_t + v \cdot \nabla_x f = Q(f,f). 
% \end{equation*}

\medskip

\noindent
{\em Step~1. Stability of the scheme.} 
Let us first assume $q<+\infty$ and prove by induction the following control 
\begin{equation}\label{eq:induction}
  \forall \, n \ge 0, \ \fa t \ge 0, \quad B_n(t) := 
  \left( \frac1q \, \Nt h^n _t \Nt_{\EE^q} ^q + K \,  \int_0^t
    \Nt h^n _\tau \Nt_{\EE^q _\nu} ^q \dd \tau  \right) \le \eps^q
\end{equation}
under a smallness condition on $\var$. 

The case $n=0$ follows from Theorem~\ref{theo:LIBE3} and the
fact that $\Nt \init{h} \Nt_{\EE^q} ^q \le (\eps/2)^q$: 
$$
\sup_{t \ge 0} \left( \Nt h^0_t \Nt_{\EE^q} ^q + K\, \int_0^t \Nt h^0_\tau
\Nt_{\EE^q _\nu} ^q \dd \tau \right) \le \eps^q.
$$

Let us now assume that \eqref{eq:induction} is satisfied at rank $n$
and let us prove it for $n+1$. A similar computation as in
Lemma~\ref{lem:apriori} yields
\begin{multline*}
\frac{{\rm d}}{{\rm d}t} \left( \frac1q \, \Nt h^{n+1} _t \Nt_{\EE^q} ^q
\right) + K \, \Nt h^{n+1} _t \Nt_{\EE^q _\nu} ^q \le \\ C \, \left( \Nt h^n
  _t\Nt_{\EE^q} \, \Nt h^{n+1} _t \Nt_{\EE^q _\nu} + \Nt h^{n+1} _t\Nt_{\EE^q} \, \Nt
  h^{n} _t \Nt_{\EE^q _\nu}\right) \, \Nt h^{n+1} _t \Nt_{\EE^q _\nu}
^{q-1}
\end{multline*}
for some constants $C,K>0$. Hence by H\"older's inequality we get
\begin{multline*}
  B_{n+1}(t) = \frac1q \, \Nt h^{n+1}_t \Nt_{\EE^q} ^q + K \, \int_0^t \Nt
  h^{n+1}_\tau
  \Nt_{\EE^q _\nu} ^q \dd \tau \\
  \le \frac1q \, \Nt \init{h} \Nt_{\EE^q} ^q + C \, \left( \sup_{\tau \ge
      0} \Nt h^n_\tau \Nt_{\EE^q} \right) \, \left( \int_0 ^t
    \Nt h^{n+1} _\tau \Nt_{\EE^q _\nu} ^q \dd \tau \right) \\ + 
  C \, \left( \int_0 ^t \Nt h^n_\tau \Nt_{\EE^q _\nu}^q \dd \tau \right)^{1/q}
  \, \left( \sup_{\tau \ge 0} \Nt h^{n+1} _\tau \Nt_{\EE^q _\nu} \dd \tau
  \right) \, \left( \int_0 ^t \Nt h^{n+1} _\tau \Nt_{\EE^q _\nu} ^q \dd \tau \right)^{1-1/q} 
  \\
  \le \frac1q \, \Nt \init{h} \Nt_{\EE^q} ^q + \left( \min \left\{ C,
      \frac{C}{K^{1/q}} \right\} \right) \, B_n ^{1/q} \, B_{n+1}(t) \\
  \le \frac1q \, \Nt \init{h} \Nt_{\EE^q} ^q + \left( \min \left\{ C,
      \frac{C}{K^{1/q}} \right\} \right) \, \var \, B_{n+1}(t)
% \\
% \left( \frac{q
%       \var^2}{K} \right)^{1/q} \, \left( \int_0^t \Nt h^{n+1}_\tau
%     \Nt_{\EE_\nu} ^q
%     \dd \tau \right)^{1-1/q} \\
%   \le \frac{\eps}2 + {C \over K} \, q^{1/q} \, \eps^{2/q} \, B_{n+1}
%   ^{1-1/q} \le \eps,
\end{multline*}
from which it follows 
\begin{equation*}
  \fa t \ge 0, \quad  B_{n+1}(t) \le \frac{2}{q} \, \Nt \init{h} \Nt_{\EE^q} ^q \le \var^q
\end{equation*}
as soon as
\[
\left( \min \left\{ C,
      \frac{C}{K^{1/q}} \right\} \right) \, \var \le \frac12. 
\]
The induction is proven. 

Passing to the limit $q \to +\infty$ it holds
\begin{equation*}
  \sup_{t \ge 0} \Nt h_t \Nt_{\EE^\infty} \le \var
\end{equation*}
assuming that the initial data satisfies $\Nt \init{h} \Nt_{\EE^\infty} \le
\var/2$. Observe that the smallness condition on $\var$ is uniform as
$q \to +\infty$, which is crucial in this limiting process. 

 \medskip

\noindent 
{\em Step~2. Convergence of the scheme.}
Let us now denote by $d^n := h^{n+1} - h^n$. % and by $s^n := h^{n+1} +
% h^n$ for $n \ge 0$ They satisfy % {\color{red}{MG: $s_n$ is never
    % used. Do we need to define it ?}} 
It satisfies
$$
\forall \, n \ge 0, \quad \partial_t d^{n+1} = \LL d ^{n+1} +
Q(d^{n+1},h^{n+1}) + Q(h^{n+1}, d^n)
$$
and 
\[
\partial_t d^{0} = \LL d ^{0} + Q(h^1,h^0).
\]

Let us denote by
\[
A^n (t) := \left( \frac1q \, \Nt d^{n} _t \Nt_{\EE^q} ^q + K \, \int_0
  ^{t} \Nt d^{n} _\tau\Nt_{\EE^q _\nu} ^q \dd \tau \right)
\]
and let us prove by induction that 
\[
\forall \, t \ge 0, \ \forall \, n \ge 0, \quad A^n(t) \le ( \bar C \, \eps)^{qn}
\]
for some constant $\bar C>0$ uniform as $\eps$ goes to zero and as $q$
goes to infinity. 

The case $n=0$ is obtained by using the estimate
\begin{multline*}
\frac{{\rm d}}{{\rm d}t} \left( \frac1q \, \Nt d^0 _t \Nt_{\EE^q} ^q
\right)  + K \, \Nt d^0 _t
\Nt_{\EE^q _\nu} ^q \\ \le C \, \left( \Nt h^1 _t \Nt_{\EE^q}  \, \Nt h^0
  _t\Nt_{\EE^q _\nu} +  \Nt h^1 _t \Nt_{\EE^q _\nu}  \, \Nt h^0
  _t\Nt_{\EE} \right) \Nt h^0 _t \Nt_{\EE^q _\nu} ^{q-1}
\end{multline*}
and the previous bounds on $h^0, h^1$ to deduce
\[
\fa t \ge 0, \quad A^0 (t) = \frac1q \, \Nt d^0 _t \Nt_{\EE^q} ^q + K \,
\int_0 ^t \Nt d^0 _\tau \Nt_{\EE^q _\nu} ^q \dd \tau \le C \,
\eps^{2} \le \var
\]
for $\var$ small enough. 

% Then we denote 
% \[
% \EE_{\nu^{1/q-1}} := W_v^{\sigma,1} W_x^{s,p}(m) \cap W^{\sigma,q}_v
% W^{s,p}_x (m\nu^{1/q-1})
% \]
% if $q <+\infty$, or simply 
% \[
% \EE_{\nu^{1/q-1}} := \EE_{\nu^{-1}} = W^{\sigma,\infty}_v W_x^{s,p}
% (m\nu^{-1})
% \]
% if $q = +\infty$. 

The propagation of the induction is obtained by estimating (similarly
as before) 
\begin{multline*}
  A^{n+1}(t) = \frac1q \, \Nt d^{n+1} _t \Nt_{\EE^q} ^q + K \, \int_0 ^t
  \, \Nt d^{n+1} _\tau \Nt_{\EE^q _\nu} ^q \dd \tau
  \\
  \le \frac1q \, \Nt \init{d}^n \Nt_{\EE^q} ^q + C \,  \int_0 ^t \Big( \Nt
  d^{n+1}_\tau \Nt_{\EE^q _\nu} \, \Nt h^{n+1} _\tau \Nt_{\EE^q} + \Nt
  d^{n+1}_\tau \Nt_{\EE} \, \Nt h^{n+1} _\tau \Nt_{\EE^q _\nu} \Big) 
  \, \Nt d^{n+1} _\tau \Nt_{\EE^q _\nu} ^{q-1} \dd \tau \\ 
+ C \, \int_0 ^t \Big( \Nt
  d^{n}_\tau \Nt_{\EE^q _\nu} \, \Nt h^{n+1} _\tau \Nt_{\EE^q} + \Nt
  d^{n}_\tau \Nt_{\EE^q} \, \Nt h^{n+1} _\tau \Nt_{\EE^q _\nu} \Big) 
  \, \Nt d^{n+1} _\tau \Nt_{\EE^q _\nu} ^{q-1} \dd \tau
  \\
  \le 2 \, C \, \var \, A^{n+1}(t) + 2 \, C \, \var \, A_n
  (t) ^{1/q} \, A_{n+1}(t)^{1-1/q}
  % \le C \, \Bigg[ \sup_{\tau \ge 0} \Nt s^n
  % _\tau \Nt_{\EE^q} \, \left( \int_0 ^{t} \Nt d^n _\tau \Nt_{\EE^q _\nu} ^q \dd \tau \right)  +
  % \sup_{\tau \ge 0} \Nt d^n _\tau \Nt_{\EE} \, \left( \int_0 ^{t} \Nt
  %   s^n _\tau \Nt_{\EE^q _\nu} \Nt d^n _\tau \Nt_{\EE^q _\nu} ^{q-1} \dd \tau \right) \Bigg]
  % \\
  % \le C' \, \eps^{1/q} \, \sup_{0 \le \tau \le t} \left( \frac1q \, \Nt d^n _\tau
  %   \Nt_{\EE^q} ^q + K \, \int_0 ^{t} \Nt d^n _\tau\Nt_{\EE^q _\nu} ^q \dd \tau
  % \right) = C' \, \eps^{1/q} \, A^n(t)
\end{multline*}
where we have used $\init{d}^n \equiv 0$ for any $n \ge 0$. Using the
induction assumption on $A_n(t)$ we deduce that 
\begin{equation*}
  A_{n+1}(t) \le 2 \, C \, \var \, A^{n+1}(t) + 2 \, C \,
  \var \, \bar C^{n} \, \var^{n} \, A_{n+1}(t)^{1-1/q}
\end{equation*}
and if $\var$ is small enough so that $2 \, C \, \var \le 1/2$
we get
\begin{equation*}
  A_{n+1}(t) \le 4 \, C \, \bar C^{n} \, \var^{n+1} \, 
  A_{n+1}(t)^{1-1/q} \ \Longrightarrow \ A_{n+1}(t) \le  (4 \,
  C)^q \, \bar C^{qn} \, \var^{q(n+1)}
\end{equation*}
which concludes the proof with $\bar C = 4 \, C$. 

Hence for $\eps$ small enough, the series $\sum_{n \ge 0} A^n(t)$ is
summable for any $t \ge 0$, and the sequence $h^n$ has the Cauchy
property in $L^\infty_t(\EE)$, which proves the convergence of the
iterative scheme. The limit $h$ as $n$ goes to infinity satisfies the
equation in the strong sense when the norm $\EE$ involves enough
derivatives, or else in the mild sense. 

Finally observe again that the smallness condition on $\var$ is
uniform as $q \to +\infty$, and by passing to the limit one gets by
induction
\begin{equation*}
  \sup_{t \ge 0} \Nt d_t ^n \Nt_{\EE^\infty}  \le (\bar C \, \var)^n
\end{equation*}
which shows again that the sequence $h^n$ is Cauchy in
$L^\infty_t(\EE^\infty)$. This proves the convergence of the iterative
scheme.

\medskip

\noindent {\em Step~3. Rate of decay.} We now consider the solution
$h$ constructed so far, first in the case $q<+\infty$. From Step 1 we
take the limit $n\to \infty$ in the stability estimate and get
\[
\sup_{t \ge 0} \left( \frac1q \, \Nt h_t \Nt_{\EE^q} ^q  + K \,  \int_0^t
\Nt h_\tau \Nt_{\EE^q _\nu}^q \dd \tau \right)  \le \eps^q.
\]
We can then apply Lemma~\ref{lem:apriori} to the solution $h_t$: 
\begin{multline*}
\frac{{\rm d}}{{\rm d}t} \left( \frac1q \, \Nt h_t \Nt_{\EE^q} ^q \right) \le \left( C \, \Nt h_t \Nt_{\EE^q} - K
\right) \, \Nt h_t \Nt_{\EE^q _\nu} ^q  \\
   \le \left( C \, q^{1/q} \, \eps - K \right) \, \Nt
h_t \Nt_{\EE^q _\nu} ^q \le \frac{\left( C \, q^{1/q} \, \eps - K \right)}{\nu_0 ^q} \, \Nt
h_t \Nt_{\EE^q} ^q,
\end{multline*}
where we have used the previous stability bound. This implies that 
\[
\Nt h_t  \Nt_{\EE^q}  \le e^{-\frac{K}{2 \nu_0^q} \, t} \, \Nt \init{h}
\Nt_{\EE^q} 
\]
under the smallness condition $C \, q^{1/q} \, \var - K \le -K/2$ on
$\var$. Moreover since $\Nt h_t \Nt_{\EE^q}$ converges to zero as $t
\to +\infty$, we integrate the previous a priori estimate from $t$ to
$+\infty$ to get
\[
\frac{K}{2} \, \int_t ^{+\infty} \Nt h_\tau \Nt_{\EE^q _\nu} ^q \dd \tau
\le \frac1q \, \Nt h_t \Nt_{\EE^q} ^q \le \frac{e^{-\frac{q \,
      K}{2\nu_0^q} \, t}}{q} \, \Nt \init{h} \Nt_{\EE^q} ^q,
\]
which implies 
\begin{equation}\label{eq:decay-int}
\int_t ^{+\infty} \left\| h_\tau \right\|_{\EE^q _\nu} ^q \dd \tau \le
\frac{2}{K \eta} \, \left\| h_t \right\|_{\EE^q} ^q \le C \, e^{-\frac{q
    \, K}{2\nu_0^q} \, t} \, \left\| \init{h}
\right\|_{\EE^q} ^q.
\end{equation}

We shall now perform a bootstrap argument in order to ensure that the
solution $h_t$ enjoys to same optimal decay rate $O(e^{-\lambda \,
  t})$ as the linearized semigroup in Theorem~\ref{theo:LIBE2}.
Assume that the solution is known to decay as
\begin{equation}\label{eq:decay-0}
\left\| h_t \right\|_{\EE^q} \le C \, e^{-\lambda_0 \, t}
\end{equation}
for some constant $C>0$, and let us prove that it indeed decays like
\[
\left\| h_t \right\|_{\EE^q} \le C' \, e^{-\lambda_1 \, t}
\]
with  $\lambda_1  =  \min  \{ \lambda_0  +  K/(4\nu_0^q),  \lambda\}$,
possibly  for some  other larger  constant $C'>0$.  Hence in  a finite
number of  steps, it proves  the desired decay rate  $O(e^{-\lambda \,
  t})$.

Assume \eqref{eq:decay-0} and write a Duhamel formulation:
\[
h_t = S_\LL(t) \init{h} + \int_0 ^t S_\LL (t-\tau) Q(h_\tau,h_\tau) \,
\dd \tau. 
\]
We go back to the original norm and we deduce from
Theorem~\ref{theo:LIBE2} and Lemma~\ref{lem:bilinear}
\[
\left\| h_t \right\|_{\EE^q} \le C \, e^{-\lambda \, t} \left\| \init{h}
\right\|_{\EE^q} + C \, \int_0 ^t e^{-\lambda \, (t-\tau)} \, \left\| h_\tau
\right\|_{\EE^q} \, \left\| h_\tau \right\|_{\EE^q _\nu} \dd \tau. 
\]

Assume $\lambda_0 < \lambda$ and denote $\lambda_1 = \min \{ \lambda_0
+ K/(4\nu_0^q), \lambda\}$. We simply estimate
\begin{multline*}
\int_0 ^t e^{-\lambda \, (t-\tau)} \, \left\| h_\tau \right\|_{\EE^q} \,
\left\| h_\tau \right\|_{\EE^q _\nu} \dd \tau 
  \\ 
\le \int_0 ^t e^{-\lambda_1 \, (t-\tau)} \, \left\| h_\tau \right\|_{\EE^q} \,
\left\| h_\tau \right\|_{\EE^q _\nu} \dd \tau 
  \\
\le C \, e^{-\lambda_1 \, t} \,
\left( \int_0 ^t e^{(\lambda_1-\lambda_0) \, \tau} \, \left\| h_\tau
\right\|_{\EE^q _\nu} \dd \tau \right) \, \left\| \init{h} \right\|_{\EE^q}
\end{multline*}
and then by integration by parts
\begin{multline*}
  \int_0 ^t e^{(\lambda_1-\lambda_0) \, \tau} \, \left\| h_\tau
  \right\|_{\EE^q _\nu} \dd \tau 
  \\
   \le \int_0 ^t \left\| h_\tau
  \right\|_{\EE^q _\nu} \dd \tau + (\lambda_1-\lambda_0) \, \int_0 ^t
  e^{(\lambda-\lambda_0) \, \tau} \, \left( \int_\tau ^t \left\| h_{\tau'}
    \right\|_{\EE^q _\nu} \dd \tau' \right) \dd \tau 
  \\
  \le C \, \left\| \init{h} \right\|_{\EE^q} + (\lambda_1-\lambda_0) \,
  \left( \int_0 ^t (t-\tau)^{1-1/q} \, e^{(\lambda_1-\lambda_0-K/(2\nu_0^q)) \, \tau} \dd \tau \right)
  \, \left\| \init{h} \right\|_{\EE^q} \\ \le C \, \left\| \init{h} \right\|_{\EE^q}
\end{multline*}
where in the last line we have used \eqref{eq:decay-int}. All in all
we deduce 
\[
\left\| h_t \right\|_{\EE^q} \le C \, e^{-\lambda_1 \, t} \left\| \init{h}
\right\|_{\EE^q}.
\]
This proves the claim and concludes the proof of the estimate 
\begin{equation*}
  \left\| h_t \right\|_{\EE^q} \le C \, e^{-\lambda \, t} \, \left\| \init{h} \right\|_{\EE^q} 
\end{equation*}
in the case $q<+\infty$, where $\lambda  = \lambda(q) >0$ is the sharp rate of the
linearized semigroup in Theorem~\ref{theo:LIBE2}, and the constant $C$ is
uniform as $q \to +\infty$.
\medskip

\noindent {\em Step~4. The case $q = +\infty$.} It is obtained by
passing to the limit in the previous estimate and using that
$\lambda(q) \to \lambda(\infty)>0$ under our assumptions, thanks to
Theorem~\ref{theo:LIBE2}.  One gets
\begin{equation*}
  \left\| h_t \right\|_{\EE^\infty} \le C \, e^{-\lambda \, t} \, \left\| \init{h} \right\|_{\EE^\infty} 
\end{equation*}
with again the sharp rate $\lambda>0$ of the linearized semigroup.

\subsection{Proof of Theorem~\ref{theo:NLBE+}}
\label{sec:proof-theor-refth-1}

We now consider the weakly inhomogeneous solutions. We split the proof
into three steps. 
\medskip

\noindent 
{\em Step~1. The spatially homogeneous evolution.} 
Consider the spatially homogeneous initial datum $\init{g} \in L^1_v
L^\infty_x(1+|v|^k)$, $k \ge 2$. From \cite[Theorem~1.2]{Mcmp} we know
that it gives rise to a unique conservative spatially homogeneous
solution $g_t \in L^1_v(1+|v|^2)$ which satisfies
\[
\left\| g_t - \mu \right\|_{L^1_v (1+|v|^k)} \le C \, e^{-\lambda t}
\]
with explicit and optimal exponential rate.
\medskip

\noindent {\em Step~2. Local-in-time stability.} 
We consider the estimate in $L^1_v L^\infty_x (1+|v|^k)$, $k
>2$.  We want to construct a solution $f_t$ that is $L^1_v L^\infty_x
(1+|v|^k)$-close to the spatially homogeneous solution $g_t$
previously considered on a finite time interval.  

Arguing as before we have the a priori bound
\begin{equation*}
  \fa t \ge 0, \quad f_t(x,v) \ge \init{f}(v) \, e^{-C \, (1+|v|) \,
    t}
\end{equation*}
where $C>0$ depends on the $L^\infty_{t,loc} L^1_v L^\infty_x (1+|v|)$
norm of the solution. 

Since $\init{f}$ is close to a non-zero spatially homogeneous solution
$\init{g}(v)$, choosing if necessary $\epsilon(M_k)$ small enough we
have
\begin{multline*}
 \fa t \ge 0, \quad  \int_{\R^3} f_t(x,v_*) \, |v-v_*| \dd v_* \ge
 e^{-C' \, t} \int_{|v_*|\le R} \init{f}(x,v_*) \, |v-v_*| \dd v_*  \\
 \ge e^{-C' \, t} \left( \int_{\R^3} \init{g}(v_*) \, |v-v_*| \dd v_*
   - \epsilon \, (1+|v|) \right) \\ \ge e^{-C' \, t} \left(
   K_{\init{g}} - \epsilon \right) \, (1+|v|) \ge K  \, e^{-C' \, t} \,
   (1+|v|)
\end{multline*}
for some constants $C', K>0$. We have used 
\begin{equation*}
  \int_{\R^3} \init{g}(v_*) \, |v-v_*| \dd v_* \ge K_{\init{g}} \, (1+|v|)
\end{equation*}
which follows from the inequalities
% the following
% facts on $g_t$:
% \begin{equation*}
%   \fa t \ge 0, \quad g_t(v) \ge \init{g}(v) \, e^{-C \, t} 
% \end{equation*}
% by the Duhamel principle and the energy bound, and 
% \begin{equation*}
%  \fa t \ge 0, \quad  \int_{\R^3} g_t(v_*) \, |v-v_*| \dd v_* \ge
%  K_{\init{g}} \, e^{-C\,t} \, (1+|v|) 
% \end{equation*}
% for some constants, where we have used
\begin{equation*}
  \int_{\R^3} \init{g}(v_*) \, |v-v_*| \dd v_* \ge |v| 
\end{equation*}
(by convexity) and 
\begin{equation*}
  \int_{\R^3} \init{g}(v_*) \, |v_*| \dd v_* >0. 
\end{equation*}

\begin{rem}
  The constant $K_{\init{g}}$ depends in general on the mass, energy
  and on the shape of $\init{g}$, more precisely on how it
  concentrates at zero velocity (recall that the momentum is
  normalized to zero). This is illustrated by the following
  counter-example
  \begin{equation*}
    g_n(v) := \left( 1- \frac{1}{n^2} \right) \,\varphi_0 + \left( \frac{1}{n^2} \right) \, \frac{\varphi_{-n} +
      \varphi_n}{2}
  \end{equation*}
  where $\varphi_0$ approximates $\delta_0$ and
    $\varphi_{\pm n}$ approximates $\delta_{\pm n}$ as $n \to 0$,
    which satisfies as $n \to \infty$
\begin{equation*}
\int_{\R^3} g_n \dd v = 1, \quad 
\int_{\R^3} g_n \, |v|^2 \dd v \sim 1, \quad 
\int_{\R^3} g_n \, |v| \dd v \sim 0.   
\end{equation*}

However, when a moment $k>2$ is assumed on $\init{g}$, it is easy to
give a bound on $K_{\init{g}}$ based on the higher moments estimates
since
\begin{equation*}
  \int_{\R^3} \init{g}(v_*) \, |v_*| \dd v_* \ge \frac{\left( \int_{\R^3}
    \init{g}(v_*) \, |v_*|^2 \dd v_* \right)^{(k-1)/(k-2)}}{\left( \int_{\R^3}
    \init{g}(v_*) \, |v_*|^k \dd v_* \right)^{1/(k-2)}}.
\end{equation*}
\end{rem}

% We first construct the solution as in the previous step on some
% $[0,T_0]$ where $T_0$ is small enough so that 
% \begin{equation*}
%   \fa t \in [0,T_0], \quad \left\{ 
% \begin{array}{l}\ds
%   \left| \n{f_t}_{L^1_v L^\infty_x (1+|v|^k)} -
%     \n{\init{f}}_{L^1_v L^\infty_x (1+|v|^k)}\right| \le \frac{\epsilon}{2},
%   \vs \\ \ds
%   \left| \n{g_t}_{L^1_v(1+|v|^k)} -
%     \n{\init{g}}_{L^1_v(1+|v|^k)}\right| \le \frac{\epsilon}{2}.
% \end{array}
% \right.
% \end{equation*}
% Observe that this $T_0$ depends only on the norms and moments of
% $\init{g}$ and $\init{d}$. 
We then consider $k' \in (2,k)$ and we define the difference $d_t:=
f_t - g_t$ and the sum $s_t := f_t+g_t$. % The difference satisfies
% \begin{equation*}
%   \n{d_t}_{L^1_v L^\infty_x (1+|v|^k)} \le
%   \n{\init{d}}_{L^1_v L^\infty_x (1+|v|^k)} + \epsilon \le 2 \,
%   \epsilon. 
% \end{equation*}
We then write the evolution equation
\[
\partial_t d_t  + v\cdot \nabla_x d_t = 2 \, Q(d_t,d_t) + 2 \, Q(g_t,d_t){\color{red}{,}} %= \mathcal P(d_t)
\]
from which we deduce the following a priori estimate arguing as in the
previous section
\begin{multline*}
  \dt \n{d_t}_{L^1_v L^\infty_x(1+|v|^{k'})} \le C_1 \, \n{d_t}_{L^1_v
    L^\infty_x(1+|v|^2)} \, \n{d_t}_{L^1_v L^\infty_x(1+|v|^{k'+1})} \\ + C_2
  \, \n{g_t}_{L^1_v L^\infty_x(1+|v|^{k'+1})} \,  \n{d_t}_{L^1_v
    L^\infty_x(1+|v|^{k'})} - K \, e^{-C' \,t} \, \n{d_t}_{L^1_v L^\infty_x(1+|v|^{k'+1})}
\end{multline*}
for some constants $C_1, C_2>0$. Observe however that here we have
to keep track of the time-dependence of the constant in the negative
part of the right hand side. Under the following a priori smallness
assumption
\begin{equation}\label{eq:smallness}
  C_1 \, \n{d_t}_{L^1_v L^\infty_x (1+|v|^2)} \le K \, e^{-C' \, t}
\end{equation}
we have
\begin{multline*}
\n{d_t}_{L^1_v L^\infty_x(1+|v|^{k'})} \le 2 \,
\epsilon \, \exp\left( C_2 \, \int_0 ^t \n{g_\tau}_{L^1_v
     (1+|v|^{k'+1})} \dd \tau\right) \\
 \le 2 \,\epsilon \, \exp\left( C_2 \, C_g \, \int_0
   ^t \min \{ 1, t^{-\beta} \} \dd \tau\right) \le e^{C'_g \,
   t} \, 2 \, \epsilon 
\end{multline*}
for some $\beta <1$. 

We then define 
\[
T_1 = T_1(\epsilon) = \frac{- \log \epsilon}{Q C'_g} \in (0,+\infty)
\]
for $Q$ to be chosen later, which yields 
\[
\fa t \in \left[0,T_1 \right], \quad e^{C'_g \,t} \, 2 \, \epsilon \le
2 \, \epsilon^{1-1/Q} \quad \mbox{and} \quad K \, e^{-C \, t} \ge K \,
 \epsilon^{\frac{C}{C'_g Q}}. 
\]
We then choose $\epsilon$ small enough so that
\begin{multline*}
\fa t \in [0,T_1], \quad C_1 \, \n{d_t}_{L^1_v L^\infty_x (1+|v|^2)}
\le C_1 \, e^{C'_g \,
   t} \, 2 \, \epsilon \\ \le 2 \, \epsilon^{1-1/Q} \le K \,
 \epsilon^{\frac{C}{C'_g Q}} \le K \, e^{-C \, t}
\end{multline*}
which is always possible as soon as 
\[
1 - \frac{1}{Q} > \frac{C}{C'_g Q}
\]
which can be ensured (uniformly as $\epsilon$ goes to zero) by taking
$Q$ large enough. This then implies the smallness condition
\eqref{eq:smallness} and thus justifies the a priori estimate. We
deduce the a priori bound
\[
\forall \, t \in [0,T_1], \quad \n{d_t}_{L^1_v L^\infty_x(1+|v|^{k'})}
\le 2 \, \epsilon^{1-1/Q}.
\]
Observe that $T_1(\epsilon) \to +\infty$ as $\eps \to 0$. The actual
construction and uniqueness of these solutions relies on the part {\bf
  (I)} of Theorem~\ref{theo:NLBE}: one uses the continuity of the
flow~\eqref{eq:estimate-presque-energy} and the scheme
\begin{equation*}
  \dt d^{n+1}_t + v \cdot \nabla_x d^{n+1}_t = 
  2 \, Q(d^n_t,d^{n+1}_t) + 2 \, Q(g_t, d^n_t). 
\end{equation*}
We skip the details of these standard arguments. 
  \medskip

\noindent
{\em Step~3. The trapping mechanism.} 
Let $\delta$ be the smallness constant of the stability neightborhood
in the part {\bf (II)} of Theorem~\ref{theo:NLBE} in $L^1_v
L^\infty_x(1+|v|^{k'})$. Then from the step 1 we know that there is
some time $T_2 >0$ depending on $\init{g}$ such that
\[
\forall \, t \ge T_2, \quad \| g_t - \mu \|_{L^1_v(1+|v|^{k'})} \le
\frac{\delta}{2}.
\]

We then choose $\epsilon$ small enough such that
% \[
% \| \init{f} - \init{g}
%  \|_\EE \le \epsilon \Rightarrow \| \mu_g - \mu_f
% \|_\EE \le \frac{\delta}{3}
% \]
% where $\mu_f$ is the global Maxwellian equilibrium associated with
% $\init{f}$, and such that
$T_1(\eps) \ge T_2(M)$ and thus
\[
\n{f_{T_2} - g_{T_2}}_{L^1_vL^\infty_x(1+|v|^{k'})} \le \frac{\delta}{2}, 
\]
from the step 3. 

It holds
\[
\n{f_{T_2} -\mu}_{L^1_vL^\infty_x(1+|v|^{k'})} \le \n{f_{T_2} -
  g_{T_2}}_{L^1_vL^\infty_x(1+|v|^{k'})} 
+ \n{g_{T_2} -  \mu}_{L^1_v(1+|v|^{k'})} \le \delta
\]
and we can therefore use the perturbative Theorem~\ref{theo:NLBE} for
$t \ge T_2$ which concludes the proof. 

%%%%%%%%%%%%%%%%%%%%%%%%%%%%%%%%%%%%%%

\subsection{Proof of Theorem~\ref{theo:expEB}}
\label{sec:proof-theor-refth-2}

We now turn to the proof of the exponential $H$-theorem. Let us first
recall existing results for polynomially decaying solutions of the
nonlinear equation:
\begin{theo}[\cite{Desvillettes-Villani-EB}]\label{theo:DV}
  Let $(f_t)_{t \ge 0}$ be a non-negative non-zero smooth solution of
  \eqref{eq:NLBE} such that for $k,s \ge 0$ big enough
\[
\sup_{t \ge 0} \left( \| f_t \|_{H^s(\T^3 \times \R^3)} + \| f_t
  \|_{L^1(1+|v|^{k+1})} \right) \le C < +\infty
\]
with initial data satisfying the lower bound~\eqref{eq:hyp-lower}.

Then for $k' \in (2,k)$, there exists an explicit polynomial function
$\varphi=\varphi(t)$ which goes to zero as $t$ goes to infinity such
that
\[
\forall \, t \ge 0, \quad \| f_t - \mu \|_{L^1_v L^\infty_x(1+|v|^{k'})} \le
\varphi(t)
\]
where $\mu$ is the global Maxwellian equilibrium associated with $f$
(same mass, momentum and temperature).
\end{theo}

\begin{proof}[Proof of Theorem~\ref{theo:DV}]
  This theorem is a consequence of
  \cite[Theorem~2]{Desvillettes-Villani-EB} about convergence to
  equilibrium for a priori smooth solutions with bounded moments and
  satisfying a Gaussian lower bound, and of part {\bf (I)} of
  Theorem~\ref{theo:NLBE} where we indeed establish such lower bounds.
  Note that the convergence in
  \cite[Theorem~2]{Desvillettes-Villani-EB} is measured in relative
  entropy, but it is a simple computation based on the
  Csisz\'ar-Kullback-Pinsker inequality (see for instance
  \cite[Chapter~9]{VillaniTOT}) and some interpolation to
  translate it into stronger norms such as the one we propose in the
  statement.
\end{proof}

Finally, combining all the previous results we can prove
Theorem~\ref{theo:expEB} as follows: we use Theorem~\ref{theo:DV} for
initial times and Theorem~\ref{theo:NLBE} for large times. The former
theorem provides an explicit time for the solution to enter the
trapping neightborhood in $L^1_v L^\infty_x(1+|v|^{k'}))$ norm of the
latter theorem. Then we write
\begin{multline*}
\int_{\T^d \times \R^3} f_t \, \log \frac{f_t}{\mu} \dd x \dd v =
\int_{\T^d \times \R^3} \left( \frac{f_t}{\mu} \, \log \frac{f_t}{\mu}
  - \frac{f_t}{\mu} +1 \right) \,\mu \dd x \dd v 
  \\
%\le \int_{\T^d \times \R^3} \left| \frac{f_t}{\mu} \, \log \frac{f_t}{\mu}
%  - \frac{f_t}{\mu} +1 \right| \,\mu \dd x \dd v \\ 
  \le 
\int_{\T^d \times \R^3} \left|  \log \frac{g_t}{\mu}
\right| \, \left| {f_t  \over \mu}-1 \right| \,\mu \dd x \dd v
\end{multline*}
for some 
\[
g_t=g_t(x,v)\in \left[\min\{\mu(v);f_t(x,v)\},\max\{\mu(v);f_t(x,v)\}\right]
\]
from the mean-value theorem.  On the one hand, if $f_t(x,v) \ge \mu(v)$
then
\begin{multline*}
  \left|\log \frac{g_t(x,v)}{\mu(v)} \right| \le \log
    \frac{f_t(x,v)}{\mu(v)}
  \le \log \left(1 + \frac{\| h_t \|_{L^\infty}}{\mu(v)}\right) 
  \\
  \le \max\left\{1,\sup_{t \ge 0} \|h_t \|_{L^\infty}\right\} \, \log
  \mu(v)^{-1} = K_1 \, (1 + |v|^2).
\end{multline*}
Moreover, if $f_t(x,v) \le \mu(v)$ one can use the exponential lower
bound $f_t(x,v) \ge A \, e^{- a \, |v|^2}$, $a > 1/2$, to get
\[
\left|\log \frac{g_t(x,v)}{\mu(v)}\right| \le   \log
\frac{\mu(v)}{f_t(x,v)} \le  K_2 \, (1 + |v|^2).
\]

Using the
bounds on the solution we hence finally deduce 
\[
\int_{\T^d \times \R^3} f_t \, \log \frac{f_t}{\mu} \dd x \dd v \le C
\, 
\int_{\T^d \times \R^3} \left| f_t - \mu \right| \, (1+|v|^2) \dd x \,
\dd v
\]
and we conclude the proof using the estimate of convergence in $L^1_v
L^\infty_x(1+|v|^2)$. 

\subsection{Structure of singularities for the nonlinear flow}
\label{sec:struct-sing-nonl}

Let us now study the singularity structure of the nonlinear flow
provided by the perturbative theorems~\ref{theo:NLBE}
and~\ref{theo:NLBE+}. We shall prove the following two properties as
we did for the linearized flow: first we show that the dominant part
of the flow in the asymptotic behavior is as regular as wanted.
Second, we prove that its worst singularities are supported by the
free motion characteristics.

\subsubsection{Asymptotic amplitude of the singularities}
\label{sec:asympt-ampl-sing}

Let us consider for instance the space $L^1_v L^\infty_x (1+|v|^k)$,
$k>2$ (other spaces satisfying the assumption of the perturbative
theorems could be used obviously). We consider some initial data
$\init{f} = \mu + \init{h} \ge 0$ in this space and assume without
loss of generality that $\Pi \init{h} =0$ (which implies $\Pi h_t =0$
for any later time). 

We start from the decomposition of the semigroup 
\begin{equation*}
  S_\LL(t) \, \init{h} =  S^s_\LL(t) \, \init{h} + S^r_\LL(t) \, \init{h}
\end{equation*}
we have introduced in Subsection~\ref{sec:study-sing-struct}. Then we
write a Duhamel formulation
\begin{multline*}
  h_t = S_\LL (t) \, \init{h} + \int_0 ^t S_\LL (t-\tau) \,
  Q(h_\tau,h_\tau) \dd \tau \\
  = \left( S^s_\LL(t) \, \init{h} + \int_0 ^t S^s _\LL (t-\tau) \,
    Q(h_\tau,h_\tau) \dd \tau \right) \\ + \left( S^r_\LL(t) \,
    \init{h} + \int_0 ^t S^r _\LL (t-\tau) \, Q(h_\tau,h_\tau) \dd
    \tau \right) =: \NN^s (t) + \NN^r(t)
\end{multline*}
(we have used here that $\Pi Q(h_\tau,h_\tau)=0$). Since 
\[
\left\{ 
  \begin{array}{l} \ds \n{ S^s_\LL(t) \, h}_{H^s_{x,v}(\mu^{-1/2})} \le
    C \, \n{h}_{L^1_{x,v}(1+|v|^k)} \, e^{-\lambda \, t} \vs \\ \ds
    \n{ S^r_\LL(t) \, h}_{L^1_{x,v}(1+|v|^k)} \le C \, \n{h}_{L^1_{x,v}(1+|v|^k)}
    \, e^{-(\nu_0-\var) \, t},
\end{array}
\right.
\]
and the nonlinear flow is known to have the decay 
\begin{equation*}
  \n{h_t}_{L^1_{x,v}(1+|v|^k)} \le C \, e^{-\lambda \, t},
\end{equation*}
we deduce that 
\begin{equation*}
  \left\{ 
\begin{array}{l} \ds
\n{ \NN^s(t)}_{H^s_{x,v} (\mu^{-1/2})} \le C \, e^{-\lambda \, t} \vs \\ \ds 
\n{ \NN^r(t)}_{L^1_{x,v}(1+|v|^k)} \le C \, e^{- \min\{ \nu_0 - \var; \, 2
  \lambda \} \, t}
\end{array}
\right.
\end{equation*}
(the factor $2$ in the exponent of the second inequality comes from
the quadratic nature of the nonlinearity). 

Then one can perform a boostrap argument in order to deduce finally
\begin{equation*}
  h_t = \bar \NN^s (t) + \bar \NN^r(t)
\end{equation*}
with 
\begin{equation*}
  \left\{ 
\begin{array}{l} \ds
\n{ \bar \NN^s(t)}_{H^s_{x,v} (\mu^{-1/2})} \le C \, e^{-\lambda \, t} \vs \\ \ds 
\n{ \bar \NN^r(t)}_{L^1_{x,v}(m)} \le C \, e^{- (\nu_0 - \var) \, t}.
\end{array}
\right.
\end{equation*}

Let us sketch the bootstrap argument. If $2 \lambda \ge \nu_0-\var$ we
are done. Suppose therefore that $2 \lambda < \nu_0 - \var$. Then plug
the decomposition $h_t = \NN^s (t) + \NN^r(t)$ into the
Duhamel formulation: 
\begin{multline*}
  h_t = S_\LL (t) \, \init{h} + \int_0 ^t S_\LL (t-\tau) \,
  Q(h_\tau,h_\tau) \dd \tau \\
  = \left( S^s_\LL(t) \, \init{h} + \int_0 ^t S^s _\LL (t-\tau) \,
    Q(h_\tau,h_\tau) \dd \tau \right) \\ + S^r_\LL(t) \, \init{h} +
  \int_0 ^t S^r _\LL (t-\tau) \, Q\left(\NN^r
    (\tau),\NN^r (\tau)\right) \dd \tau  \\
  + \int_0 ^t S^r _\LL (t-\tau) \, Q\left(\NN^r (\tau),\NN^s
    (\tau)\right) \dd \tau+ \int_0 ^t S^r _\LL (t-\tau) \,
  Q\left(\NN^s (\tau),\NN^s (\tau)\right) \dd \tau.
\end{multline*}
Then observe that in the decomposition of the linearized flow one has
\[
\n{ S^r_\LL(t) \, h}_{H^s_{x,v}(\mu^{-1/2})} \le C \, \n{h}_{H^s_{x,v}(\mu^{-1/2})}
    \, e^{- \lambda \, t}.
\]
Therefore if one defines 
\begin{multline*}
\tilde \NN^s(t) := \left( S^s_\LL(t) \, \init{h} + \int_0 ^t S^s _\LL (t-\tau) \,
    Q(h_\tau,h_\tau) \dd \tau \right) \\ + \int_0 ^t S^r _\LL (t-\tau) \,
  Q\left(\NN^s (\tau),\NN^s (\tau)\right) \dd \tau
\end{multline*}
and 
\begin{multline*}
\tilde \NN^r(t) := S^r_\LL(t) \, \init{h} +
  \int_0 ^t S^r _\LL (t-\tau) \, Q\left(\NN^r
    (\tau),\NN^r (\tau)\right) \dd \tau  \\
  + \int_0 ^t S^r _\LL (t-\tau) \, Q\left(\NN^r (\tau),\NN^s
    (\tau)\right) \dd \tau,
\end{multline*}
one checks that 
\begin{equation*}
  \left\{ 
\begin{array}{l} \ds
\n{ \tilde \NN^s(t)}_{H^s_{x,v} (\mu^{-1/2})} \le C \, e^{-\lambda \, t} \vs \\ \ds 
\n{ \tilde \NN^r(t)}_{L^1_{x,v}(1+|v|^k)} \le C \, e^{- \min\{ \nu_0 - \var; \, 3
  \lambda \} \, t}
\end{array}
\right.
\end{equation*}
(notice the factor $3$ in argument of the exponential). Hence by
iterating this argument a finite number of times, one gets the
conclusion.

In a way similar to the linear setting, the nonlinear flow splits in two
parts. The first one has the following properties: (1) it is as smooth
as wanted , (2) has Gaussian decay in the small linearization space,
(3) the exponential time decay rate is sharp. The second part of the
solution decays exponentially in time with a rate as close as wanted
to $\nu_0$, the onset of the continuous spectrum, and carries all the
singularities.

\subsubsection{Localization of the $L^2$ singularities}
\label{sec:local-l2-sing}

We consider now the space $L^\infty_{x,v}(1+|v|^k)$, $k >6$ (again
other spaces could be considered). We know that the solution $h_t$ to
the nonlinear equation remains uniformly bounded in this space along
time and decays exponentially fast to zero as time goes to
infinity. We start again from the Duhamel formula. In
Subsection~\ref{sec:study-sing-struct} we showed the following
decomposition of the linearized semigroup
\[
S_\LL(t) \, \init{h} \in (\mbox{Id} - \Pi_{\LL,0}) \, \left(
  e^{-\nu(v) \, t} \, \init{h}(x-vt,v) \right)  + O(t^{-\theta}) \, H^\alpha_{x,v,loc}
\]
for some small $\alpha>0$ and some $\theta>0$. We can then prove arguing
exactly as in \cite{MR1798557} that
\begin{equation*}
  \int_0 ^t S_\LL (t-\tau) \, Q(h_\tau,h_\tau) \dd \tau \in H^\alpha_{x,v,loc}
\end{equation*}
for some small $\alpha>0$, due to the velocity-averaging nature of the
bilinear collision operator. This proves finally that the nonlinear
solution satisfies 
\[
h_t \in (\mbox{Id} - \Pi_{\LL,0}) \, \left(
  e^{-\nu(v) \, t} \, \init{h}(x-vt,v) \right) +
O(t^{-\theta}) \, H^\alpha_{x,v,loc}
\]
which captures the localization of the $L^2$ singularities. 
  
\subsection{Open questions}

A first natural question is whether our methods could be extended to
the case of Boltzmann equations with \emph{long-range
  interactions}. In the case of non-cutoff hard and moderately soft
potentials, the linearized operator has a spectral gap
\cite{MR2322149,MR2784329} and we expect our factorization method to
be applicable in this case by using a different decomposition of the
linearized collision operator, such as the one used in
\cite{MR2254617} in order to quantify the spectral gap in velocity
only. In the case of very soft potentials, the linearized collision
operator does not have a spectral gap anymore and the expected time
decay rate is a stretched exponential. It is an interesting question
to investigate whether our factorization method could be used when
\emph{generalized coercivity estimates} replace spectral gap
estimates. Another direction opened by this work is the question of
obtaining spectral gap estimates in physical space for kinetic
equations in the whole space confined by a potential (a work is in
progress in the case of the kinetic Fokker-Planck equation in the
whole space).

We end up with what seems to us the most interesting open question
suggested by this study. In contrast with many dispersive or fluid
PDE's, the Boltzmann equation (and kinetic equations in general) does
not seem to have a clear notion of \emph{critical space}, and it has
been debatted whether such a notion would indeed apply to it. Our
perturbative study proves that the space $L^1_v L^\infty_x
(1+|v|^{2+0})$ is supercritical. But what is more interesting is that
as far as the velocity variable is concerned the space $L^1_v
(1+|v|^2)$ \emph{is} critical, as shown by the studies
\cite{MR1697562,MR1716814} in the spatially homogeneous
case. Therefore we can now focus on the spatial variable \emph{only}
in order to identify a critical space ``below'' $L^\infty_x$. A first
step in this direction would be to use averaging lemma on the
nonlinear flow in order to prove perturbative well-posedness in $L^1_v
L^p_x (1+|v|^{2+0})$ for some $p <+\infty$ possibly large but not
infinite. A natural conjecture is then to ask for the critical space
in the variable $x$ to be compatible with the incompressible
hydrodynamic limit (which is ``blind'' to the functional space used in
the velocity variable roughly speaking) and therefore to be
$L^3_x(\T^3)$ as for the three-dimensional incompressible
Navier-Stokes equations.

\bigskip

%\newpage 
\bibliographystyle{acm}
%\bibliography{biblio,mybiblio}
\bibliography{./spectreFP}

%\begin{thebibliography}{100}

\bigskip
\bigskip
\signmpg \signsm \signcm 

\end{document}